\documentclass[aos]{imsart}

\RequirePackage{amsthm,amsmath,amsfonts,amssymb}
\RequirePackage[numbers,sort&compress]{natbib}
\RequirePackage[colorlinks,citecolor=blue,urlcolor=blue]{hyperref}
\RequirePackage{graphicx}

\startlocaldefs
\theoremstyle{plain}

\theoremstyle{definition}

\usepackage[colorlinks]{hyperref}
\usepackage[dvipsnames]{xcolor}

\usepackage{graphicx}

\usepackage{color, latexsym}
\usepackage{booktabs}

\graphicspath{{figures/}}

\usepackage{booktabs}
\usepackage{multirow}

\usepackage{bbm}

\theoremstyle{plain} 
\newtheorem{thm}{Theorem}
\newtheorem{Corollary}{Corollary}
\newtheorem{prop}{Proposition}
\newtheorem{lem}{Lemma}
\newtheorem{Remark}{Remark}

\newcommand{\R}{\mathbb{R}} 
\newcommand{\cD}{\mathcal{D}} 
 
\newcommand{\cO}{\mathcal{O}} 
\newcommand{\cC}{\mathcal{C}} 
\newcommand{\cN}{\mathcal{N}} 
\newcommand{\cF}{\mathcal{F}} 
\newcommand{\cK}{\mathcal{K}}

\usepackage{float}
\usepackage{subcaption}

\def\1{1\!{\rm l}}

\newcommand{\leqa}{\lesssim}

\newcommand{\EM}{\ensuremath}

\newcommand{\al}{\alpha}
\newcommand{\be}{\beta}

\newcommand{\ga}{\gamma}

\newcommand{\la}{\lambda}

\newcommand{\si}{\sigma}

\newcommand{\ta}{\tau}
\newcommand{\veps}{\varepsilon}
\newcommand{\vphi}{\varphi}

\newcommand{\ze}{\zeta}

\newcommand{\cA}{\EM{\mathcal{A}}}
\newcommand{\cB}{\EM{\mathcal{B}}}
\newcommand{\cE}{\EM{\mathcal{E}}}
\newcommand{\cG}{\EM{\mathcal{G}}}
\newcommand{\cH}{\EM{\mathcal{H}}}
\newcommand{\cI}{\EM{\mathcal{I}}}
\newcommand{\cJ}{\EM{\mathcal{J}}}

\newcommand{\cS}{\EM{\mathcal{S}}}
\newcommand{\cT}{\EM{\mathcal{T}}}
\newcommand{\cU}{\EM{\mathcal{U}}}

\newcommand{\cX}{\EM{\mathcal{X}}}

\newcommand{\cZ}{\EM{\mathcal{Z}}}

\newcommand{\psg}{{\langle}}
\newcommand{\psd}{{\rangle}}

\DeclareMathOperator*{\argmin}{arg\,min}

\definecolor{blendedblue}{rgb}{0.2,0.2,0.7}

\newcommand{\given}{\,|\,}
\newcommand{\mh}{\mathbb{H}}

\definecolor{commentpaul}{rgb}{0,0.6,0}

\definecolor{modif}{rgb}{0.55,0,0.8}

\definecolor{commentrick}{rgb}{0.1,0.2,0.8} 

\endlocaldefs

\begin{document}

\begin{frontmatter}
\title{Robust Gaussian process regression: posterior contraction}
\runtitle{Robust GPs}

\begin{aug}
\author[A]{\fnms{Isma\"el}~\snm{Castillo} \ead[label=e1]{ismael.castillo@sorbonne-universite.fr}}
\author[A,C]{\fnms{Paul}~\snm{Egels}\ead[label=e2]{paul.egels@unibocconi.it}}
\author[B]{\fnms{Riccardo}~\snm{Lazzarini}\ead[label=e3]{riccardo.lazzarini@phd.unibocconi.it}}
\author[B,C]{\fnms{Botond}~\snm{Szab\'o}\ead[label=e4]{botond.szabo@unibocconi.it}}

\address[A]{Laboratoire de Probabilit\'es, Statistique et Mod\'elisation (LPSM), Sorbonne Universit\'e, 4 place Jussieu 75005 Paris (France),}
\address[B]{Bocconi Institute for Data Science and Analytics, Bocconi University, via Rontgen 1, 20136, Milano (Italy),}
\address[C]{Department of Decision Sciences, Bocconi University, via Rontgen 1, 20136, Milano (Italy)\\
\vspace{0.2cm}
\printead[presep={ }]{e1},\,\,\printead[presep={ }]{e2},\\ \printead[presep={ }]{e3},\,\,\printead[presep={ }]{e4}}

\end{aug}

\begin{abstract}
Gaussian processes form a popular and flexible class of prior distributions with principled uncertainty quantification, but
the Gaussian likelihood on which the corresponding posterior relies makes them highly sensitive to
outliers and contaminated observations. We study Gaussian process regression when among the $n$ observations, $m$ of these are possibly contaminated. 
A generic framework to obtain posterior contraction rates in this setting is considered, from a frequentist perspective as the sample size~$n$ and the number of corrupted observations~$m$ grow. We adopt a generalized Bayes approach in which the log--likelihood is replaced by a robustified version; we discuss various choices thereof, such as one recently proposed 
by Altamirano et al. \cite{pmlr-v235-altamirano24a} that preserves conjugacy. We investigate the minimax rate of estimation under adversarial contamination in terms of the Sobolev smoothness $\be$ of the regression function and $m$. The question of building adaptive procedures to both $\be$ and $m$ is investigated, and we show that a polynomial loss for adaptation must be paid in general. We instantiate our results for Mat\'ern
Gaussian process and deep ReLU network priors. Numerical experiments for existing and new robust algorithms are provided, illustrating the practical applicability of our results.

%
\end{abstract}

\begin{keyword}[class=MSC]
\kwdgroup[type=primary]{\kwd{62G05}}
\kwdgroup[type=secondary]{\kwd{62G20}}
\end{keyword}

\begin{keyword}
\kwd{Gaussian process}
\kwd{Bayesian nonparametrics}
\kwd{Robustness}
\kwd{Minimax rates}
\kwd{Adaptation}
\end{keyword}

\end{frontmatter}




\section{Introduction}

\subsection{Motivation and contributions}
\label{sec:motivation-contributions}

%

Gaussian processes (GPs) are a central tool in modern statistics, with applications ranging from spatial statistics \cite{cressiebook} to machine learning \cite{rasmussenwilliamsbook}, among their many practical uses. Their popularity stems from their ability to model complex
functions while representing uncertainty through a posterior distribution and
associated credible sets. In the broadly used Gaussian nonparametric regression model, GPs are particularly attractive computationally, since they form a conjugate family of prior distributions, so that the associated posterior distribution is available in explicit form. Most modern sampling techniques for GPs often use this conjugacy as a key ingredient, in particular ones that involve some form of approximation such as variational or Vecchia approximations.

In applications, however, the data-generating mechanism may deviate from the assumed statistical model. While for physical reasons it may still be reasonable in the regression setting to assume  that typical observations are subject to Gaussian noise, the available dataset can often be subject to outliers and contaminated observations. In the sequel we will assume that in the sample of $n$ data pairs $(X,Y)$ (for a covariate $X$ and an output $Y$), a certain number `$m$' of these pairs are possibly outliers or `contaminated' data. Possible sources of contamination can be, for instance, measurement errors or recording faults, which are often assumed to occur independently from the original, `clean', data--generating process, as in the classical Huber contamination model \cite{huber64}. Another, more general, form of contamination is {\em adversarial contamination}  in which the corrupted observations
may be chosen with full knowledge of the clean sample. This is the most demanding
notion of robustness 
as it requires no assumption on how the outliers are
generated. Such worst-case scenario is becoming a standard benchmark in machine
learning, where models are increasingly trained on large and unverified datasets
that may be subject to data poisoning or systematic corruption. We will mainly focus on this more general, adversarial contamination, but our results hold for any weaker form of data corruption as well. 


The problem of the presence of outliers is particularly acute with a Gaussian likelihood, whose quadratic loss allows extreme responses to have potentially large influence on posterior inference and prediction. To address this issue, a number of proposals have recently appeared; we discuss some of these below; our approach in the present work covers among others the proposal by Altamirano et al. \cite{pmlr-v235-altamirano24a}: the method is based, as for many robust Bayesian algorithms, on a modification of the likelihood function in Bayes' formula, down-weighting outliers. A key point though is that this modification {\em still preserves} overall conjugacy, thus keeping straightforward computability (the simulations in Section \ref{sec: numerics} below illustrate this paragraph's claims).
Indeed, when the imposed weights on the additive terms in the log-likelihood do not depend on the function $f$ being estimated,
the weighted loss remains quadratic in $f$; under a GP prior the generalized
posterior is then again Gaussian and available in closed form, exactly as in
standard GP regression \cite{rasmussenwilliamsbook}. Robustness is thus achieved without giving up the
computational tractability of the Gaussian model, and in particular without
resorting to MCMC or variational approximations.



One main goal of this work is to provide a principled framework under which
Gaussian process regression is made robust to (the just-described) data contamination settings while still keeping attractive computational properties of standard GPs. While a substantial literature studies the frequentist
asymptotic properties of GP posteriors under a clean sample (including \cite{vandervaartgpcontraction, JMLRvdvvz}, see also \cite{ghosal_van_der_vaart_2017, icstf} for general accounts of the theory), theoretical understanding in the contaminated case is still limited. The work \cite{pmlr-v235-altamirano24a} proposes a control of the so--called Posterior influence function (PIF), that controls by how much a single corrupted observation perturbs the posterior. However, it remains unclear if the robustified posterior distribution converges and if it does, at which rate. In particular, one would expect that if the number $m$ of outliers is quite small, one can still achieve standard optimal minimax contraction rates.

Our contributions are as follows:
\begin{enumerate}
\item We propose a theoretical framework to derive contraction rates for robustified posteriors in nonparametric regression. Robustness is guaranteed under the adversarial contamination scheme, and is achieved by replacing the likelihood of each observation by a `reweighted' version. 
Our results are not specific to a class of priors, although due to their broad use in practice we mostly focus on Gaussian process priors (neural network priors are also considered). Under prior mass type conditions, we show that this approach can achieve a contraction rate of order $\veps_n(\be)+\sqrt{m/n}$, where $\veps_n(\be)$ is the usual nonparametric rate for estimating a $\be$--smooth function without contamination and $m$ the number of outliers in the dataset. This rate can be achieved both by deterministic and hierarchical (random) tuning of the priors.
\item Our theory covers the recent proposal by Altamirano et al. \cite{pmlr-v235-altamirano24a}, that implicitly corrects the `bias' induced by the weight function using a score-matching type idea while retaining conjugacy (which enables fast posterior computations); it also covers a variety of different weights, from simple ones (thresholding weights) to more elaborated ones (weights depending on both outcomes and covariates).
\item We investigate optimality from the minimax point of view in this context. When both the regularity $\be$ of the true regression function and number $m$ of outliers are given, we determine the minimax estimation rate  in terms of the quadratic loss 
(up to a logarithmic factor). 
In the high contamination setting, the rate is strictly faster than the rate derived by the generic approach in point 1. We provide a novel generalized Gaussian Process approach, based on the trimming procedure introduced in \cite{robust_mean_estimation}, reaching this minimax lower bound. This procedure requires the knowledge of $m$.
\item We prove the at first sight surprising result that, without the knowledge of $\be$ and $m$ (that is, the most relevant practical situation), the rate $\veps_n(\be)+\sqrt{m/n}$ achieved by Gaussian process--based methods (once made adaptive to $\be$) in fact cannot be improved. This can be explained by the fact that simultaneous adaptation to the smoothness $\be$ and the number of outliers $m$ is impossible, and that there is a polynomial--in--$n$ loss for adaptation.
\item We carry out a simulation study for a number of weight functions suggested by the previous theoretical insights. Our findings illustrate the robustness of the proposed algorithms. The excellent empirical behaviour of the debiased weights approach observed in \cite{pmlr-v235-altamirano24a} is confirmed. We also introduce a novel trimming based GP method. We propose a number of refinements  -- whose optimality are justified from the obtained theory -- to improve the overall performance of both methods. 
\end{enumerate}

\subsection{Related work}
\label{sec:related-work}

{\em Heavy-tailed likelihoods.}
The most common Bayesian strategy for robustness replaces the Gaussian
observation noise with a heavy-tailed distribution, most notably a Student-$t$
likelihood \cite{JMLR:v12:jylanki11a}. The theoretical appeal of this idea rests
on a classical body of work on automatic outlier rejection: for suitably
heavy-tailed models, as an observation is sent to infinity the posterior reverts
to the one obtained by discarding that observation, so that extreme points are
ignored without being modelled explicitly \cite{dawid1973posterior,
ohagan1979outlier,west1984outlier,andrade2006bayesian,desgagne2015robustness,
ohagan2012bayesian}. This behaviour is appealing, but it is intrinsically
asymptotic in the magnitude of the outlier and is tailored to stochastic
heavy-tailed noise. Such heavy-tailed likelihoods break 
the conjugacy of GP regression, so that posterior computation generally requires Laplace approximations, variational inference, or MCMC. Moreover, such method lacks rigorous theoretical underpinning, especially in the adversarial contamination setting.

{\em Mixture and contamination likelihoods.}
A second approach introduces latent outlier indicators and models the noise as a
mixture of a nominal and a contaminating component
\cite{boxBayesianrob, LuFangTian}. Such models can be effective when the
contamination mechanism is well specified, but the resulting robustness
guarantee depends on the chosen contaminating distribution. For GP priors, the
closed-form conditioning formulas are again typically lost, even when
conditionally conjugate latent-variable representations are available, resulting in a substantially higher computational complexity.

{\em Outlier detection for Gaussian processes.}
Several robust GP methods explicitly detect, trim, or downweight suspicious
observations, including trimmed marginal-likelihood procedures and Robust
Gaussian Processes via Relevance Pursuit
\cite{pmlr-v216-andrade23a,relevancePursuitAment}. Our methods are similar in spirit,
but embed the downweighting in a single generalized-posterior construction and come with theoretical contraction and minimax guarantees.

{\em Generalized and Gibbs posteriors.}
Generalized Bayes replaces the negative log-likelihood by a general loss, thereby
decoupling Bayesian updating from a fully specified sampling model
\cite{bissiriWalker,safeBayesian,SyringMartin}. In this spirit, inspired by the $\rho$--estimation method they previously introduced, Baraud and Birgé \cite{baraudbirge20, baraud24} introduce general approaches that lead to robustness for two specific generalized Bayes posterior distributions; while leading to very broad results, the computability of these approaches is yet little understood (a first step in this direction is the recent work \cite{khribchalquier26} that explores a variational approximation version thereof). 
A number of works in the machine learning literature consider generalized Bayes approaches, in particular 
the robust and conjugate GP construction of
\cite{pmlr-v235-altamirano24a}, that is particularly appealing for its computability, and that we cover as one of the examples of our theory. 
Related robust generalized-posterior methods include
MMD-Bayes \cite{pmlr-v118-cherief-abdellatif20a}, robust median-of-means posteriors
\cite{minsker2014scalable, minsker2017robust, yao2025generalized}, and parametric M-posteriors
\cite{marusic2025theoretical}, which establish, in various settings, robustness to
misspecification or outliers, concentration and Bernstein--von Mises results, and
posterior analogues of classical robustness measures. PAC-Bayesian analyses of
MMD-losses offer a further route to finite-sample guarantees
under contamination \cite{pmlr-v118-cherief-abdellatif20a}. A distinct, fully Bayesian route to robustness is the coarsened posterior of
\cite{Miller03072019}, which conditions on the event that the data lie within a
neighbourhood of the model rather than reproducing it exactly; although not a
generalized-Bayes construction, it connects to these ideas because, for
relative-entropy neighbourhoods, it coincides with a tempered (power) posterior.

{\em Frequentist robust statistics and minimax theory.}
Robust estimation under contamination has a long history, from the classical work
of Huber and Tukey \cite{Huber_Ronchetti_2009} to the modern computational and minimax theory
\cite{Diakonikolas2019RecentAI,chen_robust_cov, robust_mean_estimation, liu2019density}. Also in a regression setting, \cite{chinot} studies robust estimation under maliciously corrupted labels via regularized empirical risk minimization with a robust convex loss and an RKHS-norm penalty. This is loosely connected to
our approach, since such a penalized estimator can be viewed as the MAP of a generalized posterior under a GP prior, though the robust loss they employ differs from ours. Most closely related to the present work, \cite{zhao2024robust} is, to our knowledge, the first to carry out a minimax analysis of nonparametric regression under an adversarial poisoning attack. 
 Although the latter work can be seen as an interesting first step, it is limited in several aspects: regularities are confined to Lipschitz and the considered estimator depends on $m$.  
The work \cite{liu2019density} considers density estimation under pointwise loss and the Huber contamination model, investigating minimax rates and showing that adaptation to both $\beta$ and $m$ is impossible. However, the loss for adaptation is not quantified and the implementation aspect of the procedures is not considered.

{\em Bayesian nonparametric contraction theory.}
Our proofs build on the general theory of posterior contraction and its
extensions to fractional, tempered, and Gibbs posteriors
\cite{ghosal_van_der_vaart_2017,temperedPosteriorBattha, JMLR:aliceIsmaelFractional, SyringMartin}. We adapt these tools to weighted
losses and contaminated samples, obtaining rates that depend explicitly on the
number of corrupted observations.

\subsection{Outline}
\label{sec:outline}

In Section \ref{sec: prelim}, we introduce the adversarial contamination setting in the nonparametric regression model and the generalized Bayesian approach, which is based on reweighting the terms in the likelihood function. In Section \ref{sec : w(y)}, we derive posterior contraction guarantees for general nonparametric priors and weights that solely depend on the response variable. We begin with a simple thresholding weight function and then extend the analysis to more general weights, including examples from the machine learning literature. We also verify that Gaussian process (GP) and deep neural network (DNN) priors satisfy the conditions of our general theorems.  
In Section \ref{sec : more complex modifications}, 
 we derive contraction rate results for more general weight functions, which may also depend on the covariates and the regression function.
In Section \ref{sec: minimax}, we derive minimax lower bounds for adversarial contamination and propose a trimming-based frequentist approach and a generalized GP approach that attain this lower bound. We also show that when the regularity of the regression function and the contamination level are unknown, attaining this lower bound is impossible, and a polynomial loss in the rate occurs. We show that the hierarchical version of the reweighted GP procedure attains the feasible adaptive rate. 
In Section \ref{sec: numerics}, we illustrate the performance of the robust approaches studied under various contamination scenarios using synthetic datasets. Finally, in Section \ref{sec:discussion}, we briefly summarize our results and discuss extensions and future research directions. Proofs and additional numerical analyses are deferred to the Supplementary Material. 

\section{Setting}
\label{sec: prelim}

\textit{Nonparametric regression.} Consider the random design Gaussian regression model, where the i.i.d. 
sample $\cD_n = ((x_i,y_i) \in \R^d \times \R, 1 \le i \le n)$ is generated as 
	\begin{equation}\label{def : model}
		y_i = f_0(x_i) + \sigma \xi_i,\qquad x_i\stackrel{\text{i.i.d.}}{\sim}P_X\, \perp\, \xi_i\stackrel{\text{i.i.d.}}{\sim}\cN(0,1),\qquad i=1,...,n,
	\end{equation}
    where  $f_0 \in \cF $ denotes the unknown regression function of interest  and $\cF$ is a subset of continuous functions in $L^2(P_X)$ (square--integrable functions with respect to $P_X$).
    Let us denote by ${P_{f_0}}$ the distribution of a single observation $(x_i,y_i)$, and by $P_{f_0}^n$ and $E_{f_0}$ the joint distribution and expectation of the $n$-sample $\cD_n$, respectively. Our goal is to estimate the ``true'' function $f_0$ using the population loss, for $f \in L_2(P_X)$,
\[ \lVert f-f_0\rVert_{2,P_X}^2 : = \int (f-f_0)^2 \, dP_X.\]
Throughout this paper we  assume that there exists some (unknown) $F_0 > 0$, such that 
\begin{equation}\label{eq : bound on f_0}
    \|f_0\|_\infty := \sup_{x \in \R^d} |f_0(x)| \leq F_0.
\end{equation}

\noindent
\textit{Generalized posterior.} Given a prior distribution $\Pi$ on $\cF$ and a weight function $w : \cF \times \R^d \times \R \to \mathbb{R}_+$, we consider the \emph{generalized posterior}, defined as 
 \begin{equation}\label{def : genpost}
 	dG(f \, | \, \cD_n) := \frac1Z \exp \left\{ - \frac1{2\sigma^2} \sum_{i =1}^n w(f,x_i,y_i) (f(x_i) - y_i)^2 \right\} \, d\Pi(f),
\end{equation}
where $Z$ is the normalizing constant (so that the last display integrates to $1$). For any subset of indices $\cI \subseteq \{1, \dots , n\}$, we denote the weighted negative log--likelihood summed on $\cI$ as
\begin{equation}\label{def : negloglik}
    L_{\cI,w}(f) := \frac1{2\sigma^2} \sum_{i\in \cI} w(f, x_i,y_i) (f(x_i) - y_i)^2.
\end{equation}
For any $m \in \{1, \dots,n\}$, we denote $[m] := \{1, \dots, m\}$ so that
$dG(f \, | \, \cD_n) \propto e^{-L_{[n],w}(f)} \, d\Pi(f).$ Whenever the weight function $w\equiv \lambda$, for some positive constant $\lambda$, we keep this notation and define the \emph{$\lambda$-tempered posterior} 
\begin{equation}\label{def : lambda-post}
    d\Pi_{\lambda}[f \, | \, \cD_n] : \propto \exp\left\{ - \frac{\lambda}{2\sigma^2} \sum_{i =1}^n  (f(x_i) - y_i)^2 \right\} \, d\Pi(f) \propto e^{-L_{[n], \lambda}(f)} \, d\Pi(f).
\end{equation}
With this definition, the usual Bayesian posterior is $d\Pi_{1}[\, \cdot \, | \, \cD_n] \propto e^{-L_{[n],1}(f)}d\Pi(f)$.


\medskip
\noindent
\textit{Adversarial contamination.} Throughout the paper, the statistician does not have access to the clean dataset $\cD_n$, but only to a corrupted version
\[\widetilde{\cD}_n = ((\tilde{x}_1,\tilde{y}_1), \cdots,( \tilde{x}_n,\tilde{y}_n)).\]
We assume that there are $m$ possibly contaminated observations with $m/n\leq \eta<1$. Both the indices and the values of the corrupted observations may depend on the full clean dataset $\cD_n$ and on some additional randomness. More precisely, let $U\sim\nu$ be an auxiliary source of randomness, independent of $\cD_n$. We denote by $\mathcal T_m$ the class of measurable maps $T$ such that for arbitrary realizations $d$  of the clean data $\cD_n$ and $u$ of the auxiliary random variable $U$ respectively, $d$ and $T(d,u)$ differ in at most $m$ coordinates. The adversarially corrupted dataset is written as
\begin{equation}
\label{def : contamination map}
    \widetilde{\cD}_n := T(\cD_n,U),
    \qquad T\in\mathcal T_m.
\end{equation}
Thus the adversary may use the entire clean sample when deciding which observations to modify and what values to replace them with. The statistician is not assumed to know either the corrupted indices or the map $T$. 

The classical Huber contamination model \cite{Huber_Ronchetti_2009}
can be viewed as a milder  
form of contamination than the adversarial setting considered here. In Huber's model, each observation is replaced, independently and with probability $\veps$, by an independent draw from an arbitrary distribution $H$, so that the observations $x_i^{\text{Hub}}, y_i^{\text{Hub}}$ in Huber's model are 
\[
    (x_i^{\text{Hub}},y_i^{\text{Hub}})
    \overset{\mathrm{i.i.d.}}{\sim}
    (1-\veps)P_{f_0}+\veps H.
\]
Although the number of replacements is random, choosing $\veps=(1 - \delta) m/n$ for some fixed $\delta \in (0,1)$ ensures that at most $m$ observations are replaced
with probability tending to one as $n\to\infty$. Thus, with probability tending to one, the Huber-contaminated
sample coincides with $T(\cD_n,U)$ for a suitable contamination
map $T\in\cT_m$. We make this comparison precise in Section~\ref{sec : lower bound}, where we use it in particular to transfer some minimax lower bounds from the Huber model to the adversarial setting. 

The additional freedom in our adversarial setting is that both
the indices selected for corruption and their replacement values
may depend on the entire clean dataset $\cD_n$. In Huber's model,
these are generated independently of the clean sample.
This distinction is not merely formal: since the set of corrupted
indices 
\begin{equation}\label{def:outliers}
\cO:=\{i\in[n]:(\tilde{x}_i,\tilde{y}_i)\neq(x_i,y_i)\}
\end{equation}
may depend on $\cD_n$, the retained observations
$\{(x_i,y_i)\}_{i\in\cO^c}$ typically do not form an i.i.d.\ sample from
$P_{f_0}$. For example, an adversary may modify precisely the $m$ largest observations $y_i$, thereby altering the distribution of the remaining uncorrupted observations. 

\medskip
\noindent
\textit{Contraction under contamination.} To describe probabilities of events depending on the corrupted data, we first introduce $Q_{f_0}:= P_{f_0}^n \otimes \nu$ the joint law of $(\cD_n,U)$ and its expectation $E_{Q_{f_0}}$. The law of the corrupted dataset $T(\cD_n,U)$ is then expressed as the pushforward
\[ \widetilde{P}_{f_0} : = T_{\#}Q_{f_0} = \left( P_{f_0}^n \otimes \nu \right) \circ T^{-1}.\]
We denote by $\widetilde{E}_{f_0}$ the corresponding expectation.
We generally suppress the dependence on $T$ and its associated
auxiliary law $\nu$, unless it is explicitly needed.
Equivalently, by change of variables, we have for any nonnegative or integrable function $\Phi$ of the contaminated data, 
\[ \widetilde{E}_{f_0}\big[ \Phi(\widetilde{\cD}_n) \big] = E_{Q_{f_0}} \left[ \Phi(T(\cD_n , U))\right].\]
The representation on the right-hand side is useful for our proofs as it allows quantities depending on the clean and corrupted  samples to be studied on the same probability space. In particular, we remark that for an event $\cA$ solely depending on the clean data, we have $Q_{f_0}(\cA) = P_{f_0}^n(\cA).$ As a general rule, the quantities carrying a tilde involve contaminated data, for instance for any $\cI \subset [n]$, we denote
\[\widetilde{L}_{\cI,w}(f) := \frac1{2\si^2} \sum_{i\in \cI} w(f, \tilde{x}_i,\tilde{y}_i) (f(\tilde{x}_i) -\tilde{y}_i)^2.\]
The generalized posterior based on the observed corrupted dataset is therefore
\begin{equation}\label{def : conta genpost}
    dG(f\, | \, \widetilde{\cD}_n) : \propto \exp \left\{- \widetilde{L}_{[n],w}(f) \right\} \, d\Pi(f).
\end{equation}
 For a loss function $\ell :\cF \times \cF \to \R_+$, we say that the generalized posterior contracts around $f_0$ at rate $\zeta_{n,m}^2$ under contamination if, for any $M_n \to \infty$, as $n \to \infty$,
\[ \sup_{T \in \cT_m}\widetilde{E}_{f_0} \, G\left( f \, : \, \ell(f,f_0) \ge M_n \zeta_{n,m}^2  \, | \, \widetilde{\cD}_n \right) \to 0. \]
To simplify the notation, we often omit the supremum over $\cT_m$,
but the convergence is always understood uniformly over this class.
Section \ref{sec : on the map formulation of adv cont} of the Supplement discusses this uniformity and shows that the above definition is equivalent to maximizing over admissible corruptions separately for each realization of the clean sample (see Lemma \ref{lem : worst case and maps}).

\medskip
\textit{Proof techniques.} The main challenge in dealing with corrupted data is separating the effects of contamination from terms ensuring appropriate posterior contraction. We achieve this by first controlling uniformly over $T$ the contribution of the (at most) $m$ corrupted observations to the generalized likelihood. The resulting bound imposes a contamination penalty that depends on the contamination only through $m$. Once this contribution is isolated, the remaining term only depends on the clean data and can be controlled under $P_{f_0}^n$ using more standard Bayesian nonparametrics arguments. The role of the weight function is then twofold: it sufficiently reduces the contribution of the outliers, while keeping the modified likelihood close to the original one for uncorrupted data, ensuring posterior contraction.

\section{{Output-dependent weights}}\label{sec : w(y)}

In this section we focus on weight functions $w$ depending only on the (contaminated) data $\tilde{y} \in \R$ and derive theoretical guarantees for the generalized posterior under general nonparametric priors. 
A key advantage of such simple weights, already pointed out by \cite{pmlr-v235-altamirano24a}, is that in case of Gaussian Process priors, the robust posterior is conjugate, allowing for easier computation.

We start our analysis with a specific, hard thresholding approach and then extend it to smoother weight functions, covering cases introduced in the machine learning literature. We also apply our general theorem for the popular Gaussian Process and deep neural network priors as specific examples.

\subsection{A simple thresholding weight function}\label{sec : thresh function}
As a first example we consider the weight 
\begin{equation}\label{def : thresh weight}
    w : y \mapsto \frac{\lambda}{1 \vee (y/\tau_n)^2}.
\end{equation} 
The idea behind this ``hard thresholding'' function is to heavily penalize values $y$ that are too large compared to $\tau_n$, while simply tempering the log-likelihood by $\lambda$ whenever $|y| \le \tau_n$. Taking $\tau_n$ as an upper bound for the largest observation coming from the uncontaminated model \eqref{def : model} with overwhelming probability allows us to detect and penalize outliers, while keeping the real data essentially untouched. 
\begin{thm}\label{thm : contconc thresh}
    Let us consider the weight function \eqref{def : thresh weight} with
    $ \tau_n \ge 2\sigma \sqrt{\log n}$ and $\lambda \in (0,1)$, and let $(\varepsilon_n)$ be a positive sequence such that $\varepsilon_n \to 0$ and $ n\varepsilon_n^2 \to \infty$ as $n \to \infty$. Furthermore, assume the prior $\Pi$ satisfies
	\begin{equation}\label{eq : prior mass condition}
	    \Pi[ \{ f \, : \, ||f - f_0||_{\infty} \leq \varepsilon_n\}] \geq e^{-  n \varepsilon_n^2}.
	\end{equation} 
   Then for $\zeta_{n,m}^2 := (m/n) \tau_n^2 + \varepsilon_n^2$, with $m/n \le \eta <1$, and a large enough constant $M>0$,
   \[ \widetilde{E}_{f_0} G\Big( f \, : \, D_{\lambda}(P_f,P_{f_0}) \geq  \frac{M}{1- \lambda}\zeta_{n,m}^2\, | \, \widetilde{\cD}_n\Big) \to 0,\]
    where $D_{\lambda}$ denotes the standard R\'enyi divergence (see e.g. Appendix \ref{app : renyi}). 
    
\end{thm}
The proof is deferred to Section \ref{proof : contconc thresh} in the supplement. Theorem \ref{thm : contconc thresh} derives a contraction rate  for the robust posterior based on the thresholding weight function in \eqref{def : thresh weight} provided the parameter $\tau_n$ is chosen large enough to control the contribution of the outliers, see Remark \ref{rmq : tau_n} below. In this case the robust posterior behaves like a $\lambda$--tempered posterior on the clean data, for which the prior mass condition \eqref{eq : prior mass condition} ensures contraction in terms of the R\'enyi divergence, see e.g.  
\cite{JMLR:aliceIsmaelFractional, temperedPosteriorBattha}. Based on the rate we can distinguish two regimes. In the low contamination regime with  $m\lesssim n\varepsilon_n^2$  (up to a $\tau_n^2$ factor), the generalized posterior contracts with the standard uncontaminated rate $\varepsilon_n^2$. On the other hand, in the high contamination regime $m\gtrsim n\varepsilon_n^2$, the rate is (up to a $\tau_n^2$ multiplicative factor) $m/n$, i.e. the contaminated proportion.

\begin{Remark}[About the choice of threshold $\tau_n$]\label{rmq : tau_n}
The thresholding parameter $\tau_n$ plays a crucial role in detecting outliers. Setting it too high prevents  from down-weighting outliers, which in turn results in a sub-optimal contraction rate in view of Theorem \ref{thm : contconc thresh}. In principle any choice satisfying $\tau_n \geq F_0 + \si \sqrt{c \log n}$, for some $c \geq 2$, is appropriate, but would require the knowledge of $F_0$. Hence we choose a more conservative value, which asymptotically dominates this threshold without increasing the rate. For small sample sizes a data driven tuning can provide a more accurate procedure, see also Section \ref{sec : recenter}.
\end{Remark}

\begin{Remark}
The contraction rate in Theorem \ref{thm : contconc thresh} is stated with respect to the R\'enyi divergence.
However, by assuming or imposing boundedness  on the posterior draws,
the R\'enyi-divergence can be seen to dominate the $L_2$ distance, i.e. $D_\lambda(P_f,P_{f_0}) \gtrsim \|f-f_0\|_{2,P_X}^2$, see Section \ref{app : renyi} for details.
Hence {under such boundedness condition} the same contraction rate holds with respect to the squared $L_2(P_X)$-norm as well, see Section \ref{app : bounded prior}. {Up to a slightly more technical proof, the boundedness assumption can often be removed, see Section \ref{sec : gen weight}}.
\end{Remark}

{In the sequel, to avoid repetitions, a sequence $(\varepsilon_n)$ as in \eqref{eq : prior mass condition} will always be assumed to verify the conditions $\varepsilon_n=o(1)$ and $n\varepsilon_n^2\to\infty$, when $n \to \infty$. }

\subsection{Application to GP and DNN priors}
\label{section applications}
In this section we provide examples of commonly used priors in machine learning and derive small deviation bounds \eqref{eq : prior mass condition} for them. These bounds together with Theorem \ref{thm : contconc thresh} readily provide contraction rates for the associated generalized posterior distributions. Furthermore, the prior small ball probability condition appears in later sections as well when considering various extensions of the weight function in the generalized Bayes method.

\subsubsection{Gaussian processes}\label{subsec : gp example} A GP $W = (W(x) \, : x \in \R^d)$ is a stochastic process whose finite dimensional marginals are multivariate Gaussian distributed. A GP prior is the law of such a process and is completely determined by its mean function $\mu : x \mapsto EW(x)$ and its covariance kernel $k :(x,x') \mapsto \text{Cov}(W(x),W(x'))$. In Gaussian regression \eqref{def : model} the posterior is also a GP and can be expressed explicitly with a mean and covariance functions given by
\begin{equation}
\begin{aligned}
\label{classic posterior mean and covariance}
    x &\mapsto \mu(x) + k(x,X)^\top(K + \si^2 I)^{-1}(\mathbf{y} - \mu(X)), \\
    (x,x') &\mapsto k(x,x') -  k(x,X)^\top(K + \si^2 I)^{-1} k(x',X),
\end{aligned}
\end{equation}
where $\mathbf{y} := (y_1,\dots,y_n)^\top$, $K =[k(x_i,x_j)]_{1 \le i,j \le n}$, $k(x,X) := (k(x,x_1), \dots ,k(x,x_n))^\top$, $\mu(X) = ( \mu(x_1),...,\mu(x_n))^\top$ and $I$ denotes the $n\times n$ identity matrix. For weights of the form $w_i = w(y_i)$, the generalized log--likelihood is
\[ L_{[n],w}(f) = \frac{1}{2 \si^2}\sum_{i=1}^n w(y_i)(y_i-f(x_i))^2 = \frac12 (\mathbf{y} - \mathbf{f})^\top J_w^{-1}(\mathbf{y}-\mathbf{f}),\]
where $\mathbf{f}^\top :=( f(x_1),\dots,f(x_n)) $ and $J_w :=\text{diag}( \si^2w(y_1)^{-1} , \dots , \si^2w(y_n)^{-1})\in\mathbb{R}^{n\times n}$. Then, by standard Gaussian calculations one can show that the generalized posterior $G$, see \eqref{def : genpost}, is also a GP with mean and covariance functions given by
\begin{equation}
\begin{aligned}
\label{posterior mean and covariance gp}
    x &\mapsto \mu(x) + k(x,X)^\top(K +  J_w)^{-1}(\mathbf{y} - \mu(X)), \\
    (x,x') &\mapsto k(x,x') -  k(x,X)^\top(K + J_w)^{-1} k(x',X).
\end{aligned}
\end{equation}

\textit{Mat\'ern covariance kernel.}
We focus on the popular centered Mat\'ern GP priors. A Mat\'ern process $W^{\gamma_n} = (W^{\gamma_n}(x) : x \in \mathbb{R}^d)$ is a centered GP with covariance function
\begin{equation}
\label{eq: matern_rescaled}
    k_{\gamma_n}(x,x') := \int_{\mathbb{R}^d} e^{i s^\top (x - x')}\, \rho_{\alpha,\gamma_n}(s)\, ds,
    \quad 
    \rho_{\alpha,\gamma_n}(s) := C_{\alpha,d}\, \gamma_n^{2\alpha}
    \bigl( c_{\alpha,d}\, \gamma_n^{2} + \| s \|_2^2 \bigr)^{-(\alpha + d/2)},
\end{equation}
where  $\alpha > 0$ and $\gamma_n>0$ denote the regularity and scale tuning parameters, respectively, and $c_{\alpha,d}, C_{\alpha,d} > 0$ are constants. The corresponding generalized posterior is also a GP, with mean and covariance functions given in \eqref{posterior mean and covariance gp}, where $k$ and $K$ are replaced by $k_{\gamma_n}$ and $K_{\gamma_n} := [k_{\gamma_n}(x_i,x_j)]_{1\le i,j\le n}$, respectively.

Let us assume that the true regression function $f_0$ in \eqref{def : model} belongs to a $\beta$--H\"older ball of radius $R>0$, defined as
\begin{equation}\label{def : holder ball}
    \cH^\beta(R) := \bigg\{ f:\R^d \to \R \, , \, \sum_{|\mathbf{k}| \le \lfloor \beta \rfloor }\| \partial^{\mathbf{k}}f\|_\infty + \sum_{|\mathbf{k}| = \lfloor \be \rfloor} \sup_{x \neq x'} \frac{| \partial^{\mathbf{k}}f(x) - \partial^{\mathbf{k}}f(x')|}{\|x-x'\|_\infty^{\beta - \lfloor \beta \rfloor}} \le R \bigg\} ,
\end{equation} 
we also denote $\cH^\be := \cup_{R>0} \cH^\be(R)$ the corresponding space.
For $\be$--Hölder functions the $L_2(P_X)$ minimax rate of estimation is of order
\begin{equation}\label{def : minimax rate}
    \varepsilon_n^* := n^{-\be / \{2 \be +d \}}.
\end{equation}
The proposition below provides the small ball probability condition \eqref{eq : prior mass condition}  with the rate $ \varepsilon_n^*$ for the Mat\'ern process if either the regularity parameter is set $\alpha=\beta$ or in case $\alpha>\beta$ by appropriately rescaling it. The $d=1$ case is given in \cite{adaptationBotond}, while for the closely related Vecchia Mat\'ern processes in \cite{szabo2026vecchia}. For completeness we provide the proof in Section \ref{proof : prior mass GP}.

\begin{prop}\label{prop : prior mass GP}
    Let $f_0 : [0,1]^d \mapsto \R$ belong to $\cH^\be(R)$ for some $\be,R>0$. Consider the rescaled Mat\'ern prior $\Pi^{(\gamma_n)}$ in \eqref{eq: matern_rescaled}, with $\alpha \geq \beta$, $\alpha + d/2 \in \mathbb{N}$ and 
    \begin{equation}
    \label{rescaling}
        \gamma_n = n^{\frac{\alpha - \beta}{ \alpha ( 2 \beta + d)}}.
    \end{equation}
    Then the prior small ball condition \eqref{eq : prior mass condition} holds for  $\veps_n \ge M\veps_n^*$, for some constant $M >0$ and large enough $n$.
\end{prop}

Note that the optimal scale $\gamma_n$ in Proposition~\ref{prop : prior mass GP} depends on the typically unknown smoothness parameter $\beta>0$. Therefore, in practice the scale parameter is often endowed with another layer of prior, i.e. we take a Lebesgue density supported on
$[1,\infty)$ such that,
\begin{equation}\label{eq:hyperprior}
  p_{\alpha,n}(\gamma)\geq c_1\exp\bigl\{-c_2\,n^{\frac{d}{2\alpha+d}}\,\gamma^{\frac{2\alpha d}{2\alpha+d}}\bigr\}, \qquad \ga \ge 1,
\end{equation}
where $c_1, c_2$ are positive constants. This results in a two-level, hierarchical prior 
\begin{align}
\Pi(\cdot) := \int_{1}^{\infty} \Pi^{(\gamma)}(\cdot)\, p_{\alpha,n}(\gamma)\, d\gamma.\label{hierarchical prior}
\end{align}
This hierarchical structure can be represented equivalently as $f \mid \Gamma = \gamma \sim \Pi^{(\gamma)}$, $\Gamma \sim p_{\alpha,n}$. Then, we write $W^{\Gamma} := (W(\Gamma x) : x \in [0,1]^d)$, with $W$ independent of $\Gamma$, for the rescaled Mat\'ern process with random scale $\Gamma$. The following proposition, verified in Section \ref{proof : prior mass hierarchical GP}, shows that this hierarchical prior satisfies the prior mass condition \eqref{eq : prior mass condition} with the sequence $\varepsilon_n^*$.
\begin{prop}\label{prop:prior mass hierarchical GP}
    Let $f_0 \in \mathcal{H}^{\beta}(R)$ for some $\beta, R  > 0$. Consider the hierarchical Mat\'ern prior $\Pi$ in \eqref{hierarchical prior}, with $\alpha \geq \beta$ and $\alpha + d/2 \in \mathbb{N}$. Then condition \eqref{eq : prior mass condition} holds with $\varepsilon_n \ge M \varepsilon_n^*$, for a constant $M > 0$ and large enough $n$.
\end{prop}

\begin{Remark}\label{rem:GP:Sobolev}
In the paper {for simplicity} we consider the prior small ball condition \eqref{eq : prior mass condition} with respect to the supremum norm. We conjecture that for Gaussian Processes the results can be extended for the weaker $L_2$-norm small ball probabilities, i.e.
	\begin{equation}\label{eq : prior mass condition_L2}
	    \Pi[ \{ f \, : \, ||f - f_0||_{2,P_X} \leq \varepsilon_n\}] \geq e^{-  n \varepsilon_n^2}.
	\end{equation} 
This technical extension requires to relate the empirical $L_2$-norm (the intrinsic metric in the regression model) to the population $L^2(P_X)$-norm and to derive supremum norm consistency results for the GP posterior.

Condition \eqref{eq : prior mass condition_L2} is more naturally linked to Sobolev regularity classes.  Let us define the $\be$-Hilbert-Sobolev ball,
    \begin{equation}\label{def : sobo ball}
    \mathcal{S}^\beta(R) := \left\{f : \R^d \mapsto \R \, : \, \int (1 + \|s\|_2^2)^\beta |\hat{f}(s)|^2 \, ds \le R^2\right\},
\end{equation}
where $\hat{f} : s \mapsto \int e^{i s^\top u} f(u)\, du$ is the multivariate Fourier transform of $f$, we also denote the corresponding Sobolev space $\cS^\be := \cup_{R > 0} \cS^\be(R)$. The $L_2(P_X)$--minimax rate of estimation over Sobolev balls is of same order $\veps_n^*$ as in $\eqref{def : minimax rate}$. The lower bound results in Section \ref{sec : lower bound} are derived with respect to this regularity class.

In the literature, condition \eqref{eq : prior mass condition_L2}  with the rate $\veps_n^*$ for $f_0\in \mathcal{S}^\beta(R)$ was derived for a wide range of appropriately tuned GP priors, including (rescaled) infinite series priors (p 337-p 339 of \cite{ghosal_van_der_vaart_2017}) or sieve priors (Lemma 3.1-3.3 of \cite{rousseau:szabo:2017}). The small ball probability condition was also derived for hierarchical GP priors, where either the scale, the regularity or the truncation parameters were endowed with a hyper-prior distribution, see for instance \cite{arbel2013bayesian}.
\end{Remark}

\subsubsection{Deep ReLU neural networks priors}
 Let $\rho : x \mapsto \max(0,x)$ be the ReLU activation function, applied coordinate-wise to vectors. For integers $N \ge 1$ and $L\ge 2$, consider a fully-connected feedforward neural network of width $N$ and depth $L$, represented by  
\begin{equation}\label{def : dnn}
    f : x \mapsto A_L \circ \rho \circ A_{L-1} \circ \rho \circ \cdots \circ \rho \circ A_0(x), \qquad x\in[0,1]^d,
\end{equation}
where $A_\ell : x \mapsto W_\ell x + b_\ell$ are affine transformations with
\[
    W_0 \in \R^{N \times d}, \qquad
    W_\ell \in \R^{N \times N} \ \ (1 \le \ell \le L-1), \qquad
    W_L \in \R^{1 \times N},
\]
\[
    b_\ell \in \R^{N} \ \ (0 \le \ell \le L-1), \qquad b_L \in \R.
\]
We endow the entries of the matrices and bias vectors with Gaussian priors,
\begin{equation}\label{def:prior:DNN}
W_\ell^{(i,j)}
\overset{\mathrm{i.i.d.}}{\sim}
\mathcal N(0,\sigma_W^2),
\qquad
b_\ell^{(k)}
\overset{\mathrm{i.i.d.}}{\sim}
\mathcal N(0,\sigma_b^2),
\end{equation}
where $\sigma_W,\sigma_b>0$ are fixed constants. This induces a prior $\Pi$ on functions of the form \eqref{def : dnn}. The training of the neural network is then performed by updating $\Pi$ using the reweighted likelihood and computing the generalized posterior in \eqref{def : genpost}. Unlike Gaussian process priors, neural network priors are not conjugate, and posterior computation generally requires approximation methods, such as variational inference, we refer to \cite{arbel2026primer} for additional discussion on the computational aspects.

\begin{prop}\label{prop : prior mass DNN}
    Let $f_0:[0,1]^d \to \R$ 
    be in $ \cH^\be(R)$ for some $\be,R >0$. For any $\delta >0$, set
    \[ N := \big\lceil n^{\frac{d}{2(2\be +d)}} \big\rceil \qquad \text{and} \qquad L:= \lceil \log^{1 + \delta}n \rceil .\]
   Let us consider the induced neural network prior with independent Gaussian parameters, see \eqref{def:prior:DNN}, where $\sigma_W,\sigma_b>0$ do not depend on $n$. Then, for large enough $n$,  \eqref{eq : prior mass condition} holds for any $\veps_n \ge \veps_n^* \log^{2+ \delta} n$, with $\veps_n^*$ denoting the minimax rate given in \eqref{def : minimax rate}.
\end{prop}
 The proof is deferred to Section \ref{proof : prior mass DNN}. The proposition establishes the prior small ball condition with the minimax rate corresponding to $\mathcal{H}^\be$. Note that the width of the network depends on the true regularity parameter $\beta$. In case this information is not available, one can instead consider an overparameterized architecture with $N \asymp \sqrt{n}$ and heavier--tailed coefficients, see \cite{JMLR:v26:24-0894}.

\begin{Remark}\label{rmk : Sobolev space}
In the case $\be > d/2$, the overparameterized DNN prior with heavy-tailed coefficients of \cite{JMLR:v26:24-0894} also satisfies the prior mass condition \eqref{eq : prior mass condition} with $\veps_n^*$ for {Sobolev truth} $f_0 \in \cS^\be$ given in \eqref{def : sobo ball} (see a slight modification of Theorem 10 therein with $p =2$).
\end{Remark}

\subsection{More general weight functions}\label{sec : gen weight}
In this section we {study a first extension of} the class of weight functions considered in the generalized posterior. We  assume 
\begin{equation}
\label{eq : cond w cont}
    0 \le w(y) \le \lambda \qquad \text{and} \qquad y^2w(y) \leq \lambda \tau_n^2\qquad \forall y\in \R,
\end{equation}
where $\tau_n >0$ is an $n$-dependent thresholding parameter and  $\lambda\in(0,1)$ is a scale parameter. 
Recalling that $F_0$ is the bound on $f_0$ in \eqref{eq : bound on f_0}, suppose also that there exists $w_0 > 0$, such that 
\begin{equation}\label{eq : cond lowerb w}
    \inf_{ |u| \leq F_0} w(u)(1-w(u)) \ge w_0(1-w_0),
\end{equation}
 and that there exist sequences $\rho_n,\Delta_n >0$ such that, for any $s \in [-F_0,F_0],$
\begin{equation} \label{conddiffw}
\sup_{|u| \leq \rho_n} |w(s+u)-w(s)| \le \Delta_n.
\end{equation}
Condition \eqref{eq : cond w cont} will be used to control the contribution of the outliers in the reweighted loglikelihood. The two other conditions are used to show contraction of the generalized posterior on clean data. Condition \eqref{eq : cond lowerb w} is mild and ensures that $w$ is uniformly bounded away from $0$ or $1$ on the segment $[-F_0,F_0]$. The second condition \eqref{conddiffw} controls the variations of $w$ in the range of $f_0$. 
Note that the hard thresholding weight \eqref{def : thresh weight} satisfies the above conditions.

Next, we provide contraction rates for the generalized posteriors under this more general set of conditions, with respect to the empirical $L_2$ norm of the inliers $x_i\in\cO^c$; see \eqref{def:outliers}. The proof can be found in Section \ref{proof : genweight fix}.

\begin{thm}\label{thm : genweight fix}
   Suppose the prior $\Pi$ satisfies the prior mass condition \eqref{eq : prior mass condition} and the weight $w$ satisfies $\eqref{eq : cond w cont}$ with $\tau_n > 2 \si \sqrt{\log n}$, \eqref{eq : cond lowerb w} and \eqref{conddiffw} with $\rho_n \ge 2 \sigma \sqrt{\log n}$  and $\Delta_n =o(1)$.
    
    Then, for the sequence $\zeta_{n,m}^2=  (m/n)\tau_n^2 + \varepsilon_n^2+ (\Delta_n \rho_n)^2$, 
    \[ \widetilde{E}_{f_0} G\left (f \, : \, \frac{1}{n-m} \sum_{i \notin \cO}|f(x_i) - f_0(x_i)|^2 \ge M \zeta_{n,m}^2 \,\big|\, \widetilde{\cD}_n \right ) {\to} 0,\]
    for $M >0$ large enough, as $n\rightarrow\infty$.
\end{thm}

The above contraction rate result holds for general priors satisfying the prior small ball probability condition, see Section \ref{section applications} for examples. The rate consists of three terms. The first term is equal to the contamination proportion $m/n$ up to a $O(\tau_n^2)$ multiplier, while the second term is the classical (fractional) Bayes contraction rate based on the clean data, following from the prior small ball probability condition. The third term describes the bias imposed by the weight function. We consider specific choices of the weight function below and derive upper bounds for the corresponding bias terms. We show that for appropriately chosen weights this bias is negligible compared to the other two factors. Finally, the rates are given with respect to the empirical $L_2$-norm on the inliers. One can extend this result by proving contraction with respect to the population $L_2(P_X)$-norm in case of bounded posteriors, see Section \ref{app : bounded prior} for details.

Below we focus on Gaussian Process priors and show that without restricting the posterior samples to bounded functions, the contraction can be extended to the more standard and practically more relevant $L_2(P_X)$-norm. The proof is deferred to Section \ref{proof : thml2}.

\begin{thm} \label{thml2}
Let $f_0 \in \cH^\beta(R)$ and assume $P_X$ has a density supported on $[0,1]^d$ and bounded below by $c >0$. Consider a Gaussian process prior putting mass one on $\al$--H\"older functions for some $\al>d/2$ and satisfying condition \eqref{eq : prior mass condition}.
Suppose furthermore that the weight function $w$ satisfies $\eqref{eq : cond w cont}$ with $\tau_n > 2 \si \sqrt{\log n}$, \eqref{eq : cond lowerb w}, \eqref{conddiffw} with $\rho_n \ge 2 \sigma \sqrt{\log n}$ and $\Delta_n = o(1)$. Setting $\zeta_{n,m}^2 = (m/n)\tau_n^2 + \varepsilon_n^2 + (\Delta_n\rho_n)^2$, further assume that, for some $\vartheta < \alpha$,
\begin{equation} \label{extram}
n \zeta_{n,m}^2 \le n^{2\vartheta/(2\vartheta+d)}. 
\end{equation} 
Then for $M>0$ sufficiently large, as $n\to\infty$,
\[ \widetilde{E}_{f_0} G\left[f:\ \|f-f_0\|_{2,P_X}^2 \ge M\zeta_{n,m}^2 \given \widetilde{\cD}_n \right] \to 0. \]
\end{thm}

We note that condition \eqref{extram} is purely technical and is typically satisfied for smooth enough 
{priors }or alternatively not overly contaminated settings. This condition also appears in  Section \ref{sec : recenter}, and we provide a more detailed discussion after Proposition \ref{prop : GP prior mass and sieve}.

Next we discuss weights satisfying the conditions of Theorems \ref{thm : genweight fix} and \ref{thml2}. Given $\tau = \tau_n > 0$, $\lambda \in (0,1)$ and a non-negative scaling function $h : \R \to \R_{+}$ consider functions $w$ of the form
\begin{equation}\label{def : expl weigh}
    w : y \mapsto \frac{\lambda}{1 + \frac{y^2}{\tau^2}h(y)}.
\end{equation} 
Any such function readily satisfies condition \eqref{eq : cond lowerb w}. Additionally, these functions also verify condition \eqref{eq : cond w cont}, provided $h \ge 0$ and
\begin{equation}\label{eq : cond on scaling function h}
    h(y) \ge 1 - (\tau/y)^2.
\end{equation}
A first natural choice to satisfy \eqref{eq : cond on scaling function h} is $h = 1$, which corresponds to the inverse multiquadratic (IMQ) weight introduced in \cite{pmlr-v235-altamirano24a}. A second  very natural choice 
is $h(y) = \exp\{- \tau^2/y^2\}$. We also note that the condition holds for the hard thresholding weight, given in \eqref{def : thresh weight}, as well. Figure \ref{fig : weight functions} below illustrates how these different weight functions behave around the origin.

The following proposition checks the remaining condition \eqref{conddiffw} for these choices and computes the associated rate $\zeta_{n,m}^2$. The proof can be found in  Section \ref{proof : expl weight} of the supplement.

\begin{prop}\label{prop : expl weight}
    Consider the weight function $w$ in \eqref{def : expl weigh} with scaling functions a) $h(y) = \exp\{ - \tau^2 / y^2 \}$ and b) $h\equiv 1$. The conditions of Theorem \ref{thm : genweight fix} and \ref{thml2} are satisfied with
     \begin{enumerate}
    \item[a)]  $\rho_n = 2\si\sqrt{\log n}$ and $\tau_n = 2 \si\log^{1/2 +\delta} n$ for some $\delta > 1$. One can choose $\Delta_n$ such that $\rho_n \Delta_n\leq \varepsilon_n$ for any $\varepsilon_n \to 0$ with $n \varepsilon_n^2 \to \infty$ and thus 
    \[\zeta_{n,m}^2 =(m/n)\log^{1 + 2\delta}n + \varepsilon_n^2.\]
    \item[b)] $\rho_n = 2\si\sqrt{\log n}$ and for large enough $n$, one can choose $\Delta_n \lesssim (\rho_n/\tau_n)^2$ so that
    \[\zeta_{n,m}^2 = (m/n)\tau_n^2 + \varepsilon_n^2 + \rho_n^6/\tau_n^4,\] 
    a bound optimised for the $m$--dependent choice    $\tau_n = \rho_n ({n/m})^{1/6}$.
    \end{enumerate}
    The thresholding weight function in \eqref{def : thresh weight} also satisfies these conditions with $\rho_n = 2\si\sqrt{\log n}$ and $\tau_n =3 \si \sqrt{ \log n}$. For $n$ large enough $\Delta_n =0$ and one recovers the rate  as in Theorem \ref{thm : contconc thresh}.
\end{prop}
Note that for the thresholding weight \eqref{def : thresh weight} and the exponential weight $h= \exp(- \tau^2/y^2)$ the optimal choice of the threshold $\tau_n$ does not depend on $m$ and hence the rate $ (m/n)\, \tau_n^2 +\varepsilon_n^2$ is attained for all contamination levels. In the case $h\equiv 1$ the optimal choice of $\tau_n$ to obtain the same rate depends on $m$ in the above upper bound.

\begin{figure}[H]
    \centering
    \includegraphics[width=0.5\linewidth]{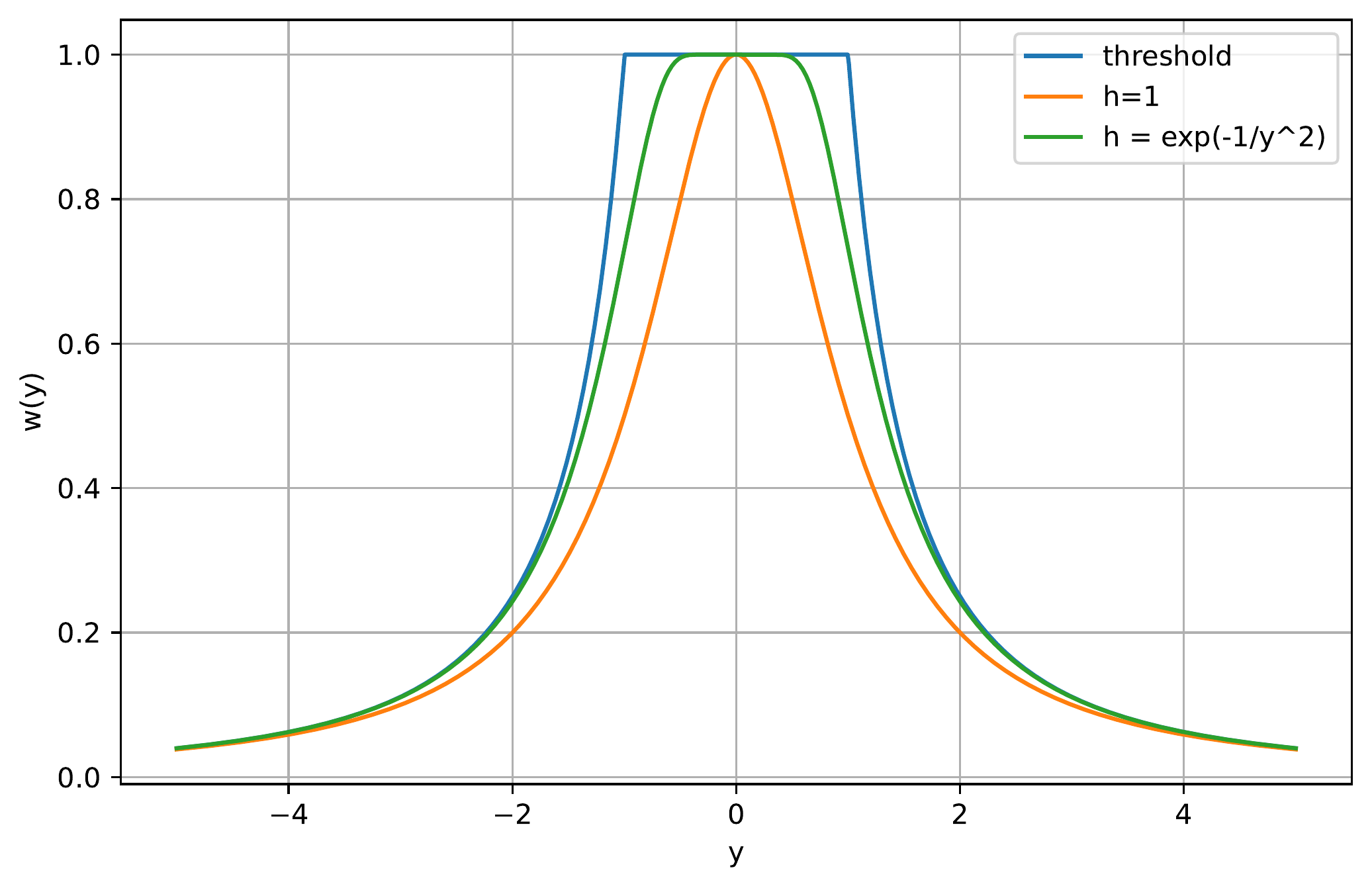}
    \caption{Example of weight functions as in \eqref{def : expl weigh} $\lambda = 1$, $\tau =1$.}
    \label{fig : weight functions}
\end{figure}

\section{{Debiased generalised likelihood and recentered weights}}
\label{sec : more complex modifications}

In this section, we consider two families of more complex robust generalized likelihoods. The first combines the weighting with a well--chosen shifted version of the original likelihood and corresponds to the approach in \cite{pmlr-v235-altamirano24a}; the second lets the weight depend on the output $y_i$, the covariate $x_i$ and possibly the functional parameter $f$.

\subsection{A debiased generalized likelihood}\label{sec : debiased}
The choice of the threshold $\tau_n$ involves a trade-off between robustness and statistical accuracy: a small $\tau_n$ is desirable because it strongly downweights large outliers, but a too small threshold also extensively downweights genuine observations
(inliers) drawn from the model \eqref{def : model}. In this section, we consider a modification of the weighted likelihood that allows a constant threshold $\tau_n = \tau$ independent of the sample size $n$. We write the generalized posterior in the form
\[
    dG(f \mid \cD_n) \propto \exp \Big\{- \tfrac{1}{2 \sigma^2} \textstyle\sum_{i=1}^n
    \ell_f(x_i,y_i) \Big\}\, d \Pi(f),
\]
for some loss function $\ell_f : \R^d \times \R \to \mathbb{R}_+$. Consider the corresponding risk and its minimizer
\begin{equation}
\label{eq : risk minimizer}
    f^{*} := \argmin_{f \in\cF} R(f,f_0),
    \qquad \mbox{with} \qquad R(f,f_0) := E_{f_0}\big[\ell_f(x,y)\big].
\end{equation}
In Section \ref{sec : w(y)}, the loss was taken as a simple reweighting of the log--likelihood $\ell_f(x,y) = w(y)(y-f(x))^2$. Note that in this case the minimizer $f^*$ of the risk  may be different from $f_0$. To correct this possible bias,
we replace the response $y$ by a transformed version $t(y)$, and consider the loss
\begin{equation}\label{eq : general loss}
    \ell_f(x,y) = w(y)\,\big(t(y) - f(x)\big)^2.
\end{equation}
The following lemma, proven in Section \ref{proof : compute risk minimizer}, identifies a choice of $t$ for which the risk is minimized at $f_0$.

\begin{lem}\label{lem : compute risk minimizer}
    Let $w$ be a differentiable, positive weight function satisfying \eqref{eq : cond w cont}. For $y = f_0(x) + \si \xi$, $\xi \sim \cN(0,1)$, assume that $P_X$-almost surely, one has $E_\xi [w(y)]\ge c>0$ as well as Stein's conditions $E_\xi [w(y)t(y)^2] \le C$ and $E_\xi[|w'(y)|] < \infty$. Further assuming that $\cF$ is a closed convex subset of $L_2(P_X)$, the risk minimizer $f^*$ in \eqref{eq : risk minimizer} with loss $\ell_f$ in \eqref{eq : general loss} satisfies 
    \[ f^*(x) =  \mathfrak{p}_\cF \left[ \frac{E_\xi[w(y)t(y)]}{E_\xi[w(y)]} \right],\]
    where $\mathfrak{p}_\cF$ is the Hilbert projection onto $\cF$ with respect to the $L^2(A_{\cdot} \cdot P_X)$-norm, with $A_x := E_{\xi}[w(f_0(x) + \sigma \xi)]$. 
    In particular, we have $f^* = f_0$ with 
    \begin{equation}\label{score matching alta loss}
    t(y) := y - \sigma^2 \frac{w'(y)}{w(y)}.
\end{equation}
\end{lem}
Let us comment on the case $t(y) = y$ considered in the previous section. A consequence of the above lemma is that, under mild conditions the minimizer   in this case   is 
 \[
    f^*(x) = f_0(x) + \sigma \frac{E_{\xi}[ \xi w(y)]}{E_{\xi}[w(y)]} = f_0(x) + \sigma^2 \frac{E_{\xi}[w'(y)]}{E_{\xi}[w(y)]}.
 \]
 In the previous sections, we had $\lim_{n\rightarrow \infty} w(y)= \lambda$ for $\tau_n \to \infty$. This together with the symmetry of  $\xi$ around zero implies that the bias $f^*(x) - f_0(x)$ tends to zero as $n$ increases. The advantage of the correction \eqref{score matching alta loss} is that this bias vanishes even for constant $\tau$. 
 
To establish posterior contraction we introduce the following conditions on the weight function.
We say that a weight $w$ is \emph{admissible} if it is differentiable, positive,
satisfies \eqref{eq : cond w cont} with a constant threshold $\tau_n = \tau > 0$ and if there exists a constant $\kappa_2 >0$, such that
\begin{equation}
\label{extra conditions on the weight}
    \sup_{y \in \R}\, \sigma^2 \,\frac{|w'(y)|}{w(y)} \le \kappa_2 .
\end{equation}
These conditions are very similar to those made in Section \ref{sec : gen weight}.
Conditions \eqref{eq : cond w cont} and \eqref{extra conditions on the weight} together can be related to our previous condition \eqref{conddiffw} and provide a uniform bound on the bias. Indeed, the mean--value theorem ensures \eqref{conddiffw} for all $\Delta_n$ such that $ \Delta_n \geq \si^{-2}\kappa_2\lambda \rho_n $. In this section we consider the correction \eqref{score matching alta loss} so that we can take $\kappa_2$ constant in $n$. However, note that taking $\kappa_2 = \kappa_2(n)$ vanishing sufficiently fast, one can control the bias term $\Delta_n \rho_n$ appearing in Theorem \ref{thm : genweight fix} (in the spirit of Proposition \ref{prop : expl weight}).

In Lemma \ref{lem : assumptions IMQ} of the supplement we show that the IMQ weight given in \eqref{def : expl weigh} with $h\equiv1$ is an admissible weight. Furthermore, we note that admissible weights satisfy Stein's conditions of Lemma \ref{lem : compute risk minimizer}. More concretely, the inequality $\sup_y w(y)t(y)^2 \le 2\lambda(\tau^2 + \kappa_2^2)$ is verified in the proof of the following theorem, while  $\sup_{y} |w'(y) | \leq \lambda \kappa_2 \sigma^{-2}$ follows immediately from the definition of admissibility.

\begin{thm}
\label{thm : debiased}
 Let $w$ be admissible and $\ell_f$ be the corresponding debiased loss
    \eqref{eq : general loss}--\eqref{score matching alta loss}.
    Suppose $\Pi$ satisfies the prior mass condition \eqref{eq : prior mass condition}. 
    Further assume that there exists $F>F_0$ such that the prior assigns probability one to $\{f \, : \|f \|_\infty \le F \}$. 
    Then, for a small enough constant  $\lambda>0$, any sequence $\zeta_{n,m}^2 \ge m/n + \varepsilon_n^2$ and $M > 0$ sufficiently large, as $n \to \infty$,
    \[
        \widetilde{E}_{f_0}\, G\left(  f : \|f-f_0\|_{2,P_X}^2 \ge M\,\zeta_{n,m}^2 
        \,\big|\, \widetilde{\cD}_n \right) \to 0.
    \]
\end{thm}
 Theorem \ref{thm : debiased}, proved in Section \ref{proof : debiased}, shows that the debiased posterior achieves the $L_2(P_X)$--squared rate $m/n + \varepsilon_n^2$ while allowing $\tau$ to remain fixed. This removes, first, the rate penalty associated with a growing threshold, but also removes the bias term $(\Delta_n \rho_n)^2$ in Theorems \ref{thm : genweight fix} and  \ref{thml2}.

Additionally, the debiasing correction preserves the conjugacy for Gaussian
process priors, retaining the computational tractability noted in Section~\ref{sec : gen
weight} and in \cite{pmlr-v235-altamirano24a}. Indeed, the generalized posterior is also a GP, with the same covariance as in \eqref{posterior mean and covariance gp} and mean function
\[
    x \mapsto \mu(x) + k(x,X)^\top (K + J_{w})^{-1} (t(\mathbf{y}) - \mu(X)), 
\]
with $t(\mathbf{y}) := (t(y_1),\dots,t(y_n))^\top$, and $\mu(X) ,K, k(x,X)$ and $J_w$ defined as in Section~\ref{section applications}.


\begin{Remark}
    In Theorem \ref{thm : debiased} we require that the prior draws are uniformly bounded functions. This can be achieved by truncation with rejection, see Section \ref{app : bounded prior} for a description of the approach. The boundedness and small $\lambda$ conditions are required for controlling the numerator in Bayes' formula using the approach developed in \cite{SyringMartin}. We conjecture that these conditions can be relaxed by a technically more involved argument. In the numerical analysis, we demonstrate this conjecture using synthetic datasets. We also note that, for convenience, we use the boundedness assumption to control the contribution of outliers. Without this condition Lemma \ref{lem : getting rid of the outliers} (with constant $\tau$) can be used instead at the cost of a logarithmic multiplicative factor.
\end{Remark}

\subsection{Recentering the weight function}\label{sec : recenter} 
Up until now the weight function $w$ depended only on the observed value $y$. As noted in Remark \ref{rmq : tau_n}, a sufficient threshold $\tau_n$ to detect outliers is $F_0 + \si \sqrt{c \log n}$ for $c \geq 2$. Since the bound $F_0$ is unknown, in Section \ref{sec : w(y)} we chose the conservative threshold $2 \si \sqrt{\log n}$, which exceeds the previous bound when $n$ is exponentially larger than the signal to noise ratio $(F_0/\si)^2$, a condition that is hard to meet in finite samples. Moreover, the weights that are even functions in $y$, such as the ones in \eqref{def : expl weigh}, will not detect the adversarial corruption $y \mapsto -y$, a flaw that negatively impact the robustness performance of the method as discussed in Section \ref{sec: numerics}.

Both these issues suggest reweighting the likelihood with respect to residuals $y-f_0(x)$, instead of responses $y$. We consider two ways to replace the unknown center $f_0(x)$. For simplicity, we consider here recentering of the thresholding weight function of Section \ref{sec : thresh function}. As the numerical analysis of Section \ref{sec: numerics} suggests, we believe the theoretical study below can also be made for the more general weights we consider in Section \ref{sec : gen weight} and \ref{sec : debiased}.

\medskip
\emph{Recentering at an estimator}. Let $\widehat f_n$ be computed from the contaminated data and consider
\begin{equation}\label{def : weight estimator recentering}
 w(\widehat f_n,x,y) := \frac{\lambda}{1\vee (y-\widehat f_n(x))^2/\tau_n^2},
 \qquad 0<\lambda<1.
\end{equation}
If $\widehat f_n$ is a robust, $L_\infty$--consistent estimator of $f_0$, then for clean observations $y-\widehat f_n(x)=\sigma\xi+o_{P_{f_0}^n}(1)$ uniformly in $x$, allowing $\tau_n$ to be chosen according to the noise scale rather than $F_0$. Since $w$ does not depend on $f$, the reweighted loglikelihood is still quadratic in $f$ and GP conjugacy is preserved. The need to recenter the weight function was already noted in \cite{pmlr-v235-altamirano24a}, where the prior mean of the GP was used. This choice of recentering, however, may be distorted by the outliers, hence our proposal to recenter at a robust point estimator. 
\begin{thm}\label{thm : recenter thresh}
   {Consider the generalized posterior based on the weight function $w$ as in \eqref{def : weight estimator recentering}}. 
    Suppose the estimator $\widehat{f}_n$ satisfies for some constant $C > 0$ and two sequences $\chi_{n,m},\delta_n \to 0$, as $n \to \infty$, 
    \begin{equation}\label{eq : suff cond estimator}
        \widetilde{P}_{f_0}\bigl( \| \hat{f}_n - f_0 \|_{\infty} > C \chi_{n,m} \bigr)
        \leq \delta_{n} .
    \end{equation}
    Suppose $\varepsilon_n$ satisfies the prior mass condition \eqref{eq : prior mass condition}. Let  $\tau_n > C\chi_{n,m} + \sigma\sqrt{3\log n}$ in \eqref{def : weight estimator recentering}. 
    For $\zeta_{n,m}^2 :=( m/n )\tau_n^2 + \varepsilon_n^2$ and a constant large enough $M > 0$, 

    \[ \widetilde{E}_{f_0} G( f \, : \, (1 - \lambda)  D_{\lambda}(P_f,P_{f_0}) \geq  M \zeta_{n,m}^2\, | \, \widetilde{\cD}_n\Big) \to 0, \qquad \text{as } n\to \infty.\]
\end{thm}
The proof is available in Section \ref{proof : recenter thresh}. It requires the sharper thresholding condition $\tau_n > C\chi_{n,m} + \sigma\sqrt{3\log n}$. Since $\chi_{n,m} \to 0$, this condition does not involve $F_0$ directly; moreover, the choice $\tau_n \geq 2\sigma\sqrt{\log n}$ satisfies it as soon as $C\chi_{n,m} < 0.26\,\sigma\sqrt{\log n}$, which typically requires a much smaller $n$ than previously. A robust trimmed--wavelet estimator satisfying the required supremum--norm consistency, under additional smoothness assumptions, is given in Section \ref{sec : sup-norm trimmed}. This frequentist construction requires tuning in terms of $m$. To obtain a $m$--independent and $L_\infty$-consistent estimator of $f_0$, one can use draws from the GP posterior of Theorem \ref{thml2}, which is shown to be $L_\infty$-consistent in Theorem \ref{thm : sup norm consistency} of the Supplement.

\medskip
\emph{Recentering at a candidate function.}
Instead of pre-computing a robust estimator $\widehat{f}_n$, we now consider
\begin{equation}\label{def : cliplik weight}
    w(f,x,y) := \frac{\lambda}{1 \vee ( y-f(x))^2/\tau_n^2} , \qquad 0 < \la < 1.
\end{equation}
Although GP conjugacy is lost with this weight, a preliminary robust estimation is not required for recentering.
To state the contraction result, we additionally require a sufficient control on the regularity of the functions $f$ drawn from the prior.
Let us denote by 
\[
    B_{1}^{\infty} := \{ f \in \mathcal{F}: \| f \|_{\infty } \leq 1\}
\]
 the unit ball in $L_{\infty}([0,1]^d)$.
\begin{thm}[Contraction under a sieve oscillation condition]\label{thm : cliplikconc gen}
    Assume uniform design $P_X = \mathcal{U}([0,1]^d)$ and that $\varepsilon_n$ satisfies  the prior mass condition \eqref{eq : prior mass condition}. Let \eqref{def : cliplik weight} hold with $\tau_n > 8 \si\sqrt{\log n}$. For $\zeta_{n,m}^2 := (m/n)\tau_n^2 + \varepsilon_n^2$, set
    $h_n := (4\zeta_{n,m}^2)^{1/d}$. Suppose, for some $\vartheta > 0$, that there are sets $H_n \subset \mathcal{H}^{\vartheta}(L_{n}) +  B_{1}^{\infty}$ with $L_{n} h_{n}^{\vartheta} = O(1)$ and $\Pi(H_n^c) \le e^{-C n \zeta_{n,m}^2}$, where  $C > 1 + \lambda(1+\lambda)/\si^2$. Then, for $M>0$ large enough, as $n \to \infty$,
    \[
    \widetilde{E}_{f_0} G\Big( f \, : \, (1-\la) D_{\lambda}(P_f,P_{f_0}) \geq  M\zeta_{n,m}^2 \log n\, | \, \widetilde{\cD}_n\Big) \to 0.
    \]
\end{thm}

Next we show that Gaussian Processes with Mat\'ern covariance kernel satisfy the assumptions of the previous theorem and can provide minimax contraction rates assuming that the contamination level is not too large and the regularity $\beta$ is high enough.  

\begin{prop}\label{prop : GP prior mass and sieve}
    For $f_0 \in \cH^\beta$ and $\Pi$ as in Proposition~\ref{prop : prior mass GP} with $\gamma_n = 1$ (that is $\alpha = \beta$), the conditions of Theorem \ref{thm : cliplikconc gen} are satisfied for any $\varepsilon_n \ge M \varepsilon_n^*$ defined in \eqref{def : minimax rate} under the condition:
    \begin{equation}\label{eq : GP modulus condition}
        n\,\ze_{n,m}^{2} \lesssim n^{\frac{2\vartheta}{d+2\vartheta}}, \qquad \text{for some } \vartheta \in (0,\be).
    \end{equation}
\end{prop}

Proofs of Theorem \ref{thm : cliplikconc gen} and Proposition \ref{prop : GP prior mass and sieve} are done in Section \ref{proof : cliplikconc} and \ref{proof : GP prior mass and sieve}, respectively. 
    Taking $\varepsilon_n \asymp \varepsilon_n^* = n^{-\be/(2\be +d)}$ and $\tau_n \asymp \log n$, two
    regimes arise.
    If $m \lesssim n\veps_n^2/\tau_n^2$, the contamination is negligible, the generalized posterior contracts at the uncontaminated rate $\varepsilon_n$. In this
    regime condition \eqref{eq : GP modulus condition} reads
    $n^{d/(2\be+d)}\le n^{2\vartheta/(d+2\vartheta)}$ for some $\vartheta < \beta$. The latter holds if and only if $\beta > d/2$. 
    If instead $n\veps_n^2/\tau_n^2 \lesssim m$, the contamination term dominates, and the contraction rate becomes $\sqrt{(m/n)\log n}$. Condition \eqref{eq : GP modulus condition} now reads $m\lesssim n^{2\be /(2\beta+d) - \kappa} / \log n$ for any $\kappa > 0$, which is the largest number of adversarial contaminations the argument can accommodate, {a mild condition 
for large enough smoothness}.

\section{Minimax rates and adaptation for adversarial contamination}
\label{sec: minimax}

In the previous sections we obtained an upper bound on the contraction rate in \(L_2(P_X)\) norm, and other related metrics, for the proposed generalized posteriors. Recall that, for $f_0\in \mathcal{S}^{\beta}$ with $\|f_0\|_{\infty}\leq F_0$ the rate is, up to logarithmic factors, $\sqrt{m/n} + n^{-\beta/(2\beta+d)}$. Note that for $m \lesssim n^{d/(2\beta+d)}$ the second term dominates the first, and the contraction rate matches the standard, uncontaminated risk. Hence, the generalized Bayes methods studied above provide optimal inference. In the high contamination regime $m \gtrsim n^{d/(2\beta+d)}$ the first term dominates the second, hence the contraction rate is \(\sqrt{m/n}\), which matches the rate of the kernel ridge regression estimator derived by Chinot~\cite{chinot}. The question arises whether this is the optimal rate in this regime.

We provide a negative answer to this question by showing that for given $m, \be$ with $m \gtrsim n^{d/(2\beta+d)}$ the above rate may be suboptimal. We prove that for $f_0\in \mathcal{S}^{\beta}$   
with 
$\|f_0\|_{\infty}\leq F_0$, the minimax rate under  adversarial contamination is,
up to logarithmic factors in \(n\) and \(m\),
\[
    ( m/n )^{
        \frac{2\beta}{2\beta+d} \vee \frac{1}{2}
    }
    + n^{-\frac{\beta}{2\beta+d}} .
\]
This leads to two distinct regimes in this high contamination setting:
\begin{itemize}
\item If \(\beta>d/2\), the minimax contamination term is    $({m}/{n})^{{2\beta}/({2\beta+d})}\ll \sqrt{m/n}$. Hence, in this smoother regime, the rate obtained for the proposed generalized posterior is suboptimal. In particular, the
phase transition in the minimax rate occurs at \(m\gtrsim\sqrt{n}\), in contrast to \(m\gtrsim n^{d/(2\beta+d)}\) appearing in the posterior contraction rate above.

\item If \(\beta\le d/2\), the minimax contamination term is \(\sqrt{m/n}\). Then the minimax rate coincides with the contraction rate obtained for the proposed generalized posteriors, up to logarithmic factors. 
\end{itemize}

\subsection{Lower bound}
\label{sec : lower bound}
 In the literature the lower bound results typically focused on the Huber-contamination model using modulus of continuity techniques, see \cite{chen_robust_cov}. We extend these arguments to the adversarial contamination scenario with fixed number of contaminated data $m$, see Section \ref{sec : proof lower bound} for the exact statement and proof.  This extension to the adversarial setting is of interest on its own right and can be used for other nonparametric models as well. The main contribution is, however, to derive minimax lower bounds for the nonparametric regression model. Previous results in the literature either provided suboptimal lower bounds \cite{chaoHuberRegression} or considered more restrictive regularity classes \cite{zhao2024robust}.

For $M > 0$, we denote by $\mathcal{S}_{\text{per}}^{\beta}(M)$ the (periodic) Sobolev ball on $[0,1]^d$ of regularity $\beta$ and radius $M$, defined as
\begin{equation}
\label{sobolev ball}
\mathcal{S}_{\text{per}}^{\beta}(M)
:=
\left\{
f = \sum_{k \in \mathbb{Z}^d} f_{k}\,\varphi_k :
\sum_{k \in \mathbb{Z}^d}(1 + \|k\|_2)^{2\beta}\, f_k^2 \le M^2
\right\},
\end{equation}
where $\{\varphi_k\}_{k \in \mathbb{Z}^d}$ is the trigonometric (Fourier) basis on $[0,1]^d$, obtained as the tensor product
\begin{equation*}
\varphi_k(x) = \prod_{\ell=1}^d \varphi_{k_\ell}(x_\ell),
\qquad x = (x_1,\dots,x_d) \in [0,1]^d,
\end{equation*}
of the univariate trigonometric basis $\{\varphi_{j}\}_{j \geq 0}$ on $[0,1]$, namely
\begin{equation*}
\varphi_0(x) = 1, \qquad
\varphi_{2j-1}(x) = \sqrt{2}\,\cos(2\pi j x), \qquad
\varphi_{2j}(x) = \sqrt{2}\,\sin(2\pi j x), \qquad j \geq 1.
\end{equation*}

\begin{Remark} 
Let $\mathcal S_{\mathrm{per}}^\beta
:= \bigcup_{M>0}\mathcal S_{\mathrm{per}}^\beta(M)$. We identify $f \in \cS_\mathrm{per}^\be$ with the restriction on $[0,1]^d$ of a function in the Sobolev space $\cS^\be$, given in \eqref{def : sobo ball}. Conversely, if $f$ is supported on a compact subset of $(0,1)^d$ and is the restriction on $[0,1]^d$ of a function in $\cS^\be$, then its periodic extension belongs to $\cS_\mathrm{per}^\be$ (see  Lemma A.3 of \cite{tsybakov2009} for $\beta \in \mathbb{N}$ and $d=1$, and Chapter 4 of \cite{GineNickl2016} for additional discussion in the more general case).
\end{Remark}

\begin{thm}
\label{lower bound nonparametric regression}
Consider the nonparametric regression model with design $x_i \overset{\text{i.i.d.}}{\sim} P_X = \mathcal{U}([0,1]^d)$. Let $M, \beta, F_0 > 0$. We have that
\begin{equation}\label{def:minimax:cont}
    \inf_{\hat{f}_n} \sup_{f_0 \in S_{\text{per}}^{\beta}(M) \; : \; \| f_0 \|_{\infty} \leq F_0} \sup_{T \in \cT_m} \widetilde{E}_{f_0, T} \left [   \| \hat{f}_n - f_0 \|_{L_2(P_X)} \right ] \gtrsim  
       \left (\frac{m}{n} \right)^{\tfrac{2 \beta}{2 \beta +d} \vee \tfrac{1}{2}}  \vee n^{-\tfrac{  \beta}{2 \beta +d}},
 \end{equation}
where $\widetilde{E}_{f_0,T}$ denotes the expectation under $\widetilde{P}_{f_0}=\widetilde{P}_{f_0,T} := (P_{f_0}^n \otimes \nu) \circ T^{-1} $.
\end{thm}

The proof is given in Section~\ref{sec : proof lower bound}. Since the generalized GP approaches considered in the previous sections do not attain this information-theoretic lower bound, we derive a novel generalized GP approach based on trimming that achieves the minimax lower bound stated above. First, we extend the trimming-based frequentist 
estimator of \cite{robust_mean_estimation}, originally derived for the finite dimensional mean estimation problem, to the nonparametric regression setting. Building on this extension, we then propose a generalized GP approach. Both methods rely on knowing $m$. Finally, we investigate the practically relevant scenario in which $m$ is unknown and show that this implies a polynomial penalty relative to the minimax lower bound derived in the theorem above.

\subsection{Upper bound}
\label{sec : upper bound}
We derive a robust method by estimating each (Fourier) coefficient of the regression function $f_0$ using the trimmed-mean
estimator of \cite{robust_mean_estimation}. Note that for uniform design in $[0,1]^d$, each Fourier coefficient of $f_0$ can be written as $f_{0,k} := \langle f_0,\varphi_k\rangle_{L^2}= E_{f_0}[y_i\varphi_k(x_i)]$. Furthermore, for $a \le b$, define the trimming function
\begin{equation}
\label{clipping function}
\phi_{a,b}(x):=(x\vee a)\wedge b.
\end{equation}
Let $\widetilde{a}_k, \widetilde{b}_k$ be the $qn$-th and $(n - nq)$-th order statistics of $\{\widetilde y_i \varphi_k(\widetilde x_i)\}_{i=1}^n$, respectively, with 
\begin{equation}
q = q(m,n,\delta) =  8\frac{m}{n} + 12 \frac{\log(8/\delta)}{n}\label{def:q}
\end{equation}
for $\delta \in (0,1)$. We then define the \textit{trimmed projection estimator} as follows 
\begin{equation}
\label{trimmed estimator aggregated}
    \hat f_n:= \sum_{\|k\|_2 \le K_{n,m}} \hat f_{n,k}\, \varphi_k, \qquad
    \hat f_{n,k} := \frac{1}{n} \sum_{i=1}^n \phi_{\widetilde a_k, \widetilde b_k}\big( \widetilde y_i \varphi_k(\widetilde x_i) \big),
\end{equation}
for some integer $K_{n,m}$ to be specified below. We note, that unlike \cite{robust_mean_estimation}, the proposed estimator does not involve a sample splitting.

\begin{prop}
\label{trimmed projection rate}
Let  $B_{n,m} = (m/n)^2 \log( n/m) \vee \log (n)/n$, $K_{n,m} = \lfloor B_{n,m}^{-1/(2\beta + d)} \rfloor$ and $\delta = \delta(n,m) = 8 \exp(-3 n B_{n,m})$. Assume uniform design in $[0,1]^d$, $\beta > d/2$ and $q \leq 1/4$ given in \eqref{def:q}. Then there exists $C > 0$ such that, under adversarial contamination and uniformly over $f_0 \in \mathcal{S}_{\text{per}}^\beta(M)$, with probability at least $1 - \exp( - 2 n B_{n,m} )$, for $n$ large enough,
\[
    \|\hat f_n - f_0\|_{L^2}^2 \le C\left[ \Big( \frac{m}{n} \sqrt{\log \tfrac{n}{m}} \Big)^{\frac{4\beta}{2\beta + d}} + \Big( \frac{n}{\log n} \Big)^{-\frac{2\beta}{2\beta + d}} \right].
\]
\end{prop}

The proof is given in Section~\ref{trimmed projection rate proof}. For $\beta > d/2$, Proposition~\ref{trimmed projection rate} shows that the trimmed
projection estimator attains the lower bound derived in Theorem~\ref{lower bound nonparametric regression}
up to logarithmic factors. Unlike the generalized Bayesian approach, however, it requires knowledge of both $m$ and $\beta$. Next we derive its generalized Bayesian counterpart.

\subsection{Generalized Bayes trimmed projection posterior}
\label{sec : generalized bayes trimmed}

We now construct a conjugate generalized Bayes procedure attaining the minimax rate up to logarithmic factors when $m$ and $\beta$ are known. 
For each $k$ with $\| k \|_2 \leq n$, compute the trimmed coefficient $\hat f_{n,k}$ in \eqref{trimmed estimator aggregated} with confidence level $\delta = n^{-(d +2)}$. Then, let $\Pi$ denote the law of the Gaussian process $f = \sum_{\|k \|_2 \leq n}  \sigma_k Z_k \varphi_k$ with $\sigma_k = (1 \vee \|k \|_2)^{ - \alpha - d/2}$ for $\alpha > 0$ and $Z_k \overset{\text{i.i.d.}}{\sim} \mathcal{N}(0,1)$. Let $B_{n,m}$ be as in the statement of Proposition \ref{trimmed projection rate} and define the learning rate $\lambda_{n,m} := 1/B_{n,m}$. We define the \textit{Generalized Bayes Trimmed Projection Posterior} as follows:
\begin{equation}
\label{bayes trimmed}
    dG(f \mid \widetilde{\cD}_n) : \propto \exp\left\{
-\frac{\lambda_{n,m}}2
\sum_{\|k\|_2\le n}(\langle f , \varphi_k \rangle_2-\hat f_{n,k})^2
\right\} d \Pi(f).
\end{equation}

The construction of the above generalized GP posterior is motivated by the asymptotic equivalence of the nonparametric regression and Gaussian white noise (GWN) models \cite{reiss}. The latter can be equivalently represented with the help of the noisy Fourier coefficients $z_k=f_k+n^{-1/2}\gamma_k$, $k=1,...,$ where $f_k=\langle f,\varphi_k\rangle_2$ and $\gamma_k\stackrel{\text{i.i.d.}}{\sim}\mathcal{N}(0,1)$. The contribution of the $k$th term of the sequence to the likelihood is proportional to $\exp\{-n(z_k-f_k)^2/2\}$. We then replace the noisy Fourier coefficient $z_k$ in the GWN model with the trimmed-estimator $\hat{f}_{n,k}$ of the Fourier coefficient in the regression model.
Moreover, the squared robust estimation error bound in Corollary \ref{corollary mendelson} (Section \ref{sec : proofstrimm}) suggests replacing the variance $1/n$ by $1/\lambda_{n,m}$, with
\[
\lambda_{n,m}\asymp
\frac{1}{(m/n)^2\log(n/m)\vee(\log n)/n},
\]
which takes into account the error due to contamination as well as sampling variability.

We note that, by conjugacy and the simultaneous factorization of the prior and loss, the generalized posterior admits the closed-form expression 
\begin{equation}\label{posterior mean and variance}
    f_k \mid \widetilde{\cD}_n \stackrel{\text{ind}}{\sim} \mathcal{N}(\bar{f}_{n,k}, \sigma^2_{k,n})\qquad\text{with}\quad 
    \bar{f}_{n,k} := \frac{ \lambda_{n,m} \hat{f}_{n,k}}{\lambda_{n,m}+ \sigma_k^{-2}}, \quad \sigma^2_{k,n} := \frac{1}{\lambda_{n,m} + \sigma_k^{-2}},
\end{equation}
for $ \| k \|_2 \leq n$. Note that at the function-space level $f \mid \widetilde{\cD}_n$ is a Gaussian process with posterior mean and covariance functions given by
\begin{align*}
x  \mapsto \sum_{\|k\|_2\le n} \frac{\lambda_{n,m}\,\hat f_{n,k}}{\lambda_{n,m}+\sigma_k^{-2}}\,\varphi_k(x),\qquad
(x,x') \mapsto \sum_{\|k\|_2\le n} \frac{\varphi_k(x)\varphi_k(x') }{\lambda_{n,m}+\sigma_k^{-2}},\qquad x,x'\in[0,1]^d.
\end{align*}
The next theorem derives a contraction rate for this procedure.
 
\begin{thm}
\label{thm:contraction adversarial trimmed projection}
Let $f_0\in\mathcal S_{\text{per}}^\beta(M)$ with $\beta>d/2$, $x_i\stackrel{\text{i.i.d.}}{\sim}\mathcal{U}([0,1]^d)$, $i=1,..,n$, $\sigma^2 = 1$ and $m = o(n)$.
Then, under adversarial contamination, for every sequence
$M_n\to\infty$,
\[
\widetilde E_{f_0}
G\!\left(
f:\|f-f_0\|_{L^2}> M_n \Big(  \frac{m^2}{n^2} \log \left (\frac{n}{m}\right) \vee \frac{\log n}{n}\Big)^{\frac{\beta \wedge \alpha}{2\alpha+d}}
\,\middle|\,\widetilde{\mathcal D}_n
\right)\to 0,  \quad \text{ as } n \to \infty.
\]
\end{thm}

The proof is given in
Section~\ref{proof : contraction adversarial trimmed projection}.
For the oracle choice $\alpha = \beta$, the contraction rate in the previous theorem matches the minimax lower bound of
Theorem~\ref{lower bound nonparametric regression}
up to logarithmic factors.

\subsection{Impossibility of $(m,\be)$--adaptation}\label{sec : noadapt}
The estimators described above explicitly use the contamination ratio $m/n$ and the regularity $\beta$ of the regression function $f_0$; however, this information is often unavailable in practice. This raises the question of whether an adaptive procedure can achieve the minimax risk \eqref{def:minimax:cont} without knowing the regularity and contamination level. We give a negative answer to this question. In fact, we show that one necessarily loses a polynomial factor in the $\|\cdot \|_{2,P_X}$ estimation rate when simultaneously adapting to $\beta$ and $m$.

  In our analysis, for simplicity, we focus on the $\be\ge d/2$ regime, as the $\beta<d/2$ case is already covered by the generalized GP procedure.  In particular, the result exhibits pairs $(m_1,\be_1), (m_2,\be_2)$ such that the loss incurred by {\em any} procedure with respect to the non-adaptive minimax rate $\veps_n(\be,m)$ given in \eqref{def:minimax:cont} is at least polynomial in $n$. 
\begin{thm} \label{thmnoadapt}
Suppose the density of the design points is bounded away from $0$ and infinity. Let $M_1\ge \mu_1$ and $M_2\ge \mu_2$, where $\mu_1,\mu_2$ are universal constants depending only on the dimension $d$ and the bounds on the density of the design points.  Set, for any given $\be\in(d/2,\infty)$ and $m=o(n)$,
\[ (\be_1,m_1)=(\be,m),\qquad (\be_2,m_2)=(d/2,0). \]
Then there exists a constant $c\in(0,\infty)$ such that
\[\rho_{n,m}:=\inf_{\hat{f}} \max_{i=1,2} \sup_{f_i\in \mathcal{S}_{\text{per}}^{\be_i}(M_i)} \sup_{T_i\in \cT_{m_i}} 
\frac{\widetilde{E}_{f_i,T_i} \| \hat{f}-f_i\|_{2,P_X}^2}{\veps_n(\be_i,m_i)^2}
 \ge  c\cdot \left\{\frac{\sqrt{m/n}}{\veps_n(\be,m)\vee \veps_n(d/2,0)}\right\}^2,
 \]
and the lower bound above is attained for functions $f_{(1)}\in \mathcal{S}_{\text{per}}^{\be_1}(M_1),f_{(2)}\in \mathcal{S}_{\text{per}}^{\be_2}(M_2)$ and adversarial maps $T_{(1)}\in \cT_{m_1}, T_{(2)}\in \cT_{m_2}$. 
\end{thm}

The loss for adaptation is made more explicit in the following corollary.
\begin{Corollary} \label{corimpo}
(i) Let $m=n^r$ for some $r\in (0,1)$. Then the loss for adaptation $\rho_{n,m}$ as in Theorem \ref{thmnoadapt} verifies
\[ \rho_{n,m} \ge n^a,\ \text{with }\quad a=(1-r)\left\{\frac{2\be-d}{2\be+d}\right\} \wedge (r-1/2). \]
In particular, a polynomial loss is incurred for adaptation as soon as $r>1/2$.

(ii) Suppose the pair $(\be,m)$ verifies $(m/n)^{2\be/(2\be+d)}=n^{-1/4}$, then
\[ 
\inf_{\hat{f}} \left\{\,
\widetilde{E}_{f_{(1)},T_{(1)}} \| \hat{f}-f_{(1)}\|_{2,P_X}^2\ \vee\  
\widetilde{E}_{f_{(2)},T_{(2)}} \| \hat{f}-f_{(2)}\|_{2,P_X}^2 \,\right\}
 \ge  c\cdot \left\{\sqrt{m/n}\right\}^2,
\]
where the functions $f_{(1)},f_{(2)}$ and adversarial maps $T_{(1)}, T_{(2)}$ are as in Theorem \ref{thmnoadapt}.
\end{Corollary} 
From the minimax point of view, Corollary \ref{corimpo}, (ii), says that for any $m=n^{r}$ for $1/2<r<3/4$, one can find a regularity $\be>d/2$ such that the (adaptive) minimax rate under contamination must be at least $(\sqrt{m/n})^2$. Methods from the first part  of the paper do not depend on $m$ and achieve this rate. In particular the hierarchical Bayesian method also adapts to the regularity, see Remark \ref{rem:GP:Sobolev}. Hence from the minimax point of view, the performance of such a procedure is minimax optimal if one requires adaptation simultaneously over large enough `rectangles' in terms of coordinates $(\be,\log m/ \log n)$: that is, the (adaptive) GP procedure from the first part of the paper is optimal if one requires simultaneous adaptation for (at least) all  $\be\in[d/2+\epsilon,\beta_{max}]$ and all $m$ with $\log m/ \log n \in [1/2+ \epsilon, 3/4 - \epsilon]$, for large enough $\be_{max}$ and small $\epsilon>0$. The proof of Theorem \ref{thmnoadapt} can be found in Section 
\ref{sec : proofnoadapt}.

\section{Numerical analysis}
\label{sec: numerics}
In this section, we compare the numerical performance of the different generalized Bayesian methods considered above over synthetic data sets, illustrating our theoretical findings. We note that \cite{pmlr-v235-altamirano24a} provide an extensive numerical comparison of the weighted approach with the IMQ weight and the debiased loss using GP priors, also covered by our theoretical and numerical analysis, against several competing robust methods, both on synthetic and real world data sets. Therefore, for comparison with other robust methods we refer the reader to this paper. While our theory applies to general priors, we restrict our experiments to Gaussian process priors, which preserve conjugacy. Since we focus on conjugate methods, we omit from these comparisons the weight function introduced in \eqref{def : cliplik weight}. 

\begin{figure}[H]
    \centering
    \includegraphics[width=\linewidth]{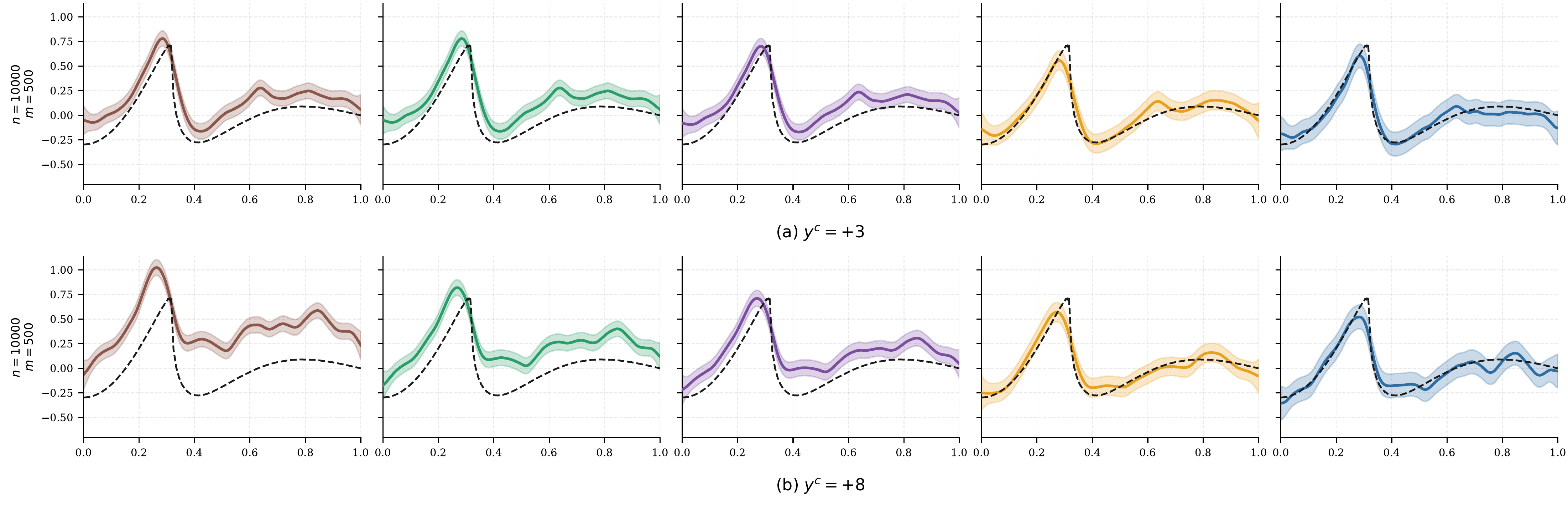}
    \caption{Scenario~1: random contamination with $n=10\,000$ and
    $\eta=1/20$. Top: $y^c=3$; bottom: $y^c=8$.
    Methods: I (brown), II (green), III (purple), IV (yellow), and V (blue).
    Black dashed curves show $f_0$. }
    \label{fig:diffused main}
\end{figure}

We consider the regression model \eqref{def : model} with uniform design on the unit interval, $x_i \overset{\mathrm{iid}}{\sim} \mathcal{U}([0,1])$, we fix $\sigma = 1$, and we take the underlying true regression function to be
\begin{equation}
\label{eq:true-function}
    f_0(x) = \sum_{j \geq 1} f_{0,j}\, \sqrt{2}\, \cos\!\left( \pi \left( j - \tfrac{1}{2} \right) x \right),
\end{equation}
with coefficients $f_{0,j} = 1.5\, \sin(j)\, j^{-1/2 - \beta}$, where $\beta = 1$, for $j \geq 4$, and $f_{0,j} = 0$ for $j \leq 3$. Note that $f_0\in \mathcal{S}^{\gamma}(M)$ for any $\gamma < \beta$ {and that $f_0$ has a localised `spike', see Figure \ref{fig:diffused main}}.

 We aim to assess how the various methods perform for finite sample sizes at different contamination levels including the `clean data' case. Hence, after generating the synthetic data, we contaminate $m = \eta n$ observations, with $\eta \in \{0, 1/40, 1/20\}$, the first choice corresponding to the uncontaminated case. We then consider the following three contamination scenarios:
\begin{enumerate}
    \item \emph{Random contamination:} We select $m$ indices at random and set the corresponding responses $\tilde{y}_i=y_i + y^c$, for $y^c \in \{3, 8\}$.
    \item \emph{Localized contamination:} We select the $m$ covariates closest to the spike of $f_0$ (where the function attains its maximum) and set the corresponding responses equal to $y_i - y^c$.
    \item \emph{Offset and random contamination:} We shift the above function by a constant offset $V=2$, so that the {true} regression  function becomes $V + f_0$. The magnitude of the response is then of order $V$, and we introduce $m$ random contaminations, as in scenario~1, all exactly equal to $-V$. In this way, those weight functions that are even function in $y$, cannot distinguish between genuine and contaminated responses, treating the latter as inliers.
\end{enumerate}
Random contamination resembles Huber-type contamination, whereas localized contamination results in more challenging scenarios since the spike of the signal is hidden by the adversary corruption. The two contamination magnitudes $y^c\in\{3,8\}$ are, respectively, slightly and substantially larger than the typical observed responses, making the latter easier to distinguish from inliers. The third setting is the most challenging one, see also the discussion in Section \ref{sec : recenter}.

In our numerical analysis we compare the following methods:
\begin{enumerate}
    \item[I] \emph{Standard GP}: the standard GP posterior as benchmark;
    \item[II] \emph{Hard threshold}: the generalized posterior \eqref{def : genpost} with thresholding weight given in \eqref{def : thresh weight}, and parameters $\tau_n = 2\sigma\sqrt{\log n}$ and $\lambda = 1$;
    \item[III] \emph{IMQ without debiasing}: the generalized posterior \eqref{def : genpost} with the IMQ weight \eqref{def : expl weigh} ($h = 1$), and parameters $\tau_n = 2\sigma\sqrt{\log n}$ and $\lambda = 1$;
    \item[IV] \emph{IMQ with debiasing}: the debiased generalized likelihood given in \eqref{eq : general loss} and  \eqref{score matching alta loss}  with the IMQ weight and a fixed threshold $\tau$ equal to the empirical quantile of order $0.7$ of $\{|y_i|\}_{i=1}^n$ (the choice of $\tau$ proposed in \cite{pmlr-v235-altamirano24a}).
    \item[V] \emph{Bayes-trimmed}: the generalized Bayes trimmed projection posterior introduced in Section \ref{sec : generalized bayes trimmed}, with trimming level $q = 3.5\,\eta + 5\,\tfrac{\log n}{n}$ and learning rate $\lambda_{n,m} = n$. 
\end{enumerate}

\begin{Remark}
Calibrating the learning rate in generalized Bayes is a well-known practical problem in the literature; see \cite{learningRate, safeBayesian}. In Method V we deviate from the parametrization used for the asymptotic theoretical analysis in Section \ref{sec : generalized bayes trimmed} to improve the practical performance for small to moderate sample sizes.
\end{Remark}

In Methods I-IV we use a Mat\'ern GP prior with regularity parameter $\alpha = 3/2$ and scale parameter $\gamma_n=n^{\frac{\alpha-\beta}{\alpha(2\beta+1)}}$, see \eqref{rescaling}, to compensate for the mismatch between the regularities of the prior and the true regression function $f_0$. To make the comparison fair, we set the regularity $\alpha = 3/2$ in Method V and rescale the prior optimally, replacing $\sigma_k^2$ defined in Section \eqref{sec : generalized bayes trimmed} by $\sigma_k^2 = t_n^2\, k^{-2\alpha - 1}$, with $t_n = n^{\frac{\alpha - \beta}{ 2\beta+1}}$. This optimal rescaling $t_n$ is derived for instance in \cite{knapik2011bayesian}. 

For every simulation scenario we run $100$ independent replications at sample size $n = 5000$, and report the average integrated squared errors
$\|\hat f_n - f_0\|_{2,P_X}^2$, with $\hat{f}_n$ the (generalized) posterior means, together with their standard deviations. To visualize the behaviour of the procedures, we also display, for a representative run, the posterior mean of each method (solid lines in the plots) together with the corresponding $95\%$ pointwise credible bands (shaded area in the plots), compared against the true regression function $f_0$. 

\begin{table}[H]
\centering
\setlength{\tabcolsep}{4pt}
\resizebox{\linewidth}{!}{%
\begin{tabular}{|ll|ccccc|}
\hline
$y^c$ & $\eta$ & Standard GP & Hard threshold & IMQ (no debias.) & IMQ (debias.) & Bayes-trimmed \\
\hline
--- & $0$ & \textbf{0.0051}\,{\scriptsize(0.0010)} & \textbf{0.0051}\,{\scriptsize(0.0010)} & 0.0053\,{\scriptsize(0.0010)} & 0.0076\,{\scriptsize(0.0020)} & 0.0057\,{\scriptsize(0.0014)} \\
\hline
\multirow{2}{*}{$+3$}
 & $1/20$ & 0.0286\,{\scriptsize(0.0046)} & 0.0286\,{\scriptsize(0.0046)} & 0.0197\,{\scriptsize(0.0035)} & 0.0085\,{\scriptsize(0.0027)} & \textbf{0.0061}\,{\scriptsize(0.0016)} \\
 & $1/40$ & 0.0113\,{\scriptsize(0.0026)} & 0.0113\,{\scriptsize(0.0026)} & 0.0090\,{\scriptsize(0.0020)} & 0.0078\,{\scriptsize(0.0022)} & \textbf{0.0057}\,{\scriptsize(0.0015)} \\
\hline
\multirow{2}{*}{$+8$}
 & $1/20$ & 0.1713\,{\scriptsize(0.0122)} & 0.0568\,{\scriptsize(0.0067)} & 0.0283\,{\scriptsize(0.0044)} & 0.0076\,{\scriptsize(0.0021)} & \textbf{0.0060}\,{\scriptsize(0.0016)} \\
 & $1/40$ & 0.0484\,{\scriptsize(0.0064)} & 0.0182\,{\scriptsize(0.0035)} & 0.0111\,{\scriptsize(0.0024)} & 0.0076\,{\scriptsize(0.0020)} & \textbf{0.0057}\,{\scriptsize(0.0015)} \\
\hline
\end{tabular}%
}
\caption{Scenario~1: random contamination. Mean integrated squared
$L_2$ error (standard deviation) over $100$ replications, $n=5000$.
Best mean in bold.}
\label{tab:L2-diffuse}
\end{table}

In the first scenario considering random contamination one can observe from Table \ref{tab:L2-diffuse} that the standard GP is vulnerable to data contamination and almost all methods substantially outperform its behavior. For the clean data the standard Bayesian approach achieves the highest accuracy, but the other methods provide comparable performance, only slightly increasing the MISE. The best performing method is the trimmed-Bayesian approach (Method V) proposed in Section \ref{sec : generalized bayes trimmed} with the de-biased IMQ weighted method (Method IV) being close second. These behaviours are also visualized in {Figure \ref{fig:diffused main}}.

\begin{table}[H]
\centering
\setlength{\tabcolsep}{4pt}
\resizebox{\linewidth}{!}{%
\begin{tabular}{|ll|ccccc|}
\hline
$y^c$ & $\eta$ & Standard GP & Hard threshold & IMQ (no debias.) & IMQ (debias.) & Bayes-trimmed \\
\hline
\multirow{2}{*}{$-3$}
 & $1/20$ & 0.2114\,{\scriptsize(0.0193)} & 0.2113\,{\scriptsize(0.0193)} & 0.1733\,{\scriptsize(0.0158)} & \textbf{0.0560}\,{\scriptsize(0.0089)} & 0.1007\,{\scriptsize(0.0070)} \\
 & $1/40$ & 0.0634\,{\scriptsize(0.0083)} & 0.0634\,{\scriptsize(0.0083)} & 0.0523\,{\scriptsize(0.0066)} & \textbf{0.0200}\,{\scriptsize(0.0042)} & 0.0418\,{\scriptsize(0.0053)} \\
\hline
\multirow{2}{*}{$-8$}
 & $1/20$ & 1.4262\,{\scriptsize(0.0988)} & 0.8438\,{\scriptsize(0.0634)} & 0.4956\,{\scriptsize(0.0428)} & \textbf{0.0188}\,{\scriptsize(0.0046)} & 0.1124\,{\scriptsize(0.0074)} \\
 & $1/40$ & 0.3927\,{\scriptsize(0.0337)} & 0.2043\,{\scriptsize(0.0183)} & 0.1073\,{\scriptsize(0.0113)} & \textbf{0.0108}\,{\scriptsize(0.0026)} & 0.0506\,{\scriptsize(0.0055)} \\
\hline
\end{tabular}%
}
\caption{Scenario~2: localized contamination. Mean integrated squared
$L_2$ error (standard deviation) over $100$ replications, $n=5000$.
Best mean in bold.}
\label{tab:L2-localized}
\end{table}

\begin{figure}[H]
    \centering
    \includegraphics[width=\linewidth]{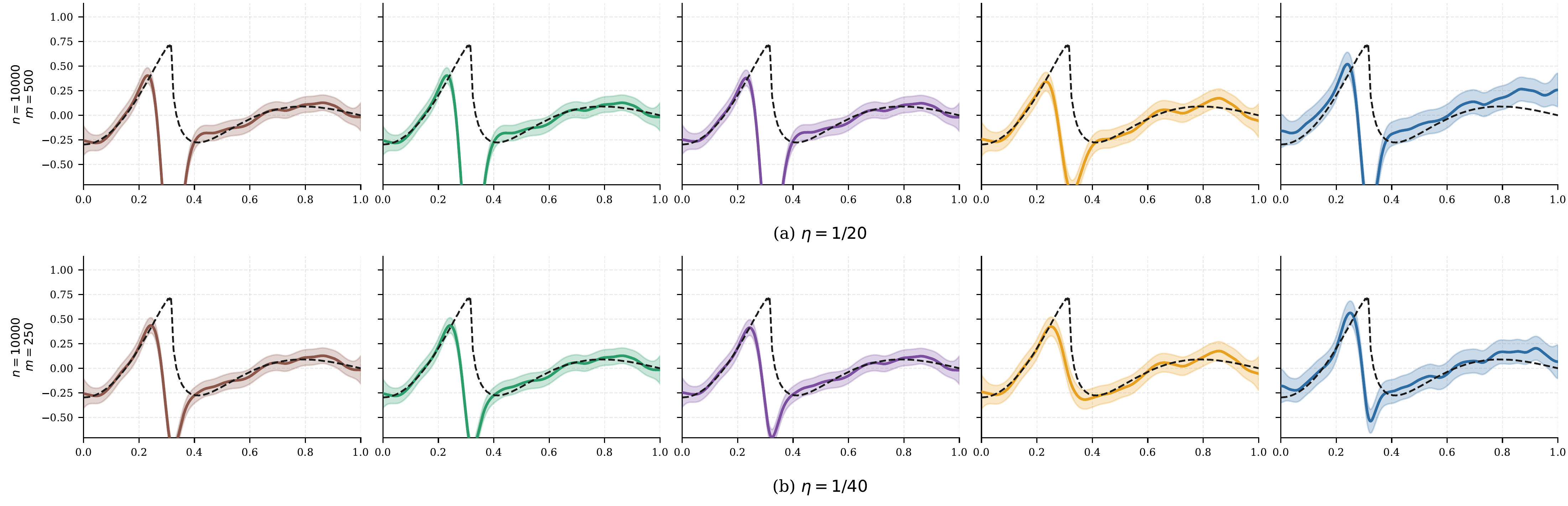}
    \caption{Scenario~2: localized contamination with shift $-3$ and
    $n=10\,000$. Top: $\eta=1/20$; bottom: $\eta=1/40$.
    Methods: I (brown), II (green), III (purple), IV (yellow), and V (blue). Black dashed curves show $f_0$.}
    \label{fig:localized-minus3}
\end{figure}

In the localized contamination case (Scenario 2), in view of Table \ref{tab:L2-localized}, the debiased IMQ weight procedure (Method IV) performs the best by far for contaminated data, while the trimmed GP (Method V) is the second best. The standard Bayesian approach, together with the hard thresholding method, is highly sensitive for data corruption in this case. Figure \ref{fig:localized-minus3} demonstrates that none of the methods can reliably recover the spike in the function due to the adversarial masking. Nevertheless, it is also clear that Method IV is substantially more accurate than the others, while Method V being the second best. We provide additional figures and tables for both Scenario 1 and 2 in the supplement.

\begin{table}[H]
\centering
\setlength{\tabcolsep}{4pt}
\resizebox{\linewidth}{!}{%
\begin{tabular}{|ll|ccccc|}
\hline
\shortstack{Outlier\\value} & $\eta$
 & Standard GP
 & \shortstack{Hard threshold\\(pilot-rec.)}
 & IMQ (debias.)
 & \shortstack{IMQ (debias.)\\(pilot-rec.)}
 & \shortstack{Bayes-trimmed\\(demeaned)} \\
\hline
\multirow{2}{*}{$-2$}
 & $1/20$ & 0.0466\,{\scriptsize(0.0057)} & 0.0466\,{\scriptsize(0.0057)} & 0.0561\,{\scriptsize(0.0066)} & \textbf{0.0080}\,{\scriptsize(0.0021)} & 0.0136\,{\scriptsize(0.0034)} \\
 & $1/40$ & 0.0159\,{\scriptsize(0.0030)} & 0.0159\,{\scriptsize(0.0030)} & 0.0198\,{\scriptsize(0.0038)} & \textbf{0.0077}\,{\scriptsize(0.0020)} & 0.0080\,{\scriptsize(0.0022)} \\
\hline
\end{tabular}%
}
\caption{Scenario~3: offset $V=2$, with diffuse outliers at $-2$.
Mean integrated squared $L_2$ error (standard deviation) over
$100$ replications, $n=5000$. Best mean in bold;
pilot-rec.\ denotes recentering at a pilot fit.}
\label{tab:L2-offset}
\end{table}

\begin{figure}[H]
    \centering
    \includegraphics[width=\linewidth]{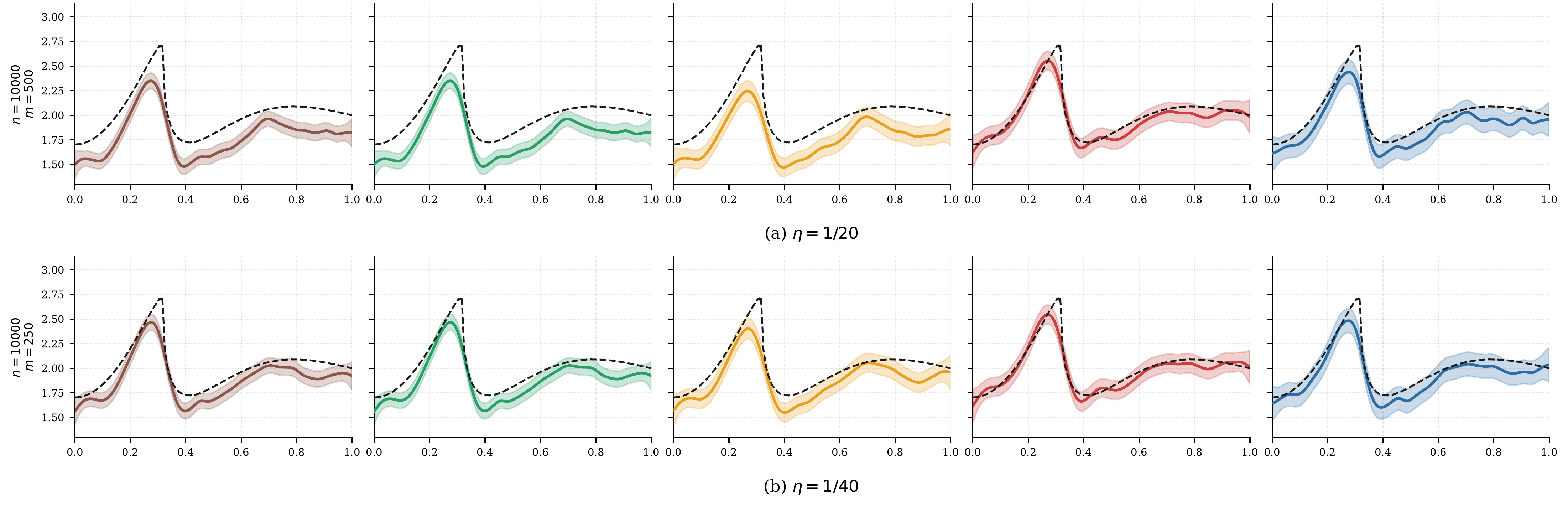}
    \caption{Scenario~3: offset and diffuse contamination with $n=10\,000$. Top: $\eta=1/20$; bottom: $\eta=1/40$. Methods: I (brown), II (green), IV (yellow), IV' (red), and V' (blue). Black dashed curves represent $f_0+2$.}
    \label{fig:shifted-outliers}
\end{figure}

Finally, we investigate Scenario 3. In this case we compare Methods I, II and IV with a demeaned version of Method~V and the centered IMQ method with debiasing. We update Method V by first estimating the coefficient $f_{0,0}$ corresponding to a constant shift and then subtracting it from the responses before estimating the remaining coefficients. In the end we add back the estimated constant shift. We denote this procedure by Method V'. This modification in principle reduces the supremum norm of the signal (by subtracting the additive constant term) and hence reduces the error in the trimmed estimator, which depends on the size of the signal. However, it also introduces additional uncertainty into the procedure by subtracting a random shift, see the supplement for additional numerical analysis of this procedure in Scenarios 1 and 2. Therefore, we advise to apply this modification if from the (contaminated) data it is clear that a constant shift is present in the signal, e.g. if all outcome variables are positive (or negative).

We also consider a modified version of Method IV following the idea introduced in Section \ref{sec : recenter}, by re-centering the weight function using a preliminary fit $\hat{f}_n$. We take $\hat{f}_n$ to be the trimmed projection estimator of Section~\ref{sec : upper bound}, computed with the same demeaning step, and set $\tau$ to the empirical $0.7$-quantile of $\{|y_i-\hat f_n(x_i)|\}_{i=1}^n$. We denote this approach by Method IV'. One can see from Table \ref{tab:L2-offset} that Method IV' performs the best and Method V' is close second. The other approaches are substantially worse, not even improving on the standard GP approach. These findings are also illustrated in Figure \ref{fig:shifted-outliers}.

In the supplement we provide additional numerical analysis. First we consider heavy tailed noises and investigate the performance of the robust GP methods in this setting. Then we investigate adaptation to the typically unknown regularity parameter $\beta$ by using the generalized hierarchical GP approach introduced in Section \ref{subsec : gp example}.

\section{Discussion}\label{sec:discussion}
We proposed various generalized Bayesian methods to deal with adversarial contamination in the context of the nonparametric regression model. We derived theoretical guarantees but also limitations for these procedures. We have considered various approaches to re-weight the terms in the factorized likelihood and derived contraction rates for the corresponding posteriors. Our general theorems were applied for Gaussian Process priors with deterministic and random rescaling and deep neural network priors. The advantage of these methods is that they do not depend on the contamination level $m$. We have also derived a novel trimming based GP approach which achieves the minimax estimation rate established in the paper, under the assumption that both the regularity of the true function and the contamination level are known. We also show that if this information is not available, then the minimax rate is polynomially higher and an appropriately reweighted generalized hierarchical GP prior can reach the corresponding information theoretic limits. In our study we cover several methods proposed in the literature but lacking theoretical underpinning, but also propose novel approaches with optimal theoretical guarantees. We also provide theory driven guidelines for the tuning of such approaches to achieve numerical performance while retaining those guarantees.

In the regression model we considered the noise variance $\sigma>0$ to be known. In practice one typically estimates it from the data or endows it with a hyper-prior distribution. We believe that our results can be extended to cover this setting as well. In Section \ref{subsec : gp example} we considered endowing the scale parameter of the GP prior with a hyper-prior. Alternatively one can also consider a plug-in estimator for the hyper-parameter. Furthermore, in practice one might want to consider a data driven choice for the threshold parameter $\tau_n$, which can outperform the perhaps overly conservative deterministic choices provided in the theoretical investigation. We leave these and other extensions for future work.

Another natural extension is to consider more general nonparametric models, including density estimation and nonparametric classification, under adversarial contamination. In these models the weight function enters the likelihood in a different manner. Nevertheless, the idea of separating the contaminated and clean data, together with the statistical and information theoretic tools developed in our paper could provide a natural basis for solving such problems. 

Another interesting and practically relevant research direction is to investigate the robustness properties of heavy-tailed likelihoods introduced in the Bayesian model. While this approach is substantially different from the reweighting techniques considered in our paper, the statistical tools developed herein can serve as a starting point for studying the corresponding rates.
\begin{acks}[Acknowledgments]
The authors would like to thank Chao Gao for insightful comments.

Large Language Models were solely used for double checking the proofs and for detecting notational inconsistencies. The statements and proofs were entirely derived by the authors and they take full responsibility for their validity.
\end{acks}

\begin{funding}
The research leading to these results was cofunded by the ANR (BACKUP project ANR-23-CE40-0018-01) and the European Union (ERC, BigBayesUQ, project number: 101041064). Views and opinions expressed are, however, those of the author(s) only and do not necessarily reflect those of the European Union or the European Research Council. Neither the European Union nor the granting authority can be held responsible for them.
\end{funding}


\bibliographystyle{imsart-number}
\bibliography{bib/bibRGPR}

\pagebreak


\appendix
\section{Additional numerical analysis}

First we collect the figures associated to the different scenarios studied in Section \ref{sec: numerics} which due to space limitation could not be displayed in the main paper. Then we investigate the behaviour of the robust Bayesian approaches in case of heavy tailed noise. Finally, we study the generalized hierarchical GP methods for adapting to the unknown regularity of the regression function.

\subsection{Additional figures and tables}

\begin{figure}[H]
    \centering
    \includegraphics[width=\linewidth]{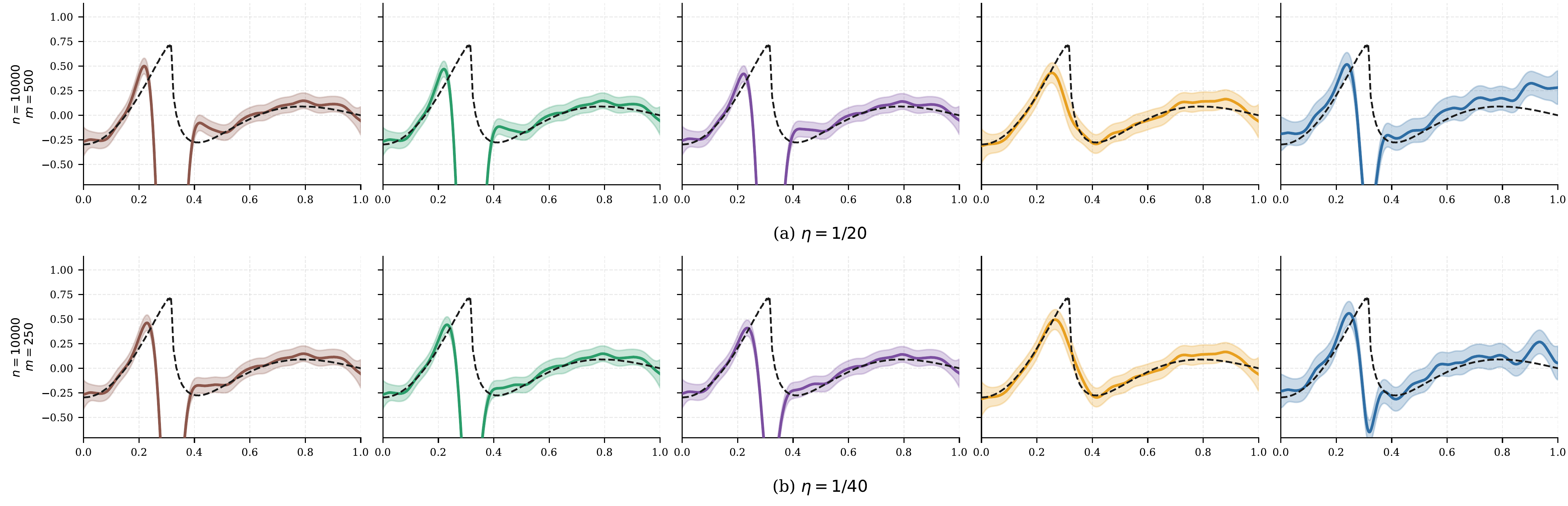}
    \caption{Scenario~2: localized contamination with shift $-8$ and $n=10\,000$. Top: $\eta=1/20$; bottom: $\eta=1/40$. Methods: I (brown), II (green), III (purple), IV (yellow), and V (blue). Black dashed curves show $f_0$.}
    \label{fig:localized-minus8}
\end{figure}

Figure \ref{fig:localized-minus8} shows the posterior means and credible sets resulting from the robust Bayesian procedures in Scenario 2 with a shift $y^c=8$  in the contamination. This figure is the counterpart of Figure \ref{fig:localized-minus3} with $y^c=3$. As can be seen in Table \ref{tab:L2-localized}, Method IV’s superior performance over the other approaches is even more pronounced in this setting.

We also investigate the performance of Methods IV' and V' in Scenarios 1 and 2. One can notice that in these Scenarios Methods IV' and V' perform a bit worse than their counterparts, i.e. Methods IV and V, respectively. However, these modifications provide additional protection against highly adversarial settings, like the one studied in Scenario 3. 

\begin{table}[H]
\centering
\setlength{\tabcolsep}{4pt}
\resizebox{\linewidth}{!}{%
\begin{tabular}{|ll|ccccccc|}
\hline
$y^c$ & $\eta$ & Standard GP & Hard threshold & IMQ (no debias.) & IMQ (debias.) & Bayes-trimmed & IMQ (debiased, pilot-rec.) & Bayes-trimmed (demean) \\
\hline
--- & $0$ & \textbf{0.0051}\,{\scriptsize(0.0010)} & \textbf{0.0051}\,{\scriptsize(0.0010)} & 0.0053\,{\scriptsize(0.0010)} & 0.0076\,{\scriptsize(0.0020)} & 0.0057\,{\scriptsize(0.0014)} & 0.0076\,{\scriptsize(0.0021)} & 0.0059\,{\scriptsize(0.0015)} \\
\hline
\multirow{2}{*}{$+3$}
 & $1/20$ & 0.0286\,{\scriptsize(0.0046)} & 0.0286\,{\scriptsize(0.0046)} & 0.0197\,{\scriptsize(0.0035)} & 0.0085\,{\scriptsize(0.0027)} & \textbf{0.0061}\,{\scriptsize(0.0016)} & 0.0086\,{\scriptsize(0.0028)} & 0.0139\,{\scriptsize(0.0035)} \\
 & $1/40$ & 0.0113\,{\scriptsize(0.0026)} & 0.0113\,{\scriptsize(0.0026)} & 0.0090\,{\scriptsize(0.0020)} & 0.0078\,{\scriptsize(0.0022)} & \textbf{0.0057}\,{\scriptsize(0.0015)} & 0.0078\,{\scriptsize(0.0023)} & 0.0085\,{\scriptsize(0.0026)} \\
\hline
\multirow{2}{*}{$+8$}
 & $1/20$ & 0.1713\,{\scriptsize(0.0122)} & 0.0568\,{\scriptsize(0.0067)} & 0.0283\,{\scriptsize(0.0044)} & 0.0076\,{\scriptsize(0.0021)} & \textbf{0.0060}\,{\scriptsize(0.0016)} & 0.0076\,{\scriptsize(0.0022)} & 0.0142\,{\scriptsize(0.0035)} \\
 & $1/40$ & 0.0484\,{\scriptsize(0.0064)} & 0.0182\,{\scriptsize(0.0035)} & 0.0111\,{\scriptsize(0.0024)} & 0.0076\,{\scriptsize(0.0020)} & \textbf{0.0057}\,{\scriptsize(0.0015)} & 0.0076\,{\scriptsize(0.0022)} & 0.0087\,{\scriptsize(0.0026)} \\
\hline
\end{tabular}%
}
\caption{Uncontaminated case and Scenario~1. Mean integrated squared
$L_2$ error (standard deviation) over $100$ replications, $n=5000$.
Best mean in bold.}
\end{table}

\begin{table}[H]
\centering
\setlength{\tabcolsep}{4pt}
\resizebox{\linewidth}{!}{%
\begin{tabular}{|ll|ccccccc|}
\hline
$y^c$ & $\eta$ & Standard GP & Hard threshold & IMQ (no debias.) & IMQ (debias.) & Bayes-trimmed & IMQ (debiased, pilot-rec.) & Bayes-trimmed (demean) \\
\hline
\multirow{2}{*}{$-3$}
 & $1/20$ & 0.2114\,{\scriptsize(0.0193)} & 0.2113\,{\scriptsize(0.0193)} & 0.1733\,{\scriptsize(0.0158)} & \textbf{0.0560}\,{\scriptsize(0.0089)} & 0.1007\,{\scriptsize(0.0070)} & 0.1325\,{\scriptsize(0.0158)} & 0.1198\,{\scriptsize(0.0090)} \\
 & $1/40$ & 0.0634\,{\scriptsize(0.0083)} & 0.0634\,{\scriptsize(0.0083)} & 0.0523\,{\scriptsize(0.0066)} & \textbf{0.0200}\,{\scriptsize(0.0042)} & 0.0418\,{\scriptsize(0.0053)} & 0.0332\,{\scriptsize(0.0071)} & 0.0463\,{\scriptsize(0.0062)} \\
\hline
\multirow{2}{*}{$-8$}
 & $1/20$ & 1.4262\,{\scriptsize(0.0988)} & 0.8438\,{\scriptsize(0.0634)} & 0.4956\,{\scriptsize(0.0428)} & \textbf{0.0188}\,{\scriptsize(0.0046)} & 0.1124\,{\scriptsize(0.0074)} & 0.0251\,{\scriptsize(0.0060)} & 0.1357\,{\scriptsize(0.0090)} \\
 & $1/40$ & 0.3927\,{\scriptsize(0.0337)} & 0.2043\,{\scriptsize(0.0183)} & 0.1073\,{\scriptsize(0.0113)} & \textbf{0.0108}\,{\scriptsize(0.0026)} & 0.0506\,{\scriptsize(0.0055)} & 0.0119\,{\scriptsize(0.0031)} & 0.0563\,{\scriptsize(0.0067)} \\
\hline
\end{tabular}%
}
\caption{Scenario~2. Mean integrated squared
$L_2$ error (standard deviation) over $100$ replications, $n=5000$.
Best mean in bold.}
\end{table}

\subsection{Heavy-tailed noise}
In this subsection, we investigate the robustness of the proposed procedures against heavy-tailed noise. Although this setting is not covered by our theoretical analysis, we assess the numerical performance of the methods under such setting as well.

We compare Methods I--V from Section \ref{sec: numerics} using the same true regression function $f_0$ defined in \eqref{eq:true-function}, changing only the noise distribution. Specifically, we consider a Student's $t$ distribution with $3$ degrees of freedom, rescaled by $1/\sqrt{3}$; a Laplace distribution with location $0$ and scale $1/\sqrt{2}$; and a Cauchy distribution with location $0$ and scale $1$. The first two distributions have variance $1$, matching the noise variance in the Gaussian experiment. The Cauchy distribution, however, has neither a finite mean nor a finite variance. In this case, $y_i \mid x_i$ follows a Cauchy distribution with location $f_0(x_i)$ and scale $1$, so $f_0(x_i)$ is the conditional median, while the conditional mean is undefined. Since we change the entire noise distribution, there is no prescribed number of outliers $m$. We therefore use a conservative contamination proportion of $1/20$ for the Bayes Trimmed Projection Posterior under all three noise distributions.

\begin{table}[H]
\centering
\setlength{\tabcolsep}{4pt}
\resizebox{\linewidth}{!}{%
\begin{tabular}{|l|ccccc|}
\hline
Noise & Standard GP & Hard threshold & IMQ (no debias.) & IMQ (debias.) & Bayes-trimmed \\
\hline
Student's $t_3$
 & 0.0051\,{\scriptsize(0.0008)}
 & 0.0048\,{\scriptsize(0.0007)}
 & 0.0047\,{\scriptsize(0.0007)}
 & 0.0276\,{\scriptsize(0.0059)}
 & \textbf{0.0045}\,{\scriptsize(0.0009)} \\
\hline
Laplace
 & 0.0051\,{\scriptsize(0.0010)}
 & 0.0051\,{\scriptsize(0.0010)}
 & 0.0050\,{\scriptsize(0.0009)}
 & 0.0216\,{\scriptsize(0.0053)}
 & \textbf{0.0047}\,{\scriptsize(0.0010)} \\
\hline
Cauchy
 & 3318.7535\,{\scriptsize(26431.3312)}
 & 0.0133\,{\scriptsize(0.0041)}
 & 0.0108\,{\scriptsize(0.0027)}
 & \textbf{0.0080}\,{\scriptsize(0.0016)}
 & 0.0105\,{\scriptsize(0.0031)} \\
\hline
\end{tabular}%
}
\caption{Mean squared $L_2$ error under rescaled Student's $t_3$, Laplace, and Cauchy noise over $100$ replications with $n=5000$. Standard deviations are reported in parentheses. Best mean is shown in bold.}
\label{tab:L2-noise}
\end{table}

\begin{figure}[H]
    \centering
    \includegraphics[width=0.9\linewidth]{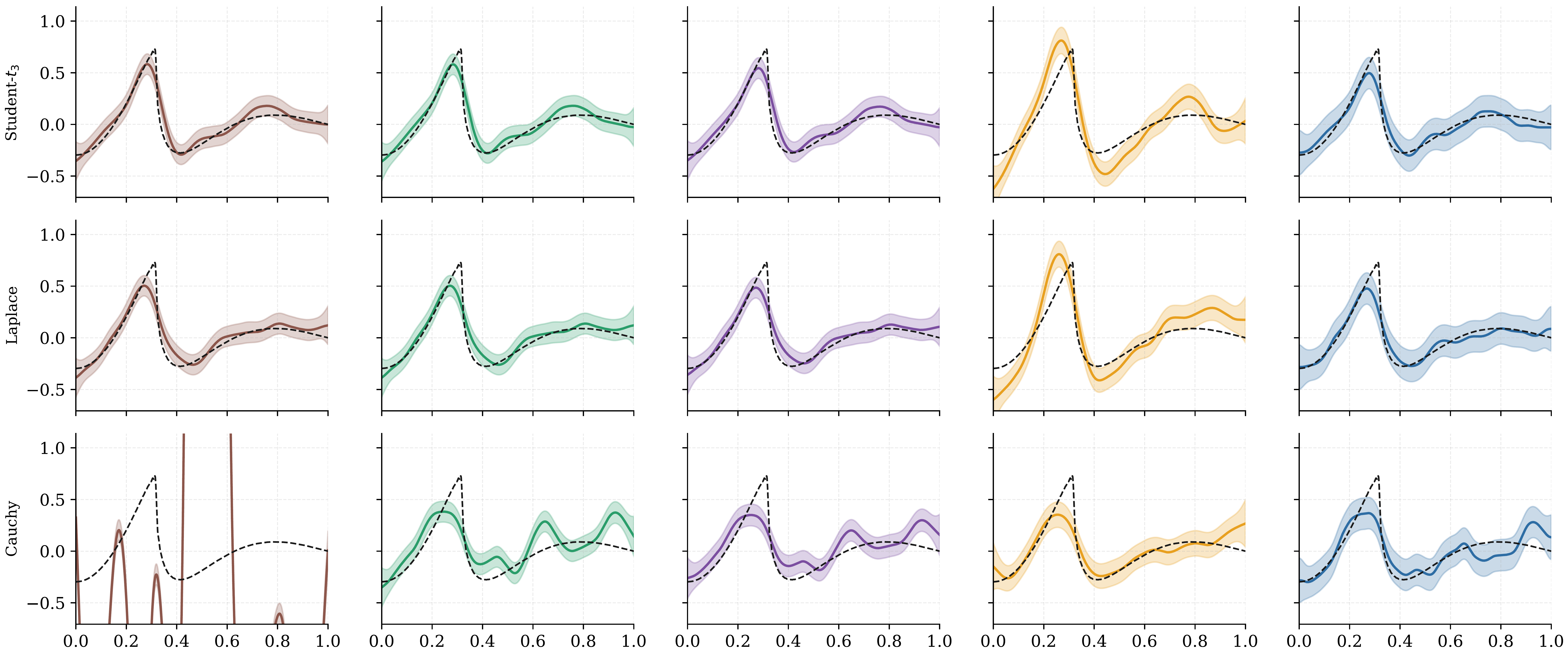}
    \caption{Comparison of methods I (brown), II (green), III (purple), IV (yellow), and V (blue) under rescaled Student's $t_3$, Laplace, and Cauchy noise, with $n=5000$. Black dashed curves show $f_0$.}
    \label{fig:diffused}
\end{figure}

In view of Table \ref{tab:L2-noise} and Figure \ref{fig:diffused} one can conclude that in the first two cases, with Student and Laplace noise, all methods except Method IV, perform similarly well, with the Bayes-trimmed approach (Method V) achieving the best performance. In contrast, in case of heavy tailed noise, the benchmark standard GP performs terribly, while the robust GP methods improve its performance substantially. The best performing method is the debiased version of the IMQ weights (Method IV).

Next we discuss this substantial difference in the performance of the standard GP method for Cauchy and lighter tailed noises.
In view of Section~4 of \cite{misspecificationInfinite}, under suitable assumptions on the true noise distribution (in Lemma~4.1 of \cite{misspecificationInfinite} the true noise random variable $\xi_0$ was assumed to satisfy $E[e^{M|\xi_0|}] < \infty$ for every $M > 0$), the GP approach (with surrogate Gaussian likelihood) contracts around the true regression function (with the optimal rate). While none of the considered noise distributions satisfy this condition, it shows that standard GPs offer a certain degree of flexibility. Indeed, in the first two cases, with  Student's $t_3$ and Laplace noise it performs well and fails only with Cauchy noise. Most of the proposed robust methods provide a slight improvement, with the Bayes-trimmed method performing the best, even using a nominal contamination proportion of $m/n = 1/20$ that does not correspond to a prescribed fraction of contaminated observations. However, the debiasing approach of the IMQ weight performs the worst among all the methods. One possible explanation is that the correction is derived using Gaussian noise and Stein's identity, as discussed in our analysis. A correction adapted to the actual noise distribution might therefore improve performance, although this would require further investigation.

The Cauchy case presents a more fundamental difficulty:
Kullback--Leibler divergence from a Cauchy distribution to a Gaussian distribution is infinite. In this case the finite-KL framework of \cite{misspecificationInfinite} does not apply. This seems to be a fundamental problem as the standard GP approach provides a nonsensible solution. The proposed robust modifications help and in this case IMQ with debiasing achieves by far the lowest mean error, followed by Bayes-trimmed.

\subsection{Adaptive methods}

In this subsection we consider hierarchical priors in the generalized posteriors with the goal of adapting to the unknown regularity of the underlying true function, see Section \ref{subsec : gp example}. In particular, we consider a Mat\'ern kernel of regularity $3/2$. Recall that the true function $f_0$ considered in the experiments above has regularity $\beta = 1$ that differs from the regularity of the Mat\'ern kernel, hence one has to rescale the prior for optimal inference. Since typically one does not have access to the regularity parameter $\beta$ in practice, the oracle rescaling given in \eqref{rescaling} (and applied so far in the numerical analysis) cannot be applied. Hence, we consider a data driven rescaling, by endowing the rescaling parameter $\gamma$ with a hyper-prior, resulting in a two-layer, hierarchical prior distribution on the regression function. We use the following hyper-prior on $\gamma$
\[
p_{\alpha,n}(\gamma)
=
\frac{2\alpha}{2\alpha+1}\,
n^{\frac{1}{2\alpha+1}}\,
\gamma^{-\frac{1}{2\alpha+1}}
\exp\!\left \{
-n^{\frac{1}{2\alpha+1}}
\left(\gamma^{\frac{2\alpha}{2\alpha+1}}-1\right)
\right \}
\mathbf{1}_{\{\gamma\ge 1\}}.
\]
Note that this prior satisfies the condition in equation \eqref{eq:hyperprior}. 

\begin{figure}[H]
    \centering
    \includegraphics[width=0.9\linewidth]{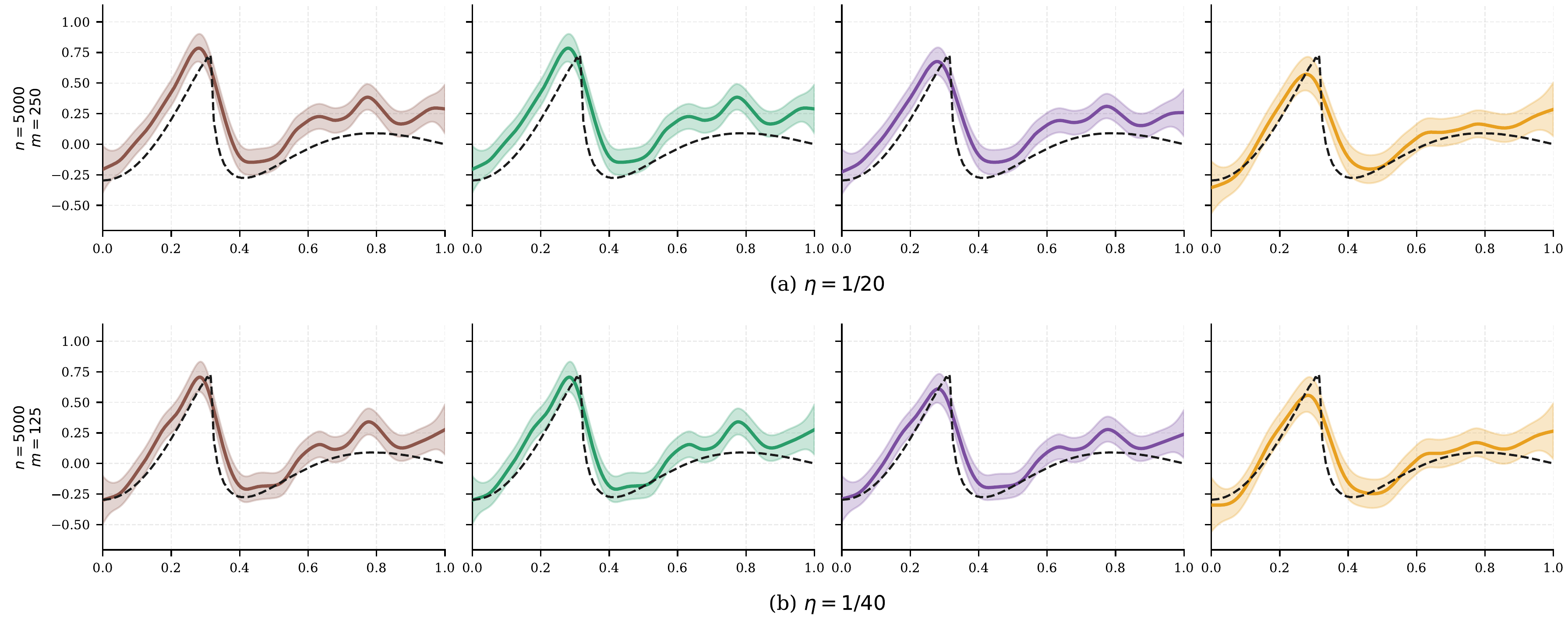}
    \caption{Scenario~1: random contamination with $n=5000$. Top: $\eta = 1/20$; bottom: $\eta = 1/40$.
    Methods: I (brown), II (green), III (purple), IV (yellow).
    Black dashed curves show $f_0$. }
    \label{fig:adapt:diffused}
\end{figure}

\begin{figure}[H]
    \centering
    \includegraphics[width=0.9\linewidth]{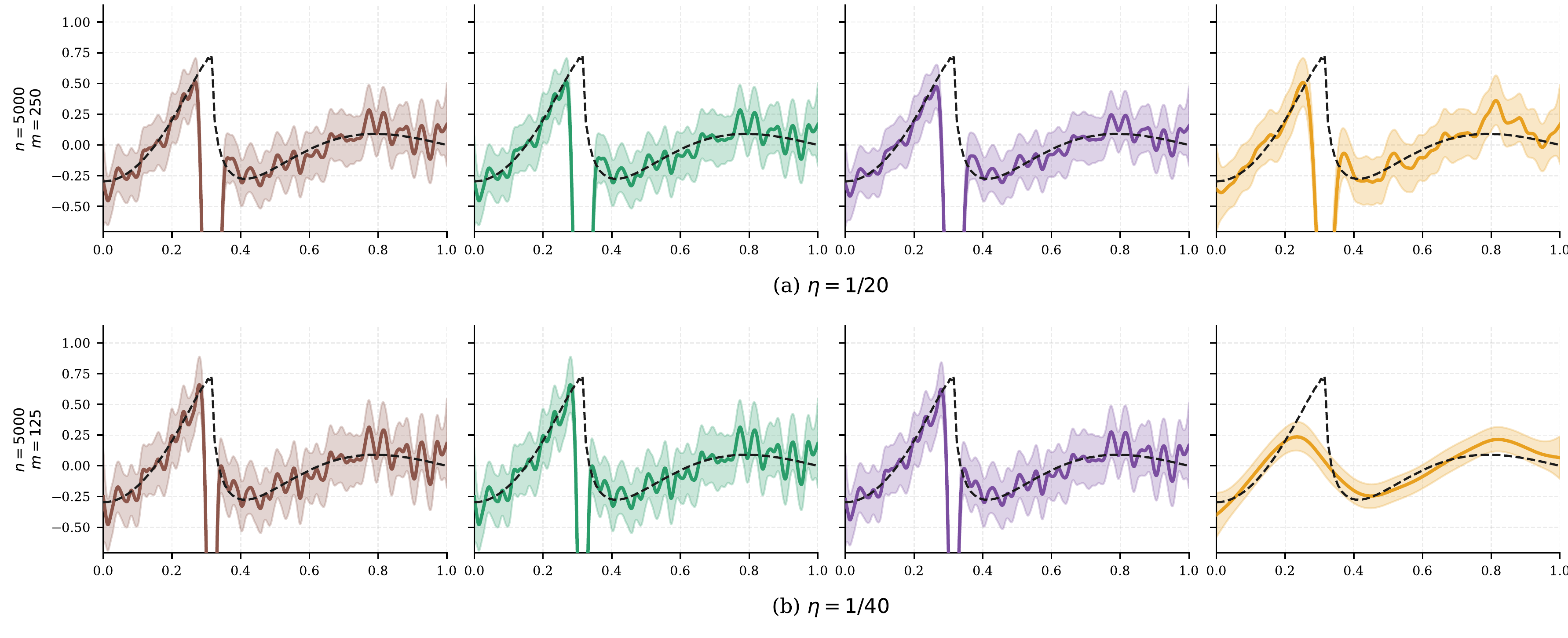}
    \caption{Scenario~2: localized contamination with shift $-3$ and
    $n=5000$. Top: $\eta=1/20$; bottom: $\eta=1/40$.
    Methods: I (brown), II (green), III (purple) and IV (yellow). Black dashed curves show $f_0$.}
    \label{fig:adapt:diffused2}
\end{figure}

We focus on Scenarios $1$ and $2$ and restrict our study to the case of $n = 5000$ and $y^c=3$ due to the bigger computational cost of the hierarchical methods. In view of Tables \ref{tab:L2-diffuse-plus3} and \ref{tab:L2-localized-minus3} and Figures \ref{fig:adapt:diffused} and \ref{fig:adapt:diffused2}, the IMQ debiased approach performs the best in the contaminated setting, sometimes providing an order of magnitude smaller MSE than the other approaches. The biggest difference is in the localized contamination setting, where the other approaches seem to overfit, while the debiased IMQ approach provides a more reasonable, smoother approximation. At the same time, in the uncontaminated setting it loses a bit from the accuracy compared to the standard GP approach, but still provides a reasonable estimate.


\begin{table}[H]
\centering
\setlength{\tabcolsep}{4pt}
\resizebox{0.8\linewidth}{!}{%
\begin{tabular}{|ll|cccc|}
\hline
$y^c$ & $\eta$ & Standard GP & Hard threshold & IMQ (no debias.) & IMQ (debias.) \\
\hline
--- & $0$ & \textbf{0.0055}\,{\scriptsize(0.0012)} & \textbf{0.0055}\,{\scriptsize(0.0012)} & 0.0061\,{\scriptsize(0.0016)} & 0.0092\,{\scriptsize(0.0034)} \\
\hline
\multirow{2}{*}{$+3$}
 & $1/20$ & 0.0288\,{\scriptsize(0.0049)} & 0.0288\,{\scriptsize(0.0049)} & 0.0204\,{\scriptsize(0.0039)} & \textbf{0.0101}\,{\scriptsize(0.0042)} \\
 & $1/40$ & 0.0117\,{\scriptsize(0.0026)} & 0.0117\,{\scriptsize(0.0026)} & 0.0098\,{\scriptsize(0.0024)} & \textbf{0.0092}\,{\scriptsize(0.0037)} \\
\hline
\end{tabular}%
}
\caption{Uncontaminated case and Scenario 1. Mean integrated squared $L_2$ error (standard deviation) over 20 replications, $n=5000$. Best mean in bold.}
\label{tab:L2-diffuse-plus3}
\end{table}

\begin{table}[H]
\centering
\setlength{\tabcolsep}{4pt}
\resizebox{0.8\linewidth}{!}{%
\begin{tabular}{|ll|cccc|}
\hline
$y^c$ & $\eta$ & Standard GP & Hard threshold & IMQ (no debias.) & IMQ (debias.) \\
\hline
\multirow{2}{*}{$-3$}
 & $1/20$ & 0.3982\,{\scriptsize(0.0339)} & 0.3981\,{\scriptsize(0.0339)} & 0.3525\,{\scriptsize(0.0311)} & \textbf{0.1958}\,{\scriptsize(0.0221)} \\
 & $1/40$ & 0.1747\,{\scriptsize(0.0229)} & 0.1746\,{\scriptsize(0.0229)} & 0.1468\,{\scriptsize(0.0200)} & \textbf{0.0213}\,{\scriptsize(0.0043)} \\
\hline
\end{tabular}%
}
\caption{Scenario 2. Mean integrated squared $L_2$ error (standard deviation) over 20 replications, $n=5000$. Best mean in bold.}
\label{tab:L2-localized-minus3}
\end{table}

\section{Additional results}

\subsection{Bounded posteriors}\label{app : bounded prior}
In this section we consider posteriors with bounded sample paths and extend the contraction rate with respect to the empirical $L_2$-norm on the inliers, derived in Theorem \ref{thm : genweight fix}, to the stronger $L_2(P_X)$ contraction. One can achieve boundedness of the posterior by using bounded priors or by imposing this constraint on the posterior by `truncation with rejection'. That is discarding samples with overly large supremum norm. 
Formally, for a probability measure $P$ on $\cF$ and a sequence $s_n >0$ (later taken to be constant or slowly diverging), such that $P(\{f \, : \|f\|_\infty\le s_n\}) >0$, we denote the \textit{truncated-by-rejection} version of $P$ as
\[ P^{[s_n]}(\cdot) := \frac{P(\, \cdot \, \cap \{f \in \cF \, : \|f\|_\infty \le s_n \})}{P(\{f \in \cF \, : \|f\|_\infty \le s_n \})}.\]
We remark that, by Bayes' formula, applying this type of truncation to the posterior or to the prior is equivalent, indeed for any $B \subseteq \cF$,
\[ G^{[s_n]}(B\given \widetilde{\cD}_n) \propto \int_{B \cap \{\|f\|_\infty \le s_n\}} e^{-\tilde L_{[n],w}(f)} \,d\Pi(f) \propto \int_{B} e^{-\tilde L_{[n],w}(f)} \,d\Pi^{[s_n]}(f). \]
In particular, employing such truncated posteriors allows one to work with bounded function draws. This can be particularly useful to simplify posterior analysis as demonstrated in Theorem \ref{thm : debiased}. The next theorem, proven in Section \ref{proof : genweight random}, leverages this boundedness to obtain a contraction rate in $L_2(P_X)$ norm.

\begin{thm}\label{thm : genweight random} Assume the prior $\Pi$ satisfies condition \eqref{eq : prior mass condition} for some $\veps_n$. Suppose the weight function $w$ satisfies conditions $\eqref{eq : cond w cont}$ with $\tau_n > 2 \si \sqrt{\log n}$, \eqref{eq : cond lowerb w}, and \eqref{conddiffw} with $\rho_n \ge 2\si \sqrt{\log n}$ and $\Delta_n = o(1)$. Let $s_n \ge F_0 + \veps_n$ and consider $G^{[s_n]}(\, \cdot \,|\, \widetilde{\cD}_n)$ the truncated-by-rejection generalized posterior. For $\zeta_{n,m}^2 := s_n^2\left[\frac{m}{n} (\tau_n^2 + s_n^2)+ \varepsilon_n^2+ (\Delta_n \rho_n)^2\right]$ and $M >0$ large enough, as $n \to \infty$,
    \[ \widetilde{E}_{f_0} G^{[s_n]}\Big(f \, : \, \|f-f_0\|_{2,P_X}^2 \ge M\zeta_{n,m}^2  \mid \widetilde{\cD}_n\Big) \to 0.
    \]
\end{thm}
If a bound $F > F_0$ is available, the truncation at $s_n =F$ in Theorem \ref{thm : genweight random} yields $L_2(P_X)$ contraction at the same rate as Theorem \ref{thml2} without assuming the GP structure (see also Theorem \ref{thm : debiased}). If no such bound is known on $F_0$, one can in principle truncate at $s_n^2 = \log n$ and inflate the rate by this logarithmic factor in our upper bound.

\subsection{A robust sup--norm-consistent estimator}
\label{sec : sup-norm trimmed}

We now show that the construction in Section \ref{sec : upper bound} can also yield an estimator that is consistent in the supremum norm. To this end, we replace the Fourier basis with a compactly supported wavelet basis. The main advantage of using wavelets is their spatial localisation, which makes it possible to control the $L^\infty$-norm of a function directly through its wavelet coefficients. 

For simplicity, we work on $[0,1]$. The argument extends in a standard way to $[0,1]^d$, at the expense of heavier notation. We consider the boundary-corrected wavelet basis introduced by Cohen, Daubechies, and Vial (see e.g. \cite{GineNickl2016}, Chapter 4 and \cite{supnormcastillo}), which will be referred to as the CDV basis and denoted by
$\left\{\psi_{\ell,k}:\ell\geq 0,\;
0\leq k\leq 2^\ell-1\right\}$.
We assume that the regularity of the basis is larger than the smoothness parameter $\beta$ introduced below. We recall the following standard properties of the CDV basis:
\begin{enumerate}
\item[(i)] $\{\psi_{\ell,k}\}_{\ell,k}$ forms an orthonormal basis of $L^2[0,1]$.
\item[(ii)] There exists a constant $C>0$, independent of $\ell$ and $k$, such that $\|\psi_{\ell,k}\|_\infty
\leq C2^{\ell/2} $.
\item[(iii)] The basis is spatially localised. More precisely, there exists a constant $C>0$ such that
\[
\sup_{x\in[0,1]}
\sum_{k=0}^{2^\ell-1}
|\psi_{\ell,k}(x)|
\leq C\,2^{\ell/2},
\qquad \ell\geq 0.
\]
\end{enumerate}
Based on this wavelet basis we define the $\beta$-regular Besov-$(\infty,\infty)$-space and the corresponding $M$-radius ball as
\begin{equation}
\label{eq : besov inf inf ball}
B_{\infty,\infty}^{\beta}(M)
:=
\left\{
f\in L^2([0,1]):
\sup_{\ell\geq 0}
\max_{0\leq k\leq 2^\ell-1}
2^{\ell(\beta+1/2)}
\left|\langle f,\psi_{\ell,k}\rangle_2\right|
\leq M
\right\},
\end{equation}
where $\beta>0$ and $M>0$ (the whole Besov space $B^{\beta}_{\infty,\infty}$ can be written as $\cup_{0 < M < \infty} B^{\beta}_{\infty,\infty}(M)$). In the rest of the subsection we assume that $f_0\in B_{\infty,\infty}^{\beta}(M)$. 

Since the basis has a double index, we denote the estimator of the wavelet coefficient introduced in \eqref{trimmed estimator aggregated} by $\hat{f}_{\ell,k}$ in this case, 
\[
    \hat{f}_{\ell,k} := \frac{1}{n} \sum_{i=1}^n \phi_{\widetilde{a}_{\ell,k}, \widetilde{b}_{\ell,k}}( \widetilde y_i \psi_{\ell,k}(\widetilde x_i) ) ,
\]
where we choose the trimming level $q$ as follows
\[
 q = 8 \frac{m}{n} +  12 \frac{\log(16 n^2)}{n},
\]
and, given a suitable truncation level $L_{n,m}$, we define the \textit{trimmed projection estimator} in the wavelet basis as 
\[
    \hat{f}_n = \sum_{\ell = 0}^{L_{n,m}} \sum_{k=0}^{2^\ell -1} \hat{f}_{\ell,k} \psi_{\ell,k}
\]
Let us also define the following set
\[
    \mathcal{I}_{n,m} = \{ (\ell,k) : 0 \leq \ell \leq L_{n,m}, 0 \leq k \leq 2^\ell -1 \}
\]
and notice that the total number of coefficients used to compute $\hat{f}_n$ is $| \mathcal{I}_{n,m} | = 2^{L_{n,m} + 1} - 1$.

\begin{thm}
\label{sup norm trimmed estimator}
Assume the regression model with uniform design on $[0,1]$, and suppose that the sample has at most $m$ adversarial contaminations with $1 \leq m = o(n)$. For $\beta > 0$,
let the truncation level be
\[
    L_{n,m} \;:=\; \left\lfloor \log_2 \min\left\{
    \left (q \sqrt{ \log(1/q)} \right )^{-\frac{1}{\beta+1}},\;
    \left( \frac{n}{\log n} \right )^{\frac{1}{2\beta+1}} \right\} \right\rfloor .
\]
Then there exists a constant $C > 0$, such that, uniformly over $f_0 \in B_{\infty,\infty}^\beta(M)$ and over the
adversarial actions, and for all $n$ sufficiently large,
\[
    \bigl\| \hat{f}_n - f_0 \bigr\|_\infty
    \;\leq\; C \left[ \left ( \frac{m}{n} \sqrt{  \log{\left (\frac{n}{m} \right )}  } \right )^{\frac{\beta}{\beta+1}} 
    \;+\; \Bigl( \frac{\log n}{n} \Bigr)^{\frac{\beta}{2\beta+1}} \right]
\]
with $\widetilde{P}_{f_0}$-probability at least $1 - 1/n$.
\end{thm}

The proof is deferred to Section \ref{sup norm trimmed projection estimator proof}. Another approach for supremum norm consistency, directly based on the generalised posterior distribution, is considered in Section \ref{sec : consistency sup norm} below.

\section{Proofs of the main results}\label{sec : proof main}
\subsection{Additional notation and a fundamental lemma}
We begin this proof section with a few additional notations. For any measurable subset $S \subseteq \cF$ and any subset of indices $J\subseteq\{1,\dots,n\}$, we define the (generalized) posterior integrals 
\begin{equation}\label{def : post integral}
     I_{J}^w(S) := \int_S \exp \left\{ L_{J,w}(f_0)- L_{J,w}(f) \right\} \, d \Pi(f),
\end{equation}
and their contaminated version

\[\widetilde{I}_{J}^w(S) := \int_S \exp \left\{ \widetilde{L}_{J,w}(f_0)- \widetilde{L}_{J,w}(f) \right\} \, d \Pi(f).\]
We also denote $\cJ_{n,m}$ the set of index sets with at most $m$ indices,
\begin{equation}\label{def : jnm}
    \cJ_{n,m} := \left\{ J\subset\{1,\dots,n\} \, : \, |J|\le m \right\}.
\end{equation}
We recall that the corrupted dataset $\widetilde{\cD}_n$ is obtained from a measurable transformation of the clean dataset $\cD_n$ with the possible use of an auxiliary variable $U$, as
\[ \widetilde{\cD}_n = T(\cD_n,U).\]
As introduced in \eqref{def:outliers}, we denote the set of corrupted indices $\cO$ as
\begin{equation}\label{def : O}
    \cO = \cO(\cD_n,U) : = \left\{ i\in \{1,\dots,n\} \, : \, (x_i,y_i) \neq (\tilde{x}_i,\tilde{y}_i) \right\},
\end{equation}
in particular, we have $\cO \in \cJ_{n,m}$, since the adversary is allowed to modify at most $m$ observations.

We now state a lemma that we will use throughout the proof section to handle the presence of outliers in the data. 
\begin{lem}\label{lem : getting rid of the outliers}
For any nonnegative sequence $\veps_{n} \to 0$, define the subset
    \[ \cF_n := \{ f \, : \,  \| f-f_0\|_{\infty} \leq \veps_{n} \} \subset \cF.\]
Assume the weight function $w$ satisfies \eqref{eq : cond w cont}.
Let $S \subseteq \cF$ be any measurable set and $\cA$ be any event depending only on the clean dataset $\cD_n$. For $n$ large enough, we have
\[ \widetilde{E}_{f_0} G(S\, | \, \widetilde{\cD}_n) \le e^{m \frac{\lambda}{\sigma^2}(\tau_n^2 + 4F_0^2)}  E_{f_0} \left [ \max_{J \in \cJ_{n,m}}I_{J^c}^w(S) / \min_{J \in \cJ_{n,m}} I_{J^c}^w(\cF_n) \, \mathbf{1}_\cA \right ] +{P}_{f_0}^n(\cA^c).\]
\end{lem}
\begin{proof}
    Recall that $Q_{f_0} = P_{f_0}^n \otimes \nu$ denotes the joint law of $(\cD_n,U)$, so that
    \[ \widetilde{E}_{f_0} G(S\, | \, \widetilde{\cD}_n) = E_{Q_{f_0}} G(S\, | \, T(\cD_n,U)).\]
    For $\cA$ depending only on $\cD_n$, we have $Q_{f_0}(\cA^c) = P_{f_0}^n(\cA^c)$. Since $G(S\, | \, T(\cD_n,U)) \le 1$, 
    \[ \widetilde{E}_{f_0} G(S\, | \, \widetilde{\cD}_n) \le E_{Q_{f_0}} \left[ G(S\, | \, T(\cD_n,U)) \mathbf{1}_\cA \right] + P_{f_0}^n(\cA^c). \]
    Using $(a+b)^2 \le 2a^2+2b^2$, for any weight $w$ satisfying \eqref{eq : cond w cont} and any $(x,y) \in \R^d\times \R$, we get
    \[ 0 \le w(y)(y-f(x))^2 \le 2 \lambda(\tau_n^2 + f(x)^2).\]
    Splitting $\widetilde{L}_{[n],w}(f)$ over contaminated indices in $\cO$ and clean ones in $\cO^c$ and using the previously displayed inequalities on the corrupted part, one obtains
    \[ L_{\cO^c,w} (f) \le \widetilde{L}_{[n],w}(f) \le L_{\cO^c,w} (f) + \frac{\la}{\si^2}\sum_{i \in \cO} (\tau_n^2 + f(\tilde{x}_i)^2).\]
    Using these inequalities to bound both numerator and denominator of the posterior yields
    \[E_{Q_{f_0}} [ G(S \given T(\cD_n,U)) \mathbf{1}_\cA] \le E_{Q_{f_0}}\left[ \frac{\int_S e^{- L_{\cO^c,w} (f) } d\Pi(f)}{\int e^{- \frac{\la}{\si^2}\sum_{i \in \cO} (\tau_n^2 + f(\tilde{x}_i)^2)}e^{- L_{\cO^c,w} (f) } d\Pi(f)} \mathbf{1}_\cA\right].\]
    We lower bound the denominator by restricting the integral on $\cF_n=\{ f \, : \|f-f_0\|_\infty \le \veps_n \}$, on such a set and for $n$ large enough, we have
    \[ f(\tilde{x}_i)^2 \le \|f\|_\infty^2 \le 2\|f_0\|_\infty^2 + 2\|f -f_0\|_\infty^2 \le 2F_0^2 + 2\veps_n^2 \le 4F_0^2,\]
    where we used $\|f_0\|_\infty \le F_0$ as well as $\veps_n \to 0$ in the last step. Since $|\cO| \le m$, we get
    \[E_{Q_{f_0}} [ G(S \given T(\cD_n,U)) \mathbf{1}_\cA] \le e^{\frac{\la}{\si^2}m(\tau_n^2 + 4F_0^2)}E_{Q_{f_0}}\left[ \frac{\int_S e^{- L_{\cO^c,w} (f) } d\Pi(f)}{\int_{\cF_n} e^{- L_{\cO^c,w} (f) } d\Pi(f)} \mathbf{1}_\cA\right].\]
    Recalling the definition of $I_{\cO^c}^w$ in \eqref{def : post integral}, and multiplying both the numerator and the denominator by $\exp \{L_{\cO^c,w}(f_0)\} \ge 0 $, yields
    \[E_{Q_{f_0}} [ G(S \given T(\cD_n,U)) \mathbf{1}_\cA] \le e^{\frac{\la}{\si^2}m(\tau_n^2 + 4F_0^2)}E_{Q_{f_0}}\left[ {I_{\cO^c}^w(S)}/{I_{\cO^c}^w(\cF_n)} \mathbf{1}_\cA\right].\]
    Since $\cO \in \cJ_{n,m}$, we make this bound uniform and independent of the contamination by taking both the minimum and the maximum over $\cJ_{n,m}$. The expectation with respect to $Q_{f_0}$ of this quantity is then equivalent to the usual $E_{f_0}$ expectation, so that
    \[E_{Q_{f_0}} [ G(S \given T(\cD_n,U)) \mathbf{1}_\cA] \le e^{\frac{\la}{\si^2}m(\tau_n^2 + 4F_0^2)}E_{{f_0}}\left[ \max_{J \in \cJ_{n,m}}{I_{J^c}^w(S)}/\min_{J \in \cJ_{n,m}}{I_{J^c}^w(\cF_n)} \mathbf{1}_\cA\right].\]
\end{proof}

\begin{Remark}\label{rmk : J dependent S}
    We note that Lemma \ref{lem : getting rid of the outliers} can  {in particular} be applied to a family of event $\{S = S(J^c) \, : J \in \cJ_{n,m}\}$ that only depends on the clean data $\cD_n$, with the idea that we are interested in what happens to the posterior on the retained (uncontaminated) set $S(\cO^c)$. In this case the inequality becomes
     \[ \widetilde{E}_{f_0} G(S(\cO^c)\, | \, \widetilde{\cD}_n) \le e^{m \frac{\lambda}{\sigma^2}(\tau_n^2 + 4F_0^2)}  E_{f_0} \left [ \max_{J \in \cJ_{n,m}}I_{J^c}^w(S(J^c)) / \min_{J \in \cJ_{n,m}} I_{J^c}^w(\cF_n) \, \mathbf{1}_\cA \right ] +{P}_{f_0}^n(\cA^c).\]
     This will be useful for proving Theorem \ref{thm : genweight fix} as well as Theorem \ref{thml2} below.
\end{Remark}

\subsection{Proofs for Section \ref{sec : thresh function}}

\subsubsection{Proof of Theorem \ref{thm : contconc thresh}}\label{proof : contconc thresh}
\begin{proof}
Consider the uncontaminated data $y_i = f_0(x_i) + \sigma \xi_i$ as in model \eqref{def : model} and define the event
\[
    \cB_n := \left \{ \max_{i=1,\dots,n} | \xi_i | \leq \sqrt{3 \log n}      \right \}.
\]
By Lemma \ref{lem : gaussian max}, we have that $ P_{f_0}^n(\cB^c_n) = o( 1/\sqrt{n})$. 
Note that, on this event, since $\| f_0 \|_{\infty} \leq F_0$, we have, for all $i=1,\dots,n$,
\begin{equation}
\label{bound obs in Bn}
    |y_i| \leq F_0 + \sigma |\xi_i | \leq F_0 + \sigma \max_{j=1,\dots,n} |\xi_j| \leq F_0 + \sigma \sqrt{3 \log n} \leq \tau_n,
\end{equation}
where the last inequality holds for $n$ large enough. Recall the notation $\cJ_{n,m} = \{ J \subset\{1,\dots,n\} \, : \, |J| \leq m \}$ and  the definition of the posterior integrals $I_{J^c}^w$ and $I_{J^c}^\lambda$ (for $w \equiv \lambda$) given in \eqref{def : post integral}, and define the event 
\[\cC_{n,m}:= \left\{ \min_{J \in \cJ_{n,m}} I_{J^c}^\lambda(\cF_n) > e^{-2m\log n - (1 + \frac{\lambda(1+\lambda)}{\si^2}) n \varepsilon_n^2} \right\},\]
with $\cF_n = \{ \| f - f_0 \|_\infty \le \veps_{n} \}$. According to Lemma \ref{lem : bound denom}, we have
\[ {P}_{f_0}^n(\cC_{n,m})\ge 1 - e^{-\frac{\lambda(1+\lambda)}{2 \si^2} n \veps_n^2}.\]

Define $ \cA := \cB_{n} \cap \cC_{n,m}$, this event depends only on the clean data and and  $P_{f_0}^n(\cA^c) =  o(1)$ as $n \to \infty$.
Applying Lemma \ref{lem : getting rid of the outliers} with $\cA$, we have, for $n$ large enough,
\[ \widetilde{E}_{f_0}G(S\, | \, \widetilde{\cD}_n) \le e^{m \frac{\lambda}{\si^2}(\tau_n^2 + 4F_0^2)} E_{f_0} \left [ \max_{J \in \cJ_{n,m}}I_{J^c}^w(S) / \min_{J \in \cJ_{n,m}} I_{J^c}^w(\cF_n) \, \mathbf{1}_\cA \right ] + o(1) .\]
Furthermore, for any set $B \subset \cF$, any $J \in \cJ_{n,m}$ and any $i \in J^c$, on the event $\cA$, since $|y_i|  \le \tau_n$ by \eqref{bound obs in Bn}, $w(y_i) = \lambda$ and thus $I_{J^c}^w(B) = I_{J^c}^\lambda(B)$ for all $J \in \cJ_{n,m}$. It follows that, on $\cA$, 
\[\min_{J \in \cJ_{n,m}}I_{J^c}^w(B) = \min_{J \in \cJ_{n,m}}I_{J^c}^\lambda(B) \qquad \text{and}\qquad \max_{J \in \cJ_{n,m}}I_{J^c}^w(B) = \max_{J \in \cJ_{n,m}}I_{J^c}^\la(B).\]
Using the lower bound for $\min_{J \in \cJ_{n,m}}I_{J^c}^\lambda(\cF_n)$ on the event $\cC_{n,m}$, one gets
\[ \widetilde{E}_{f_0}G(S\, | \, \widetilde{\cD}_n) \le e^{m[ 2\log n + \frac{\lambda}{\si^2}(\tau_n^2 + 4F_0^2)] + ( 1 + \frac{\lambda(1+\lambda)}{\si^2})n \varepsilon_n^2}E_{f_0} \left [\max_{J \in \cJ_{n,m}}I_{J^c}^\la(S) \, \mathbf{1}_{\cA} \right ] + o(1) .\]
Let $J \in \cJ_{n,m}$, using the tensorization $D_{\lambda}(P_f^{n-|J|},P_{f_0}^{n-|J|}) =(n-|J|)D_{\lambda}(P_f,P_{f_0})$ of the R\'enyi divergence  (recalled in Appendix \ref{app : renyi}), Fubini's theorem and $\mathbf{1}_\cA \le 1$, one obtains
\[E_{f_0}[ I_{J^c}^{\lambda}(S) \mathbf{1}_\cA] \le \int_S \int \left( \frac{dP_f^{n-|J|}}{dP_{f_0}^{n-|J|}}\right)^{\lambda} \, dP_{f_0}^{n-|J|} \, d\Pi(f) = \int_S e^{-(1-\lambda)(n-|J|)D_{\lambda}(P_f,P_{f_0})} d\Pi(f).\]
Taking $S := \{ f : (1- \eta)(1-\lambda)D_{\lambda}(P_f,P_{f_0}) \ge M  \zeta_{n,m}^2 \}$ and noting that $n-|J| \ge n(1- \eta)$, 
\[ E_{f_0}[ I_{J^c}^{\lambda}(S)\mathbf{1}_\cA] \le e^{-M \frac{(n-|J|)}{1 - \eta} \zeta_{n,m}^2} \le e^{-M n\zeta_{n,m}^2}.\]
This bound is independent of the index set $J$ and, since $|\cJ_{n,m}| =  \sum_{k=0}^m {n \choose k} \le e^{2 m \log n}$, 
\[ E_{f_0} \left [\max_{J \in \cJ_{n,m}}I_{J^c}^\la(S) \, \mathbf{1}_{\cA} \right ] \le \sum_{J \in \cJ_{n,m}}E_{f_0} [I_{J^c}^\la(S) \, \mathbf{1}_{\cA}] \le e^{2 m \log n  -M n\zeta_{n,m}^2}. \]
Using $4 \log n \le \tau_n^2/\si^2$ as well as $ 4 F_0^2 \le \tau_n^2$ for $n$ large enough,

\[ \widetilde{E}_{f_0} G(S\, | \, \widetilde{\cD}_n) \le e^{ \frac{1 + 2 \la}{\si^2}m\tau_n^2 + (1 + \frac{\lambda(1+ \lambda)}{\si^2})n \varepsilon_n^2} e^{-M n \zeta_{n,m}^2} + o(1) .\]
Taking $\zeta_{n,m}^2 = \frac{m}{n} \tau_n^2 + \varepsilon_n^2 $ and $M > \frac{1+2\lambda}{\si^2} \vee (1 + \frac{\lambda(1+ \lambda)}{\si^2}) $, we obtain the result.
\end{proof}

\begin{lem}\label{lem : bound denom}
    Let $\Pi$ be a prior on $\cF$ satisfying
    \[ \Pi(\cF_n) = \Pi[\{ f \, : \| f - f_0\|_\infty \le \varepsilon_{n} \}] \ge e^{-n \varepsilon_n^2} .\] 
    Setting $\cJ_{n,m} := \{ J \subset \{1,\dots, n\} \, : \, |J| \leq m\}$, with $P_{f_0}^n$--probability at least $1 - e^{-\frac{\lambda(1+\lambda)}{2 \si^2} n \veps_n^2}$,
    \[ \min_{J \in \cJ_{n,m}}I_{J^c}^{\lambda}(\cF_n) \ge  e^{- 2m \log n-\left( 1 + \frac{\lambda(1 + \lambda)}{\si^2}\right) n\veps_n^2}.\]
   
\end{lem}
\begin{proof}
     In case $m = 0$, we use the following conventions
    \[
        \cJ_{n,0} = \{ \emptyset \}, \qquad | \cJ_{n,0} | = 1, \qquad \min_{J \in \cJ_{n,0}} I^{w}_{J^c}(\cF_n) = I^{w}_{[n]}(\cF_n).
    \]
    Setting
    \[ \gamma_{n,m} : = 2 m \log n + \Big( 1 + \frac{\lambda(1+\lambda)}{ \si^2}\Big) n \varepsilon_n^2,\]
    a union bounds yields
    \[ P_{f_0}^n \left (\min_{J \in \cJ_{n,m}} I_{J^c}^\lambda(\cF_n) \le e^{-\gamma_{n,m}} \right ) \le \sum_{J\in \cJ_{n,m}} P_{f_0}^n ( I_{J^c}^\lambda(\cF_n) \le e^{-\gamma_{n,m}}). \]
    Let $J \in \cJ_{n,m}$, by the prior assumption, we have  $\Pi(\cF_n) \ge e^{-n \veps_{n}^2}>0$. Setting $\overline{\Pi} := \Pi|_{\cF_n}$, 
    \[ I_{J^c}^{\lambda}(\cF_n) \ge e^{- n \varepsilon_{n}^2} \int \exp\{ L_{J^c, \lambda}(f_0) - L_{J^c,\lambda}(f) \}\, d \overline{\Pi}(f) .\]
    Applying Jensen's inequality to $\overline{\Pi}$ with the convex map $x \mapsto x^{-1}$ and then Markov's inequality, 
    \begin{align*}
        P_{f_0}^n (I_{J^c}^\lambda(\cF_n) \le e^{-\gamma_{n,m}}) &\le  P_{f_0}^n \left(\int \exp\{ L_{J^c, \lambda}(f_0) - L_{J^c,\lambda}(f) \}\, d \overline{\Pi}(f) \le e^{n\varepsilon_n^2-\gamma_{n,m}} \right)\\
        & \le P_{f_0}^n \left(\int \exp\{ - L_{J^c,\lambda}(f_0) + L_{J^c,\lambda}(f) \}\, d \overline{\Pi}(f) \geq e^{-n\varepsilon_n^2+\gamma_{n,m}} \right)\\
        & \le e^{ n \veps_{n}^2 - \gamma_{n,m}} E_{f_0}\int \exp\{ -L_{J^c,\lambda}(f_0) + L_{J^c,\lambda}(f) \}\, d \overline{\Pi}(f).
    \end{align*} 
    On the generic set $J^c$ the data are i.i.d, by definition of the negative log--likelihood $L_{J^c, \lambda}(f)$,
    \[ \exp\{ -L_{J^c, \lambda}(f_0) + L_{J^c,\lambda}(f) \} = \left( \frac{dP_{f_0}^{n-|J|}}{dP_{f}^{n-|J|}} \right)^{\lambda} .\]
    Combining Fubini's theorem and the definition of the R\'enyi divergence \eqref{def : renyi} 
    yields
    \begin{align*} 
    E_{f_0} \int \left( \frac{dP_{f_0}^{n-|J|}}{dP_{f}^{n-|J|}} \right)^{\lambda}\, d \overline\Pi(f) 
    &= \int \int \left(\frac{dP_{f_0}^{n-|J|}}{dP_{f}^{n-|J|}} \right)^{\lambda} dP_{f_0}^{n-|J|} \, d\overline{\Pi}(f)\\
    &=\int \left ( \frac{dP_{f_0}^{n - |J|}}{d P_{f}^{n - |J|}} \right )^{\lambda + 1 } dP_{f}^{n- |J|}\, d\overline{\Pi}(f)\\
    &= \int e^{\lambda D_{1 + \lambda}(P_{f_0}^{n-|J|}, P_f^{n-|J|})}\, d \overline{\Pi}(f).
    \end{align*}
    Using the tensorization property of the R\'enyi divergence, and Lemma \ref{lem : L2 from renyi}, on the set $\cF_n$,
    \[ D_{1 + \lambda}(P_{f_0}^{n-|J|}, P_f^{n-|J|}) = (n-|J|)D_{1 + \lambda}(P_{f_0}, P_f) \le \frac{1 +\lambda}{2 \si^2}(n-|J|) \varepsilon_n^2 \leq  \frac{1 +\lambda}{2 \si^2} n \varepsilon_n^2. \]
    The expectation in the last but one display is upper bounded by $e^{{\lambda(1+\lambda)}({2 \si^2})^{-1} n \veps_{n}^2}$ and thus
     \[  P_{f_0}^n (I_{J^c}^\lambda(\cF_n) \le e^{-\gamma_{n,m}})\le e^{-\gamma_{n,m} +  n \veps_{n}^2 + \frac{\lambda(1+\lambda)}{2 \si^2} n \veps_{n}^2 } \le e^{- 2m \log n-\frac{\lambda(1+\lambda)}{2 \si^2} n \veps_{n}^2}. \]
    This holds for any $J \in \cJ_{n,m}$, so from the union bound above, one obtains
    \[ P_{f_0}^n \left (\min_{J \in \cJ_{n,m}} I_{J^c}^\lambda(\cF_n) \le e^{-\gamma_{n,m}} \right )\le |\cJ_{n,m}| e^{- 2 m \log n-\frac{\lambda(1+\lambda)}{2 \si^2} n \veps_{n}^2} \le e^{-\frac{\lambda(1+\lambda)}{2 \si^2} n \veps_{n}^2}, \]
    where we used $|\cJ_{n,m}| = \sum_{k=0}^m {n \choose k} \le  e^{2m\log n}$ in the last step, 
     completing the proof.
\end{proof}

\subsection{Proofs for Section \ref{section applications}}
\subsubsection{Proof of Proposition \ref{prop : prior mass GP}} \label{proof : prior mass GP}

\begin{proof}
The proof is based on the non-centered small ball inequality (Lemma $11.19$ of \cite{ghosal_van_der_vaart_2017}),
\[ \Pi^{(\gamma)}[ \|f - f_0\|_\infty \le \veps_n^*] \ge e^{-\varphi_{f_0}^{\gamma}(\veps_n^*/2)},\]
where the concentration function $\varphi^{\gamma}_{f_0}(\cdot)$ is defined as
\begin{equation}
\label{concentration function rescaled gp}
 \varphi_{f_0}^{\gamma}(\veps) := \inf_{h \in \mathbb{H}^{\gamma} \, : \|f_0-h\|_\infty \le \veps}\|h\|_{\mathbb{H}^{\gamma}}^2 - \log\Pi^{(\gamma)}[ \|f\|_\infty \le \veps].
\end{equation}
Here $\mathbb{H}^{\gamma}$ denotes the RKHS of the prior $\Pi^{(\gamma)}$. Using the results in Lemma \ref{small ball rescaled matern} and Lemma \ref{bias conc func rescaled matern}, we obtain 
\begin{equation}
\label{bound concentration function}
    \varphi_{f_0}^{\gamma}(\varepsilon) \lesssim \gamma^{-2 \alpha} \varepsilon^{- \frac{2 \alpha + d - 2 \beta}{\beta}} + \gamma^d \varepsilon^{- d / \alpha},
\end{equation}
where the condition $\gamma \geq 1$ is satisfied because $\alpha \geq \beta$ by assumption. Replacing in the above display the choices of $\gamma_n$ and $\varepsilon_n$, we have that
\[
    \varphi_{f_0}^{\gamma_n}(\varepsilon_n) \leq n \varepsilon_n^2 
\]
provided the constant $M$ in $\varepsilon_n$ is taken sufficiently large and $n$ is large enough.
\end{proof}

The following three lemmas describe the Reproducing Kernel Hilbert Space, the small deviation bound and the decentralization term associated to the rescaled Mat\'ern process. These results were derived for the 1-dimensional case in \cite{adaptationBotond}. We provide the proof for the general $d$-dimensional case below using the same arguments for completeness.

\begin{lem}[RKHS of the time-rescaled Mat\'ern process]\label{lem:rkhs}
Let $\alpha>0$ and $\gamma>0$, and let $W^{\gamma}$ be the time-rescaled
Mat\'ern process, with covariance function \eqref{eq: matern_rescaled}. The RKHS
$\mathbb{H}^{\gamma}$ of $W^{\gamma}$ is the set of functions $h:[0,1]^d\to\mathbb{R}$ of the form
$h(x)=g(\gamma x)$ for $g\in \cS^{\alpha+d/2}([0,\gamma]^d)$, with
\begin{equation}
\label{equivalence of norms}
    \|h\|_{\mathbb{H}^{\gamma}}\;\asymp\;\|g\|_{\cS^{\alpha+d/2}([0,\gamma]^d)},
\end{equation}
up to multiplicative constants depending on $\alpha$ and $d$ only. When $C_{\alpha,d}=c_{\alpha,d}=1$, the relation in  \eqref{equivalence of norms} is a (norm) identity. 
\end{lem}

\begin{proof}
 The argument here follows the one in Lemma $13$ in \cite{adaptationBotond}.  We take, for simplicity, $C_{\alpha,d}=c_{\alpha,d}=1$. By Lemma $11.35$ of \cite{ghosal_van_der_vaart_2017}, $\mathbb{H}^{\gamma}$ is the set of functions
\[
    h(x)=\int_{\mathbb{R}^d}e^{i\langle s,x\rangle}\psi(\lambda)\rho_{\alpha,\gamma}(s)\,ds,
    \qquad x\in[0,1]^d,
\]
with $\int|\psi|^2\rho_{\alpha,\gamma}<\infty$, the square norm $\|h\|^2_{\mathbb{H}^\gamma}$ being
the infimum of the latter integral over all $\psi$ that give the same function $h$ on
$[0,1]^d$. Since $\rho_{\alpha,\gamma}(s)=\tau^{-d}(1+\|s/\gamma\|^2)^{-\alpha-d/2}$,
the substitution $s=\tau\mu$,
together with $\phi(\mu):=\psi(\gamma\mu)(1+\|\mu\|^2)^{-\alpha-d/2}$, gives
\[
    h(x)=\int_{\mathbb{R}^d}e^{i\langle\mu,\gamma x\rangle}\phi(\mu)\,d\mu=g(\gamma x),
    \qquad
    \int_{\mathbb{R}^d}|\psi|^2\rho_{\alpha,\gamma}\,ds
    =\int_{\mathbb{R}^d}|\phi(\mu)|^2(1+\|\mu\|^2)^{\alpha+d/2}\,d\mu,
\]
where $g(y)=\int e^{i\langle\mu,y\rangle}\phi(\mu)\,d\mu$. When defined on all of $\mathbb{R}^d$, the functions $g$ arising in this way
are exactly the elements of $\cS^{\alpha+d/2}(\mathbb{R}^d)$, since $\phi$ is the Fourier
transform of $g$ up to normalisation and the required condition on $\phi$ is precisely the one in the definition of that (Sobolev) space.

Note that two functions $\psi_1$, $\psi_2$ provide the same $h$ on $[0,1]^d$ if and
only if the corresponding functions $\phi_1$, $\phi_2$ define $g_1$ and $g_2$, respectively, satisfying  $g_1 1_{[0,\gamma]^d}=g_2 1_{[0,\gamma]^d}$, and $\phi_j$, $j=1,2$ is the Fourier transform of exactly one $g_j\in \cS^{\alpha+d/2}(\mathbb{R}^d)$. Hence
\[
    \|h\|^2_{\mathbb{H}^\gamma}
    =\inf\bigl\{\|G\|^2_{\cS^{\alpha+d/2}(\mathbb{R}^d)}:\ G|_{[0,\gamma]^d}=g\bigr\}
    =\|g\|^2_{\cS^{\alpha+d/2}([0,\gamma]^d)},
\]
the infimum being over all extensions of $g$, which is the asserted norm identity.
\end{proof}

\begin{lem}{(Centered small ball probability)}
\label{small ball rescaled matern}
Let $W^{\gamma}$ be the time-rescaled Matérn process with covariance function \eqref{eq: matern_rescaled} with $\alpha + d/2 \in \mathbb{N}$, restricted to $[0,1]^d$. Then, for $\varepsilon < 1$, $\gamma \geq 1$ and a multiplicative constant depending only on $\alpha$ and $d$,
\[
   - \log \Pi^{(\gamma)}  \left ( f: \| f \|_{\infty} \leq \varepsilon  \right ) \lesssim \gamma ^{d} \varepsilon^{- d/ \alpha}.
\]
\end{lem}

\begin{proof}
    Since $\gamma \geq 1$, by definition of the spectral density of the process we have that
    \[
        \rho_{\alpha,\gamma}(s ) \leq C_{\alpha,d} \gamma^{2 \alpha}( c_{\alpha,d} + \| s\|^2)^{- \alpha - d/2} =: \rho_{\alpha,\gamma}^*(s),
    \]
    and $\rho_{\alpha,\gamma}^*(s)$ is $\gamma^{2 \alpha}$ times the spectral density of the Matérn process with $\gamma = 1$, i.e. the spectral density of the space-rescaled process $x \mapsto \gamma^{\alpha} W^{1}(x)$. By Lemma $11.36$ of \cite{ghosal_van_der_vaart_2017}, 
    \begin{equation}
    \label{space-rescaled small ball}
        - \log P \left( \sup_{x \in [0,1]^d} | \gamma^{\alpha} W^{1}(x) | \leq \varepsilon \right )  \lesssim (\varepsilon / \gamma^{\alpha})^{- d/\alpha} = \gamma^{d} \varepsilon^{-d / \alpha}.
    \end{equation}
    Since $\rho^{*}_{\alpha,\gamma} - \rho_{\alpha,\gamma} \geq 0$, by Bochner's theorem,  this difference function is the spectral density of a stationary Gaussian process $Z$, and $\gamma^{\alpha} W^{1} \overset{d}{=} W^{\gamma} + Z$ with $W^{\gamma}$ independent of $Z$. Anderson's lemma then gives $P( \sup_{x\in[0,1]^d} | \gamma^{\alpha}W^{1}(x) | \leq \varepsilon) \leq P(\sup_{x\in[0,1]^d} |W^{\gamma}(x)| \leq \varepsilon)$, which together with \eqref{space-rescaled small ball} proves the claim. 
\end{proof}

\begin{lem}
\label{bias conc func rescaled matern}
    Let $f_0 \in \mathcal{H}^{\beta}([0,1]^d)$ with $0 < \beta \leq \alpha$ and $\alpha + d/2 \in \mathbb{N}$. Then, for $\varepsilon \in (0,1)$, $\gamma \geq 1$ and $\gamma \lesssim \varepsilon^{- 1 / \beta}$,
    \[
        \inf_{h \in \mathbb{H}^{\gamma}: \| f_0 - h \|_{\infty} \leq \varepsilon} \| h \|_{\mathbb{H}^{\gamma}}^2 \lesssim \gamma^{d} + \gamma^{-2 \alpha} \varepsilon^{- \frac{2 \alpha + d - 2 \beta}{\beta}}.
    \]
\end{lem}

\begin{proof}
    By Whitney's theorem, a function $f_0\in \mathcal{H}^{\beta}([0,1]^d)$ can be extended to a function on
    $\mathbb{R}^d$ with compact support and H\"older norm bounded by a multiple of
$\|f_0\|_{\mathcal{H}^{\beta}([0,1]^d)}$, the multiple depending only on $\beta$ and $d$. By an abuse of notation, we will denote this function always by $f_0$. Let $\psi : \mathbb{R}^d \to \mathbb{R}$ be a kernel of order $k := \alpha + d/2$, i.e. $\int \psi = 1$, $\int s^{a} \psi(s) ds = 0$ for $1 \leq |a| \leq k$ and $\int \|s\|^{k+1} | \psi(s) | ds < \infty$. We use the following notation $\psi_{\sigma}(\cdot) = \sigma^{-d} \psi( \cdot / \sigma)$. The $d$-dimensional version of Lemma 11.31 of \cite{ghosal_van_der_vaart_2017} gives, for $\sigma \in (0,1)$ and all $a$ such that $|a| \leq k$, 
\begin{equation}
\label{approximation bounds}
    \sup_{x \in \mathbb{R}^d} | f_0 - f_0 \ast \psi_{\sigma} |(x) \lesssim \sigma^{\beta}, \quad \sup_{x \in \mathbb{R}^d} | D^{a}(f_0 \ast \psi_\sigma)(x) | \lesssim \left (  \frac{1}{\sigma} \right )^{( | a| - \beta) \vee 0}.
\end{equation}
Take $\sigma = \varepsilon^{1 / \beta}$ and define $h(x) = (f_0 \ast \psi_{\sigma})(x)$ for $x \in [0,1]^d$. We have $h = g(\gamma x)$ with $g(x) = (f_0 \ast \psi_{\sigma})(x / \gamma)$, $x \in [0,\gamma]^d$. Then $\| f_0 - h \|_{\infty} \lesssim \varepsilon$ and $h \in \mathbb{H}^{\gamma}$ by Lemma \ref{lem:rkhs}.

Changing variables $s = \gamma u$, 
\begin{equation}
\begin{aligned}
    \| h \|_{\mathbb{H}^{\gamma}}^2 & = \| g \|_{\cS^k([0,\gamma]^d)}^2 \leq \| f_0 \ast \psi_{\sigma}(\cdot / \gamma) \|^2_{\cS^k(\mathbb{R}^d)} = \\ & \sum_{| a| \leq k} \gamma^{d - 2 |a|} \int_{\mathbb{R}^d} | D^{a}(f_0 \ast \psi_{\sigma})(u) |^2 du \lesssim \sum_{\ell = 0}^k \gamma^{d - 2 \ell} \left ( \frac{1}{\sigma}\right )^{2 ( (\ell - \beta) \vee 0)},
\end{aligned}
\end{equation}
where in the last step we used the bound on the right hand side of \eqref{approximation bounds}. Let $b = \lfloor \beta \rfloor$. We split the sum of the previous display as follows
\[
    \sum_{\ell = 0}^k \gamma^{d - 2 \ell} \left ( \frac{1}{\sigma}\right )^{2 ( (\ell - \beta) \vee 0)} = \sum_{\ell \leq b} \gamma^{d - 2 \ell}  + \sum_{\ell > b} \gamma^{d - 2 \ell} \left ( \frac{1}{\sigma}\right )^{2 ( (\ell - \beta) \vee 0)},
\]
The first sum is bounded by $\gamma^d$ for $\gamma \geq 1$. In the second sum we have $\ell > b$ and consecutive terms have ratio $(\gamma \sigma)^{-2}$, which, under the assumption $\gamma \lesssim 1 /\sigma$, satisfies $(\gamma \sigma)^{-2} \gtrsim  1$. Thus, this sum is bounded by a multiple of its last term
\[
    \gamma^{d - 2 k} \sigma^{-2(k - \beta)} = \gamma^{- 2 \alpha} \sigma^{- (2 \alpha + d - 2 \beta)}.
\]
Substituting $\sigma = \varepsilon^{1 / \beta}$ gives the assertion.

\end{proof}

\subsubsection{Proof of Proposition \ref{prop:prior mass hierarchical GP}}\label{proof : prior mass hierarchical GP}

\begin{proof}
    The proof follows an argument similar to the proof of $(29)$ in \cite{adaptationBotond}. Define the sequence
    \[
        \gamma_n = n^{\frac{\alpha - \beta}{\alpha(2\beta + d)}},
    \]
    and note that $\gamma_n \geq 1$, since $\alpha \geq \beta$, so that the interval
    $[\gamma_n, 2\gamma_n]$ is contained in the support of $p_{\alpha,n}$. By Lemma $11.19$ of \cite{ghosal_van_der_vaart_2017},
    \begin{equation}
    \label{conditional lower bound}
        P\bigl( \| W^{\Gamma} - f_0 \|_{\infty} \leq \varepsilon_n \mid \Gamma = \gamma \bigr)
        \geq e^{-\varphi_{f_0}^{\gamma}(\varepsilon_n/2)},
    \end{equation}
    where we recall that $\varphi_{f_0}^{\gamma}(\cdot)$ is the concentration function of the
    time-rescaled process, defined in \eqref{concentration function rescaled gp}. Let
    $\gamma \in [\gamma_n, 2\gamma_n]$. Since $(\alpha - \beta)/\alpha < 1$, we have
    $2\gamma_n \leq (\varepsilon_n/2)^{-1/\beta}$ for $n$ large enough, so the
    condition $\gamma \lesssim \varepsilon^{-1/\beta}$ required by Lemma \ref{bias conc func rescaled matern} is satisfied. The bound in 
    \eqref{bound concentration function}, yields 
    \[
        \varphi_{f_0}^{\gamma}(\varepsilon_n/2)
        \lesssim\ \gamma_n^{-2\alpha}\varepsilon_n^{-\frac{2\alpha + d - 2\beta}{\beta}}
        + \gamma_n^{d}\varepsilon_n^{-d/\alpha}
        \asymp n^{\frac{d}{2\beta + d}}.
    \]
    The right hand side is bounded by $n\varepsilon_n^2/2$ once the constant
    $M$ in the definition of $\varepsilon_n$ is taken large enough. Combining the definition of the randomly time-rescaled process
    with \eqref{conditional lower bound} and the upper bound in the previous display, we obtain
    \begin{equation}
    \begin{aligned}
    \label{final prior mass}
        \Pi\bigl( f : \| f - f_0 \|_{\infty} \leq \varepsilon_n\bigr)
        & = P\bigl( \| W^{\Gamma} - f_0 \|_{\infty} \leq \varepsilon_n \bigr) \\
        & = \int_{1}^{\infty}
            P\bigl( \| W^{\Gamma} - f_0 \|_{\infty} \leq \varepsilon_n \mid \Gamma = \gamma \bigr)
            g_{\alpha,n}(\gamma)\, d\gamma \\
        & \geq \int_{\gamma_n}^{2\gamma_n} e^{-\varphi_{f_0}^{\gamma}(\varepsilon_n/2)}
            g_{\alpha,n}(\gamma)\, d\gamma \\
        & \geq e^{-n\varepsilon_n^2/2}\, c_1 \gamma_n
            \exp\left\{ -d_1 n^{\frac{d}{2\alpha + d}} (2\gamma_n)^{\frac{2\alpha d}{2\alpha + d}} \right\},
    \end{aligned}
    \end{equation}
    where the last step uses condition \eqref{eq:hyperprior}. Taking minus the
    logarithm of the right hand side of \eqref{final prior mass}, and using
    \[
        n^{\frac{d}{2\alpha + d}} \gamma_n^{\frac{2\alpha d}{2\alpha + d}}
        = n^{\frac{d}{2\beta + d}} = \frac{n\varepsilon_n^2}{M^2},
        \qquad
        \log\frac{1}{c_1\gamma_n} \leq \log\frac{1}{c_1},
    \]
    the latter because $\gamma_n \geq 1$, we obtain
    \[
        -\log \Pi\bigl[\{ f : \| f - f_0 \|_{\infty} \leq \varepsilon_n \}\bigr]
        \leq \frac{n\varepsilon_n^2}{2} + \log\frac{1}{c_1}
        + 2^{\frac{2\alpha d}{2\alpha + d}} d_1 \frac{n\varepsilon_n^2}{M^2}
        \leq n\varepsilon_n^2,
    \]
    for $n$ large enough, provided $M$ is chosen large enough. This concludes
    the proof.
\end{proof}

\subsubsection{Proof of Proposition \ref{prop : prior mass DNN}}\label{proof : prior mass DNN}
\begin{proof}
    We employ Theorem 1 of \cite{kongkim} (see also \cite{kohler_full} for the original construction) with $M_n = (n/ \log n)^{1/\{2(2\beta +d)\}}$, up to increasing the depth and width of the network by adding zeroes in the coefficients of the matrices and vectors and adding identity layers. There exists a DNN $\tilde{f}_0$ of width $N = \lceil n^{d/\{2(2\beta +d)\}}  \rceil$ and depth $L = \lceil \log^{1+\delta} n\rceil$, such that, for some constant $c_0$ and $n$ large enough.
    \[ \| f_0 - \tilde{f}_0 \|_\infty \le c_0 M_n^{-2\beta} = c_0 \varepsilon_n^* \log ^{\frac{\be}{2 \be + d}}n. \]
    Furthermore, all the coefficients of $\tilde{f}_0$ are bounded in absolute value by a constant $C >0$. Thus, for any $\varepsilon_n \ge \varepsilon_n^* \log n$, as $n$ gets large enough, $\| f_0 - \tilde{f}_0 \|_\infty \le \varepsilon_n/2$, the triangle inequality yields
    \[  \Pi[ \{ f \, : \|f - f_0\|_\infty \le \varepsilon_n \}] \ge \Pi[ \{ f \, : \|f - \tilde{f}_0\|_\infty \le \varepsilon_n/2 \}].\]
    Under the prior $\Pi$, the function $f$ is a DNN with same width and depth as $\tilde{f}_0$, hence in view of Lemma \ref{lem : propag dnn} below we obtain, by independence of the coefficients under the prior,
     \[  \Pi[ \{ f \, : \|f - f_0\|_\infty \le \varepsilon_n \}] \ge \prod_{\ell \ge 0} \prod_{i,j,k}\Pi \left( |W_\ell^{(i,j)} - \tilde{W}_\ell^{(i,j)}| \le \varepsilon_n/V \, , \, |b_\ell^{(k)} - \tilde{b}_\ell^{(k)}| \le \varepsilon_n/V\right),\]
     where $\tilde{W},\tilde{b}$ are the matrix/vectors of $\tilde{f}_0$ and
     \[ V : = 2(d+1)(N+1)^{L-1}(C\vee1)^L(L+1).\]
     In the remaining of the proof we focus on the coefficients $W_\ell^{(i,j)}$. The biases $b_\ell^{(k)}$ can be handled analogously. Also, and without loss of generality, assume $C>1$; notice that, because of the choice of $N$ and $L$, as $n$ gets large enough, we have
     \[ V \le N^L (C \vee 1)^L = (CN)^L.\]
     Therefore, (again only looking at $W_\ell^{(ij)}$) we have
     \[  \prod_{\ell \ge 0} \prod_{i,j}\Pi \left( |W_\ell^{(i,j)} - \tilde{W}_\ell^{(i,j)}| \le \varepsilon_n/V \right) \ge \prod_{\ell \ge 0} \prod_{i,j}\Pi \left( |W_\ell^{(i,j)} - \tilde{W}_\ell^{(i,j)}| \le \varepsilon_n(CN)^{-L}\right).\]
     The number of terms in the product above is $dN + N^2(L-1) + N \le N^2L$, as $n$ gets large enough. We obtain
     \[  \prod_{\ell \ge 0} \prod_{i,j}\Pi \left( |W_\ell^{(i,j)} - \tilde{W}_\ell^{(i,j)}| \le \varepsilon_n(CN)^{-L}\right)\ge \left[P \left( |w-\tilde{w}| \le \varepsilon_n(CN)^{-L}\right) \right]^{N^2L},\]
     where $w \sim \cN(0,\si_W^2)$ and $|\tilde{w}| \le C.$ Writing $\phi$ for the standard Gaussian density we have
     \[P \left( |w-\tilde{w}| \le \varepsilon_n(CN)^{-L}\right) = \int_{\tilde{w}-\varepsilon_n(CN)^{-L}}^{\tilde{w}+\varepsilon_n(CN)^{-L}} \si_W^{-1} \phi(w/\si_W) \, dw . \]
     The integrand being symmetric, we can assume $\tilde{w} >0$ without loss of generality and obtain the bound
     \[P \left( |w-\tilde{w}| \le \varepsilon_n(CN)^{-L}\right) \ge \frac{2 \varepsilon_n}{\si_W (CN)^L} \exp\bigg\{-\frac1{2 \si_W^2}(\tilde{w} + \varepsilon_n(CN)^{-L})^2 \bigg\}. \]
     Notice that as $n$ gets large enough we both have $\varepsilon_n(CN)^{-L} \le 1$ and $\varepsilon_n \ge e^{-L}$, along with $\tilde{w} \le C$ this shows existence of a constant $C'>C>1$  such that
     \[ P \left( |w-\tilde{w}| \le \varepsilon_n(CN)^{-L}\right) \ge (C'N)^{-L} \ge N^{-2L}.\]
     Plugging this into the product form of the bound yields, for some $C'' >0$
     \[ \prod_{\ell \ge 0} \prod_{i,j}\Pi \left( |W_\ell^{(i,j)} - \tilde{W}_\ell^{(i,j)}| \le \varepsilon_n(CN)^{-L}\right) \ge N^{-2 N^2L^2} \ge e^{-C'' n (\varepsilon_n^*)^2 \log^{2(1+\delta)+1} n} \ge e^{-\frac12 n \varepsilon_n^2},\]
     where we used $N^2 \lesssim n (\varepsilon_n^*)^2$, $L^2\lesssim \log^{2(1+\delta)}n$ as well as $\varepsilon_n \ge \varepsilon_n^* \log^{2+\delta}n$ in the last step. The exact same bounds can be obtained for the product over the bias vectors and thus
     \[ \Pi(\{ f\, : \, \|f -f_0\|_\infty \le \veps_n \}) \ge e^{-n\veps_n^2}. \qedhere\]
\end{proof}

\vspace{.1cm}

\begin{lem}[Lemma 24 in \cite{JMLR:v26:24-0894}]\label{lem : propag dnn}
Let $f,\tilde{f}$ be two DNNs on $[0,1]^d$ as in \eqref{def : dnn} with width $N$ and depth $L$. Suppose that all their coefficients are bounded by $C>0$ and that for all $\ell \ge 0$ and $(i,j,k)$,
\[ |W_\ell^{(i,j)} - \tilde{W}_\ell^{(i,j)}| \le \delta \qquad \text{and} \qquad|b_\ell^{(k)} - \tilde{b}_\ell^{(k)}| \le \delta. \]
Then, 
\[ \|f-f^*\|_\infty \le \delta (d+1)(N+1)^{L-1}(C\vee1)^L(L+1).\]
\end{lem}

\subsection{Proofs for Section \ref{sec : gen weight}}

\subsubsection{Proof of Theorem \ref{thm : genweight fix}}\label{proof : genweight fix} 
\begin{proof}
We follow the steps of the proof of Theorem \ref{thm : contconc thresh}. For $J \in \cJ_{n,m}$ defined in \eqref{def : jnm}, set
\[ S_{n,m}(J^c) := \left\{ f :  \frac{1}{n - |J|} \sum_{i \in J^c} |f(x_i)-f_0(x_i)|^2 \ge M \frac{\zeta_{n,m}^2}{1 - \eta}\right\} ,\]
and recall $m/n \le \eta<1$.
We aim to show that the expected posterior mass of $S_{n,m}(\cO^c)$ vanishes. For an event $\cA$ depending only on the clean data, to be chosen below, an application of Lemma \ref{lem : getting rid of the outliers} yields (see also remark \ref{rmk : J dependent S}), for $\cF_n = \{ \| f - f_0 \|_\infty \le \veps_n \}$, as $n$ gets large enough,
\begin{equation}\label{eq : proof gen w first step}
    \widetilde{E}_{f_0}G(S_{n,m}(\cO^c)\, | \, \widetilde{\cD}_n) \le e^{m \frac{\lambda}{\si^2}(\tau_n^2 + 4F_0^2)} E_{f_0} \left [ \max_{J \in \cJ_{n,m}}I_{J^c}^w(S_{n,m}(J^c)) / \min_{J \in \cJ_{n,m}} I_{J^c}^w(\cF_n) \, \mathbf{1}_\cA \right ] + {P}_{f_0}^n(\cA^c).
\end{equation}
Here we recall that for any subset of indices $\cI \subset \{1, \dots ,n\},$
\[L_{\cI,w}(f_0) - L_{\cI,w}(f) = \frac{1}{2 \si^2} \sum_{i \in \cI} w(y_i) [(f_0(x_i) - y_i)^2 -(f(x_i) - y_i)^2].\]
For any $x \in \R^d$, denote 
\begin{align}
a(x) = a(f,f_0,x):= \sigma^{-1}\{ f(x) - f_0(x)\}\label{def:a(x)}
\end{align}
and notice
\[\sigma^{-2}[(f_0(x) - y)^2 -(f(x) - y)^2] = -a(x)^2 + 2a(x) \sigma^{-1}(y-f_0(x)).\]
Therefore, for any $S \subset \cF$, with $y_i-f_0(x_i) = \si\xi_i$, we have
\[I_{\cI}^w(S) = \int_S \exp\left\{ -\frac12 \sum_{i\in \cI}w(y_i)\left[a(x_i)^2 - 2  a(x_i)\xi_i\right] \right\} \, d\Pi(f). \]
Applying Lemma \ref{lem : denom prob proof gen weight}, there exists an event $\cC_{n,m}$ with $P_{f_0}^{n}(\cC_{n,m}) \to 1$, such that on $\cC_{n,m}$,
\[  \min_{J \in \cJ_{n,m}} I_{J^c}^w(\cF_n) \ge e^{- 2 m \log n-\left(\frac{\lambda(1+\la)}{\sigma^2}  + 2\right)n\varepsilon_{n}^2 - \frac{\Delta_n \rho_n}{\sigma^2}n \varepsilon_{n}}.\]
Recall that $\rho_n > 2\si\sqrt{\log n}$ and set the event $\cA$ to be chosen in \eqref{eq : proof gen w first step} as
\[\cA := \cC_{n,m} \cap \left \{ \max_{1 \le i \le n}|\xi_i| \le \rho_n/\si \right \}, \qquad \text{so that} \qquad  P_{f_0}^{n}(\cA^c) \to 0,\]
where the right hand side of the previous display follows from Lemma \ref{lem : gaussian max} and a union bound. 

This yields, since $4 F_0^2 \leq \tau_n^2 / 2$ for $n$ large enough,
\[
\widetilde{E}_{f_0}G(S_{n,m}(\cO^c) | \, \widetilde{\cD}_n) \le e^{m \frac{3\lambda}{ 2\si^2}\tau_n^2 } e^{2 m \log n} e^{(\frac{\lambda(1+\la)}{\sigma^2}  + 2)n\varepsilon_{n}^2 + \frac{\Delta_n \rho_n}{\sigma^2}n \varepsilon_{n}}E_{f_0} \left [ \max_{J \in \cJ_{n,m}} {I_{J^c}^w(S_{n,m}(J^c))} \, \mathbf{1}_{\cA} \right] + o(1).
\]
For any $J \in \cJ_{n,m}$, by Fubini's theorem, we have
\[
    E_{f_0} \left [I_{J^c}^w(S_{n,m}(J^c)) \, \mathbf{1}_{\cA} \right] \leq E_{X} \int_{S_{n,m}(J^c)} E_{Y \mid X} \exp\left\{ -\frac12 \sum_{i\in J^c}w(y_i)\left[a(x_i)^2 - 2  a(x_i)\xi_i\right] \right\} \mathbf{1}_{\cA} \, d\Pi(f).
\]
We apply the first inequality of Lemma \ref{lem : proof gen weight bound on E} with $\cI =  J^c$, $b_i = f_0(x_i)$ and $a_i = a(x_i)$ to the inner expectation of the previous display to obtain
\begin{equation}\label{proof : genweight fixed design step}
    E_{f_0}[ I_{J^c}^w(S_{n,m}(J^c)) \mathbf{1}_{\cA}] \leq E_{X} \int_{S_{n,m}(J^c)}\exp \left\{- \frac{ w_0(1-w_0)}{4} \sum_{i \in J^c} a_i^2 + \frac{\Delta_n \rho_n}{\sigma} \sum_{i \in J^c} |a_i| \right\} \, d\Pi(f).
\end{equation}
Using Cauchy-Schwarz inequality, we obtain
\[ E_{f_0}[ I_{J^c}^w(S_{n,m}(J^c)) \mathbf{1}_{\cA}] \leq E_X \int_{S_{n,m}(J^c)}\exp\left\{- \frac{ w_0(1-w_0)}{4} \sum_{i\in J^c} a_i^2 + \sqrt{n-|J|} \,\frac{\Delta_n \rho_n}{\sigma} \left (\sum_{i \in J^c} a_i^2\right )^{\frac12} \right\} \, d\Pi(f).\]
Let $\al_w := w_0(1-w_0)/4$, with the inequality $\sqrt{xy} \leq \al_w x/2 + y/(2\al_w) $, we further bound
\[\exp\left\{- \al_w \sum_{i \in J^c} a_i^2 + \sqrt{n-|J|} \,\frac{\Delta_n \rho_n}{\sigma} \left (\sum_{i \in J^c} a_i^2 \right)^{\frac12} \right\} \le \exp\left\{- \frac{ \al_w}{2} \sum_{i\in J^c} a_i^2 + (n-|J|) \frac{(\Delta_n \rho_n)^2}{2\al_w \sigma^2} \right\} .\]
Recalling $a_i := \sigma^{-1}\{f(x_i) - f_0(x_i)\}$ as well as $|J|/n \le m/n \leq  \eta < 1$ and 
\[ S_{n,m}(J^c) =  \left \{ f \, : \, \frac{1}{n-|J|} \sum_{i \in J^c} |f(x_i) - f_0(x_i)|^2 \geq  M \frac{\zeta_{n,m}^2}{1 - \eta} \right \} ,\]
we get, {with $n-|J|\ge n-m=(1-\eta)n$,}
\[ E_{f_0}[ I_{J^c}^w(S_{n,m}(J^c)) \mathbf{1}_{\cA}] \leq e^{(n-|J|) \frac{(\Delta_n \rho_n)^2}{2\al_w \sigma^2}} e^{- \frac{\al_w M}{2\si^2} n \zeta_{n,m}^2} \le e^{n\frac{(\Delta_n \rho_n)^2}{2\al_w \sigma^2}} e^{- \frac{\al_w M}{2\si^2} n \zeta_{n,m}^2}.\]
Moreover, as in the proof of Theorem \ref{thm : contconc thresh}, bounding the maximum by the sum,
\[
    E_{f_0} \left [ \max_{J \in \cJ_{n,m}} {I_{J^c}^w(S_{n,m}(J^c))} \, \mathbf{1}_{\cA} \right] \leq \sum_{J \in \cJ_{n,m}} E_{f_0}[ I_{J^c}^w(S_{n,m}(J^c)) \mathbf{1}_{\cA} ].
\]
Gathering all the previous bounds together and using $|\cJ_{n,m} | \leq e^{2 m \log n}$ yields
\[ \widetilde{E}_{f_0}G(S_{n,m}(\cO^c) \, | \, \widetilde{\cD}_n) \le e^{m( 4 \log n + \frac{3\lambda}{2\si^2}\tau_n^2)} e^{(\frac{\lambda(1+\la)}{\sigma^2} + 2) n\veps_n^2 + \frac{\Delta_n \rho_n}{\si^2}n\veps_n +  n\frac{(\Delta_n\rho_n)^2}{2\al_w \sigma^2 }} e^{- \frac{\al_w M}{2 \si^2} n \zeta_{n,m}^2} + o(1).\]
Using $\tau_n^2/\si^2 \ge 4 \log n$ and $\Delta_n \rho_n \veps_n/\si^2 \le 1/2 \cdot (\Delta_n\rho_n/\si)^2 + \veps_n^2/(2\si^2)$, we get
\[ \widetilde{E}_{f_0}G(S_{n,m}(\cO^c) \, | \, \widetilde{\cD}_n) \le e^{  \frac{2+3\lambda}{2\si^2}m\tau_n^2} e^{(\frac{1/2+\lambda(1+\la)}{\sigma^2} + 2) n\veps_n^2 +  (\frac{1}{2} + \frac{1}{2\al_w })\si^{-2} n(\Delta_n \rho_n)^2} e^{- \frac{\al_w M}{2 \si^2} n \zeta_{n,m}^2} + o(1).\]
Taking $M$ large enough such that $\frac{ \al_w M}{2 \sigma^2 } > ( \frac{2 +3\la}{2\si^2})  \vee (2 + \frac{1/2+\lambda(1+\la)}{ \sigma^2}) \vee (\frac{1}{2} + \frac{1}{2 \al_w})\sigma^{-2}$ and setting $\zeta_{n,m}^2 =  \frac{m}{n} \tau_n^2 + \varepsilon_n^2 + (\Delta_n\rho_n)^2$, the last display is a $o(1)$, as desired.
\end{proof}
\begin{lem}\label{lem : denom prob proof gen weight}
    Let $\cF_n := \{ f \, : \, \|f - f_0\|_\infty \le \veps_n \}$ and consider a prior $\Pi$ such that $\Pi(\cF_n) \ge e^{- n \veps_n^2 }$. Assume $w$ is a weight function satisfying \eqref{eq : cond w cont}--\eqref{eq : cond lowerb w}--\eqref{conddiffw} with $\Delta_n = o(1)$ and $\rho_n \ge 2 \si\sqrt{\log n}$. Define the set $\cJ_{n,m}:= \{ J \subset [n] \ : \ |J| \leq m\}$. Then, as $n \to \infty$,
    \[P_{f_0}^n\left[ \min_{J \in \cJ_{n,m}} I_{J^{c}}^w(\cF_{n}) \geq e^{- 2 m \log n-\left(\frac{\lambda(1+\la)}{\sigma^2}  + 2\right)n\varepsilon_{n}^2 - \frac{\Delta_n \rho_n}{\sigma^2}n \varepsilon_{n}}\right] \to 1.\]
\end{lem}

\begin{proof}
    This proof is similar to that of Lemma \ref{lem : bound denom} with the difference that now the likelihood is weighted by $w$ and not simply tempered by $\lambda$. In case $m = 0$, we use the conventions
    \[
        \cJ_{n,0} = \{ \emptyset \}, \quad | \cJ_{n,0} | = 1, \quad \min_{J \in \cJ_{n,0}} I^{w}_{J^c}(\cF_n) = I^{w}_{[n]}(\cF_n).
    \]
    Let us set
    \[ \gamma_{n} := 2 m \log n + \left(\frac{\lambda(1+\la)}{\sigma^2} + 2\right)n\varepsilon_{n}^2 + \frac{\Delta_n \rho_n}{\sigma^2}n \varepsilon_{n} \quad \text{and} \quad \cA_{n,m}:= \left \{ \min_{J \in \cJ_{n,m}} I_{J^c}^w(\cF_n) \ge e^{-\gamma_n} \right \}.\]
    We aim to show $P_{f_0}^{n}(\cA_{n,m}) \to 1$ as $n \to \infty$,
    Recall that by assumption we have $\Pi(\cF_n) \geq e^{-n\varepsilon_{n}^2} > 0$. Denoting $\overline{\Pi} := \Pi|_{\cF_n}$, we get 
    \[ I_{J^c}^w(\cF_n) \geq e^{-n \varepsilon_n^2} \int_{\cF} \exp\left\{ -\frac12 \sum_{i \in J^c}w(y_i)\left[a(x_i)^2 - 2  a(x_i)\xi_i\right] \right\}\, d\overline{\Pi}(f), \]
    where $a(x)$ was defined in \eqref{def:a(x)}. As in the proof of Lemma \ref{lem : bound denom}, combining the previous display with a union bound and applying Jensen's inequality to $x \mapsto x^{-1}$, we obtain
    \begin{equation}
    \label{to apply markov}
         P_{f_0}^n(\cA_{n,m}^c) \le \sum_{J \in \cJ_{n,m}}P_{f_0}^{n}\left[\int_{\cF_n} \exp\left\{ \frac12 \sum_{i \in J^c} w(y_i)\left[a(x_i)^2 - 2  a(x_i)\xi_i\right] \right\}\, d\overline{\Pi}(f) \ge e^{\gamma_n - n \varepsilon_n^2} \right].
    \end{equation} 
    We bound the following quantity uniformly over all $J \in \cJ_{n,m}$
    \[ \cE_{den} :=E_{f_0} \int_{\cF_n} \exp\left\{ \frac12 \sum_{i \in J^c} w(y_i)\left[a(x_i)^2 - 2  a(x_i)\xi_i\right] \right\}\, d\overline{\Pi}(f). \]
    By Fubini's theorem, we obtain 
    \[ \cE_{den} = \int_{\cF_n} E_X E_{Y \mid X} \exp\left\{\frac12 \sum_{i \in J^c} w(f_0(x_i) + \si \xi_i)\left[a(x_i)^2 - 2a(x_i)\xi_i \right]\right\}\, d\overline{\Pi}(f) . \]
    We apply the second inequality of Lemma \ref{lem : proof gen weight bound on E} (conditionally on the design points), with $ \cI = J^c$, $b_i = f_0(x_i)$ and $a_i = a(x_i)$, to get
    \[ \cE_{den} \leq \int_{\cF_n} E_X  \left \{ e^{{\la(1+\la)}  \sum_{i \in J^c} |a(x_i)|^2} \prod_{i \in J^c} \left[ e^{\frac{\Delta_n \rho_n}{\si} |a(x_i)|} + 2e^{-\frac{\rho_n^2}{4\si^2}}\right]  \right \} \, d \overline{\Pi}(f).\]
    Restricted on $\cF_{n}$, we have $|a(x_i)| = \si^{-1}|f(x_i) - f_0(x_i)|\le \varepsilon_{n}/\sigma$, thus
     \[ \cE_{den} \leq e^{\frac{\la(1+\la)}{\si^2} (n-|J|)\varepsilon_{n}^2}  \left[ e^{ \frac{\Delta_n \rho_n}{\si^2} \veps_n} + 2e^{-\rho_n^2/(4\si^2)} \right]^{n-|J|}.\]
     Using $e^{ \frac{\Delta_n \rho_n}{\si^2} \varepsilon_{n}}\ge1$, along with the inequality $\log(1+x) \leq x$, for $x>0$, this yields
     \[\left[ e^{ \frac{\Delta_n \rho_n}{\si^2} \varepsilon_{n}} + 2e^{-\rho_n^2/(4\si^2)} \right]^{n-|J|} \leq e^{\frac{\Delta_n \rho_n}{\si^2} (n-|J|) \varepsilon_{n}}\left[1 + 2e^{-\frac{\rho_n^2}{4\si^2}}\right]^{n-|J|}\leq e^{ \frac{\Delta_n \rho_n}{\si^2}n \varepsilon_{n}}\times e^{2ne^{-\rho_n^2/(4\si^2)}}.\]
     Since $\rho_n \ge 2\si\sqrt{\log n}$, we have
$2n e^{-\rho_n^2/(4\si^2)} \le 2$. Hence, uniformly over
$J \in \cJ_{n,m}$, using $n-|J|\le n$,
\[
\cE_{den}
\le e^2 \exp\left\{
\frac{\lambda(1+\la)}{\si^2}n\varepsilon_n^2
+\frac{\Delta_n\rho_n}{\si^2}n\varepsilon_n
\right\}.
\]
Consequently, applying Markov's inequality to each probability on the right-hand side of \eqref{to apply markov}, using the bound in the previous display and the fact that $|\cJ_{n,m}| \leq e^{2 m \log n}$, one gets
\begin{equation*}
P_{f_0}^n(\cA_{n,m}^c)
\le |\cJ_{n,m}|e^{n\varepsilon_n^2-\gamma_n}
e^2\exp\left\{
\frac{\lambda(1+\la)}{\si^2}n\varepsilon_n^2
+\frac{\Delta_n\rho_n}{\si^2}n\varepsilon_n
\right\}\\
\le e^{2-n\varepsilon_n^2}
\to 0,
\end{equation*}
which concludes the proof.
\end{proof}

    \begin{lem} \label{lem : proof gen weight bound on E}
        Let $\si >0$ and $w$ be a weight function satisfying \eqref{eq : cond w cont}--\eqref{eq : cond lowerb w}--\eqref{conddiffw}. Let $\cI \subset \{1, \dots, n\}$ with $\Delta_n=o(1)$; for any $i \in\cI$, let $a_i \in \R$, $|b_i| \le F_0$ and $\xi_i \overset{iid}{\sim} \cN(0,1)$. Define
        \[ \psi_i := \frac12 w(b_i + \si \xi_i)(a_i^2 - 2a_i\xi_i). \]
        Then, for $n$ large enough
        \begin{align*}
            &E(  e^{- \sum_{i \in \cI} \psi_i} \mathbf{1}_{\max_i |\xi_i| \leq \rho_n/\si}) \leq  e^{ \frac{-w_0(1-w_0)}{4} \sum_{i \in \cI} a_i^2 + \frac{\Delta_n \rho_n}{\si}\sum_{i \in \cI} |a_i|}, \\  
             & E ( e^{\sum_{i \in \cI} \psi_i}) \le  e^{\la(1+\la)\sum_{i \in \cI} a_i^2} \prod_{i \in \cI} \left[ e^{\frac{\Delta_n\rho_n}{\si} |a_i|} + 2e^{-\frac{\rho_n^2}{4\si^2}}\right]. 
        \end{align*}
    \end{lem}
    \begin{proof}
        First define 
        \[ \psi^*_i := \frac12 w(b_i)(a_i^2 - 2a_i\xi_i).\]
        Thanks to $|b_i| \le F_0$ and the property \eqref{conddiffw} of $w$, we have for $|\xi_i| \le \rho_n/\si$
        \begin{equation}\label{eq : diff in psi proof gen weight}
            |\psi_i - \psi^*_i| \leq \Delta_n \left( \frac12 a_i^2 + |\xi_i| \, |a_i| \right).
        \end{equation} 
        Also, using \eqref{eq : cond lowerb w} as well as independence (and Gaussianity) of the $\xi_i$, we obtain
        \begin{align}
            \label{eq : proof sumpsi2}E\, e^{- \sum_{i \in \cI} \psi_i^*} &= e^{-\frac{1}{2} \sum_{i \in \cI} w(b_i)a_i^2} E\, e^{ \sum_{i \in \cI} w(b_i)a_i \xi_i} =e^{\frac{1}{2} \sum_{i \in \cI} a_i^2 (-w(b_i) + w(b_i)^2)} \le e^{-\frac{w_0(1-w_0)}{2} \sum_{i \in \cI} a_i^2}.
        \end{align}
       Combining inequalities \eqref{eq : proof sumpsi2} and \eqref{eq : diff in psi proof gen weight} as well as $\Delta_n = o(1)$, we obtain, for large enough $n$,
        \begin{align*}
        E(  e^{- \sum_{i \in \cI} \psi_i} \mathbf{1}_{ \max_i |\xi_i| \leq \rho_n/\si }) 
        &\leq  e^{ \frac{\Delta_n - w_0(1-w_0)}{2} \sum_{i \in \cI} a_i^2 + \frac{\Delta_n \rho_n}{\si} \sum_{i \in \cI}|a_i| }\\
        &\leq e^{ \frac{- w_0(1-w_0)}{4} \sum_{i \in \cI} a_i^2 + \frac{\Delta_n \rho_n}{\si} \sum_{i \in \cI} |a_i| }. 
        \end{align*}
        For the other bound, using independence, one gets
        \[E \, e^{\sum_{i \in \cI} \psi_i} = \prod_{i \in \cI} \left[ E\, e^{\psi_i} \mathbf{1}_{|\xi_i| \leq \rho_n/\si} 
        + E\, e^{\psi_i} \mathbf{1}_{|\xi_i| > \rho_n/\si}\right].\]
        Noticing that, thanks to \eqref{eq : cond w cont}, we have $E \, e^{\psi_i^*} \leq e^{\frac{\lambda +\lambda^2}{2} a_i^2}$, the first term in the product can be bounded, using  first \eqref{eq : diff in psi proof gen weight} and then $\Delta_n = o(1)$, by \[E\, e^{\psi_i} \mathbf{1}_{|\xi_i| \leq \rho_n/\si } \leq e^{\frac{\lambda +\lambda^2 + \Delta_n}{2}a_i^2 + \frac{\Delta_n \rho_n}{\si} |a_i|} \le e^{\la(1+\la)a_i^2 + \frac{\Delta_n \rho_n}{\si} |a_i|}.   \] 
        For the second term, notice first $\psi_i \leq \lambda a_i^2 /2 + \lambda |a_i||\xi_i|$. Then using Cauchy-Schwarz inequality, the Gaussian tail bound $P(|\xi_i | > t) \leq 2 e^{-t^2/2}$ for $t > 0$ and the bound on the moment generating function $E e^{t |\xi_i|} \leq 2 e^{t^2/2}$ ($t \in \mathbb{R}$), leads to 
        \begin{equation}
        \begin{aligned}
            E\, e^{\psi_i} \mathbf{1}_{|\xi_i| > \rho_n/\si} & \leq e^{\frac{\lambda}{2}a_i^2}E\,( e^{\lambda |a_i||\xi_i|}\mathbf{1}_{|\xi_i| > \rho_n / \sigma } ) \\ &  \leq e^{\frac{\lambda}{2} a_i^2} \sqrt{ E e^{2 \lambda |a_i| |\xi_i| }} \cdot \sqrt{E \mathbf{1}_{|\xi_i| > \rho_n/\sigma}} \\ & \leq
             e^{\frac{\lambda}{2} a_i^2} \cdot \sqrt{2} e^{\lambda^2 a_i^2} \cdot \sqrt{2} e^{- \frac{\rho_n^2}{4 \sigma^2}} = 2 e^{(\frac{\lambda}{2} + \lambda^2) a_i^2} e^{- \frac{\rho_n^2}{4 \sigma^2}}.
        \end{aligned}
        \end{equation} 
        Using $\lambda/2 \le \lambda$ and combining the previous bound we obtain
        \[E \, e^{\sum_{i \in \cI} \psi_i} \le \prod_{i \in \cI} e^{\la(1+\la)a_i^2}\left[ e^{\frac{\Delta_n \rho_n}{\si} |a_i|} + 2e^{-\frac{\rho_n^2}{4\si^2}}\right].\]
    \end{proof}

\subsubsection{Proof of Theorem \ref{thml2}} \label{proof : thml2}
Before proving the main theorem, we derive some additional technical results. Recall the definitions of the sets $\cJ_{n,m}$ and $\cO = \cO(T,U)$, see \eqref{def : jnm} and \eqref{def : O}, respectively. 
 
\begin{prop} \label{priormasslemma}
Let $\{B_n(J^c) : J \in \cJ_{n,m} \}$ be a family of sets, where, for each $J$, $B_n(J)$ is measurable with respect to the clean design points $x_i$. Suppose one works under the conditions of Theorem \ref{thml2} and that
\begin{equation} \label{pmformula}  
\max_{J \in \cJ_{n,m}} E_{X}\Pi[B_n(J^c)] = o\left( e^{-d_1 m\ta_n^2 - d_2 n\veps_n^2- d_3n(\Delta_n\rho_n)^2} \right),
\end{equation}
with $d_1= \frac{1 + 3 \lambda/2}{\sigma^2}, d_2=\frac{1/2+\lambda(1+\la)}{ \sigma^2} + 2 $ and $d_3=\frac{1}{2 \alpha_w \sigma^2} + \frac{1}{2\sigma^2}$, where we recall that $\alpha_w : = w_0( 1 - w_0)/4$.
Then, as $n\to\infty$,
\[ \widetilde{E}_{f_0} G\left(B_n(\cO^c) \given \widetilde{\cD}_n\right)=o(1). \]
\end{prop}
\begin{proof}  
We follow the proof of Theorem \ref{thm : genweight fix}. First, the proof of Lemma \ref{lem : getting rid of the outliers} gives (see also Remark \ref{rmk : J dependent S}), for any event $\cA$ depending only on the clean data, 
\[ \widetilde{E}_{f_0} [G(B_n(\cO^c) \given \widetilde{\cD}_n)] \le e^{m \frac{\lambda}{\si^2}(\tau_n^2 + 4F_0^2)}E_{f_0} \left [  \max_{J \in \cJ_{n,m}} I_{J^c}^{w}(B_n(J^c)) / \min_{J \in \cJ_{n,m}} I_{J^c}^{w}(\cF_n) \mathbf{1}_{\cA}  \right ] + {P}_{f_0}^n(\cA^c).\]
Next, in view of Lemma \ref{lem : denom prob proof gen weight}, on an event $\cC_{n,m}$, with $P_{f_0}^n(\cC_{n,m}^c) = o(1)$, 
\begin{align*}
\min_{J \in \cJ_{n,m}} I_{J^c}^w(\cF_{n}) &\ge e^{- 2 m \log n-\left(\frac{\lambda(1+\la)}{\sigma^2}  + 2 \right)n\varepsilon_{n}^2 - \frac{ \Delta_n\rho_n}{\si^2}n\veps_n}\\
&\geq  e^{- 2 m \log n-\left(\frac{1/2+\lambda(1+\la)}{\sigma^2}  + 2 \right)n\varepsilon_{n}^2 -\frac{1}{2\si^2}n (\Delta_n\rho_n)^2},
\end{align*}
where we used
$\Delta_n\rho_n \varepsilon_n/\si^2 \le 1/2 \cdot  (\Delta_n \rho_n/\si)^2 + \varepsilon_n^2/(2\si^2)$.
Then, consider the event
\[\cA := \cC_{n,m} \cap \left \{ \max_{1 \le i \le n}|\xi_i| \le \rho_n/\si\right \}, \qquad \text{satisfying} \qquad P_{f_0}^{n}(\cA) \to 1.\]
By combining the above displays, for $n$ large enough and a constant $D_1>0$,
\begin{equation}
\label{inter step bound sup} \widetilde{E}_{f_0}G(B_n(\cO^c) \, | \, \widetilde{\cD}_n) \le e^{D_1(m \tau_n^2  + n\varepsilon_{n}^2 + n(\Delta_n\rho_n)^2)} E_{f_0} \left [ \max_{J \in \cJ_{n,m}}{I_{J^c}^w(B_n(J^c)} \, \mathbf{1}_{\cA} \right ] + o(1).
\end{equation}

One then proceeds along the proof of Theorem \ref{thm : genweight fix}, and, for each $J \in \cJ_{n,m}$, obtains
\begin{align*} E_{Y\given X}\left[  I_{J^c}^w(B_n(J^c)) \mathbf{1}_{\cA}\right] 
&\le \int_{B_n(J^c)} \exp\left\{ -\frac{\alpha_w}{2} \sum_{i \in J^c} a_i^2+ \frac{n(\Delta_n \rho_n)^2}{2\al_w \sigma^2}\right\}d\Pi(f)\\
&\le \Pi[B_n(J^c)]e^{n(\Delta_n \rho_n)^2/(2\al_w \sigma^2)},
\end{align*}
with $a_i=\sigma^{-1}[f(x_i)-f_0(x_i)]$ and $\alpha_w=w_0(1-w_0)/4$, where we have bounded the negative term in the last exponential by $0$. 
Furthermore, one can bound the maximum in \eqref{inter step bound sup} by the sum over $J \in \cJ_{n,m}$. 
Then taking the expectation of the previous display with respect to $E_X$ and next bounding the sum by the maximum of the terms 
times the number of subsets  $J \in \cJ_{n,m}$, one obtains, for constants $d_1,d_2,d_3> 0$ as in the statement, 
\[  \widetilde{E}_{f_0}G(B_n(\cO^c) \, | \, \widetilde{\cD}_n) \le e^{d_1m\ta_n^2+ d_2n\veps_n^2 + d_3n (\Delta_n\rho_n)^2} \max_{J \in \cJ_{n,m}} E_X\Pi[B_n(J^c)] + o(1),  \]
which using the assumption of the proposition concludes the proof.
\end{proof}

Suppose the prior $\Pi$ is a zero--mean Gaussian process taking values in $(\cC^0([0,1]^d),\|\cdot\|_\infty)$, the Banach space of continuous functions on $[0,1]^d$, with RKHS $\mh$ and concentration function at $f_0$ denoted by $\vphi_{f_0}(\cdot)$ given in \eqref{concentration function rescaled gp}.  Consider, for all $J \in \cJ_{n,m}$, the empirical $L_2$--norm at design points indexed by $J^c$,
\[
    \| h \|_{n, J^c}^2 := \frac{1}{n - |J|} \sum_{i \in J^c} h(x_i)^2.
\]
The following is partly inspired by the proof of Theorem 2 in \cite{JMLRvdvvz} (the proof below slightly revisits the argument to ensure one applies Gaussian process arguments from \cite{van2008reproducing} within a separable Banach space).

\begin{prop} \label{interl2}
Assume the density of the design is bounded below by $c>0$ on $[0,1]^d$. Suppose the Gaussian process prior $\Pi$ as above gives probability $1$ to $\alpha$--H\"older functions $\cH^\alpha$ on $[0,1]^d$ and suppose $\vphi_{f_0}(\veps_n/2)\le n\veps_n^2$ holds. 

For $J \in \cJ_{n,m}$, $M \geq 2$ and $1\geq v_n \geq \varepsilon_n$, let $B_n(J^c)$ be the (design--dependent) set 
\[ B_n(J^c) :=\left\{  f\in \cC^0([0,1]^d):\ \|f-f_0\|_{2,P_X} > 2Mv_n\ge 8 \|f-f_0\|_{n, J^c} \right\}. \]
Then for any $\ga <\al$, there exists a constant $C > 0$ such that, for any $r_n\ge 1 \vee \sqrt{n} \varepsilon_n $  and $L\geq 1$  with $L r_n \geq M v_n $, one has
\[ \max_{J \in \cJ_{n,m}} E_X\Pi\left[ B_n(J^c) \right] \le  \exp\left\{-C (n-m) 
((M/L)v_n/r_n)^{2d/(2\ga+d)}  \right\} + \exp\left\{-2 L^2r_n^2\right\}. \]
\end{prop}

\begin{proof} 
Let us fix $\ga <\alpha$ and suppose that $\ga\notin\mathbb{N}$. In that case, the H\"older space $\cH^\ga$ coincides with the Besov space $B^{\gamma}_{\infty,\infty}$ (see Equation \eqref{eq : besov inf inf ball}) with equivalent norms. The case $\ga\in\mathbb{N}$ is treated at the end of the proof. By definition of the concentration function $\vphi_{f_0}$, there exists $g_n\in\mh$ such that 
\begin{equation}\label{construction approximation gn} 
\|g_n-f_0\|_\infty \le \veps_n/2,\qquad \|g_n\|_\mh^2 \le n\veps_n^2. 
\end{equation}

We first state a number of properties used throughout the proof. First, we shall use the fact that the norm $\|f\|_\ga$ is sub-Gaussian, where $\|\cdot\|_\ga$ denotes  the norm in $\cH^\ga$, for $f$ drawn from the Gaussian process as in the statement. More precisely, there exists a large enough constant $K_\gamma>0$  such that for large enough $x>0$ (the constant $2$ is for numerical convenience)
\begin{equation} \label{borellga}
\Pi\left[ \|f\|_{\ga} > K_\gamma x\right] \le e^{-2x^2}. 
\end{equation}

 Second,  for the function $g_n$ chosen in \eqref{construction approximation gn}, we shall use the bound
\begin{equation} \label{rkhsga}
 \|g_n\|_{\ga} \leqa \sqrt{n}\veps_n,
\end{equation}
so that $g_n\in \cH^{\ga}$. The verification that \eqref{borellga} and \eqref{rkhsga} hold is postponed at the end of the proof.
 Third, we will also need the interpolation-type inequality, for any $f\in \cH^\ga([0,1]^d)$,
 \begin{equation} \label{interpga}
  \|f\|_\infty \leqa \|f\|_\ga^{d/(2\ga+d)} \|f\|_2^{2\ga/(2\ga+d)}. 
 \end{equation}
This is Lemma 15 in  \cite{JMLRvdvvz}, which holds for covariates $X$ with density bounded from below by $c$, up to the fact that we have replaced the norm $\|\cdot\|_{\ga\given \infty}$ therein (which by definition is the norm of the Besov space $B^{\gamma}_{\infty,\infty}$) by the $\|\cdot\|_\ga$ norm, which one can do up to a constant, since as noted above both norms are equivalent when $\ga\notin\mathbb{N}$. 

Let us  introduce the set, for $L\ge 1$, 
\begin{equation} \cG_n := \left\{f:\ \|f\|_\ga \le K_\gamma L r_n \right\},\label{def:Gn}
\end{equation}
where $K_\gamma$ is the constant in \eqref{borellga} above. Using \eqref{borellga} with $x=Lr_n$ yields
\begin{align}
 \Pi[ \cG_n^c] & \le e^{-2 L^2 r_n^2}. \label{sievehoelder} 
\end{align}

Since $v_n\ge \veps_n$, the bounds in \eqref{construction approximation gn} also hold if $\veps_n$ is replaced by $v_n$. 
Next by the triangle inequality and writing $\|\cdot\|_2$ as a shorthand for $\|\cdot\|_{2,P_X}$, one has $\|f-f_0\|_2\le \|f-g_n\|_2+v_n/2$ and $\|f-f_0\|_{n,J^c}\ge \|f-g_n\|_{n,J^c}-v_n/2$, so that, assuming $M\ge 2$, for $J \in \cJ_{n,m}$
\begin{align*} 
 \Pi[B_n(J^c)]  & = \Pi[ \|f-f_0\|_2> 2Mv_n\,,\, 2\|f-f_0\|_{n,J^c}\le Mv_n/2 ] \\
 &\le \Pi[ \|f-f_0\|_2> (M+1/2)v_n\,,\, 2\|f-f_0\|_{n,J^c} \le (M-1)v_n] \\
 & \le \Pi[ \|f-g_n\|_2> Mv_n\ge  2\|f-g_n\|_{n,J^c}].
\end{align*} 
By definition of $B_n(J^c)$ as in the statement, and $\cG_n$ as in \eqref{def:Gn}, the following inclusion holds
\[ B_n(J^c) \subset \{f\in \cG_n:\ \|f-g_n\|_2>M v_n\, ,\, \|f-g_n\|_2\ge 2\|f-g_n\|_{n,J^c}\} \cup \cG_n^c. \] 
Taking first the expectation under $\Pi$ followed by the expectation under $E_{X}$ and next using Fubini's theorem gives, for each $J \in \cJ_{n,m}$
\begin{align*}    
E_{X} \Pi[B_n(J^c)]  & \le \int_{f\in \cG_n:\ \|f-g_n\|_2>Mv_n} P^n_{X}\left[\|f-g_n\|_2\ge 2\|f-g_n\|_{n, J^c} \right]d\Pi(f)
+  \Pi[ \cG_n^c] \\
 &\le \int_{f\in\cG_n:\ \|f-g_n\|_2>Mv_n} \exp\left\{-\frac{n - m}5 \left(\frac{\|f-g_n\|_2^2}{\|f-g_n\|_\infty^2} \right) \right\} d\Pi(f)+  \Pi[ \cG_n^c],
 \end{align*}
where the last line follows from Bernstein's inequality, applied to the  variables $(f-g_n)(X_i)^2$ which are independent, each with expectation $\|f-g_n\|_2^2$ and variance bounded from above by $E_X[(f-g_n)(X_1)^4]\le \|f-g_n\|_2^2 \|f-g_n\|_\infty^2$. Now we apply the interpolation inequality \eqref{interpga} to control the supremum norm. Since $g_n$ has finite $\|\cdot\|_\ga$--norm by \eqref{rkhsga}, one can use \eqref{interpga} with the function $f-g_n$ to get, for a constant $c_1=c_1(\ga)$,
\[ \|f-g_n\|_\infty \le c_1 \|f-g_n\|_\ga^{d/(2\ga+d)} \|f-g_n\|_2^{2\ga/(2\ga+d)}. \]
Writing  $C_\gamma$ for the constant in \eqref{rkhsga}, we have that $\| f - g_n \|_\gamma \leq \| f \|_\gamma + \| g_n \|_\gamma \leq (K_\gamma L + C_\gamma) r_n \leq (K_\gamma + C_\gamma) L r_n$ which holds since $L \geq 1$. 
Combining with the above bounds shows that $E_{X} \Pi[B_n(J^c)]$ is bounded from above by
\[ \int_{f\in\cG_n:\ \|f-g_n\|_2>Mv_n} \exp\left\{-c_2 (n-m) (Lr_n)^{-2d/(2\ga+d)} 
\|f-g_n\|_2^{2d/(2\ga+d)}  \right\} d\Pi(f)+  \Pi[ \cG_n^c],\]
for some constant $c_2=c_2(\ga) > 0$.
Bounding the exponential in the last display using $\|f-g_n\|_2>Mv_n$ yields the result.

It now remains to prove \eqref{borellga}--\eqref{rkhsga}. Starting with \eqref{borellga}, we can use the same argument as in the proof of Proposition 2 in (the supplement of) \cite{cr15}, slightly adapted  to deal with $[0,1]^d$ (instead of $[0,1]$). We include it here for completeness. 
By assumption the Gaussian prior $\Pi$ sits on $\cH^\ga$ for any $\ga<\al$. On the other hand, for  $\delta>0$ we have the continuous embeddings $\cH^{\ga+2\delta} \hookrightarrow B^{\ga + \delta}_{d/\delta,1} \hookrightarrow B^{\ga}_{\infty,\infty}$ and as noted above $B^{\ga}_{\infty,\infty}=\cH^\ga$ since $\ga\notin\mathbb{N}$. Since $\Pi$ sits on $\cH^{\ga+2\delta}$ for $\delta$ small enough, $\Pi$ can be seen as the law of a Gaussian random element in the separable Banach space $B^{\ga + \delta}_{d / \delta,1}$. By Borell's inequality (e.g. Corollary 5.1 in \cite{van2008reproducing}), inequality \eqref{borellga} holds with the norm $\|\cdot\|_{\ga+\delta,d / \delta,1}$ of $B^{\ga + \delta}_{d / \delta,1}$. Since the latter norm dominates the $\|\cdot\|_\ga$--norm by the embeddings noted above, it also holds as stated with $\|\cdot\|_\ga$ (with a possibly different constant $K_\gamma$), which proves \eqref{borellga}.

Turning to \eqref{rkhsga}, let  $(\mathbb{H}',\|\cdot\|_{\mathbb{H}'})$ denote the RKHS of $\Pi$ viewed as the law of a Gaussian random element in $(\mathbb{B}',\|\cdot\|_{\mathbb{B}'})$ with $\mathbb{B}':=B^{\ga + \delta}_{d/\delta,1}$. Recall that above we denoted $\mathbb{H}$ for the RKHS of $\Pi$ viewed as the law of a Gaussian random element in $\cC^0$. By Lemma 8.1 in \cite{van2008reproducing}  both RKHS coincide $\mathbb{H}=\mathbb{H}'$ as well as their norms $\|\cdot\|_{\mathbb{H}}=\|\cdot\|_{\mathbb{H}'}$ (by Lemma 7.1 in \cite{van2008reproducing}). On the other hand,  the RKHS norm $\|\cdot\|_{\mathbb{H}'}$ always dominates the norm of the corresponding (separable) Banach space the Gaussian process lives on ($\|\cdot\|_{\mathbb{B}'}$ in this case), see Eq. (2.5) in \cite{van2008reproducing}. 
Combining  these bounds gives
\[ \|g_n\|_{\mathbb{B}'} \leqa \|g_n\|_{\mathbb{H}'} = \|g_n\|_{\mathbb{H}}\le \sqrt{n}\veps_n,\]
where the last inequality follows from \eqref{construction approximation gn}. Now we use again $\mathbb{B}':=B^{\ga +\delta}_{d/\delta,1} \hookrightarrow B^{\ga}_{\infty,\infty}=\cH^\ga$ to get $\|g_n\|_{\ga}\leqa \|g_n\|_{\mathbb{B}'}$, which together with the last display yields \eqref{rkhsga}.

Finally, we deal with the case that $\ga\in\mathbb{N}$. In that case, one does the reasoning above with $\ga'\notin\mathbb{N}$ instead of $\ga$, with $\ga<\ga'<\al$. Since $(M/L)v_n/r_n \leq 1$ by the assumed bounds on the quantities involved, the obtained upper bound in terms of $\ga'$ in the last display of the statement is smaller than the same bound with $\ga'$ replaced by $\ga$, hence the result.
\end{proof}

\begin{Corollary} \label{corl2}
Setting $v_n^2 =  \veps_n^2+\ta_n^2(m/n) + (\Delta_n\rho_n)^2$ and  $r_n^2=nv_n^2$ 
 in Proposition \ref{interl2}, one gets for $n$ large enough and for any $0 < \gamma < \alpha$,
\begin{equation*}
\begin{aligned}
\max_{J \in \cJ_{n,m}} E_X \Pi\left[B_n(J^c)\right]
&\le \exp\left\{-C (1-\eta) M^{d/(2\ga+d)} n^{2\ga/(2\ga+d)}\right\} \\
&\quad \quad  + \exp\left\{-2M\left(m\tau_n^2+n\veps_n^2+n(\Delta_n\rho_n)^2\right)\right\}.
\end{aligned}
\end{equation*}
\end{Corollary}
We are now equipped to prove the main theorem on $L^2(P_X)$-contraction.
\begin{proof}[Proof of Theorem \ref{thml2}]
Define 
$v_n^2 := \varepsilon_n^2 + \tau_n^2 (m/n) + (\Delta_n\rho_n)^2$.
By Theorem \ref{thm : genweight fix} with $\zeta_{n,m}^2 = v_n^2 $, for suitably large $M_1>0$, we have $\widetilde{E}_{f_0} G [\|f-f_0\|_{n,\cO^c}>M_1 v_n  \given \widetilde{\cD}_n ] \to 0$. 

Note that the sequence $\veps_n$ verifies $\vphi_{f_0}(\veps_n / 2 )\le n\veps_n^2$. Hence, we can apply Proposition \ref{interl2} and in particular Corollary \ref{corl2}, to obtain, for any $\vartheta < \alpha$,
\begin{equation} 
\begin{aligned}
\label{remaining mass condition}
\max_{J \in \cJ_{n,m}} &E_{X}\Pi\left[  \|f-f_0\|_{2,P_X} > 2Mv_n\ge 8 \|f-f_0\|_{n, J^c}  \right] \leq \\ &  \exp\left\{-C  M^{d /(2 \vartheta +d)} ( 1 - \eta) n^{2\vartheta/(2\vartheta +d)} \right\} + \exp\left\{-2M n \zeta_{n,m}^2\right\}, 
\end{aligned}
\end{equation}
where $M \geq 2$.
By assumption we have that $ n \zeta_{n,m}^2 \le n^{2\vartheta /(2\vartheta+d)}$ for some $\vartheta < \alpha$, so that the last display is bounded by $2\exp\{-D n \zeta_{n,m}^2)\}$, where we can take the constant $D > 0$ larger than the constants $d_1,d_2,d_3$ in the statement of Proposition \ref{priormasslemma} provided $M$ is taken large enough.  Thus, applying Proposition \ref{priormasslemma}, one obtains
\[  \widetilde{E}_{f_0} G\left[ f:\ \|f-f_0\|_{2,P_X} > 2Mv_n\ge 8 \|f-f_0\|_{n,\cO^c}  \given \widetilde{\cD}_n \right]
=o(1). \] 
On the other hand, for large enough $M>0$,  the  application of Theorem \ref{thm : genweight fix} at the beginning of the proof gives
\[\widetilde{E}_{f_0} G\left[ f:\  2Mv_n\ge 8 \|f-f_0\|_{n,\cO^c}  \given \widetilde{\cD}_n \right]
=1-o(1).\]
Since $\zeta_{n,m}^2 = v_n^2$, putting the last two displays together and changing the value of the constant $M$ yields $\widetilde{E}_{f_0} G\left[f:\ \|f-f_0\|_{2,P_X} \ge M\zeta_{n,m} \given \widetilde{\cD}_n \right] \to 0$ as desired.
\end{proof} 

\subsubsection{Supremum-norm consistency} \label{sec : consistency sup norm}

\begin{thm} \label{thm : sup norm consistency}
Under the notation and conditions of Theorem \ref{thml2}, 
 for $L_n\ge 1$, $\vartheta <\al$, $\vartheta \notin \mathbb{N}$, let us set 
 \begin{equation} \label{snr}
 R_n = R_n(\be, \vartheta):= 2^{L_n d/2}\zeta_{n,m} + 2^{-\be L_n} + \sqrt{n}\zeta_{n,m}2^{-\vartheta L_n}.
 \end{equation}
 In case $\be\in\mathbb{N}$, suppose further that $f_0\in B^{\beta}_{\infty,\infty}$. 
 Then, for large enough $M>0$, as $n\to \infty$,
  \[ \widetilde{E}_{f_0} G\left[f:\ \|f-f_0\|_\infty \ge M R_n \given \widetilde{\cD}_n \right] \to 0. \]
\end{thm}
\begin{proof}
The proof uses a similar argument as the one of the proof of Proposition 2 in \cite{cr15}, adapted to the present contaminated setting. We include the full argument for completeness. For a continuous function $g$, let us denote by $g_{lk}=\psg g,\psi_{lk} \psd_2$ its wavelet coefficients into a CDV basis of $[0,1]^d$ (see Section \ref{sec : sup-norm trimmed} for the case $d=1$, see also Chapter 4 of \cite{GineNickl2016})
Under the assumed conditions, since $\cH^\be=B^{\beta}_{\infty,\infty}$ when $\be\notin\mathbb{N}$, we always have $f_0\in B^{\beta}_{\infty,\infty}$, and this is assumed when $\be\in\mathbb{N}$. In particular, $|f_{0,lk}|\le L 2^{-l(d/2+\be)}$ for large enough $L$ by definition of the $\|\cdot\|_{\be,\infty,\infty}$ norm. For $f$ a (smooth enough) function with coefficients $(f_{lk})$,
\begin{align*}
\|f& -f_0\|_\infty  = \| \sum_{l} \sum_{k} (f_{lk}-f_{0,lk})\psi_{lk} \|_\infty \\
& \leqa \sum_{l\le L_n} 2^{ld/2} \max_k |f_{lk}-f_{0,lk}| + \sum_{l> L_n} 2^{ld/2} \max_k |f_{lk}| 
+ \sum_{l> L_n} 2^{ld/2} \max_k |f_{0,lk}|\\
& \quad =: \qquad  \qquad(I)  \qquad \quad\ \ + \qquad \qquad (II)  \quad \quad+ \qquad \qquad (III). 
\end{align*}
By using the smoothness condition on $f_0$ in terms of wavelet coefficients, we have $(III)\leqa 2^{-\be L_n}$.
The term (I) is bounded using Cauchy-Schwarz inequality and the fact that a maximum of nonnegative terms is always bounded by its sum:
\[ (I)\leqa 2^{L_n d/2} \Big\{\sum_{l\le L_n}\sum_k (f_{lk}-f_{0,lk})^2\Big\}^{1/2} \leqa 2^{L_n d/2} \|f-f_0\|_2\leqa 2^{L_n d /2} \|f-f_0\|_{2,P_X},\]
where we have used that the standard $L_2$-norm is bounded from above by the $\|\cdot\|_{2,P_X}$-norm using the condition on the design. Then, by Theorem \ref{thml2}, the $\|f-f_0\|_{2,P_X}$ term in the last display is bounded by $M\zeta_{n,m}$ under the generalised posterior for large enough $M>0$. 

It now remains to bound the term (II).  Similarly as for term (III), we have 
\[ (II)\leqa 2^{-L_n\vartheta} \max_{l,k} 2^{l(\vartheta+d/2)} |f_{lk}| \leqa 2^{-L_n\vartheta} \|f\|_{\vartheta},\]
where the last inequality uses that the norms in $B^{\vartheta}_{\infty,\infty}$ and $\mathcal{H}^\vartheta$ are equivalent for $\vartheta\notin\mathbb{N}$. To conclude it suffices to show that the generalised posterior concentrates in balls for the $\|\cdot\|_\vartheta$-norm of radius of order $\sqrt{n}\zeta_{n,m}$. Recall from the proof of Proposition \ref{interl2} that a Borell-type inequality holds for $\|f\|_\vartheta$ if $f$ is drawn according to $\Pi$, see \eqref{borellga} therein. Using \eqref{borellga} with $x=L\sqrt{n}\zeta_{n,m}$ for large enough $L$ and combining with Proposition \ref{priormasslemma} (in the case where the sets do not depend on $J$, see also remark \ref{rmk : J dependent S})  yields $\widetilde{E}_{f_0} G\left[f:\ \|f\|_\vartheta \ge M' \sqrt{n}\zeta_{n,m} \given \widetilde{\cD}_n \right] \to 0$ for large enough $M'$, as desired, which concludes the proof of Theorem \ref{thm : sup norm consistency}.
\end{proof}

\subsubsection{Proof of Proposition \ref{prop : expl weight}}\label{proof : expl weight}
\begin{proof}
 Notice that in the first two cases $h$ is differentiable. Hence, by the mean--value theorem, to obtain \eqref{conddiffw} it suffices to show that, for $|y| \le \rho_n$, $|s| \le F_0$ and some $z \in [s,s+y],$
    \[ |w(s+y) - w(s)| \le \rho_n|w'(z)| \le \Delta_n.\]
    \begin{enumerate}
        \item Case $h(y) = \exp\{- \tau^2 / y^2 \}.$ Set $u := z/\tau$, direct computations yield
        \begin{equation*}
        \begin{aligned}
        |w'(z)| = \frac{\frac{2 \lambda}{\tau} |u+\frac{1}{u}| e^{-1/u^2}}{(1 + u^2 e^{-1/u^2})^2} & \le \frac{2\lambda}{\tau} \left|\frac{z}{\tau} + \frac{\tau}{z}\right| e^{- \frac{\tau^2}{z^2}} \\ & \le  2 \lambda \left(1 + \frac{4\rho_n^2}{\tau^2}\right) \frac{e^{- \tau^2/z^2}}{|z|} \\ & \le \frac{2 \lambda}{\tau}\left(1 + \frac{4\rho_n^2}{\tau^2}\right) e^{- \frac{\tau}{2 \rho_n}},
        \end{aligned}
        \end{equation*}
        where we used $|z| \le  |s| \vee |s + y| \le F_0 + \rho_n \le 2\rho_n$ for large enough $\rho_n$, and $e^{-1/x^2}/|x| \le e^{-1/|x|}$ in the last step. Taking $\rho_n = 2\si\sqrt{\log n}$ and $\tau_n = 2 \si\log^{1/2 +\delta} n = \rho_n \log^{\delta}n$ for some $\delta >1$, one has, for large enough $n$
        \[ \rho_n |w'(z)|  \leq \frac{2 \lambda}{\log^{\delta}n} \left ( 1 + \frac{4}{\log^{2 \delta}n} \right ) e^{-\log^{\delta} (n)/2} =: \Delta_n.
        \]
        Then $ \rho_n \Delta_n = o(n^{-1/2})$.
    \item Case $h = 1$. We obtain $w'(z) = -2 \lambda z \tau^2 / (\tau^2 + z^2)^2$ and thus when $\rho_n$ is large enough
        \[ \rho_n|w'(z)| \le \frac{2\lambda}{\tau^2}\rho_n |z| \leq \frac{2\lambda}{\tau^2}\rho_n (\rho_n + F_0) \le \frac{4\lambda}{\tau^2}\rho_n^2 =: \Delta_n.\]
    \item For the thresholding weight $w(y) = \lambda/[1\vee (y/\tau_n)^2]$ for $|y| \le \rho_n = 2 \si \sqrt{\log n}$ and $|s| \le F_0$, we have, as $n$ gets sufficiently large,
    \[ |s+y| \le 2\si\sqrt{\log n} + F_0 \le 3 \si \sqrt{\log n} = \tau_n\qquad \text{and} \qquad |s| \le F_0 \le \tau_n.\]
    Therefore, as $n$ gets large enough, $w(s+y) = w(s) = \lambda$ and thus $\Delta_n = 0$.
    \end{enumerate}
\end{proof}

\subsection{Proofs for Section \ref{sec : debiased}}
\subsubsection{Proof of Lemma \ref{lem : compute risk minimizer}}\label{proof : compute risk minimizer}
\begin{proof}
 Since the loss $\ell_f$ in \eqref{eq : general loss} depends on $f$ only through its value $f(x)$ at the point $x$,
conditioning on the design gives
\[
    R(f,f_0) = \int_{\R^d} \phi_x\big(f(x)\big)\, dP_X(x),
    \qquad \phi_x(r) := E_\xi\big[w(y)\,(t(y)-r)^2\big],
\]
with $y = f_0(x) + \sigma\xi$ and $\xi \sim \cN(0,1)$. For each fixed $x$, the map $\phi_x$
is quadratic in $r$,
\[
    \phi_x(r) = A_x\,r^2 - 2 B_x\, r + C_x,
    \quad A_x := E_\xi[w(y)],\ \ B_x := E_\xi[w(y)t(y)],\ \ C_x := E_\xi[w(y)t(y)^2],
\]
which is well defined since $A_x < \lambda$ by \eqref{eq : cond w cont}, $C_x \le C$ by assumption and $|B_x| \le (A_x C_x)^{1/2}$ by Cauchy--Schwarz. Its leading coefficient $A_x$ is strictly positive because $w > 0$. Completing the square, we obtain,
\[
\phi_x(r) = A_x\,(r - B_x/A_x)^2 + (C_x - B_x^2/A_x),
\]
so $\phi_x$ has the unique minimizer $r = B_x/A_x$. Let us denote 
\[
    g(x) := \frac{E_\xi[w(y)t(y)]}{E_\xi[w(y)]} = \frac{B_x}{A_x}.
\]
By assumption $g$ is upper bounded by $\sqrt{\lambda C}/c$ and thus belongs to $L_2(P_X)$. Integrating the last but one display with respect to $P_X$ shows that, for every $f$ with $R(f,f_0) < \infty$,
\begin{equation}\label{eq : square decomp}
    R(f,f_0) - R(g,f_0) = \int_{\R^d} A_x\,\big(f(x)-g(x)\big)^2\, dP_X(x) \ \ge\ 0 .
\end{equation}
Since $ \la \ge A_x \ge c >0$, the right-hand side of the previous display vanishes if and only if $f = g$ in a $L_2(P_X)$ sense (by equivalence of norms). Hence $g$ is the unique minimizer of $R(\cdot,f_0)$ over all of $L_2(P_X)$. Since $\cF$ is a closed convex subset of $L_2(P_X)$ and the two norms are equivalent, $\mathcal{F}$ is also a closed convex subset of $L_2(A_{\cdot} \cdot P_X)$. Hence, by the Hilbert projection theorem, we have $f^* = \mathfrak{p}_\cF (g)$.

With the choice of transform $t$ in \eqref{score matching alta loss}, we have $w(y)t(y) = w(y)y - \sigma^2 w'(y)$ and thus
\[ g(x) = f_0(x)
    + \sigma\,\frac{E_{\xi}[\xi\, w(y)]}{E_{\xi}[w(y)]}
   -  \sigma^2\,\frac{E_{\xi}[w'(y)]}{E_{\xi}[w(y)]} =: f_0(x) + (I) - (II).
\]
The correction $(\mathrm{II})$ cancels the bias $(\mathrm{I})$ exactly. Indeed, we assumed $E_{\xi}[|w'(y)|] < \infty$ and thus Stein's
identity $E_{\xi}[\xi\, \psi(\xi)] = E_{\xi}[\psi'(\xi)]$ applied to
$\psi(\xi) = w(f_0(x)+\sigma\xi)$, with $\psi'(\xi) = \sigma\, w'(f_0(x)+\sigma\xi)$, gives
$E_{\xi}[\xi\, w(y)] = \sigma\, E_{\xi}[w'(y)]$. This leads to $(\mathrm{I}) = (\mathrm{II})$ and since $f_0 \in \mathcal{F}$, we have
\[
    f^* = \mathfrak{p}_\cF (g) = \mathfrak{p}_\cF (f_0) = f_0,
\]
which concludes the proof.

\end{proof}

\subsubsection{Proof of Theorem \ref{thm : debiased}}\label{proof : debiased}

\begin{proof}

Write  $\|\cdot\|_{2} := \|\cdot\|_{2,P_X}$ and set
\[
    S_{n,m} := \{ f : \| f - f_0 \|_{2}^2 \geq M \zeta_{n,m}^2 \}, \qquad
    \mathcal{F}_n := \{ f : \| f - f_0 \|_{\infty} \leq \varepsilon_n \},
\]
and, for measurable $S \subset \mathcal{F}$, let
\[
    I_{[n]}^{w,\mathrm{de}}(S) := \int_{S} \exp\Big\{ -\frac{1}{2\si^{2}} \sum_{i=1}^n (\ell_f - \ell_{f_0})(x_i,y_i) \Big\}\, d\Pi(f)
\]
be the (generalized) integral of \eqref{def : post integral}, the superscript 'de'
indicating that we use the weighted debiased loss \eqref{eq : general loss}, \eqref{score matching alta loss}.

Since $\Pi$ assigns probability one to $\{\|f\|_\infty \le F\}$, using the bound
\eqref{bound on w t^2} {below (see Lemma \ref{lem : proof risk equivalent L2}) on the function $w\cdot t^2$} together with $w\le\lambda$ from \eqref{eq : cond w cont}, we have, for all $(x,y)$,
\[
    \ell_f(x,y) = w(y)\big(t(y) - f(x)\big)^2 \le 2 w(y) t^2(y) + 2 w(y) f^2(x)
    \le 4\lambda\tau^2 + 4 \lambda \kappa_2^2 + 2\lambda F^2 = \lambda \Gamma,
\]
with $\Gamma := 4\tau^2 + 4  \kappa_2^2 + 2 F^2$. As only the $m$ contaminated indices differ between $\widetilde{\mathcal D}_n$ and $\mathcal D_n$, and
$0 \le \ell_f \le \lambda \Gamma$ pointwise, the contaminated and clean debiased log-likelihoods satisfy
\[
     \frac{1}{2\si^{2}}\Big| \sum_{i=1}^n \ell_f(\tilde x_i,\tilde y_i) - \sum_{i=1}^n \ell_f(x_i,y_i) \Big|
    \le \frac{ m  \lambda \Gamma}{2\si^{2}}.
\]
Repeating the argument of Lemma~\ref{lem : getting rid of the adv outliers with boundedness} with $s_n = F$ and using $\Pi(S_{n,m}\setminus H_n) = 0$ with $H_n := \{\|f\|_\infty \le F\}$, we
obtain, for any event $\mathcal{A}$ depending only on the clean data, 
\begin{equation}\label{eq : thm4 reduction}
    \widetilde{E}_{f_0} G( S_{n,m} \mid \widetilde{\mathcal{D}}_n )
    \leq e^{ m\lambda \Gamma\sigma^{-2}}\, E_{f_0}\!\left[ \frac{I_{[n]}^{w,\mathrm{de}}(S_{n,m})}{I_{[n]}^{w,\mathrm{de}}(\mathcal{F}_n)}\mathbf{1}_{\mathcal{A}} \right]
    + P_{f_0}^n(\mathcal{A}^c).
\end{equation}
We apply Lemma~\ref{lem : denom prob debiased} to the prior $\Pi$ so that we obtain an event $\mathcal{A}$ such that $P_{f_0}^n(\mathcal{A}^c) = o(1)$ and, on $\cA$, the following lower bound holds
\begin{equation}\label{eq : thm4 denom}
    I_{[n]}^{w,\mathrm{de}}(\mathcal{F}_n) \geq \tfrac12 \Pi(\mathcal{F}_n) e^{-C^* n\varepsilon_n^2}
    \geq \tfrac12 e^{-(C^*+1) n\varepsilon_n^2}, 
\end{equation}
the constant $C^*$ being the one in the statement of Lemma~\ref{lem : denom prob debiased}.

For the numerator, Fubini gives
$E_{f_0}[ I_{[n]}^{w,\mathrm{de}}(S_{n,m})] = \int_{S_{n,m}} \big( E_{f_0}[ e^{-\frac{1}{2\si^2}(\ell_f - \ell_{f_0})(x,y)}] \big)^n d\Pi(f)$.
Fix $f$ with $\|f\|_\infty \le F$. The variable $-\frac{1}{2\si^2}(\ell_f - \ell_{f_0})(x,y)$ has mean
$-\frac{1}{2\si^2} E_{f_0}[\ell_f - \ell_{f_0}]$, variance $\frac{1}{4\si^4} \mathrm{Var}_{f_0}(\ell_f - \ell_{f_0})$, and, by the
uniform bound \eqref{bound excess loss} of Lemma \ref{lem : proof risk equivalent L2} , 
\[
    \big| {-}\frac{1}{2\si^2}(\ell_f - \ell_{f_0})(x,y) \big| = \frac{1}{2\si^2}\,|\ell_f(x,y) - \ell_{f_0}(x,y)| \leq \frac{\lambda B (F+F_0)}{2\si^2}  =: \lambda c .
\]
Applying Lemma~\ref{lem : bound on moment generating function} with $t=1$ to the centered variable, whose
absolute value is at most $2 \lambda c$, gives
\[
    E_{f_0}[ e^{-\frac{1}{2\si^2}(\ell_f - \ell_{f_0})(x,y)}]
    \leq \exp\left \{ -\frac{1}{2\si^2} E_{f_0}[\ell_f - \ell_{f_0}] + \frac{1}{4\si^4} \mathrm{Var}_{f_0}(\ell_f - \ell_{f_0})\, \varphi(2 \lambda c) \right \},
\]
with  $\varphi(u) := (e^u - 1 - u)/u^2$.
Inserting the mean bound \eqref{eq : risk equiv L2} and the variance bound \eqref{variance bound} of Lemma  \ref{lem : proof risk equivalent L2}, the exponent
is at most $-K\|f - f_0\|_{2}^2$ with 
\[
    K := \frac{\lambda \kappa_1}{2\sigma^2} - \frac{\lambda^2 B^2}{4\sigma^4}\,\varphi(2 \lambda c) .
\]
We show that $K > 0$ for $\lambda$ small enough. As $\varphi$ is increasing with $\varphi(u) \to \tfrac12$ when $u \to 0$, for $\lambda \le 1/(2c)$ we have
$\varphi(2\lambda c) \le \varphi(1) < 1$, and the last display is (strictly) positive whenever $\lambda < 2\sigma^2\kappa_1/B^2$.
Thus $K > 0$ for every $\lambda < \frac{1}{2c} \wedge \frac{2\sigma^2 \kappa_1}{B^2}$.

Consequently $\big( E_{f_0}[ e^{-\frac{1}{2\si^2}(\ell_f - \ell_{f_0})(x,y)}] \big)^n \leq e^{-Kn\|f - f_0\|_{2}^2}$, and since
$\|f - f_0\|_2^2 \geq M\zeta_{n,m}^2$ on $S_{n,m}$ while $\Pi$ is a probability measure,
\begin{equation}\label{eq : thm4 numerator}
    E_{f_0} \left [ I_{[n]}^{w,\mathrm{de}}(S_{n,m}) \right ] \leq \int_{S_{n,m}} e^{-Kn\|f - f_0\|_{2}^2}\, d\Pi(f)
    \leq e^{-nKM\zeta_{n,m}^2} .
\end{equation}
Combining \eqref{eq : thm4 reduction}, \eqref{eq : thm4 denom} and \eqref{eq : thm4 numerator},
\[
    \widetilde{E}_{f_0} G(S_{n,m} \mid \widetilde{\mathcal{D}}_n)
    \leq 2\exp\big\{  m\lambda \Gamma\sigma^{-2} + (C^*+1)n\varepsilon_n^2 - nKM\zeta_{n,m}^2 \big\} + o(1),
\]
and since $\zeta_{n,m}^2 \geq \frac{m}{n} + \varepsilon_n^2$, the exponent is at most
$(\lambda \Gamma\sigma^{-2} + C^* + 1 - KM)n\zeta_{n,m}^2$. We can conclude by taking $M$ large enough so that $M > (\lambda \Gamma\sigma^{-2} + C^* + 1)/K$.
\end{proof}
\begin{lem} \label{lem : proof risk equivalent L2}
    Assume the weight $w$ is admissible, namely it is positive, differentiable and satisfies
    \eqref{eq : cond w cont} and \eqref{extra conditions on the weight}. Then there exists a constant $\kappa_1 >0$, such that  
    \begin{equation}\label{eq : cond kappa1}
        \inf_{|s|\le F_0}E_\xi[w(s + \si \xi)] \ge \lambda\kappa_1.
    \end{equation}
    The debiased loss $\ell_f$ defined in \eqref{eq : general loss} and \eqref{score matching alta loss}
    satisfies
    \begin{equation}\label{eq : risk equiv L2}
        \lambda \kappa_1 \,\|f-f_0\|^2_{L^2(P_X)}
        \;\le\; E_{f_0}\big(\ell_f(x,y)-\ell_{f_0}(x,y)\big)
        \;\le\; \lambda\,\|f-f_0\|^2_{L^2(P_X)}.
    \end{equation}
    If in addition $\mathcal{F}$ consists of uniformly bounded functions, i.e.\ there
    exists $F>0$ such that $\|f\|_\infty \le F$ for all $f\in\mathcal{F}$, then, setting
    \begin{equation}\label{eq : def M}
        B :=  2\sqrt{2}\,(\tau^2+\kappa_2^2)^{1/2} + F+F_0,
    \end{equation}
    we have the variance bound
    \begin{equation}\label{variance bound}
        E_{f_0}\big(\ell_f(x,y)-\ell_{f_0}(x,y)\big)^2
        \leq \lambda^2 B^2\,\|f-f_0\|^2_{L^2(P_X)} ,
    \end{equation}
    and, for all $(x,y)\in\mathbb{R}^d\times\mathbb{R}$ and all $f\in\mathcal{F}$,
    \begin{equation}\label{bound excess loss}
        \big|\ell_f(x,y)-\ell_{f_0}(x,y)\big| \leq \lambda B(F+F_0).
    \end{equation}
\end{lem}

\begin{proof}
    Let $\bar{w} := w/\lambda$ denote the rescaled weight function, which is independent of $\lambda$.
    An admissible weight function is independent of $n$, positive and differentiable (hence continuous), we can thus find a positive constant $w_1>0$, such that
    \[ \inf_{|u|\le F_0}\bar{w}(u) \ge  w_1.\]
    Denoting $\Phi$ for the c.d.f. of the normal Gaussian distribution, for $|s| \le F_0$, we have
    \[ E_\xi[ w(s + \si \xi)] \ge \lambda w_1 P_\xi(|s+\si \xi| \le F_0) \ge  \lambda w_1[ \Phi((F_0+s)/\si) + \Phi((F_0-s)/\si) -1].\]
    The right hand side of the last display is decreasing in $|s|$, taking the infimum over $|s|\le F_0$, we obtain 
    \[ \inf_{|s|\le F_0}E_\xi[ w(s + \si \xi)] \ge \lambda w_1 (\Phi(2F_0) -1/2) =: \lambda \kappa_1 >0.\]
    For the debiased loss the population risk minimizer is the truth, $f^\ast=f_0$, so the
    square decomposition \eqref{eq : square decomp} with $g=f_0$ gives
    \[
        E_{f_0}\big(\ell_f(x,y)-\ell_{f_0}(x,y)\big)
        = \int_{\mathbb{R}^d} A_x\,\big(f(x)-f_0(x)\big)^2\, dP_X(x),
        \qquad A_x = E_\xi\big[w(f_0(x)+\sigma\xi)\big].
    \]
    By \eqref{eq : cond w cont} one has $w\le\lambda$, hence $A_x\le\lambda$, while since
    $\|f_0\|_\infty\le F_0$ equation \eqref{eq : cond kappa1}
    gives $A_x \ge \lambda \kappa_1 > 0$; therefore $\lambda \kappa_1 \le A_x\le\lambda$ for every $x$, and the previous display yields \eqref{eq : risk equiv L2}. For the remaining bounds,
    write $\Delta f := f(x)-f_0(x)$ and use the difference of squares,
    \begin{equation*}
    \begin{aligned}
        \ell_f(x,y)-\ell_{f_0}(x,y)
        & = w(y)\Big[\big(t(y)-f(x)\big)^2-\big(t(y)-f_0(x)\big)^2\Big]
        \\ & = -\,\Delta f\; w(y)\big(2t(y)-f(x)-f_0(x)\big),
    \end{aligned}
    \end{equation*}
    so it suffices to bound $w(y)\,|2t(y)-f(x)-f_0(x)|$.  We notice that, for all $y \in \mathbb{R}$, the following bound holds 
    \begin{equation}
    \begin{aligned}
    \label{bound on w t^2}
        w(y) t^2(y) & = w(y) \left ( y - \sigma^2 \frac{w'(y)}{w(y)} \right )^2  
                    \\ & \leq 2 w(y) y^2 + 2 w(y) \left ( \sigma^2 \frac{|w'(y)|}{w(y)} \right )^2  \\ & \leq 
                    2 \lambda (\tau^2 +  \kappa_2^2),
    \end{aligned}
    \end{equation}
    where in the last step we used conditions~\eqref{eq : cond w cont} and~\eqref{extra conditions on the weight} to bound the first and second terms, respectively.
    Moreover, 
    writing $w(y)|t(y)| = \sqrt{w(y)}\cdot \sqrt{w(y)}\,|t(y)|$ and using $w(y)\le\lambda$
    together with \eqref{bound on w t^2}, we obtain 
    \[
        w(y)|t(y)|
        \le \sqrt{\lambda}\,\big(w(y)t(y)^2\big)^{1/2}
        \le \sqrt{2} \lambda \,\big(\tau^2+\kappa_2^2 \big)^{1/2},
    \]
    while $w(y)|f(x)| \le \lambda F$ and $w(y)|f_0(x)| \le \lambda F_0$, so that
    \[
        w(y)\,\big|2t(y)-f(x)-f_0(x)\big|
        \le 2\sqrt{2} \lambda \,\big(\tau^2+\kappa_2^2\big)^{1/2} + \lambda(F+F_0)
        = \lambda B,
    \]
    with $B$ as in the statement. Thus, we have 
    \begin{equation}\label{eq : lip excess}
        \big|\ell_f(x,y)-\ell_{f_0}(x,y)\big| \le  \lambda B\,|\Delta f|.
    \end{equation}
    Squaring \eqref{eq : lip excess} and taking expectation under $f_0$ gives
    \[
        E_{f_0}\big(\ell_f(x,y)-\ell_{f_0}(x,y)\big)^2
        \le \lambda^2 B^2\,\|f-f_0\|^2_{L^2(P_X)},
    \]
    which is \eqref{variance bound}. Finally, since $\|f\|_\infty\le F$ and
    $\|f_0\|_\infty\le F_0$ imply $|\Delta f|\le F+F_0$ for all $x$, plugging this into
    \eqref{eq : lip excess} yields \eqref{bound excess loss}.
\end{proof}

\begin{lem}{\cite[Lemma $7.26$]{lafferty2010concentration}} \label{lem : bound on moment generating function}
      Let $X$ be a random variable with zero expectation and variance $\text{Var}(X) = \sigma^2$. Assume $| X | \le c$ for some constant $c > 0$. Then, for any $t > 0$, 
      \begin{equation*}
          E e^{t X } \leq \exp \left \{ t^2 \sigma^2 \left ( \frac{e^{tc} - 1 - tc}{(tc)^2} \right )  \right \}.
      \end{equation*}
\end{lem}

\begin{lem}\label{lem : denom prob debiased}
    Let $\cF_n := \{ f : \|f-f_0\|_\infty \le \varepsilon_n \}$ with $\varepsilon_n \to 0$
    and $n\varepsilon_n^2 \to \infty$. Let $\ell_f$ be the debiased loss
    \eqref{eq : general loss} with $t(y)$ given by \eqref{score matching alta loss} and assume $w$ admissible, see \eqref{extra conditions on the weight}.
    Then, with $C^\ast := \tfrac{\lambda(1 +B)}{2\sigma^2}$, where $B$ is defined in \eqref{eq : def M}, for any prior $\Pi$ with $\Pi(\mathcal{F}_n) > 0$,  as $n\to\infty$,
    \[
        P_{f_0}^n\!\left[\, \int_{\cF_n}
        \exp\!\Big\{ -\tfrac{1}{2\sigma^2}\textstyle\sum_{i=1}^n
        \big(\ell_f(x_i,y_i)-\ell_{f_0}(x_i,y_i)\big)\Big\}\,d\Pi(f)
        \;>\; \tfrac12\,\Pi(\cF_n)\,e^{-C^\ast n\varepsilon_n^2}\,\right] \to 1 .
    \]
\end{lem}

\begin{proof}
    We follow the same proof strategy of Lemma~$1$ in \cite{SyringMartin}. Let
    \[
        D_n := \int_{\cF_n}
        \exp\!\Big\{ - \frac{1}{2\sigma^2}\textstyle\sum_{i=1}^n
        \big(\ell_f(x_i,y_i)-\ell_{f_0}(x_i,y_i)\big)\Big\}\,d\Pi(f),    \]
    and, for $f\in\cF$, the mean and variance of the per-observation excess loss,
    \[
        \mu(f) := E_{f_0}\big[\ell_f(x,y)-\ell_{f_0}(x,y)\big],
        \qquad
        v(f) := \operatorname{Var}_{f_0}\big(\ell_f(x,y)-\ell_{f_0}(x,y)\big),
    \]
    together with the standardized empirical excess loss
    \[
        Z_n(f) := \frac{\sum_{i=1}^n\big(\ell_f(x_i,y_i)-\ell_{f_0}(x_i,y_i)\big)
                        - n\,\mu(f)}{\{n\,v(f)\}^{1/2}},
    \]
    which satisfies $E_{f_0}[Z_n(f)]=0$ and $E_{f_0}[Z_n(f)^2]=1$ since the observations are
    i.i.d.\ under $P_{f_0}^n$. Note that the above display implies the identity
    $\sum_{i=1}^n(\ell_f-\ell_{f_0})(x_i,y_i) = n\,\mu(f) + \{n\,v(f)\}^{1/2}Z_n(f)$. We define the following random subset of $\mathcal{F}_n$:
    \[
        \cZ_n := \{ f \in \cF_n : |Z_n(f)| \leq  \sqrt{n}\varepsilon_n\}.
    \]
    On $\cF_n$ the risk-equivalence and variance bounds of
    Lemma~\ref{lem : proof risk equivalent L2}, combined with
    $\|f-f_0\|_{L^2(P_X)}^2 \le \|f-f_0\|_\infty^2 \le \varepsilon_n^2$, give
    $0 \le \mu(f) \le \lambda\varepsilon_n^2$ and $v(f) \le  \lambda^2 B^2\varepsilon_n^2$, the variance
    constant $B$ being finite because every $f\in\cF_n$ satisfies
    $\|f\|_\infty \le F_0+\varepsilon_n$. Hence on $\cZ_n$,
    where $|Z_n(f)| < \sqrt{n}\varepsilon_n$,
    \begin{equation}
    \begin{aligned}
        \frac{1}{2\sigma^2}\textstyle\sum_{i=1}^n(\ell_f(x_i,y_i)-\ell_{f_0}(x_i,y_i))
        & = \frac{1}{2\sigma^2} n\,\mu(f) + \frac{1}{2\sigma^2}\{n\,v(f)\}^{1/2}Z_n(f)
        \\ & \leq \frac{\lambda(1 + B)}{2 \sigma^2} n \varepsilon_n^2 = C^* n \varepsilon_n^2,
    \end{aligned}
    \end{equation}
    so that, restricting the integral defining $D_n$ to $\cZ_n$,
    \[
        D_n \;\ge\; \int_{\cZ_n} e^{-\frac{1}{2\sigma^2}\sum_i(\ell_f-\ell_{f_0})(x_i,y_i)}\,d\Pi(f)
        \;\ge\; e^{-C^\ast n\varepsilon_n^2}\,\Pi(\cZ_n)
        = e^{-C^\ast n\varepsilon_n^2}\big(\Pi(\cF_n)-\Pi(\cZ_n^c)\big).
    \]
    From this lower bound, by Markov's inequality,
    \[
        P_{f_0}^n\!\Big[ D_n \le \tfrac12\Pi(\cF_n)e^{-C^\ast n\varepsilon_n^2}\Big]
        \le P_{f_0}^n\!\Big[ \Pi(\cZ_n^c) \ge \tfrac12\Pi(\cF_n)\Big]
        \le \frac{2\,E_{f_0}[\Pi(\cZ_n^c)]}{\Pi(\cF_n)} ,
    \]
    and Fubini's theorem simplifies the expectation to
    \[
        E_{f_0}[\Pi(\cZ_n^c)]
        = \int_{\cF_n} P_{f_0}^n\{ |Z_n(f)|\ge \sqrt{n} \varepsilon_n\}\,d\Pi(f),
    \]
    while Chebyshev's inequality, using $E_{f_0}[Z_n(f)^2]=1$, gives $P_{f_0}^n[|Z_n(f)|\ge \sqrt{n}\varepsilon_n]
    \le (n\varepsilon_n^2)^{-1}$, whence
    $E_{f_0}[\Pi(\cZ_n^c)] \le (n\varepsilon_n^2)^{-1}\Pi(\cF_n)$ and
    \[
        P_{f_0}^n\!\Big[ D_n \le \tfrac12\Pi(\cF_n)e^{-C^\ast n\varepsilon_n^2}\Big]
        \le \frac{2}{\Pi(\cF_n)}\cdot\frac{\Pi(\cF_n)}{n\varepsilon_n^2}
        = \frac{2}{n\varepsilon_n^2}.
    \]
    Since $n\varepsilon_n^2\to\infty$, the right-hand side vanishes, completing the proof.
\end{proof}
\subsubsection{The IMQ weight is admissible}\label{proof : assumptions IMQ}

\begin{lem}
\label{lem : assumptions IMQ}
    The weight function in \eqref{def : expl weigh} with $h=1$ and $\lambda \in (0,1]$ is positive and infinitely differentiable. Moreover, it satisfies \eqref{eq : cond w cont} and \eqref{extra conditions on the weight} with the constant $\kappa_2 := \sigma^2/\tau$.
\end{lem}
\begin{proof}
For $h =1$, the weight is
\[
    w(y)=\frac{\lambda}{1+(y/\tau)^2}
        =\frac{\lambda\tau^2}{\tau^2+y^2},
    \qquad \lambda,\tau>0.
\]
Since its denominator is strictly positive on $\mathbb{R}$,
we have $w\in C^\infty(\mathbb{R})$. Moreover, a direct computation gives
\[
    w'(y)
    =- \frac{2 y \lambda  \tau^2}{(\tau^2 + y^2)^2},
\]
and hence
\[
    \sigma^2\frac{|w'(y)|}{w(y)}
    = \frac{2 \sigma^2 |y|}{\tau^2 + y^2} = \frac{\sigma^2}{\tau} \cdot \frac{2 |y| \tau}{\tau^2 + y^2} \leq \frac{\sigma^2}{\tau},
\]
where in the last step we used $2 |y| \tau \leq y^2 + \tau^2$. 
\end{proof}

\subsection{Proofs for Section \ref{sec : recenter}}
\subsubsection{Proof of Theorem \ref{thm : recenter thresh}}\label{proof : recenter thresh}

\begin{proof}
    We follow the proof of Theorem~\ref{thm : contconc thresh} and only detail the two points
    where the recentering affects the argument. 

    First we apply a slightly modified version of Lemma~\ref{lem : getting rid of the outliers} to the posterior with recentered
    weights. 
    The first condition $0 \leq w(x,y) \leq \lambda$ follows directly from the definition of $w$. For the second one, fix $x \in \mathbb{R}^d$ and set $u := y - \hat{f}_n(x)$. If $|u| \leq \tau_n$, we have
    \[y^2 w(x,y) = \lambda y^2 \leq \lambda ( |\hat{f}_n(x)| + \tau_n)^2.\] 
    Otherwise, when $|u| > \tau_n$, writing $y = u + \hat{f}_n(x)$,
    \[
        y^2 w(x,y) = \lambda \frac{y^2 \tau_n^2}{u^2}
        = \lambda \tau_n^2 \Bigl( 1 + \frac{\hat{f}_n(x)}{u} \Bigr)^2
        \leq \lambda \tau_n^2 \Bigl( 1 + \frac{|\hat{f}_n(x)|}{\tau_n} \Bigr)^2
        = \lambda \bigl( |\hat{f}_n(x)| + \tau_n \bigr)^2 .
    \]
    Both cases give the same bound, so that for every $x \in \mathbb{R}^d$,
    \begin{equation}\label{condition 8 refined}
        \sup_{y \in \mathbb{R}} y^2 w(x,y)
        \leq \lambda \bigl( \| \hat{f}_n \|_{\infty} + \tau_n \bigr)^2 .
    \end{equation}
    Unlike in \eqref{eq : cond w cont}, the bound \eqref{condition 8 refined} is random because it involves $\| \hat{f}_n \|_{\infty}$, where $\hat{f}_n$ is computed from the contaminated data. Consider the event with $\widetilde{P}_{f_0}$--probability at least $1 - \delta_n$,
    \[
        \mathcal{G}_n =\mathcal{G}_n(\widetilde{\mathcal{D}}_n) := \bigl\{ \| \hat{f}_n - f_0 \|_{\infty} \leq C \chi_{n,m} \bigr\}.
    \]
    On $\mathcal{G}_n$ and for $n$ sufficiently large, we have $\| \hat{f}_n \|_{\infty} \leq  \| f_0 \|_{\infty}+\| \hat{f}_n - f_0 \|_{\infty}  \leq F_0 + C\chi_{n,m} \leq F_0 + 1$. Therefore on $\cG_n$,
    \[ \sup_{y \in \R} y^2w(x,y) \le \lambda(F_0+1+\tau_n)^2 =: \lambda t_n^2.\]
     Since $\cG_n$ depends on $\widetilde\cD_n$ we first pull it back by $T$ and replace the event $\cA$ in Lemma \ref{lem : getting rid of the outliers} by $\cA \cap T^{-1}(\cG_n)$, where $\cA$ depends only on $\cD_n$. The first step of Lemma $\ref{lem : getting rid of the outliers}$ then reads
    \[ \widetilde{E}_{f_0}G(S \given \widetilde{\cD}_n) \le E_{Q_{f_0}}[ G(S \given T(\cD_n,U))\mathbf{1}_{\cA \cap T^{-1}(\cG_n)}] + P_{f_0}^n(\cA^c) + \delta_n,\]
    where we used a union bound as well as the assumption $Q_{f_0}(T^{-1}(\cG_n))= \widetilde{P}_{f_0}(\cG_n) \ge 1 - \delta_n$. Note that, on $T^{-1}(\mathcal{G}_n)$, the estimator $\hat{f}_n = \hat{f}_n(T(\cD_n,U))$ satisfies the bound $ \| \hat{f}_n - f_0 \|_{\infty} \leq C \chi_{n,m}$. Hence, on $T^{-1}(\cG_n)$, we can follow the proof of Lemma \ref{lem : getting rid of the outliers} with $\tau_n$ replaced by $t_n$, and obtain
   \[ \widetilde{E}_{f_0} G(S\, | \, \widetilde{\cD}_n) \le e^{m \frac{\lambda}{\sigma^2}(t_n^2 + 4F_0^2)}  E_{Q_{f_0}} \left [ \max_{J \in \cJ_{n,m}}I_{J^c}^w(S) / \min_{J \in \cJ_{n,m}} I_{J^c}^w(\cF_n) \, \mathbf{1}_{\cA \cap T^{-1}(\mathcal{G}_n)} \right ] +{P}_{f_0}^n(\cA^c) + \delta_n.\]
   We have obtained a similar inequality as in the first step of the proof of Theorem \ref{thm : contconc thresh} up to the vanishing additive term $\delta_n$.  Note also the expectation is taken with respect to $Q_{f_0}$ because of the event $T^{-1}(\cG_n)$ and the weight $w$, which still depends on $\hat{f}_n$.
    
    It remains to check that the weight equals $\lambda$ on the clean observations, for which we use the choice of $\tau_n$. Recall from the proof of
    Theorem~\ref{thm : contconc thresh} the event
    $\mathcal{B}_n = \{ \max_{i = 1,\dots,n} |\xi_i| \leq \sqrt{3 \log n} \}$,
    that has probability at least $1 - o(n^{-1/2})$ by Lemma~\ref{lem : gaussian max}. It follows that, on $\mathcal{B}_n \cap T^{-1}(\mathcal{G}_n)$, for every observation $(x_i,y_i)$ from the model \eqref{def : model}, we have 
    \[
        | y_i - \hat{f}_n(x_i) |
        \leq \| \hat{f}_n - f_0 \|_{\infty} + \sigma | \xi_i |
        \leq C \chi_{n,m} + \sigma \sqrt{3 \log n} \leq \tau_n  ,
    \]
    given the assumption on $\tau_n$. Consequently, on $\cB_n \cap T^{-1}(\mathcal{G}_n)$, for any $J \in \cJ_{n,m}$, we have that \[I^{w}_{J^c}(B) = I^{\lambda}_{J^c}(B), \qquad \text{for all } B \subseteq \cF.\]
    Employing this identity first by setting $\cA = \cB_n  \cap \cC_{n,m}$ as in the proof of Theorem \ref{thm : contconc thresh}, and then using $\mathbf{1}_{T^{-1}(\cG_n)}\le 1$ to remove the dependence on the law of contaminated dataset, we obtain
     \[ \widetilde{E}_{f_0} G(S\, | \, \widetilde{\cD}_n) \le e^{m \frac{\lambda}{\sigma^2}(t_n^2 + 4F_0^2)}  E_{{f_0}} \left [ \max_{J \in \cJ_{n,m}}I_{J^c}^\la(S) / \min_{J \in \cJ_{n,m}} I_{J^c}^\la(\cF_n) \, \mathbf{1}_{\cA} \right ] +{P}_{f_0}^n(\cA^c) + \delta_n.\]
    The remaining argument from Theorem~\ref{thm : contconc thresh} can therefore be applied with $\tau_n$ replaced by $t_n$. We have, for $n$ large enough, that $t_n^2 + 4F_0^2 \leq 2\tau_n^2$.  Hence,
    \[
        m\frac{\lambda }{ \sigma^2} \bigl(
        t_n^2+4F_0^2\bigr) + 4 m \log n
        \leq
        \frac{4/3 + 2 \lambda}{ \sigma^2} m \tau_n^2
        \leq
        \frac{4/3 + 2 \lambda}{\sigma^2} n \zeta_{n,m}^2.
    \]
    The resulting constants can be absorbed into $M$, provided it is sufficiently large.
\end{proof}
\subsubsection{Proof of Theorem \ref{thm : cliplikconc gen}}\label{proof : cliplikconc}
We start with a slightly stronger version of Lemma \ref{lem : getting rid of the outliers} that comes from the fact that with this recentering, the log--likelihood is also bounded in $f$.
\begin{lem}\label{lem : cliplike getting rid of the adv outliers}
    For any nonnegative sequence $\veps_{n} \to 0$, define the subset
    \[ \cF_n := \{ f \, : \,  \| f-f_0\|_{\infty} \leq \veps_{n} \} \subset \cF.\]
Consider the weight function $w$ given in \eqref{def : cliplik weight}. Then for any measurable subset $S \subset \cF$ and event $\cA$ depending only on the clean data, for $n$ large enough,
\[ \widetilde{E}_{f_0} G(S\, | \, \widetilde{\cD}_n) \le e^{ \frac{\lambda}{\sigma^2}m\tau_n^2 } E_{f_0} \left[ I_{[n]}^{w}(S)/I_{[n]}^{w}(\cF_n)  \, \mathbf{1}_{\cA}\right] + {P}_{f_0}^n(\cA^c).\]
\end{lem}
\begin{proof}
This proof follows the one of Lemma \ref{lem : getting rid of the outliers} with the only difference that, in this case, for all $x,y,f$, we have
\[ w(f,x,y)(y-f(x))^2 \le \lambda\tau_n^2.\]
 Thus, for all $f$ and all $(y,x)$, we have
\[
    - \frac{ \lambda}{ 2\sigma^2} m \tau_n^2  \le \widetilde{L}_{[n],w}(f)- L_{[n],w}(f)\le \frac{ \lambda}{2\sigma^2} m \tau_n^2.
\]
It follows that, for any $S \subset \cF$ and event $\cA$
\[
    \widetilde{E}_{f_0}[ G( S \mid \widetilde{\cD}_n)] \leq e^{\frac{\lambda}{\sigma^2} m \tau_n^2} E_{f_0} \left [ \frac{\int_S \exp \{ -L_{[n],w}(f) \} d\Pi(f)}{\int_\mathcal{F} \exp\{- L_{[n],w}(f) \} d\Pi(f)} \mathbf{1}_{\cA} \right ] + P_{f_0}^n(\cA^c),
\]
where the quantities inside the expectation on the right hand side of the previous display depend only on the uncontaminated sample.

We obtain the desired bound by multiplying both the numerator and the denominator by $\exp \{ L_{[n],w}(f_0)\}$ and restricting the denominator integral to the set $\cF_n$.
\end{proof}
\begin{proof}[Proof of Theorem \ref{thm : cliplikconc gen}]
    In view of Lemma \ref{lem : cliplike getting rid of the adv outliers}, we have, for some event $\cA$ to  be chosen below,
    \[ \widetilde{E}_{f_0} G(S\, | \, \widetilde{\cD}_n) \le e^{ \frac{\lambda}{\sigma^2}m\tau_n^2 } E_{f_0} \left[ I_{[n]}^{w}(S)/I_{[n]}^{w}(\cF_n)  \, \mathbf{1}_{\cA}\right] + P_{f_0}^n(\cA^c).\]
    As before we first deal with the denominator. Let us define the event $\cB_n$ as
    \[\cB_n : = \left \{ \max_{i = 1,\dots,n } |\xi_i| \le  2 \sqrt{\log n }  \right \}.\]
    By Lemma \ref{lem : gaussian max} we have $P_{f_0}^n(\cB_n^c) =  o(n^{-1})$.

    Under $P_{f_0}^n$, on the event $\cB_n$ and for any $f \in \cF_n = \{ f \, : \, \| f - f_0\|_\infty \le \veps_n \}$, we have, for $n$ large enough,
    \[ (y_i - f(x_i))^2 \le 2(y_i -f_0(x_i))^2 + 2(f(x_i) -f_0(x_i))^2 \le 8 \si^2 \log n + 2 \veps_n^2 \le \tau_n^2.\]
    Therefore, under those conditions, for any $f \in \cF_n$ (and in particular for $f = f_0$), we have
    \[ L_{[n],w}(f) = \frac{\lambda}{2 \si^2} \sum_{i=1}^n (y_i - f(x_i))^2 \wedge \tau_n^2 = \frac{\lambda}{2 \si^2} \sum_{i=1}^n (y_i - f(x_i))^2 =: L_{[n],\lambda }(f).\]
    Thus, on $\cB_n$, we have $I_{[n]}^w(\cF_n) = I_{[n]}^\lambda(\cF_n)$. Let $\cC_{n}$ be the event of Lemma \ref{lem : bound denom} (with $m=0$) such that $P_{f_0}^n(\cC_{n}^c) = o(1)$ and, on $\cC_n$, 
    \[
        I_{[n]}^w(\cF_n) \geq e^{-(1 + \frac{\lambda(1 + \lambda)}{\sigma^2}) n \varepsilon_n^2}.
    \]
    Setting $\cA = \cB_n \cap \cC_{n} \cap \Omega_n$ with $\Omega_n$ to be specified below, we obtain 
    \[ \widetilde{E}_{f_0} G(S\, | \, \widetilde{\cD}_n) \le e^{ \frac{\lambda}{\sigma^2}m\tau_n^2 }  e^{(1 + \frac{\lambda(1+\lambda)}{\si^2}) n \varepsilon_n^2}E_{f_0} \left[ I_{[n]}^{w}(S)  \, \mathbf{1}_{\cA}\right] + P_{f_0}^n(\cA^c).\]
    We now study the integral in the numerator. Let us denote by $\cI_f^{>}$ the set of indices for which $|f(x_i)|$ exceeds $5 \si\sqrt{\log n}$, i.e.
    \begin{equation}\label{def : index f too large}
        \cI_{f}^> := \left\{ i \in [n] \, : \, |f(x_i)| > 5 \si\sqrt{\log n} \right\}.
    \end{equation}
    Notice that under $P_{f_0}^n$ and on the event $\cB_n$ we have $|y_i| \le F_0 + 2\si\sqrt{\log n}$, so that, as $n$ gets large enough,
    \begin{align*}
        |y_i-f(x_i)| & \le |y_i| + |f(x_i)| \le F_0 + 7 \si\sqrt{\log n} \le \tau_n  \!\!\!\!&\text{if $i \notin \cI_{f}^>$}, \\
        |y_i-f(x_i)| &\ge |f(x_i)| - |y_i| \ge 3\si\sqrt{\log n} - F_0 \ge 2 \si \sqrt{\log n}   &\text{if $i \in \cI_{f}^>$}.
    \end{align*}
    Therefore, $L_{\cI_f^>,w}(f) \ge 2 \lambda |\cI_f^{>}|\log n$ and $L_{(\cI_f^>)^c,w}(f) = L_{(\cI_f^>)^c,\lambda}(f)$; this yields
    \[ L_{[n],w}(f) \ge 2 \lambda |\cI_f^{>}|\log n + L_{(\cI_f^>)^c,\lambda}(f).\]
    Plugging this into the definition of the numerator integral $I_{[n]}^w$, we obtain the bound
    \[E_{f_0} \left[ I_{[n]}^{w}(S)  \, \mathbf{1}_{\cA}\right] \le E_{f_0} \bigg[ \int_S n^{- 2 \lambda |\cI_f^>|} e^{ L_{[n],\lambda}(f_0) -L_{(\cI_f^>)^c,\lambda}(f)} \, d\Pi(f) \mathbf{1}_{\cA} \bigg], \]
    where we used that under $\cB_n$, $L_{[n],w}(f_0)=L_{[n],\lambda}(f_0)$.
    We split this integral in two parts depending on the size of $|\cI_f^{>}|$, for this we introduce the notation
    \[ V_n : = \{ f \in \cF \, : \, |\cI_{f}^{>}| \ge n \zeta_{n,m}^2 \}.\]
    
    Let us call $(I)$ (resp. $(II)$) the integral over $S \cap V_n$ (resp. $S \cap V_n^c$). We have first
    \[  (I) = E_{f_0} \bigg[\int_{S \cap V_n}n^{- 2 \lambda |\cI_f^>|}  e^{ L_{(\cI_f^>)^c,\lambda}(f_0) -L_{(\cI_f^>)^c,\lambda}(f)} \exp\{L_{\cI_f^>,\lambda}(f_0)\} \,  \mathbf{1}_{\cA}  d\Pi(f) \bigg].\]
    Note that $\mathbf{1}_{\cA} \le  \prod_{i \in \cI_f^>} \mathbf{1}_{|y_i| \le F_0 +2 \si \sqrt{ \log n}}$ and that, for each $f$, the set $\cI^{>}_f$ is fixed given the design points. Thus, we apply the chain rule $E_{f_0} = E_XE_{Y|X}$ and the preceding indicator bound inside the prior integral, so that by Fubini's theorem and conditional independence, we get 
\begin{align*}
(I) \leq E_X \int_{S \cap V_n}
&\, n^{-2\lambda |\cI_f^>|}  \left( E_{Y | X}
    e^{-\frac{\lambda}{2 \si^2}[(Y-f(X))^2 - (Y-f_0(X))^2]}
\right)^{|(\cI_f^>)^c|} \\
&\times \left( E_{Y|X}
    \big[
        e^{\frac{\lambda}{2 \si^2}(Y - f_0(X))^2}
        \mathbf{1}_{|Y| \le F_0 + 2 \si \sqrt{\log n}}
    \big]
\right)^{|\cI_f^>|}
\, d\Pi(f).
    \end{align*}
Let $Z := (Y - f_0(X)) / \sigma$, $a(X) := (f(X) - f_0(X))/\sigma$ and $C_n := F_0 + 2 \sigma \sqrt{ \log n}$. By direct Gaussian computation, we have that
\begin{equation*}
\begin{aligned}
    E_{Y | X}e^{-\frac{\lambda}{2 \si^2}[(Y-f(X))^2 - (Y-f_0(X))^2]}  & = e^{- \frac{\lambda}{2} a^2(X)} E_{Y \mid X} [  e^{\lambda a(X) Z}] \\ & = e^{- \frac{\lambda}{2} a^2(X)} e^{\frac{\lambda^2}{2} a^2(X)} \\ & = e^{- \frac{\lambda(1 - \lambda)}{2 \sigma^2}( f(X) - f_0(X))^2} \leq 1,
\end{aligned}
\end{equation*}
as well as, for $n$ large enough, using the notation $\tilde{\sigma}=\sigma/\sqrt{1-\lambda}$,
\begin{equation*}
\begin{aligned}
      E_{Y|X} \big[e^{\frac{\lambda}{2 \si^2}(Y - f_0(X))^2} \mathbf{1}_{ |Y| \le C_n} \big] & = 
      \frac{1}{\sqrt{2\pi\si^2}} \int_{-C_n}^{C_n} e^{- \frac{(1 - \lambda)}{2 \sigma^2} ( y - f_0(X))^2} dy \\
      &=\frac{1}{\sqrt{2\pi\tilde{\si}^2(1-\lambda)}} \int_{-C_n}^{C_n} e^{- \frac{1}{2 \tilde{\sigma}^2} ( y - f_0(X))^2} dy \\
      & \leq \frac{1}{\sqrt{1-\lambda}}
\end{aligned}
\end{equation*}
Plugging in these bounds into the expression of $(I)$ above, and using the definition of $V_n$, for $n$ large enough,
    \[ (I) \le E_X \int_{S \cap V_n} e^{- |\cI_f^>| \left(2 \lambda \log n +0.5 \log (1-\lambda) \right)} \, d \Pi(f)\le e^{-\lambda n\zeta_{n,m}^2 \log n}.\]

    We now turn to the control of the second integral,
    \[(II) :=  E_{f_0} \bigg[ \int_{S \cap V_n^c} n^{- 2\lambda |\cI_f^>|} e^{ L_{[n],\lambda}(f_0) -L_{(\cI_f^>)^c,\lambda}(f)} \, d\Pi(f) \mathbf{1}_{\cA} \bigg]. \]
    We have to split again this integral over the sieve set $H_n$ and its complement satisfying $\Pi(H_n^c) \le e^{-C n \zeta_{n,m}^2}$. We denote these integrals by $(II_{H})$ and $(II_{H^c})$, respectively. On $H_n^c$ the control of the integral is equivalent to that of $(I)$ and we obtain
    \[ (II_{H^c}) \le E_X \int_{S \cap V_n^c \cap H_n^c} e^{- |\cI_f^>| \left( 2 \lambda \log n  +0.5 \log (1-\lambda)\right)} \, d \Pi(f)\le \Pi(H_n^c) \le e^{- C n\zeta_{n,m}^2}.\]
    The last integral to control is $(II_H)$, where, by completing the likelihood, we get
    \[ (II_H) =  E_{f_0} \bigg[\int_{S \cap V_n^c \cap H_n} n^{-2\lambda |\cI_f^>|}e^{L_{[n],\lambda}(f_0)-L_{[n],\lambda}(f)} \exp\{ L_{\cI_f^>,\lambda}(f)\} \, d\Pi(f) \mathbf{1}_{\cA}\bigg].\]
    Notice first $n^{-2\lambda |\cI_f^>|} \le 1$; using $(y - f(x))^2 \le 2y^2 + 2f(x)^2$, we have on the event $\cB_n$,
    \[ \exp\{ L_{\cI_f^>,\lambda}(f)\} \le e^{\frac{\lambda}{\si^2}|\cI_f^>|(F_0 +2 \si \sqrt{ \log n})^2} \exp\Big\{ \frac{\lambda}{\si^2} \sum_{i \in \cI_f^>}f(x_i)^2 \Big\}.\]
    Since $(F_0+2\si\sqrt{\log n})^2 \le 5\si^2\log n$ for $n$ large enough, and $|\cI_f^>| \le n \ze_{n,m}^2$ on $V_n^c$, this provides
    \[ (II_H) \le e^{5\lambda n \zeta_{n,m}^2\log n} E_{f_0} \bigg[\int_{S \cap V_n^c \cap H_n} e^{L_{[n],\lambda}(f_0)-L_{[n],\lambda}(f)}\exp\Big\{ \frac{\lambda}{\si^2} \sum_{i \in \cI_f^>}f(x_i)^2 \Big\} \, d\Pi(f) \mathbf{1}_{\cA}\bigg].\]

We show below that there exists an event $\Omega_n$, with $P_{f_0}^n(\Omega_n^c)=o(1)$, and constant $C_{II}^{'}>0$, such that for $f\in V_n^c\cap H_n$,
\begin{equation}\label{eq:help:II_H}
        \exp \left \{   \frac{\lambda}{\sigma^2} \sum_{i \in \mathcal{I}_f^{>}} f(x_i)^2 \right \} \leq \exp \{ C_{II}' |\mathcal{I}_f^{>}| \log n \}.
\end{equation}
Then, setting $\cA = \cB_n \cap \cC_n \cap \Omega_n$, so that $P_{f_0}^n(\cA^c) = o(1)$, and recalling $V_n^c = \{ f \, : |\cI_f^>| \le n \zeta_{n,m}^2 \}$, we obtain, for $C_{II} >0$ large enough, that
    \[ (II_H) \le e^{C_{II} n \zeta_{n,m}^2\log n} E_{f_0} \bigg[\int_{S \cap V_n^c \cap H_n} e^{L_{[n],\lambda}(f_0)-L_{[n],\lambda}(f)} \, d\Pi(f) \mathbf{1}_{\cA}\bigg].\]
    Recalling $S = \{ f \, : (1-\lambda)D_\lambda(P_f,P_{f_0}) \ge M \zeta_{n,m}^2 \log n\}$, Fubini's theorem yields
    \[ (II_H) \le e^{C_{II} n \zeta_{n,m}^2\log n} \int_{S}  e^{-(1-\lambda)nD_\lambda(P_f,P_{f_0}) }\, d\Pi(f) \le e^{-(M-C_{II}) n \zeta_{n,m}^2\log n}.\]
    Combining all the previous bounds into the very first inequality of this proof, we get
    \[ \widetilde{E}_{f_0} G(S\, | \, \widetilde{\cD}_n) \le e^{ \frac{\lambda}{\sigma^2}m\tau_n^2 }e^{(1 + \frac{\lambda(1+\lambda)}{\si^2}) n \varepsilon_n^2}  \left[ \underbrace{e^{-C_I n\zeta_{n,m}^2 \log n}}_{(I)}+ \underbrace{e^{- C n\zeta_{n,m}^2}}_{(II_{H^c})} + \underbrace{e^{-(M-C_{II}) n \zeta_{n,m}^2\log n}}_{(II_H)}\right] + o(1).\]
    Provided $M>0$ is taken larger than $C_{II}$, the largest term in the above bracket is $e^{-C n \zeta_{n,m}^2}$. The last displayed bound then goes to zero for any $n\zeta_{n,m}^2 \ge m \tau_n^2 + n \veps_n^2$ and any $C > ( \lambda/\si^2) \vee (1 + \lambda(1+\lambda)/\si^2) = 1+\la(1+\la)/\si^2$, providing the statement of the theorem.

 It remains to prove \eqref{eq:help:II_H}. For each $x_i$, $i\in \mathcal{I}_f^{>}$, we control the size of $f\in V_n^c\cap H_n$ using a finite sum of its values on a carefully chosen finite subset $B(x_i)\subset (x_j:\, j\in (\mathcal{I}_f^{>})^c)$, together with an $L_nh_n^{\vartheta}=O(1)$ error term. By definition, $|f(x_j)|\lesssim \sqrt{\log n}$ for $j\in (\mathcal{I}_f^{>})^c$, so this gives an $O(\log n)$ upper bound for $f(x_i)^2$, $i\in \mathcal{I}_f^{>}$. Choosing $B(x_i)$ is somewhat subtle. For each $x_i$, we select nearby points and  bound the function evaluation using $\lfloor \vartheta\rfloor$-order polynomial interpolation on $B(x_i)$. We note that a lower-order polynomial interpolation would result in a larger error term $L_nh_n^{\vartheta}$ and hence impose stronger constraints on $\zeta_{n,m}\asymp h_n^{d/2}$. To  control the polynomial interpolation uniformly, we select the points in $B(x_i)$ near the points of a sufficiently fine grid, then combine techniques from polynomial approximation with concentration inequalities. We provide the technical details below.

 Recall from the statement of the theorem that $h_n := (4 \zeta_{n,m}^2)^{1/d}$ and define $\ell := \max( \lfloor \vartheta \rfloor, 1)$. We consider the grid of mesh $1 / \lfloor h_n^{-1} \rfloor$ of the unit hypercube. Let us denote by $C(x_i)$, for $i \in \mathcal{I}_f^{>}$, an $r_{\vartheta,d} \ell / \lfloor h_n^{-1} \rfloor$ side-length hypercube on this grid containing $x_i$, for some large enough $r_{\vartheta,d} \in \mathbb{N}$ to be specified later (note that are multiple choices for such hypercubes, but we can take an arbitrary one). Furthermore, let us denote by $x_j'$, $j = 1,...,k$, with $k = (\ell +1)^d$, the equidistant grid points in this hypercube with spacing $r_{\vartheta,d} / \lfloor h_n^{-1} \rfloor$ (we omit in the notation their dependence on $x_i$ for simplicity). Then in view of Lemma \ref{good points approx} with $s_n = n \zeta_{n,m}^2$, $B = \{ x_i : i \in \mathcal{I}_f^{>}\}$ and $\delta = 1/2$ (with $\ell$ as above and $r = r_{\vartheta,d}$), there exists an event $\Omega_n$ satisfying
    \[
        P_{f_0}^n(\Omega_n^c) \leq \frac{e^{-n \zeta_{n,m}^2/2}}{4 \zeta_{n,m}^2} = o(1),
    \]
    such that, on $\Omega_n$, for any $i \in \mathcal{I}_f^{>}$, there exist $j_1(i),...,j_k(i) \in \mathcal{I}_f^{\leq}:=(\mathcal{I}_f^{>})^c$ with $x_{j_b(i)} \in C(x_i)$ and 
    \[
        \sum_{b=1}^k \| x_{j_b(i)} - x_b' \|_2^2 \leq 4 d k \left ( \frac{4 s_n}{ n} \right )^{2/d}.
    \]
    Taking $\delta_{n} = r_{\vartheta,d} \ell / \lfloor h_n^{-1} \rfloor$, we get that
    \[
        \sum_{b=1}^k \| x_{j_b(i)} - x_b' \|_2^2 \leq \frac{4 d k}{r_{\vartheta,d}^2 \ell^2} \delta_n^2.
    \]
    Note that for large enough choice of $r_{\vartheta,d}^2$ the multiplicative factor in front of $\delta_n^2$ can be made arbitrarily small. 
    
    By assumption we have for $f\in H_n$ that $f(x) = f_1(x) + f_2(x)$ with $\| f_2 \|_{\infty} \leq 1$ and $f_1 \in \mathcal{H}^{\vartheta}(L_n)$ with $L_n h_n^\vartheta = O(1)$. The last equality further implies $L_n \delta_n^{\vartheta} = O(1)$. Hence, after translating $C(x_i)$ to $[0,\delta_n]^d$, we apply Lemma \ref{sampling app argument} to $f_1$ of H\"older regularity $\beta = \vartheta$, with $\delta = \delta_n$, and $r_{\vartheta,d}^2$ large enough, so that
    \begin{equation*}
    \begin{aligned}
        | f(x_i) | \leq 1 + |f_1(x_i) | \leq 1 + C_{\vartheta,d} \left [ \left (\sum_{b=1}^k f_1(x_{j_b(i)})^2\right )^{1/2} + L_n \delta_n^{\vartheta} \right ] \leq \tilde{C}_{\vartheta,d}\sigma \sqrt{\log n},
    \end{aligned}
    \end{equation*}
    for a constant $\tilde{C}_{\vartheta,d}$ large enough. Taking the square of both sides above and summing over $i\in \mathcal{I}_f^{>}$ we arrive at \eqref{eq:help:II_H}, concluding the proof of the theorem.
\end{proof}

    \begin{lem}\label{lem : gaussian max}
    Let $\xi_1,\dots,\xi_n$ be i.i.d.\ standard Gaussian random variables. Then, for every $t>0$,
    \[
        P\Big( \max_{i=1,\dots,n} |\xi_i| > t \Big) \leq n \frac{\sqrt{2/\pi}}{t} e^{- t^2/2}.
    \]
    In particular, taking $t=\sqrt{2c\log n}$ with $c\geq 1$ gives
    $P\big(\max_{i=1,\dots,n}|\xi_i| > \sqrt{2c\log n}\big) \le n^{1-c} / \sqrt{2 c \log n}$.
\end{lem}

\begin{proof}
    The standard bound $\bar{\Phi}(t):= 1 - \Phi(t) \leq \phi(t) / t$, implies that $P(|\xi_1|>t) = 2(1-\Phi(t)) \le 2\phi(t)/t$, and the claim follows by a union bound.
\end{proof}

Before stating the next lemma, let us introduce some additional notation. Consider
\[
    M := \left \lfloor \left ( \frac{\delta n}{ 2 s_n }  \right )^{1/d} \right \rfloor \geq 1
\]
One can partition $[0,1]^d$ into $M^d$ equivolume cubes of side length $1/M$. For an arbitrary $x \in [0,1]^d$, let us denote by $C_{R}(x) $ an arbitrary $R/M$-hypercube on this grid containing $x$, for some $R \in \mathbb{N}$ (note that this hypercube need not be unique and we can take an arbitrary one). 

\begin{lem}
\label{good points approx}
    Let $n \geq 1$, $\ell, r \in \mathbb{N} : \ell r \leq M$, $\delta \in (0,1)$ and $s_n \leq \delta n /2$. Let $(X_i)_{i=1}^n$ be independent $\mathcal{U}([0,1]^d)$ random variables. There exists an event $\Omega_n$ with 
    \[
        P(\Omega_n^c) \leq \frac{\delta n}{ 2 s_n} \exp \left \{ - \frac{(1 - \delta)^2}{\delta} s_n  \right \},
    \]
    such that, on $\Omega_n$, for any $B \subset \{X_1,...,X_n\}$ with $|B| \leq s_n$, any $x \in B$, there exist $x_1,...,x_k \in ( \{ X_1,...,X_n \} \setminus B) \cap C_{\ell r} (x)$, $k = (\ell + 1)^d$, such that
    \[
        \sum_{j= 1}^k \| x_j - x_j'(x) \|_2^2 \leq 4 d(\ell +1)^d \left ( \frac{2 s_n}{\delta n} \right )^{2 /d},
    \]
    where $x_1'(x),...,x_k'(x)$, $k = (\ell +1)^d$ denote the $r/M$-length equidistant grid points in $C_{\ell r}(x)$.
\end{lem}

\begin{proof}
    In view of Lemma \ref{lem : chernoff bound}, there exists an event $\Omega_n$, that holds with the required probability, such that, on $\Omega_n$, each $1/M$-side length hypercube on the $1/M$-grid of $[0,1]^d$ contains at least $\delta n/M^d \geq 2 s_n$
    data points $X_i$.
    Given $B \subset \{X_1,...,X_n\}$ with $|B| \leq s_n$ and $x \in B$, each hypercube of side length $1/M$ in $C_{\ell r}(x)$ contains at least $2s_n$ of the data points $(X_i)$ under $\Omega_n$. Thus, since $|B| \leq s_n$, each of these cubes contains at least $s_n$ data points that are not in $B$. 
    
    Let us assign to each $r/M$-grid points $x_j'(x)$, $j = 1,...,(\ell +1)^d$ a $1/M$ length cube containing it. Recall that this cube contains a data point $x_j \in \{ X_1,...,X_n \} \setminus B$, hence $\| x_j - x_j'(x) \|_2^2 \leq d/M^2$. Therefore, under $\Omega_n$, there exist $x_1,...x_k \in ( \{ X_1,...,X_n \} \setminus B ) \cap C_{\ell r}(x)$, such that
    \[
      \sum_{j=1}^k \| x_j - x_j' \|_2^2 \leq \frac{k d}{M^2} \leq 4 d (\ell +1)^d \left ( \frac{2 s_n}{\delta n} \right )^{2 /d} ,
    \]
    concluding the proof.
    \end{proof}

\begin{lem}
\label{sampling app argument}
    Assume $f : [0, \delta]^d \to \mathbb{R}$ satisfies $f \in \mathcal{H}^{\beta}(L)$ for $\beta >0 $ and $\delta,L > 0$. Let $q = \lfloor \beta \rfloor$, $\bar{q} := \max(q,1)$ and $k := (1 + \bar{q})^d$.
    
    We denote by $(x_j')_{j=1}^k$ the $\delta/ \bar{q}$-equidistant grid points in $[0,\delta]^d$. Let $x_1,...,x_k \in [0,\delta]^d$ satisfy  $\sum_{j=1}^k \| x_j - x_j' \|_2^2 \leq \epsilon_{\beta,d}^2 \delta^2$, for some sufficiently small $\epsilon_{\beta,d}$ depending only on $\beta,d$. Then, there exists a constant $C_{\beta,d}$, depending only on $\beta$ and $d$, such that, for all $x_0 \in [0,\delta]^d$, 
    \[
        | f(x_0) | \leq C_{\beta,d} \left [   \left ( \sum_{j=1}^k f^2(x_j)  \right )^{1/2}   + L \delta^\beta           \right ].
    \]
\end{lem}

\begin{proof}
    Let us denote by $\mathcal{P}_{q}$ the space of polynomials spanned by the basis $\{ x^{\nu} : | \nu | \leq q\}$. Then, let us consider the Taylor approximation of $f$ around $x_0$ of the form
    \[
        p_{x_0}(x) = \sum_{| \nu | \leq q} \frac{1}{\nu!} \frac{\partial^{|\nu|} f}{\partial x^{\nu}}(x_0) (x - x_0)^{\nu},
    \]
    where $\nu ! = \nu_1 ! \cdots \nu_d!$ (for $q = 0$, $p_{x_0}(x) \equiv f(x_0)$). Note that, given the regularity of $f$, we can control the Taylor remainder and obtain, for all $x,x_0 \in [0,\delta]^d$,
    \begin{equation}
    \label{error taylor}
        |f(x) - p_{x_0}(x) | \leq C_{\beta,d} L \| x - x_0 \|_2^{\beta}.
    \end{equation}
    Moreover, we have that $p_{x_0}(x_0) = f(x_0)$ and $p_{x_0}(x) \in \mathcal{P}_{q}$.
    
    Consider the grid from the statement rescaled to the unit cube $[0,1]^d$, namely $( \delta^{-1} x_j' )_{j=1}^k$, and its lower set of nodes $\mathcal{K} = \{ k / \bar{q} : k \in \mathbb{N}_0^d, |k|_1 \leq q \}$. This set of nodes forms a unisolvent set for $q$-order polynomial interpolation (see e.g. \cite{floater2017nodal}), hence the corresponding Vandermonde matrix $V_d'$ is not singular. This further implies that the minimal singular value of $V_d'$ is lower bounded by a positive constant $c_{\beta,d} > 0$ depending only on $\beta$ and $d$. Let us denote by $\mathcal{J} \subset \{ 1,...,k\}$ the subset of indices corresponding to the grid points in $\mathcal{K}$, i.e. $\forall j \in \mathcal{J}$, $\exists x \in \mathcal{K}$ such that $\delta^{-1} x_j' = x$.
    
    Next, let us denote by $V_d$ the Vandermonde matrix corresponding to the up to $q$-order monomials evaluated at $\delta^{-1} x_j$, $j \in \mathcal{J}$. Note that
    \[
        \sigma_{\text{min}}(V_d) \geq \sigma_{\text{min}}(V_d') - \| V_d - V_d' \|_{\text{op}},
    \]
    where $\sigma_{\text{min}}(A)$ denotes the smallest singular value of the matrix $A$ and $\|A\|_{\text{op}}$ its operator norm. Furthermore, using that the operator norm is bounded by the Frobenius norm (denoted by $\| \cdot \|_{\text{F}}$), we get that
    \begin{equation}
    \begin{aligned}
        \| V_d - V_d' \|_{\text{op}}^2 \leq \| V_d - V_d' \|_{\text{F}}^2 & \leq \sum_{j \in \mathcal{J}} \sum_{|\nu | \leq q} | (\delta^{-1}x_j)^\nu - (\delta^{-1} x_j')^{\nu} |^2 \\ & \leq C_{\beta,d} \sum_{j=1}^k \delta^{-2} \| x_j - x_j' \|_2^2 \leq C_{\beta,d} \epsilon_{\beta,d}^2,
    \end{aligned}
    \end{equation}
where the second inequality follows from the following arguments. Let $z_j=\delta^{-1}x_j$ and $u_j=\delta^{-1}x_j'$.
    Since $z_j,u_j\in[0,1]^d$, for every multi-index $\nu$
    with $|\nu|_1\le q$,
    \[
    \begin{aligned}
    |z_j^\nu-u_j^\nu|^2
    &\le
    \left(\sum_{i=1}^d
    |z_{j,i}^{\nu_i}-u_{j,i}^{\nu_i}|\right)^2
    \le
    \left(\sum_{i=1}^d
    \nu_i|z_{j,i}-u_{j,i}|\right)^2\\
    &\le
    \left(\sum_{i=1}^d\nu_i^2\right)
    \left(\sum_{i=1}^d|z_{j,i}-u_{j,i}|^2\right)
    \le q^2\|z_j-u_j\|_2^2
    = q^2\delta^{-2}\|x_j-x_j'\|_2^2,
    \end{aligned}
    \]
    where the first inequality follows by telescoping the products
    and bounding all remaining factors by $1$;
    the second uses
    \[
    |a^n-b^n|
    =
    |a-b|\left|\sum_{r=0}^{n-1}a^{n-1-r}b^r\right|
    \le n|a-b|,
    \qquad a,b\in[0,1],\quad n\ge1,
    \]
    while terms with $\nu_i=0$ vanish; the third inequality is Cauchy--Schwarz; and the last uses
    $\sum_i\nu_i^2\le(\sum_i\nu_i)^2\le q^2$.
    
    Hence, for $\epsilon_{\beta,d} \leq c_{\beta,d} / (2 \sqrt{C_{\beta,d}})$,
    \[
        \sigma_{\text{min}}(V_d) \geq \sigma_{\text{min}}(V_d') - \sqrt{C_{\beta,d}}\epsilon_{\beta,d} \geq c_{\beta,d} - c_{\beta,d}/2 = c_{\beta,d}/2 > 0.
    \]
    This implies that the grid of points $\{ \delta^{-1} x_j : j \in \mathcal{J} \}$ form a unisolvent set on $[0,1]^d$ for fitting a $q$-order polynomial. Rescaling this set to $[0,\delta]^d$ will again form a unisolvent set for the $q$-order polynomial interpolation. Therefore, there exists a unique polynomial of order $q$ interpolating pairs $(x_j, p_{x_0}(x_j))_{j \in \mathcal{J}}$. Since $p_{x_0} \in \mathcal{P}_q$, it coincides with this unique solution. 
    
    Let us denote by $\tilde{V}_d$ the Vandermonde matrix formed by the monomials $(x_j^r)$ for $j \in \mathcal{J}$ and $|r| \leq q$. Then the unisolvent polynomial can be written in the form
    \[
        p_{x_0}(x) = v(x) \tilde{V}_d^{-1} [ p_{x_0}(x_j)]_{j \in \mathcal{J}},
    \]
    where $v(x) = [ \phi_1(x), \phi_2(x),..., \phi_{|\mathcal{J}|}(x)]$ where $\phi_j(x)$ denotes the $j$-th monomial. Define the following diagonal matrix $D_\delta = \text{diag}( \delta^{|\nu|}: |\nu | \leq q) $. We have $\tilde{V}_d = V_d D_\delta$ and $v(x) = v(\delta^{-1}x) D_\delta$, so that,
    \[
        v(x) \tilde{V}_d^{-1} = v(\delta^{-1}x) D_\delta D_\delta^{-1} V_d^{-1} =  v(\delta^{-1} x) V_d^{-1},
    \]
    hence
    \[
        p_{x_0}(x) = v(\delta^{-1} x) V_d^{-1} [p_{x_0}(x_j)]_{j \in \mathcal{J}}.
    \]
    
    Finally, note that $\| v(\delta^{-1} x) \|_2^2 = \sum_{j=1}^{|\mathcal{J}|} \phi_j(\delta^{-1}x )^2 \leq |\mathcal{J}|=\binom{ q + d}{d}$, since $\delta^{-1} x \in [0,1]^d$ and the monomials satisfy $|\phi_j(x) | \leq 1$ in the unit cube. Therefore, by combining the above formulas with the polynomial approximation error \eqref{error taylor} and triangle inequality, we get that
    \begin{equation}
    \begin{aligned}
        |f(x_0)| = |p_{x_0}(x_0) | & \leq \frac{\| v(\delta^{-1} x_0) \|_2}{\sigma_{\text{min}}(V_d)} \left ( \sum_{j \in \mathcal{J}} p_{x_0}(x_j)^2 \right )^{1/2} \\ & \leq \frac{\sqrt{\binom{q +d }{d}}}{c_{\beta,d}/2} \left [  \left ( \sum_{j \in \mathcal{J}} f(x_j)^2\right )^{1/2} + C_{\beta,d} L \left (  \sum_{j \in \mathcal{J}} \| x_j - x_0 \|_2^{2 \beta} \right )^{1/2} \right ] \\ & \leq C_{\beta,d} \left [ \left ( \sum_{j=1}^k f(x_j)^2 \right )^{1/2} + L \delta^{\beta}\right ],
    \end{aligned}
    \end{equation}
    concluding the proof of the lemma.
\end{proof}

    \begin{lem}\label{lem : chernoff bound}
        Let $M \ge 1$ and consider $\{C_j \}_{j =1, \dots ,M^d}$ the partition of the hypercube $[0,1]^d$ into $M^d$ equivolume cubes. For $n \ge 1$ and $i \in \{ 1 , \dots , n \}$, consider $X_i \overset{iid}{\sim} \cU([0,1]^d)$. Let
        \[N_j := \sum_{i =1}^n \mathbf{1}_{X_i \in C_j},\] the number of data points in the cube $C_j$. Then for all $\delta \in (0,1)$,
        \[P \left(\forall j, \, N_j \ge \delta \frac{n}{M^d}\right) \ge 1 - M^d e^{-\frac{(1-\delta)^2}{2} \frac{n}{M^d}}\]
    \end{lem}
    \begin{proof}
        For any $j \in \{ 1 , \dots , M^d\}$ and any $i \in \{ 1 , \dots ,n\}$, $ Y_i^{(j)} :=\mathbf{1}_{X_i \in C_j}$ is a Bernoulli random variable of parameter $\text{Vol}(C_j) = 1/M^d$. A Chernoff bound provides, for all $\delta \in (0,1)$
        \[P( N_j \le \delta\frac{n}{M^d}) \le e^{- \frac{(1-\delta)^2}{2}\frac{n}{ M^d}}.\]
        We conclude the proof by applying a union bound for $\{N_j\}$.
        
    \end{proof}

\subsubsection{Proof of Proposition \ref{prop : GP prior mass and sieve}} \label{proof : GP prior mass and sieve}
\begin{proof}
  The prior mass bound is verified in Proposition~\ref{prop : prior mass GP}. Moreover its proof yields the concentration inequality
    $\varphi_{f_0}(\veps_n)\le n\veps_n^2$ for $\veps_n\ge M\veps_n^*$, where $\varphi_{f_0}$ denotes the
    concentration function of $\Pi$, see \eqref{concentration function rescaled gp}.

    Let $\mathbb H_1$ be the unit ball of the RKHS of $\Pi$, let $\mathbb B_1$ be that of
    $\big(\cC([0,1]^d),\|\cdot\|_\infty\big)$, the space of continuous functions on $[0,1]^d$ endowed with the supremum norm, and set
    \[
        H_n := M_n\mathbb H_1 + \ze_{n,m}\mathbb B_1, \qquad M_n := (\sqrt{2} + \sqrt{2 C})\sqrt n\,\ze_{n,m}.
    \]
    The centered small-ball exponent $\varphi_0(\cdot ) := -\log\Pi(\|f\|_\infty\le\,\cdot\,)$ is nonincreasing and dominated by $\varphi_{f_0}$, so
    that, since $\ze_{n,m}\ge\veps_n$, monotonicity provides 
    \[
        \varphi_0(\ze_{n,m}) \;\le\; \varphi_0(\veps_n) \;\le\; \varphi_{f_0}(\veps_n) \;\le\; n\veps_n^2 \;\le\; n\ze_{n,m}^2 .
    \]
    Borell's inequality \cite{van2008reproducing}, combined with $\Phi^{-1}(e^{-u})\ge-\sqrt{2u}$  and $1-\Phi(t)\le e^{-t^2/2}$,
    then gives
    \[
        \Pi(H_n^c) \;\le\; 1-\Phi\Big(\Phi^{-1}\big(e^{-\varphi_0(\ze_{n,m})}\big)+M_n\Big)
        \;\le\; e^{-C\,n\ze_{n,m}^2}.
    \]

    It remains to show the required embedding. The RKHS of the Mat\'ern covariance kernel is
    norm-equivalent to the Sobolev space $\cS^{\be+d/2}([0,1]^d)$ \cite{JMLRvdvvz}, and the Sobolev embedding $\cS^{\be+d/2}\hookrightarrow\cH^{\vartheta}$ (see Chapter 4 of \cite{Taylor2023PDE1}), for all $\vartheta < \beta$, provides a constant $c_\iota=c_\iota(d,\be, \vartheta)$ such that
    $\| h \|_{\cH^{\vartheta}}\le c_\iota\|h\|_{\mathbb H}$ for every $h\in\mathbb H$. It follows that
    $M_n\mathbb H_1\subset\cH^{\vartheta}(c_\iota M_n)$. Write $f\in H_n$ as $f = f_2+f_1$ with $f_2\in M_n\mathbb H_1$ and $\|f_1\|_\infty\le\ze_{n,m} \leq 1$ (for $n$ large enough). Furthermore, \eqref{eq : GP modulus condition} is equivalent to $\ze_{n,m}^2 \lesssim n^{- d/(d + 2 \vartheta)}$, so that $L_n := c_\iota M_n \lesssim \sqrt{n} \zeta_{n,m} \lesssim (\zeta_{n,m}^2)^{-\vartheta/d}$.
\end{proof}

\subsection{Proofs for Section \ref{sec : lower bound}} \label{sec : proof lower bound}

Let $X_{1:n}:=(X_1,\ldots,X_n)$ be an i.i.d.\ sample from a
distribution $P_\theta$ on $\mathcal X$, where $\theta\in\Theta$
is the unknown parameter. As in the introduction, the contaminated
sample is $\widetilde X_{1:n}=T(X_{1:n},U)$, with $T\in\cT_m$
and auxiliary randomness $U\sim\nu$ independent of $X_{1:n}$.
In this section, we make the dependence on $T$ explicit by writing
$  \widetilde{P}_{\theta,T}  := (P_\theta^n\otimes\nu)\circ T^{-1}$
for the law of the contaminated sample and $\widetilde{E}_{\theta,T}$ for the
corresponding expectation. For a pseudo-metric $L$ on $\Theta$, we define the minimax risk
under adversarial contamination as
\[
  \mathcal R^*_{n,m}
  :=
  \inf_{\hat\theta_n}
  \sup_{\theta\in\Theta}
  \sup_{T\in\cT_m}
  \widetilde{E}_{\theta,T}\!\left[L(\hat\theta_n,\theta)\right],
\]
where the infimum is over all estimators $\hat\theta_n$ based
on the contaminated sample $\widetilde X_{1:n}$.

We recall the modulus of continuity introduced by \cite{chen_robust_cov}, defined as
\[
\omega(\varepsilon,\Theta,L)
:=
\sup\left\{
L(\theta_1,\theta_2):
\mathrm{TV}(P_{\theta_1},P_{\theta_2})
\le \frac{\varepsilon}{1-\varepsilon},
\ \theta_1,\theta_2 \in \Theta
\right\},
\]
where $\mathrm{TV}(\cdot,\cdot)$ denotes the total variation distance.
The following result establishes a general way to obtain minimax lower bounds with adversarial contamination.

\begin{prop}
\label{lower bound chen gao}
Assume that there exists a sequence $\psi_n(0)$ and a constant $c_1 > 0$, such that
\[
\inf_{\hat\theta_n}
\sup_{\theta\in\Theta}
P_\theta^{n}
\left\{
L(\hat\theta_n,\theta)
\gtrsim \psi_n(0)
\right\}
\ge c_1.
\]
Then there exists a constant $c_2 > 0$ such that
\begin{equation}
\label{first result}
\inf_{\hat\theta_n}
\sup_{\theta\in\Theta}
\sup_{T\in \cT_m} \widetilde{P}_{\theta,T}
\left\{
L(\hat\theta_n,\theta)
\gtrsim \psi_n(m)
\right\}
\ge c_2,
\end{equation}
where
\[
\psi_n(m)
=
\psi_n(0)
\vee
\omega\!\left(\frac{m}{2 n},\Theta,L\right).
\]

Consequently, the adversarial minimax risk satisfies
\begin{equation}
\label{second result}
\mathcal R^*_{n,m}
\gtrsim \psi_n(m).
\end{equation}
\end{prop}

\begin{proof}[Proof of Proposition \ref{lower bound chen gao}]

When $\psi_n(0) \geq \omega(m/(2n), \Theta, L)$, we have $\psi_n(m) = \psi_n(0)$. Thus
\begin{equation*}
\begin{aligned}
     \inf_{\hat{\theta}_n} \sup_{\theta \in \Theta} \sup_{T \in \cT_m} \widetilde{P}_{\theta,T} \left \{ L(\hat{\theta}_n, \theta) \gtrsim \psi_n(m)   \right \} 
     &\geq 
    \inf_{\hat{\theta}_n} \sup_{\theta \in \Theta}  P_\theta^{n} \left \{ L(\hat{\theta}_n, \theta) \gtrsim \psi_n(m)   \right \}\\ 
     & = \inf_{\hat{\theta}_n} \sup_{\theta \in \Theta}  P_\theta^{n} \left \{ L(\hat{\theta}_n, \theta) \gtrsim \psi_n(0)   \right \} \geq c_1 > 0,
\end{aligned}    
\end{equation*}
where in the first step we used that the identity map $T(x_{1:n},u) = x_{1:n}$  belongs to $\cT_m$.

It remains to treat the case $\psi_n(0) < \omega(m/(2n), \Theta, L)$, for which it suffices to show
\[
    \inf_{\hat{\theta}_n} \sup_{\theta \in \Theta} \sup_{T \in \cT_m} \widetilde{P}_{\theta, T} \left \{  L  (\hat{\theta}_n, \theta) \geq \tfrac{1}{4}\,\omega \left (\tfrac{m}{2n}, \Theta, L \right ) \right \} \geq c_3 > 0
\]
then we take the constant of the statement as  $c_2 = c_1 \wedge c_3 > 0$. Note first that for an identifiable model $\omega(0,\Theta,L)=0\le\psi_n(0)$, so $m=0$ always satisfies the condition of the first case. Hence in this second case we necessarily have $m \geq 1$.

The following construction is adapted from the proof of Theorem~5.1 in \cite{chen_robust_cov}, with modifications to handle the adversarial model with at most $m$ corruptions. Set $\gamma := m/(2n)$. By definition of the supremum, we can find a pair $\theta_1,\theta_2 \in \Theta$ satisfying
\[
     L(\theta_1,\theta_2)\geq \frac{1}{2} \omega(\gamma, \Theta, L) \quad  \text{and} \quad \mathrm{TV}(P_{\theta_1}, P_{\theta_2}) \leq \gamma/(1 - \gamma).
\]
Then there exists $\gamma' \leq \gamma $ such that $\mathrm{TV}(P_{\theta_1}, P_{\theta_2}) = \gamma'/(1- \gamma')$.
For these $\theta_1, \theta_2 \in \Theta$, let $p_{\theta_1}, p_{\theta_2}$ denote the respective densities of  $P_{\theta_1}, P_{\theta_2}$ with respect to the common dominating measure $P_{\theta_1}+ P_{\theta_2}$. We define the distributions $Q_1, Q_2$ via
\begin{equation*}
\begin{aligned}
    & d Q_1 = \frac{(p_{\theta_2} - p_{\theta_1}) \mathbf{1}\{p_{\theta_2} \geq p_{\theta_1}\}}{\mathrm{TV}(P_{\theta_1}, P_{\theta_2})} d ( P_{\theta_1} + P_{\theta_2}), \\ &
    d Q_2 = \frac{(p_{\theta_1} - p_{\theta_2}) \mathbf{1}\{p_{\theta_1} \geq p_{\theta_2}\}}{\mathrm{TV}(P_{\theta_1}, P_{\theta_2})} d ( P_{\theta_1} + P_{\theta_2}).
\end{aligned}
\end{equation*}
These are probability measures, and the mixtures
\[
    Q_{1,\gamma'} = (1 - \gamma') P_{\theta_1} + \gamma' Q_1, \quad Q_{2,\gamma'} = (1 - \gamma') P_{\theta_2} + \gamma' Q_2
\]
satisfy $Q_{1,\gamma'} = Q_{2,\gamma'}$, since the identity $\mathrm{TV}(P_{\theta_1}, P_{\theta_2}) = \gamma'/(1-\gamma')$ makes their densities coincide (until here we basically followed Theorem~5.1 in \cite{chen_robust_cov}).

Relying on these measures, we design two randomized adversarial maps $T_1, T_2 \in \cT_m$. Generate $Y_i \overset{\text{i.i.d}}{\sim} \text{Bern}(\gamma')$ for $i=1,...,n$ and denote by $R$ their joint law, namely $R := \text{Bern} (\gamma')^{\otimes n}$. For each $j \in \{1,2\}$, consider a clean sample $X^{(j)}_{1:n} \sim P_{\theta_j}^{n}$ and, independently of it, generate variables
\[
    \quad Z^{(j)}_i \overset{\text{i.i.d}}{\sim} Q_j, \quad i=1,\dots,n.
\]
Take the random vectors
$X_{1:n}^{(1)},X_{1:n}^{(2)},Z_{1:n}^{(1)},Z_{1:n}^{(2)}$ to be mutually independent and independent of $Y_{1:n}$.
Let us define the set $\mathcal{O}_0 : = \{ i : Y_i =1\}$. If $|\mathcal{O}_0|\leq m$, we set $\mathcal{O} \equiv \mathcal{O}_0$; otherwise, choose $\mathcal{O}$ uniformly among the subsets of $\mathcal{O}_{0}$ of size $m$, using further independent randomness. 
Using the same clean sample and the same $Y_{1:n}$, we define, for each $j$, the following samples
\[
    \widetilde{X}_i^{(j)} = \begin{cases} Z^{(j)}_i,  &i \in \mathcal{O}, \\ X_i^{(j)}, & i \notin \mathcal{O},  \end{cases}, \quad \overline{X}_i^{(j)} = \begin{cases} Z^{(j)}_i,  &i \in \mathcal{O}_0, \\ X_i^{(j)}, & i \notin \mathcal{O}_0.  \end{cases}
\]
Let $U$ collect the Bernoulli variables $Y_i$, the replacement draws,
and the additional randomness used to select $\mathcal O$.
By construction, $U$ is independent of the clean samples. For $j=1,2$, the preceding construction
defines a measurable map $T_j$ such that
\[
  \widetilde X_{1:n}^{(j)}
  =T_j(X_{1:n}^{(j)},U).
\]
Since the construction modifies at most $m$ coordinates for every
input, $T_j\in\cT_m$, and the law of
$\widetilde X_{1:n}^{(j)}$ is $\widetilde{P}_{\theta_j,T_j}$. Moreover, we denote by $\overline{\mathbb{P}}_j$ the law of the second sample $(\overline{X}_i^{(j)})_{i=1}^n$. Note that
\[
    \overline{\mathbb{P}}_j = Q_{j, \gamma'}^{\otimes n}.
\]
Since $Q_{1,\gamma'} = Q_{2,\gamma'}$, we have
\[
     \overline{\mathbb{P}}_1 =  \overline{\mathbb{P}}_2 =: \overline{\mathbb{P}}.
\]
For each $j$, the joint law of the two samples constructed above is a coupling of $\widetilde{P}_{\theta_j,T_j}$ and $\overline{\mathbb{P}}$. On $\{ | \mathcal{O}_0 | \leq m\}$, the two samples coincide because $\mathcal{O} = \mathcal{O}_0$. By the coupling inequality of the total variation distance, we obtain
\[
    \mathrm{TV}(\widetilde{P}_{\theta_j,T_j}, \overline{\mathbb{P}}) \leq R( |\mathcal{O}_0| > m), \quad j =1,2.
\]
Consequently, the triangle inequality yields
\[
     \mathrm{TV}(\widetilde{P}_{\theta_1,T_1},\widetilde{P}_{\theta_2,T_2}) \leq \mathrm{TV}(\widetilde{P}_{\theta_1,T_1}, \overline{\mathbb{P}}) + \mathrm{TV}(\widetilde{P}_{\theta_2,T_2}, \overline{\mathbb{P}}) \leq 2  R( |\mathcal{O}_0| > m).
\]
Now $|\mathcal{O}_0| \sim \mathrm{Bin}(n, \gamma')$ with mean $\mu = n\gamma' \leq n\gamma = m/2$, then by the multiplicative Chernoff bound, for any $a > \mu$, $R(|\mathcal{O}_0| \geq a) \leq e^{-\mu}(e\mu/a)^a$. Note that the right-hand side of the latter expression is increasing in $\mu$ on $\mu < a$, so taking $a = m+1$ and the worst case $\mu = m/2$,
\[
    R(|\mathcal{O}_0| > m) = R(|\mathcal{O}_0| \geq m+1)
    \leq e^{-m/2}\left(\frac{e\,m}{2(m+1)}\right)^{m+1} =: \rho_m.
\]
We have that $2\rho_m \le 2\rho_1 = e^{3/2}/8 < 0.561$ for all $m \ge 1$: the continuous extension has $(\ln\rho)'(x) = \bigl(\tfrac1x - \ln(1+\tfrac1x)\bigr) - (\ln 2 - \tfrac12) < 0$ for $x \ge 2$, since $\tfrac1x - \ln(1+\tfrac1x)$ decreases and at $x=2$ already equals $\tfrac12 - \ln\tfrac32 < \ln 2 - \tfrac12$, so $\rho_m$ decreases for $m \ge 2$, while $\rho_1 > \rho_2$. Hence
\[
    \mathrm{TV}(\widetilde{P}_{\theta_1, T_1}, \widetilde{P}_{\theta_2,T_2}) \leq 2\rho_m \leq 2\rho_1 < 1 \qquad \text{for all } m \geq 1.
\]
Finally, given the assumption on $L$, the events $\{L(\hat\theta_n,\theta_j) < \tfrac12 L(\theta_1,\theta_2) \}$, $j=1,2$, are disjoint. Recall that we have $L(\theta_1, \theta_2) \geq \frac{1}{2} \omega(\gamma, \Theta, L)$, thus by monotonicity of the probability measures $\widetilde{P}_{\theta,T}$ and Le Cam's method (\cite{tsybakov2009}, Theorem $2.2$)
\begin{equation*}
\begin{aligned}
    & \inf_{\hat{\theta}_n} \sup_{\theta \in \Theta} \sup_{T \in \cT_m} \widetilde{P}_{\theta, T} \left \{  L(\hat{\theta}_n, \theta) \geq \tfrac{1}{4}\omega \left (\gamma, \Theta, L \right ) \right \} \geq \\ &
     \inf_{\hat{\theta}_n} \max_{j \in \{1,2\}} \widetilde{P}_{\theta_j, T_j} \left \{  L(\hat{\theta}_n, \theta_j) \geq \tfrac{1}{2} L(\theta_1,\theta_2) \right \}
     \geq \tfrac{1}{2}\left(1- \mathrm{TV}(\widetilde{P}_{\theta_1,T_1}, \widetilde{P}_{\theta_2,T_2} )\right) \geq \tfrac{1}{2}(1 - 2\rho_1) > 0.
\end{aligned}
\end{equation*}
This proves \eqref{first result}. The lower bound on the minimax risk \eqref{second result} follows from \eqref{first result} by Markov's inequality.
\end{proof}

\subsubsection{Proof of Theorem \ref{lower bound nonparametric regression}} \label{lower bound nonparametric regression proof}

In this case the parameter $\theta \in \Theta$ coincides with the regression function $f \in \mathcal{F} \subset L^2(P_X)$, and the observations $\{(x_i,y_i)_{i=1}^n\}$ form a sample from the distribution $P_{f_0}^n$ for some $f_0$, see \eqref{def : model}.
\begin{proof}
    Standard minimax theory for nonparametric regression model with $f_0 \in \mathcal{S}_{\text{per}}^\beta(M)$ yields $\psi_n(0) \gtrsim n^{- \frac{ \beta}{2 \beta + d}}$ (\cite{tsybakov2009}, Theorem $2.11$; see also Chapter 6 of \cite{GineNickl2016}). We need to show 
    \[
        \omega\left (\frac{m}{2n}, \mathcal{S}_{\text{per}}^{\beta}(M), \| \cdot \|_{L^2} \right ) \gtrsim \left ( \frac{m}{n} \right )^{\frac{2 \beta}{2 \beta + d} \vee \frac{1}{2}},
    \]
    so that Proposition~\ref{lower bound chen gao} implies the claim.
    
    We recall that the statistical model $P_f$ admits a density with respect to the Lebesgue measure on $\mathbb{R} \times [0,1]^d$ given by
    \[
        p_f(x,y) = p_f(y \mid x) p_X(x), 
    \]
    where $p_X(x)$ is the design density, while $p_f(y \mid x)$ is the density of a Gaussian distribution with mean $f(x)$ and variance $\sigma^2$.  As a consequence, we have that 
\begin{equation}
\begin{aligned}
\label{first bound on the Tv}
    \mathrm{TV}(P_{f_1}, P_{f_2}) &=  \frac{1}{2} \int_{\mathbb{R} \times [0,1]^d} | p_{f_1}(y,x) - p_{f_2}(y,x) | \ dy \ dx   \\ &= \frac{1}{2} \int_{[0,1]^d} \int_{\mathbb{R} } | p_{f_1}(y \mid x) - p_{f_2}(y \mid x) | \ dy \ p_{X}(x) dx  \\
    & = \int_{[0,1]^d} \mathrm{TV}\left (\mathcal{N} (f_1(x),\sigma^2), \mathcal{N}(f_2(x),\sigma^2)\right ) \ p_{X}(x) dx.
\end{aligned}
\end{equation}
For any fixed $x \in [0,1]^d$, by Pinsker's inequality, we have the following
\begin{equation}
\begin{aligned}
\label{second bound on the TV}
     \mathrm{TV}\left (\mathcal{N}(f_1(x),\sigma^2), \mathcal{N}(f_2(x),\sigma^2)\right ) 
     &\leq   \sqrt{ \frac{\mathrm{KL}\left (\mathcal{N}(f_1(x),\sigma^2), \mathcal{N}(f_2(x),\sigma^2)\right )}{2}}\\
     &=  \frac{|f_1(x) - f_2(x)|}{2 \sigma }.
\end{aligned}
\end{equation}
Combining \eqref{first bound on the Tv} with \eqref{second bound on the TV}, we get
\[
    \mathrm{TV}(P_{f_1}, P_{f_2}) \leq \frac{1}{2 \sigma } E_{X}\big( |f_1(X) - f_2(X)|\big) = \frac{1}{2 \sigma} \| f_1 - f_2 \|_{L^1([0,1]^d)},
\]
where in the last step we used that $X \sim \mathcal{U}([0,1]^d)$. 
This in turn implies, using the notation $\gamma := \frac{m}{2n}$, that
\begin{align*}
 \omega(\gamma, \mathcal F, L) 
 &\geq \sup_{f_1,f_2 \in \mathcal{F} \ : \ \| f_1 - f_2 \|_{L^1} \leq   \frac{2\sigma\gamma}{1 - \gamma}} L(f_1,f_2) 
 \geq \sup_{f_1,f_2 \in \mathcal{F} \ : \ \| f_1 - f_2 \|_{L^1} \leq  \sigma\gamma} L(f_1,f_2) = : (\text{OPT})
\end{align*}
since we are taking suprema over a  smaller set. We recall that $L(f ,g) = \| f - g \|_{L^2([0,1]^d)}$. Since $0 \in \mathcal{F}$, restricting the supremum to pairs of the form $f_1 = 0$ and $f_2 = g$ yields
\[
    (\text{OPT}) \geq \sup_{g \in \mathcal{F}: \| g \|_{L^1} \leq \sigma\gamma} \| g \|_{L^2}.
\]
We split into case $\beta > d/2$ and $\beta \leq d/2$. In the former case we consider the function
\begin{equation}
\label{example function}
    f(x) = A \, K\!\left( \frac{x - x_0}{h} \right), \qquad x_0 = (1/2,\dots,1/2) \in [0,1]^d,
\end{equation}
where the scaling $A$ and the bandwidth $h$ are chosen, for some $R  > 0$, as 
\[
    A := \frac{\gamma\sigma}{h^{d}\, \|K\|_{L^1}},
    \qquad
    h := \left( \frac{R \, \|K\|_{L^1}}{\gamma\sigma\, \|K\|_{\cS^{\beta}}} \right)^{-2/(2\beta+d)},
\]
where we write 
\[
    \|K\|_{\cS^{\beta}}:= \left ( \int_{\mathbb{R}^d} \bigl(1+\|\xi\|_2^2\bigr)^{\beta}\, |\widehat{K}(\xi)|^2\, d\xi \right )^{1/2}
\]
for the (fractional) Sobolev norm of $K$, with $\hat{K}$ the multivariate Fourier transform of $K$. We can take the function $K : \mathbb{R}^d \to [0,\infty)$ so that $K \in C^{\infty}(\mathbb{R}^d)$ (and of course $\| K \|_{L^1} > 0$). Moreover, we assume $K$ is supported on the closed cube $\{|u|_{\infty} \le 1/2\}$, so $K$ belongs to  $C_c^{\infty}(\mathbb{R}^d)$, the space of compact supported functions with infinetely many derivatives, and $K$ together with all its partial derivatives vanishes on the boundary $\{|u|_{\infty} = 1/2\}$.

We assume $\gamma$ small enough that $h < 1$. Then the support of $f$ is contained in the cube $x_0 + [-\tfrac{h}{2},\tfrac{h}{2}]^d \subseteq [0,1]^d$. Since $f$ vanishes together with all its derivatives in a neighborhood of the boundary of $[0,1]^d$, its periodic extension to the torus $\mathbb{T}^d$ is $C^{\infty}$. Consequently requiring that 
\[
    \sum_{k \in \mathbb{Z}^d} (1 + \|k\|)^{2\beta}\, f_k^2 \leq M^2
\]
is equivalent to say (\cite{tsybakov2009}, Lemma A.$3$ for $\beta \in \mathbb{N}$ and $d=1$, \cite{GineNickl2016} and \cite{Taylor2023PDE1} for the general case)
\begin{equation}
  \int_{\mathbb{R}^d} \bigl(1+\| \xi \|_2^2\bigr)^{\beta}\, |\widehat{f}(\xi)|^2\, d\xi \leq R^2\label{def:Sob:Fourier}
\end{equation}
for some $R = R(M,d,\beta) > 0$.
By the scaling of the Fourier transform, $\widehat{f}(\xi) = A\, h^{d}\, e^{ i \langle \xi, x_0 \rangle}\, \widehat{K}(h\xi)$, and the change of variables $\zeta = h\xi$ gives
\begin{equation}
\label{sobolev control of the norm}
\begin{aligned}
    \|f\|_{\cS^{\beta}}^2
    &= A^2 h^{2d} \int_{\mathbb{R}^d} \bigl(1+\|\xi\|^2\bigr)^{\beta}\, |\widehat{K}(h\xi)|^2\, d\xi
    = A^2 h^{d} \int_{\mathbb{R}^d} \bigl(1+ h^{-2}\|\zeta\|^2\bigr)^{\beta}\, |\widehat{K}(\zeta)|^2\, d\zeta \\
    &\le A^2 h^{d-2\beta} \int_{\mathbb{R}^d} \bigl(1+\|\zeta\|^2\bigr)^{\beta}\, |\widehat{K}(\zeta)|^2\, d\zeta
    = A^2 h^{d-2\beta}\, \|K\|_{\cS^{\beta}}^2,
\end{aligned}
\end{equation}
where the inequality uses $h \le 1$, so that $1 + h^{-2}|\zeta|^2 \le h^{-2}\bigl(1+|\zeta|^2\bigr)$. Substituting the choices of $A$ and $h$,
\[
    \|f\|_{\cS^{\beta}}^2
    \le A^2 h^{d-2\beta}\,\|K\|_{\cS^{\beta}}^2
    = \frac{\gamma^2\sigma^2}{h^{2d}\,\|K\|_{L^1}^2}\, h^{d-2\beta}\,\|K\|_{\cS^{\beta}}^2
    = \frac{\gamma^2\sigma^2\,\|K\|_{\cS^{\beta}}^2}{\|K\|_{L^1}^2}\, h^{-(2\beta+d)}
    = R^2,
\]
so that $f \in \mathcal{S}_{\text{per}}^{\beta}(M)$. Notice that the function satisfies the following bound
\[
    \| f \|_{\infty} \leq A \| K \|_{\infty} \lesssim \gamma^{\tfrac{2 \beta - d}{2 \beta + d}}
\]
and thus has $L^\infty$-norm bounded by $F_0$ since we are in the case $\beta > d/2$ and we can take $\gamma$ small enough. 
Moreover, we have that 
\[
    \|f\|_{L^1} = A \, h^{d}\, \|K\|_{L^1} = \sigma\gamma.
\]
By computing the $L_2$ norm of the function we get a lower bound for the modulus of continuity 
\[
    \|f\|_{L^2} = A \, h^{d/2}\, \|K\|_{L^2}
    = \frac{\gamma\sigma\, \|K\|_{L^2}}{h^{d/2}\, \|K\|_{L^1}}
    \asymp \gamma \, h^{-d/2}
    \asymp \gamma \cdot  \gamma^{-\tfrac{d}{2\beta+d}}
    = \left ( \frac{m}{n} \right )^{\tfrac{2\beta}{2\beta+d}},
\]
where we used $h \asymp \gamma^{2/(2\beta+d)}$.

Finally in the case $\beta \leq d/2$, we consider the same function as in \eqref{example function} but with the following scaling
\[
    A = \frac{F_0}{ \| K \|_{\infty}},
\]
so that $\| f_0 \|_{\infty} \leq F_0$ and the following bandwidth
\[
    h = \left (  \frac{\gamma\sigma}{ \| K \|_{L^1} A} \right )^{1/d},
\]
which implies $\| f \|_{L^1} =\sigma \gamma$ and $\| f \|_{L^2} \asymp \gamma^{1/2}$. It remains to show that $ f \in \cS_{\text{per}}^{\beta}(M)$. Similar computation to \eqref{sobolev control of the norm} yields
\[
    \| f \|_{\cS^\beta} \lesssim A h^{d/2 - \beta} \| K \|_{\cS^\beta},
\]
that is controlled by $R$ provided $\gamma$ small enough.

\end{proof}
\subsection{Proofs for Section \ref{sec : upper bound}} \label{sec : proofstrimm}

We obtain the following result by applying Theorem~\ref{lugosi and mendelson theorem}, which is a slight adaptation and extension of Theorem~1 of \cite{robust_mean_estimation} to our setting. Our proof of the theorem does not use sample splitting and does not require the clean-data distribution to be continuous (see Section \ref{sec : trimmed mean lugosi}).

Before stating the following corollary on the robust estimation of the Fourier coefficients of the regression function, we recall that the Fourier bases $(\varphi_{k})_{k\in\mathbb{N}}$ and the true function $f_0$ satisfy that $\sup_{k\in \mathbb{Z}^d} \|\varphi_{k}\|_\infty \le F_1$ and $\| f_0 \|_{\infty} \leq F_0$, respectively. Furthermore, for simplicity, assume that the design is uniform and compact, and even more conveniently we consider $x_i \stackrel{\text{i.i.d.}}{\sim} P_X = \mathcal{U}([0,1]^d)$.

\begin{Corollary}
\label{corollary mendelson}
Let $\delta \in (0,1/2]$ and assume that $m$ and $\delta$ are such that $q \leq 1/4$. Then, with $\widetilde{P}_{f_0}$-probability at least $1 - \delta$,
\[
    |\hat f_{n,k} - f_{0,k}| \le K_0\, \sigma_Z \left( \frac{m}{n} \sqrt{\log \left (\frac{n}{m} \right )} + \sqrt{\frac{\log(1/\delta)}{n}} \right),
\]
where $K_0 > 0$ is a constant and $\sigma_Z^2 := (F_0^2 + \sigma^2) F_1^2$.
\end{Corollary}

\begin{proof}

We verify that $y_i \varphi_{k}(x_i)$ has mean $f_{0,k}$ and is sub-Gaussian; the conclusion then follows by applying Theorem~\ref{lugosi and mendelson theorem} (Theorem~$1$ in \cite{robust_mean_estimation}). The mean property is immediate. For sub-Gaussianity, recall that $y_i = f_0(x_i) + \sigma \xi_i$ with $\xi_i \sim \mathcal{N}(0,1)$ independent of $x_i$; hence, for every $t \in \mathbb{R}$,
\begin{equation*}
\begin{aligned}
    E_{f_0} \left[ e^{t(y_i \varphi_{k}(x_i) - f_{0,k})} \right] 
    &= E_{f_0}\!\left[ e^{t(f_{0}(x_i)\varphi_{k}(x_i) - f_{0,k})}\, E_{f_0}\!\left[ e^{t \sigma  \xi_i \varphi_{k}(x_i)}   \mid x_i \right] \right] \\
    &= E_{f_0}\!\left[ e^{t(f_0(x_i)\varphi_{k}(x_i) - f_{0,k})}\, e^{\tfrac{t^2}{2} \sigma^2 \varphi_{k}(x_i)^2} \right] \\
    &\le e^{\tfrac{t^2}{2} \sigma^2 F_1^2}\, E_{f_0} \!\left[ e^{t(f_0(x_i)\varphi_{k}(x_i) - f_{0,k})} \right] 
    \le e^{\tfrac{t^2}{2} \sigma_Z^2},
\end{aligned}
\end{equation*}
where in the last step we used Hoeffding's lemma applied to the centered random variable $f_0(x_i)\varphi_{k}(x_i) - f_{0,k}$, with $|f_0(x_i)\varphi_{k}(x_i)| \le F_0 F_1$.

\end{proof}

\subsubsection{Proof of Proposition \ref{trimmed projection rate}} \label{trimmed projection rate proof}

\begin{proof}
Decompose
\[
    \|\hat f_n - f_0\|_{L^2}^2 = \sum_{\|k\| \le K_{n,m}} (\hat f_{n,k} - f_{0,k})^2 + \sum_{\|k\| > K_{n,m}} f_{0,k}^2  = (I) + (II).
\]

For the second term (squared bias), note that
\begin{align*}
(II)=\sum_{\|k\| > K_{n,m}} f_{0,k}^2 \leq  K_{n,m}^{-2\beta}\sum_{\|k\| > K_{n,m}} f_{0,k}^2\|k\|^{2\beta}\leq M K_{n,m}^{-2\beta},
\end{align*}
using that $f_0\in S_{\text{per}}^{\beta}(M)$.

Let $N_{K_{n,m}}$ denote the number of Fourier coefficients used to construct $\hat f_n$.
Since $N_{K_{n,m}}\lesssim K_{n,m}^d=o(n)$, using $\log n \leq nB_{n,m}$, we have for $n$ large enough,
\[
N_{K_{n,m}}\delta
\le n e^{-3nB_{n,m}}
\le e^{-2nB_{n,m}}.
\]
An application of Corollary~\ref{corollary mendelson} with
$\delta=8\exp(-3nB_{n,m})$ to each coefficient $f_{0,k}$ for $\|k\|_2 \leq K_{n,m}$ and a union bound gives, with probability at least
$1-e^{-2nB_{n,m}}$,
\begin{equation}
|\hat f_{n,k}-f_{0,k}|
\lesssim \sigma_Z\left(
  \frac{m}{n} \sqrt{\log \left( \frac{n}{m} \right )}
  +\sqrt{\frac{\log(1/\delta)}{n}}
  \right) \lesssim \sqrt{B_{n,m}}
\end{equation}
uniformly in $\|k\|\le K_{n,m}$, where we used
$\log(1/\delta)\le3nB_{n,m}$ and the uniform boundedness
of $\sigma_Z$. It follows that $(I)\lesssim K_{n,m}^d B_{n,m}$. Combining the bounds on $(I)$ and $(II)$,
\[
    \|\hat f_n - f_0\|_{L^2}^2 \lesssim K_{n,m}^d B_{n,m} + K_{n,m}^{-2\beta},
\]
and the stated rate follows by the choice $K_{n,m} \asymp B_{n,m}^{-1/(2\beta + d)}$, which balances the two terms.
\end{proof}

\subsection{Proofs for Section \ref{sec : generalized bayes trimmed}}
\subsubsection{Proof of Theorem \ref{thm:contraction adversarial trimmed projection}} \label{proof : contraction adversarial trimmed projection}

\begin{proof}
Recall that $\lambda_{n,m} = 1/ B_{n,m}$ and set 
\[
    K^*_{n,m} := \lambda_{n,m}^{1 / (2 \alpha + d)}, \quad \varepsilon_{n,m} := B_{n,m}^{(\alpha \wedge \beta)/(2 \alpha + d)}.
\]
Define the event
\[
    A_n := \left\{ \max_{\| k\|_2 \leq n} |\hat f_{n,k} - f_{0,k}| \;\lesssim\; \frac{m}{n}\sqrt{\log \left ( \frac{n}{m}   \right )} + \sqrt{\frac{\log n}{n}} \right\}.
\]
Applying Corollary~\ref{corollary mendelson} with the choice of $\delta = n^{-(d+2)}$ and using the union bound over the $|\{ k : \| k \|_2 \leq n \}| \asymp  n^d$ indices yields, as $n \to \infty$
\[
    \widetilde{P}_{f_0}(A_n^c) \lesssim n^d \cdot n^{-(d+2)} =o(1).
\]
Since $G(\cdot \mid \widetilde{\mathcal{D}}_n) \le 1$, it suffices to show that $\widetilde{E}_{f_0} \mathbf{1}_{A_n} G( \|f - f_0\|_{L^2} > M_n \varepsilon_{n,m}  \mid \widetilde{\mathcal{D}}_n) \to 0$. By Markov's inequality and Parseval's identity,
\[
    G\!\left( \|f - f_0\|_{L^2}^2 > \varepsilon_{n,m}^2 M_n^2 \,\middle|\, \widetilde{\mathcal{D}}_n \right)
    \;\le\; \frac{1}{\varepsilon_{n,m}^2 M_n^2} \sum_{k \in \mathbb{Z}^d} \int (f_k - f_{0,k})^2 \, dG(f_k \mid \widetilde{\mathcal{D}}_n).
\]
For $k$ such that $\|k \|_2 > n$, $f_k \equiv 0$ under the posterior, so $\int (f_k - f_{0,k})^2 dG = f_{0,k}^2$. For $\|k \|_2 \leq n$, the bias-variance decomposition  gives
\[
    \int (f_k - f_{0,k})^2 \, dG(f_k \mid \widetilde{\mathcal{D}}_n) \;=\; \sigma^2_{k,n} + (\bar f_{n,k} - f_{0,k})^2,
\]
where the posterior mean $\bar{f}_{n,k}$ and variance $\sigma_{k,n}^2$ are defined in Equation \eqref{posterior mean and variance} of the appendix.

It is therefore sufficient to show
\[
    \mathrm{(I)}+\mathrm{(II)}+\mathrm{(III)}:=
     \sum_{\| k \|_2 \leq n} \mathcal \sigma^2_{k,n} +
 \widetilde{E}_{f_0} \mathbf{1}_{A_n} \sum_{\|k \|_2 \leq n} (\bar f_{n,k} - f_{0,k})^2 
 +\sum_{\|k\|_2 > n} f_{0,k}^2 \lesssim
 \varepsilon_{n,m}^2.
\]

\emph{Bound on $\mathrm{(I)}$.} Splitting at $\|k\|_2 = K_{n,m}^*$,
\begin{align*}
    \mathrm{(I)} &\le\; \sum_{\|k\|_2 \le K_{n,m}^*} \frac{1}{\lambda_{n,m}} + \sum_{K_{n,m}^*<\|k\|_2\leq n  } \sigma_k^2
    \leq\; \frac{(K_{n,m}^*)^d}{\lambda_{n,m}} + 
    \sum_{K_{n,m}^*<\|k\|_2 } \|k\|_2^{-d-2\alpha} \\
    &\lesssim \lambda_{n,m}^{-\frac{2\alpha}{2\alpha + d}}+ (K_{n,m}^*)^{-2\alpha}
    \;\lesssim\; \varepsilon_{n,m}^2,
\end{align*}
where we used the definition of $K_{n,m}^*$ and
\begin{equation}
\label{bound sum by integral}
    \sum_{\|k\|_2 > K_{n,m}^*} \|k\|_2^{-d-p} \asymp \int_{\|x\|_2 > K_{n,m}^*} \|x\|_2^{-d-p}\, dx \asymp \int_{K_{n,m}^*}^\infty r^{-d-p}\, r^{d-1}\, dr \asymp (K_{n,m}^*)^{-p}, \quad p > 0,
\end{equation}
applied with $p = 2\alpha$. 

\emph{Bound on $\mathrm{(II)}$.} Writing $\bar f_{n,k} - f_{0,k} = \frac{\lambda_{n,m}}{\lambda_{n,m} + \sigma_k^{-2}}(\hat f_{n,k} - f_{0,k}) - \frac{\sigma_k^{-2}}{\lambda_{n,m} + \sigma_k^{-2}} f_{0,k}$, we obtain
\[
    (\bar f_{n,k} - f_{0,k})^2 \;\le\; 2 \frac{\lambda_{n,m}^2}{(\lambda_{n,m} + \sigma_k^{-2})^2}(\hat f_{n,k} - f_{0,k})^2 + 2 \frac{\sigma_k^{-4}}{(\lambda_{n,m} + \sigma_k^{-2})^2} f_{0,k}^2 = 2 a_k + 2 b_k.
\]
For $\sum_k b_k$, using $\frac{\sigma_k^{-4}}{(\lambda_{n,m} + \sigma_k^{-2})^2} \le \min\{ \sigma_k^{-4}/\lambda_{n,m}^2,\, 1\}$ and splitting at $\|k\|_2 = K_{n,m}^*$,
\begin{align*}
    \sum_{ \| k \|_2 \leq n} b_k 
    &\;\le\; \lambda_{n,m}^{-2} \sum_{\|k\|_2 \le K_{n,m}^*} \|k\|_2^{4\alpha+2d}\, f_{0,k}^2 \;+\; \sum_{\|k\|_2 > K_{n,m}^*} f_{0,k}^2 \\
    &\;\le\; \lambda_{n,m}^{-2}\, (K_{n,m}^*)^{(4\alpha+2d-2\beta) \vee 0} \sum_{k} \|k\|_2^{2\beta} f_{0,k}^2 + (K_{n,m}^*)^{-2\beta}\sum_{\|k\|_2 > K_{n,m}^*} \|k\|_2^{2\beta}f_{0,k}^2 \\
    &\;\lesssim\; M^2\, (K_{n,m}^*)^{-2 (\beta \wedge (2\alpha + d))}+M^2(K_{n,m}^*)^{-2\beta}  \lesssim  \varepsilon_{n,m}^2,
\end{align*}
where the second line uses $\|k\|_2^{4\alpha+2d-2\beta} \le (K_{n,m}^*)^{(4\alpha+2d-2\beta) \vee 0}$ for $\|k\|_2 \le K_{n,m}^*$, and the third uses $\lambda_{n,m} = (K_{n,m}^*)^{2\alpha+d}$ together with the Sobolev assumption $\sum_k \|k\|_2^{2\beta} f_{0,k}^2 \le M^2$.
For $\sum_k a_k$, on $A_n$ we have $|\hat f_{n,k} - f_{0,k}|^2 \lesssim (m/n)^2 \log(n/m) + \log( n )/ n \asymp B_{n,m}$ uniformly in $\|k \|_2 \leq n$. Splitting at $\|k\|_2 = K_{n,m}^*$:
\begin{align*}
    \widetilde{E}_{f_0} \mathbf{1}_{A_n} \sum_{\|k\|_2 \le K_{n,m}^*} a_k 
    & \lesssim \sum_{\|k\|_2 \le K_{n,m}^*} B_{n,m} \;\lesssim\; (K_{n,m}^*)^d\, B_{n,m} \;\asymp\; \lambda_{n,m}^{d/(2\alpha+d)}\, B_{n,m}, \\
    \widetilde{E}_{f_0} \mathbf{1}_{A_n} \sum_{K_{n,m}^* < \|k\|_2 \leq n} a_k 
    & \lesssim B_{n,m} \sum_{\|k\|_2 > K_{n,m}^*} \frac{\lambda_{n,m}^2}{(\lambda_{n,m} + \sigma_k^{-2})^2}
    \;\le\; B_{n,m}\, \lambda_{n,m}^2 \sum_{\|k\|_2 > K_{n,m}^*} \sigma_k^{4} \\
    &\;=\; B_{n,m}\, \lambda_{n,m}^2 \sum_{\|k\|_2 > K_{n,m}^*} \|k\|_2^{-d-(4\alpha+d)}
    \;\asymp\; B_{n,m}\, \lambda_{n,m}^2\, (K_{n,m}^*)^{-(4\alpha+d)} \\
    &\;\asymp\; B_{n,m}\, \lambda_{n,m}^2 \cdot \lambda_{n,m}^{-1} \cdot \lambda_{n,m}^{-\frac{2\alpha}{2\alpha+d}}
    \;=\; B_{n,m}\, \lambda_{n,m}^{d/(2\alpha+d)},
\end{align*}
where the tail sum applies \eqref{bound sum by integral} with $p = 4\alpha + d$, using $\sigma_k^{2} = \|k\|_2^{-(2\alpha+d)}$, and the final step uses the fact that $K_{n,m}^* = \lambda_{n,m}^{1/(2 \alpha +d )}$ after factoring $(K_{n,m}^*)^{-(4\alpha+d)} = (K_{n,m}^*)^{-(2\alpha+d)}\, (K_{n,m}^*)^{-2\alpha}$. Since we have $\lambda_{n,m}\, B_{n,m} = 1$, both terms are bounded by a constant multiple of $\lambda_{n,m}^{-2\alpha/(2\alpha+d)} \le \varepsilon_{n,m}^2$.

\emph{Bound on $\mathrm{(III)}$.} By the Sobolev assumption,
\[
    \mathrm{(III)} = \sum_{\|k\|_2 > n} f_{0,k}^2\leq  n^{-2\beta}\sum_{\|k\|_2 > n} f_{0,k}^2\|k\|_2^{2\beta} \;\le\; n^{-2\beta} M^2 = o(\varepsilon_{n,m}^2).
\]
Combining $\mathrm{(I)}, \mathrm{(II)}, \mathrm{(III)}$ yields the claim.
\end{proof}
\subsection{Proofs for Section \ref{sec : noadapt}} 
\label{sec : proofnoadapt}
\subsubsection{Proof of Theorem \ref{thmnoadapt}}

\begin{proof}
We construct below pairs $(f_{(i)},T_{(i)})$ for $i=1,2$. In slight abuse of notation we denote them simply $(f_1,T_1), (f_2,T_2)$, dropping the brackets. Also denote $v_1:=\veps_n(\be_1,m_1)$, $v_2:=\veps_n(\be_2,m_2)$ as shorthand.  

It is enough to establish the result working with the $\|\cdot\|_{2,P_X}$--norm, as by Cauchy--Schwarz' inequality $E\|\cdot\|_{2,P_X}\le (E\|\cdot\|_{2,P_X}^2)^{1/2}$ yielding the result for the squared--loss. We denote $L(f,g)=\|f-g\|_{2,P_X}$ the corresponding loss. 
The infimum $\rho_{n,m}$ in the statement is bounded from below by, setting $v:=\max(v_1,v_2)$ and using that the maximum is larger than the average,
\begin{align}
\sqrt{\rho_{n,m}}& \ge \frac12  \inf_{\hat{f}} \left[\frac{\widetilde{E}_{f_1,T_1}L(\hat{f},f_1)}{v_1} +  \frac{\widetilde{E}_{f_2,T_2}L(\hat{f},f_2)}{v_2}\right]\nonumber\\
& \ge (2v)^{-1} \inf_{\hat{f}} \left[\widetilde{E}_{f_1,T_1}L(\hat{f},f_1) +  \widetilde{E}_{f_2,T_2}L(\hat{f},f_2)\right].\label{eq:help:noadapt}
\end{align}
Set, for an arbitrary  function $g$ in $\mathcal{S}_{\text{per}}^\be(M_1)$ and the given contamination number $1\le m=o(n)$,
\begin{align*}
f_1=g,&\qquad m_1=m,\\
f_2=g+f, & \qquad m_2=0\,,\,\be_2=d/2,
\end{align*}
where the function $f$ is constructed in a similar way as in the proof of Theorem \ref{lower bound nonparametric regression}, with slight adaptation of the notation.  
 One sets $x_0=(1/2,\ldots,1/2)$ and for some $h,A$ to be chosen, any $x\in[0,1]^d$,
\[ f(x)= A K\left(\frac{x-x_0}{h}\right), \]
where the choice of kernel $K$ is the same as in the proof of Theorem \ref{lower bound nonparametric regression}. We further set  $\cK:=\cK(f)=[x_0\pm h/2]^d$ the support of $f$ and for some $R>0$ 
\begin{equation} \label{teh}
 A = \frac{R}{\|K\|_{\cS_{\text{per}}^{d/2}}},\qquad h =  \left\{\frac{m}{n} \cdot \frac{\|K\|_{\cS_{\text{per}}^{d/2}} }{2R\|K\|_1}\right\}^{1/d}. 
\end{equation} 
Then, for a constant $c_1>0$ depending on the uniform lower bound on the density of the design points only, we have, for $c_2=c_1 \|K\|_2\sqrt{R/\{2  \|K\|_{\cS_{\text{per}}^{d/2}} \|K\|_1\}}$, 
\begin{equation} \label{ell2b}
 \|f\|_{2,P_X}\ge c_1 A h^{d/2} \|K\|_2 =  c_2\cdot \sqrt{\frac{m}{n}}, 
\end{equation} 
using $ \|f\|_{2,P_X}\geq c_1 \| f \|_{L^2([0,1]^d)}$, together with $\| f \|_{L^2([0,1]^d)} = A h^{d/2} \|K\|_{L^2}$. 
Also, the function $f_1$ is $\be_1$--smooth with $\be_1=\be$ as $g\in \mathcal{S}_{\text{per}}^\be(M_1)$ by construction. The function $f_2$ on the other hand belongs to $\mathcal{S}_{\text{per}}^{d/2}(M_2)$ for $M_2\ge \mu_2$ large enough. Indeed, $\|f_2\|_{\cS_{\text{per}}^{d/2}}\le \|f\|_{\cS_{\text{per}}^{d/2}}+\|f_1\|_{\cS_{\text{per}}^{d/2}}$. Since $\|f_1\|_{\cS_{\text{per}}^{d/2}}\le \|f_1\|_{\cS_{\text{per}}^{\be}}$ for $\be\ge d/2$, this last term can be made smaller than any given constant $c>0$, up to replacing $g$ by $k_2\cdot g$ for a small enough constant $k_2>0$ in the definition of $g$. Also, as in the proof of Theorem \ref{lower bound nonparametric regression}, one computes $\|f\|_{\cS^{d/2}}\le R$ which in view of the equivalence given in \eqref{def:Sob:Fourier} implies $\|f\|_{\cS_{\text{per}}^{d/2}}\le M(R)$, so that $\|f_2\|_{\cS_{\text{per}}^{d/2}}\le c+M(R)\le \mu_2\le M_2$ provided $\mu_2$ is larger than $c+M(R)$, where the constant $R$ is chosen below and depends only on the dimension and the bounds of the density of the design points.  

The adversarial maps $T_1, T_2$ are constructed as follows. The map $T_2$ does not modify the data (coincides with the identity map), since we assume that $m_2=0$ and there are no contaminations in that case. The map $T_1$ is adversarial (but only $X$--adversarial, in that the outliers' set $\cO$ only depends on $x_i$'s, not $y_i$'s) and is defined as follows. Given a $n$--tuple $(x_1,\ldots,x_n)$, let $\mathcal{I}_{\cK}^X$ denote the indices of the $x_i$'s that belong to the support $\cK$ of the function $f$.  
The map $T_1$ sends the clean sample
$((x_i,y_i))_{i=1}^n$ to
$((\tilde x_i,\tilde y_i))_{i=1}^n$
according to the following rule:
\begin{itemize}
\item for  indices $j\in \mathcal{I}_{\cK}^X$ and if the cardinality of $\mathcal{I}_{\cK}^X$ is at most $m$,  one sets
\[ (\tilde{x}_j,\tilde{y}_j)=(x_j,y_j+f(x_j)); \]
if the cardinality of $\mathcal{I}_{\cK}^X$ is more than $m$, then one defines $\tilde{x}_j,\tilde{y}_j$ according to the last display for the $m$ smallest indices in $\mathcal{I}_{\cK}^X$ and sets $(\tilde{x}_j,\tilde{y}_j)=(x_j,y_j)$ otherwise;
\item for indices $j\notin \mathcal{I}_{\cK}^X$, one sets
\[ (\tilde{x}_j,\tilde{y}_j)=(x_j,y_j). \]
\end{itemize}
Let us define the event, with $|B|$ the cardinality of the set $B$,
\[ \cA = \left\{ \left|\mathcal{I}_{\cK}^X\right| \le m \right\}. \]
The random variable $\left|\mathcal{I}_{\cK}^X\right|$ has a binomial $\text{Bin}(n,\int_{\cK} f_X)$ distribution with expectation $n\int_{\cK} f_X \le c_3 n h^d$, where $c_3$ depends on the uniform upper bound of the design density. In particular, if we choose $h$ such that
\begin{equation} \label{condh} 
 2c_3 n h^d\le m,
\end{equation} 
then by Markov's inequality $P_X^{n }[\cA]=1-P_{X}^{n}[\cA^c]\ge 1- E_X[\left|\mathcal{I}_{\cK}^X\right|]/m\ge 1-1/2=1/2$. 
By taking $R$ in \eqref{teh} large enough depending on $c_3, d$ and the kernel $K$ only,  Condition \eqref{condh} is satisfied. Note that one can then also accommodate a possible constraint $\|f\|_\infty\le F_0$ as long as $F_0$ is large enough since $\|f\|_\infty\le A\|K\|_\infty=R \|K\|_\infty/\|K\|_{\mathcal{S}_{\text{per}}^{d/2}}$.

Now we claim that the distributions $(\tilde{x}_i,\tilde{y}_i)$ {\em conditional} on the event $\cA$ coincide under $P_{f_1,T_1}$ and $P_{f_2,T_2}$. Indeed, given $\cA$ and $(x_i)$, the distribution of $(\tilde{x}_i,\tilde{y}_i)$ is specified by the conditional distribution of $(\tilde{y}_i)$ given $(\tilde{x}_i), \cA$ (or equivalently $(x_i),\cA$ since $\tilde{x}_i=x_i$ under both maps). Under the first map, since $f_1(x_i)+f(x_i)=(g+f)(x_i)$ if $i\in \mathcal{I}_{\cK}^X$ and $f_1(x_i)=g(x_i)$ for $i\notin \mathcal{I}_{\cK}^X$, this conditional is the tensor product
\[ \bigotimes_{i\in \mathcal{I}_{\cK}^X} \, \mathcal{N}((g+f)(x_i),\sigma^2) \, \otimes \, \bigotimes_{i\notin \mathcal{I}_{\cK}^X} \, \mathcal{N}(g(x_i),\sigma^2). \] 
 Under the second map,  the conditional distribution at stake is the tensor product of distributions $\mathcal{N}(f_2(x_i),\sigma^2)=\mathcal{N}((g+f)(x_i),\sigma^2)$ if $i\in\mathcal{I}_{\cK}^X$ and $\mathcal{N}(f_2(x_i),\sigma^2)=\mathcal{N}(g(x_i),\sigma^2)$ if $i\notin\mathcal{I}_{\cK}^X$ since in that case $x_i$ falls outside the support of $f$, which shows that the conditionals coincide under both maps. This means that
\begin{align*}
\widetilde{E}_{f_1,T_1}L(\hat{f},f_1) \ge \widetilde{E}_{f_1,T_1}& \left[L(\hat{f},f_1)\mathbf{1}_{\cA}\right]= \widetilde{E}_{f_1,T_1}\left[ L(\hat{f},f_1) \given \cA \right] P_X^n[\cA]\\
& = \widetilde{E}_{f_1,T_1}\left[\widetilde{E}_{f_1,T_1}\left[ L(\hat{f},f_1) \given \cA, (\tilde{x}_i) \right]\given \cA \right] P_X^n[\cA].
\end{align*}
Similarly 
\begin{align*}
 \widetilde{E}_{f_2,T_2}L(\hat{f},f_2) \ge \widetilde{E}_{f_2,T_2}& \left[L(\hat{f},f_2)\mathbf{1}_{\cA}\right]= \widetilde{E}_{f_2,T_2}\left[ L(\hat{f},f_2) \given \cA \right] P_X^n[\cA] \\
 & = \widetilde{E}_{f_2,T_2}\left[\widetilde{E}_{f_2,T_2}\left[ L(\hat{f},f_2) \given \cA, (\tilde{x}_i) \right] \given \cA \right] P_X^n[\cA].
\end{align*}
Since the inner conditionals are the same, as well the distributions of $(\tilde{x}_i) \mid \cA$ (which we denote by $E_X[ \ \cdot \ \given \cA ]$ in the next display) under both maps,  using $P_X^n[\cA]\ge 1/2$ this yields
\begin{align*}
 \left[\widetilde{E}_{f_1,T_1}L(\hat{f},f_1) +  \widetilde{E}_{f_2,T_2}L(\hat{f},f_2)\right]
& \ge  (1/2) E_X \left[ \widetilde{E}_{f_1,T_1}\left[ \left\{ L(\hat{f},f_1) + L(\hat{f},f_2)\right\} \given \cA, (\tilde{x}_i) \right]  \given \cA \right] \\
& \ge  (1/2) L(f_1,f_2)=\|f_1-f_2\|_{2,P_X}/2,
\end{align*}
where the last line uses that $L$ is a distance, so verifies the triangle inequality. Since $f_2-f_1=f$, one deduces in view of \eqref{eq:help:noadapt} and using the lower bound on $\|f\|_{2,P_X}$ obtained in \eqref{ell2b},
\[ \sqrt{\rho_{n,m}} \ge \frac{\|f\|_{2,P_X}}{4v}\ge \frac{c_2\sqrt{m/n}}{4v},\]
which replacing by $v=v_1\vee v_2$ concludes the proof of the theorem.
\end{proof}   

\section{Additional proofs} \label{app : add proofs}

\subsection{Proofs for Section \ref{app : bounded prior}}
\subsubsection{Proof of Theorem \ref{thm : genweight random}}\label{proof : genweight random}
We start by providing a more general version of Lemma \ref{lem : getting rid of the outliers} by splitting the domain of the generalized posterior into a bounded set (with respect to the supremum norm) and its complement. Then a more careful analysis of separately controlling these two regimes can provide more refined upper bounds, in particular in the case of bounded or truncated posteriors one of the terms vanishes. 
\begin{lem}\label{lem : getting rid of the adv outliers with boundedness}
    For any nonnegative sequence $\varepsilon_n \to 0$, define the subset
    \[ \cF_n := \{ f \in \cF \, : \,  \| f-f_0\|_{\infty} \leq \veps_{n} \} \]  For a nonnegative sequence $s_n$, possibly going to $\infty$, consider
\[ H_n := \{f \in \cF \, : \|f\|_\infty \le s_n \}.\]
Suppose the weight function $w$ satisfies \eqref{eq : cond w cont}.
Then, for any measurable subset $S \subset \cF$, any event $\cA$ depending only on the clean data and for $n$ large enough,
\begin{align*}
\widetilde{E}_{f_0}G(S\given \widetilde{\cD}_n) 
&\le e^{\frac{\lambda}{\sigma^2}m(2\tau_n^2 
+ 4F_0^2 + s_n^2)}E_{f_0} \bigg[\frac{I_{[n]}^{w}(S\cap H_n)}{I_{[n]}^{w}(\cF_n)} \mathbf{1_{\cA}}\bigg] \\
&\qquad+ {E}_{Q_{f_0}}[G(S \setminus H_n \given T(\cD_n,U)) \mathbf{1}_{\cA}] + {P}_{f_0}^n(\cA^c).
\end{align*}
In particular, for the $s_n$ truncated-by-rejection posterior $G^{[s_n]}$, we have
\begin{align*}
\widetilde{E}_{f_0}G^{[s_n]}(S\given \widetilde{\cD}_n) 
\le e^{\frac{\lambda}{\sigma^2}m(2\tau_n^2 
+ 4F_0^2 + s_n^2)}E_{f_0} \bigg[\frac{I_{[n]}^{w}(S\cap H_n)}{I_{[n]}^{w}(\cF_n \cap H_n)} \mathbf{1_{\cA}}\bigg] + {P}_{f_0}^n(\cA^c).
\end{align*}
\end{lem}
\begin{proof}
The proof follows the one of Lemma \ref{lem : getting rid of the outliers}. For any event $\cA$ that depends only on the clean data
\[ \widetilde{E}_{f_0} G(S \given \widetilde{\cD}_n) \le E_{Q_{f_0}} [G(S \given T(\cD_n,U))\mathbf{1}_\cA ] + P_{f_0}^n(\cA^c).\]
We can split the posterior over $H_n :=\{ f\, : \, \|f\|_\infty \le s_n\}$ and its complement
\[ \widetilde{E}_{f_0} G(S \given \widetilde{\cD}_n) \le E_{Q_{f_0}} [G(S \cap H_n \given T(\cD_n,U))\mathbf{1}_\cA ] +  E_{Q_{f_0}} [G(S \setminus H_n \given T(\cD_n,U))\mathbf{1}_\cA ] + P_{f_0}^n(\cA^c).\]
Thanks to \eqref{eq : cond w cont}, we still have, for any $(x,y) \in \R^d \times \R$,
\[ 0 \le w(y)(y-f(x))^2 \leq 2\lambda( \tau_n^2 + f(x)^2).\]
Recalling the definition of contaminated log--likelihood $\widetilde{L}_{[n],w}$ in \eqref{def : conta genpost} and the set of corrupted indices \eqref{def : O}, we have 
\[ -\frac{\lambda}{\sigma^2}\sum_{i \in \cO}(\tau_n^2 + f(\tilde x_i)^2)\le \widetilde{L}_{[n],w}(f)- L_{[n],w}(f)\le \frac{\lambda}{\sigma^2}\sum_{i \in \cO}(\tau_n^2 + f(\tilde x_i)^2).\]
We can then write for any $S \subset \cF$ and event $\cA$
\[ E_{Q_{f_0}} [G(S \cap H_n \given T(\cD_n,U))\mathbf{1}_\cA ] \le {E}_{Q_{f_0}} \bigg[ \mathbf{1_{\cA}}\frac{ \int_{S \cap H_n} e^{\frac{\lambda}{\si^2} \sum_{i \in \cO}(\tau_n^2 + f(\tilde x_i)^2)} \exp\{  - L_{[n],w}(f) \} \, d\Pi(f)}{\int_\cF e^{- \frac{\lambda}{\si^2} \sum_{i \in \cO}(\tau_n^2 + f(\tilde x_i)^2)} \exp\{ - L_{[n],w}(f) \} \, d\Pi(f)} \bigg].\]
To lower bound the denominator we restrict the integral on $\cF_n =\{ f \, : \, ||f - f_0||_{\infty} \leq \varepsilon_{n}\}$, on such set and for large enough $n$,
\[ f(\tilde x_i)^2 \le \|f\|_\infty^2 \leq 2 \|f_0\|_\infty^2 + 2\|f-f_0\|_\infty^2 \le 2 F_0^2 + 2\varepsilon_{n}^2 \le 4 F_0^2,\]
where we have used $\|f_0\|_\infty \le F_0$ as well as $\veps_{n} \to 0$ in the last step. Similarly, for the numerator, restriction on $H_n$ ensures $f(\tilde x_i)^2 \le s_n^2$ for all $i \in \cO$. Recall the definition of the posterior integrals $I_{[n]}^w$ in \eqref{def : post integral}, we have
\[ \widetilde{E}_{f_0}G(S\given \widetilde{\cD}_n) \le e^{\frac{\lambda}{\sigma^2}m(2\tau_n^2 + 4F_0^2 + s_n^2)}E_{f_0} \bigg[\frac{I_{[n]}^{w}(S\cap H_n)}{I_{[n]}^{w}(\cF_n)} \mathbf{1_{\cA}}\bigg] + {E}_{Q_{f_0}}[G(S \setminus H_n \given T(\cD_n,U)) \mathbf{1}_{\cA}] + {P}_{f_0}^n(\cA^c).\]
Which provides the first inequality of the lemma. 

Next, for the truncated-by-rejection posterior $G^{[s_n]}$, by definition we have
\[ G^{[s_n]}(S \setminus H_n \given T(\cD_n,U)) = 0.\]
The rest of the proof is identical up to replacing the prior $\Pi$ with its restricted version $\Pi^{[s_n]}$ thanks to the equivalence noted in Section \ref{app : bounded prior}. This only changes $I_{[n]}^w(\cF_n)$ into $I_{[n]}^w(\cF_n \cap H_n)$ and gives 
\[ \widetilde{E}_{f_0}G^{[s_n]}(S\given \widetilde{\cD}_n) \le e^{\frac{\lambda}{\sigma^2}m(2\tau_n^2 + 4F_0^2 + s_n^2)}E_{f_0} \bigg[\frac{I_{[n]}^{w}(S\cap H_n)}{I_{[n]}^{w}(\cF_n\cap H_n)} \mathbf{1_{\cA}}\bigg] + {P}_{f_0}^n(\cA^c).\]
\end{proof}

\begin{proof}[Proof of Theorem \ref{thm : genweight random}]
    We set $H_n := \{ f : \| f \|_{\infty} \leq s_n \}$ and 
    \[ S_{n,m} = \{ f \, : \, \| f - f_0\|_{2,P_X}^2 \ge M \zeta_{n,m}^2 \}.\]
    Applying Lemma \ref{lem : getting rid of the adv outliers with boundedness}, for an event $\cA$ to be chosen below, we obtain
    \[ \widetilde{E}_{f_0} G^{[s_n]}(S_{n,m} \, | \, \widetilde{\cD}_n) \le e^{\frac{\lambda}{\sigma^2}m(2\tau_n^2 + 4F_0^2 + s_n^2)}E_{f_0} \bigg[\frac{I_{[n]}^{w}(S_{n,m} \cap H_n)}{I_{[n]}^{w}(\cF_n \cap H_n)} \mathbf{1_{\cA}}\bigg] + P_{f_0}^n(\cA^c),\]
    
    Since $s_n \ge F_0 + \veps_n$ and $\|f_0\|_\infty \le F_0$, by the triangle inequality, for sufficiently large $n$, we have $\mathcal{F}_n \subset H_n$ and $I_{[n]}^w(\cF_n \cap H_n) = I_{[n]}^w(\cF_n)$.
    Thus, we apply Lemma \ref{lem : denom prob proof gen weight} with $m = 0$ to obtain the set $\cC_{n,0}$ with $P_{f_0}^n(\cC_{n,0}^c) \to 0$ as $n \to \infty$, such that, on $\cC_{n,0}$,
    \[
        I^{w}_{[n]}(\cF_n) \geq e^{-\left(\frac{\lambda(1+\lambda)}{\sigma^2}  + 2\right)n\varepsilon_{n}^2 - \frac{\Delta_n \rho_n}{\sigma^2}n \varepsilon_{n}} \ge e^{-\left(\frac{2+\lambda(1+\lambda)}{\sigma^2}  + 2\right)n\varepsilon_{n}^2 - \frac{2}{\sigma^2}n (\Delta_n \rho_n)^2},
    \]
    where $\si^{-2}\Delta_n\rho_n\varepsilon_n \le 2\si^{-2}(\Delta_n\rho_n)^2 + 2\si^{-2}\varepsilon_n^2$ was used in the last step.
    As in the proof of Theorem \ref{thm : genweight fix}, we take
    \[
        \cA := \cC_{n,0} \cap \left \{ \max_{1 \le i \le n}|\xi_i| \le \rho_n/\si \right \} \qquad \text{such that} \qquad  P_{f_0}^{n}(\cA^c) \to 0,
    \]
    This yields, since $4F_0^2 \le \tau_n^2$ for $n$ large enough,
    \[ \widetilde{E}_{f_0} G^{[s_n]}(S_{n,m} \, | \, \widetilde{\cD}_n) \le e^{\frac{\lambda}{\sigma^2}m(3\tau_n^2 + s_n^2)} e^{\left(\frac{2+\lambda(1+\lambda)}{\sigma^2}  + 2\right)n\varepsilon_{n}^2 + \frac{2}{\sigma^2}n (\Delta_n \rho_n)^2}E_{f_0} [ I_{[n]}^{w}(S_{n,m} \cap H_n) \mathbf{1}_{\cA}] + o(1).\]
    By applying Fubini's theorem and the first inequality of Lemma \ref{lem : proof gen weight bound on E} with $\cI = [n]$, $b_i = f_0(x_i)$ and $a_i = a(x_i)$, given in \eqref{def:a(x)}, we obtain
    \[  E_{f_0}[ I_{[n]}^{w}(S_{n,m} \cap H_n) \mathbf{1}_{\cA}] \leq  \int_{S_{n,m} \cap H_n} E_{X}\left(e^{- \frac{ w_0(1-w_0)}{4}  a(X)^2 + \frac{\Delta_n \rho_n}{\si} |a(X)| }\right)^{n} \, d\Pi(f).\]
    Let $\al_w := w_0(1-w_0)/4$, using the inequality
    \[ - \al_w a^2 + (\Delta_n\rho_n/\si)|a| \le -\frac{\al_w}{2}a^2 + \frac{1}{2\al_w}(\Delta_n\rho_n/\si)^2,\]
    we obtain 
    \[  E_{f_0}[ I_{[n]}^{w}(S_{n,m} \cap H_n) \mathbf{1}_{\cA}] \leq  e^{\frac{1}{2 \al_w\si^2} n(\Delta_n\rho_n)^2}\int_{S_{n,m} \cap H_n} E_{X}\left(e^{- \frac{\al_w}{2}  a(X)^2  }\right)^{n} \, d\Pi(f).\]
    Being restricted on $H_n$ and since $s_n >F_0$, we have
    \[ |a(X)| = \si^{-1}|f(X) - f_0(X)| \le \si^{-1}(s_n + F_0) \le 2\si^{-1}s_n.\]
    Let $c := \al_w/2$ and $d := 4(s_n/\si)^2$, applying the definition of convexity to the map $x \mapsto e^{-cx}$ between $0$ and $d$, for any $t\in[0,1]$, we get
    \[ e^{-c(1-t)d} \le t + (1-t)e^{-cd} .\]
    Taking $1-t = |a(X)|^2/d \in[0,1]$, we obtain
     \[ e^{-\frac{\al_w}{2}|a(X)|^2} \le 1 - v_n |a(X)|^2, \qquad \text{where} \qquad v_n := \frac{1- e^{- \frac{2\al_ws_n^2}{\si^2}}}{4(s_n/\si)^2}.\]
    Since $s_n \ge F_0$, we have $v_n \ge (2s_n / \sigma )^{-2} (1 - e^{-2\al_wF_0^2/\si^2})$. Taking the expectation, we obtain, for some constant $C > 0$
    \[ E_{X}e^{- \frac{\al_w}{2}  a(X)^2 } \le 1 - Cs_n^{-2}E_X|a(X)|^2.\]
    Using $1 -x \le e^{-x}$ and raising to the power $n$,
    \[ E_{X}\left(e^{- \frac{\al_w}{2}  a(X)^2  }\right)^{n} \le e^{-Cn s_n^{-2} E_X |a(X)|^2}.\]
    Recalling the definition of $S_{n,m}$, we get
    \[  E_{f_0}[ I_{[n]}^{w}(S_{n,m} \cap H_n) \mathbf{1}_{\cA}] \leq  e^{\frac{1}{2 \al_w\si^2} n(\Delta_n\rho_n)^2} e^{-C Mn s_n^{-2} \zeta_{n,m}^2}.\]
    We have finally obtained
    \[ \widetilde{E}_{f_0} G^{[s_n]}(S_{n,m} \, | \, \widetilde{\cD}_n) \le e^{\frac{\lambda}{\sigma^2}m(3\tau_n^2 + s_n^2)} e^{\left(\frac{2+\lambda(1+\lambda)}{\sigma^2}  + 2\right)n\varepsilon_{n}^2 + \frac{1}{\sigma^2}(2 +\frac{1}{2\al_w})n (\Delta_n \rho_n)^2} e^{-C Mn s_n^{-2} \zeta_{n,m}^2} + o(1).\]
    Taking $\zeta_{n,m}^2 = s_n^2 \left[ (m/n)(\tau_n^2 + s_n^2) +\varepsilon_n^2 + (\Delta_n\rho_n)^2 \right]$ and $M$ large enough, the last display is $o(1)$, this concludes the proof.

\end{proof}

\subsection{Proofs for Section \ref{sec : sup-norm trimmed}}
\subsubsection{Proof of Theorem \ref{sup norm trimmed estimator}} \label{sup norm trimmed projection estimator proof}
\begin{proof}
Throughout, $C$ denotes a positive constant depending only on $\sigma$, $F_0$, $M$,
$\beta$ and on the basis, whose value may change from line to line.

As in the proof of Corollary \ref{corollary mendelson}, we apply Theorem
\ref{lugosi and mendelson theorem}  (an adapted version of Theorem 1 of \cite{robust_mean_estimation}) to each coefficient. Let
$Z_{\ell,k} = y\,\psi_{\ell,k}(x)$ with $(y,x) \sim P_{f_0}$. Then
$E_{f_0}[Z_{\ell,k}] = \langle f_0, \psi_{\ell,k} \rangle_2$, and $Z_{\ell,k}$ is
sub-Gaussian with variance proxy
$v_\ell^2 = (\sigma^2 + F_0^2)\|\psi_{\ell,k}\|_\infty^2 \lesssim 2^{\ell}$. Indeed, for
every $t \in \mathbb{R}$,
\begin{equation*}
\begin{aligned}
E_{f_0} \left[ e^{t(y \psi_{\ell,k}(x) - \langle f_0, \psi_{\ell,k} \rangle_2 )} \right]
&= E_{f_0}\!\left[ e^{t(f_{0}(x)\psi_{\ell,k}(x) - \langle f_0, \psi_{\ell,k} \rangle_2)}\,
   E_{f_0}\!\left[ e^{t \sigma \xi \psi_{\ell,k}(x)} \mid x \right] \right] \\
&= E_{f_0}\!\left[ e^{t(f_0(x)\psi_{\ell,k}(x) - \langle f_0, \psi_{\ell,k} \rangle_2)}\,
   e^{\tfrac{t^2}{2} \sigma^2 \psi_{\ell,k}(x)^2} \right] \\
&\le e^{\tfrac{t^2}{2} \sigma^2 \| \psi_{\ell,k} \|_{\infty}^2}\,
   E_{f_0}\!\left[ e^{t(f_0(x)\psi_{\ell,k}(x) - \langle f_0, \psi_{\ell,k} \rangle_2)} \right]
\le e^{\tfrac{t^2}{2} (\sigma^2 + F_0^2) \| \psi_{\ell,k} \|_{\infty}^2},
\end{aligned}
\end{equation*}
where in the last step we used Hoeffding's lemma applied to the centred random variable
$f_0(x)\psi_{\ell,k}(x) - \langle f_0, \psi_{\ell,k} \rangle_2$, $x\sim \mathcal{U}([0,1])$, whose range is controlled
by $|f_0(x)\psi_{\ell,k}(x)| \leq F_0 \|\psi_{\ell,k}\|_\infty$. Moreover, by the
orthonormality of the basis $\{\psi_{\ell,k}\}$, the variance of $Z_{\ell,k}$ is bounded by
a constant that does not depend on $(\ell,k)$: since the design is uniform and
$E_{f_0}[y^2 \mid x] = f_0(x)^2 + \sigma^2$,
\[
    \operatorname{Var}_{f_0}(Z_{\ell,k}) \leq E_{f_0} Z_{\ell,k}^2
    = \int_{[0,1]} \bigl(f_0(x)^2 + \sigma^2\bigr) \psi_{\ell,k}(x)^2 \, dx
    \leq (F_0^2 + \sigma^2) \|\psi_{\ell,k}\|_{2}^2 = F_0^2 + \sigma^2 .
\]

We apply Theorem \ref{lugosi and mendelson theorem} to each estimator $\hat{f}_{\ell,k}$, $(\ell,k) \in \mathcal{I}_{n,m}$, at the confidence level
\begin{equation}\label{eq:delta n}
    \delta_n := \frac{1}{2 n^2}  ,
\end{equation}
and note that the corresponding trimming fraction is 
\begin{equation}\label{eq:q n}
    q = 8\,\eta + 12\,\frac{\log(8/\delta_n)}{n}
      = 8\,\frac{m}{n} + 12\,\frac{\log(16 n^2) }{n} \leq 8 \frac{m}{n} + 72 \frac{\log n}{n},
\end{equation}
where in the last step we used that $\log\bigl(16 n^2) \le 6 \log n$.

Since $m= o(n)$, the condition $q\leq 1/4$ of Theorem \ref{lugosi and mendelson theorem} holds for $n$ large enough. Hence, by applying the theorem to $Z_{\ell,k}$ with
$ s^2 = s_{\ell,k}^2 := \text{Var}_{f_0}(Z_{\ell,k}) \leq F_0^2 + \sigma^2$ and $v_\ell^2 \lesssim 2^{\ell}$, for every
$(\ell,k) \in \mathcal{I}_{n,m}$ and uniformly over the adversarial actions, we get that
\[
    \widetilde{P}_{f_0}\left( \bigl| \hat{f}_{\ell,k} - \langle f_0, \psi_{\ell,k} \rangle_2 \bigr|
    > C\Big( 2^{\ell/2}\, q \sqrt{\log(1/q)}
    + s\sqrt{\tfrac{\log(1/\delta_n)}{n}} \Big) \right) \leq \delta_n .
\]
Let $\mathcal{A}_{\ell,k}$ denote the complement of the event in the last display. By
union bound and by the choice \eqref{eq:delta n},
\[
    \widetilde{P}_{f_0}\Biggl( \bigcup_{(\ell,k) \in \mathcal{I}_{n,m}} \mathcal{A}_{\ell,k}^{\,c} \Biggr)
    \;\leq\; \sum_{(\ell,k) \in \mathcal{I}_{n,m}} \widetilde{P}_{f_0}\bigl( \mathcal{A}_{\ell,k}^{\,c} \bigr)
    \;\leq\; |\mathcal{I}_{n,m}| \, \delta_n  =   |\mathcal{I}_{n,m}| \cdot \frac{1}{2 n^2}
    \leq \frac{1}{n},
\]
where in the last step we used the fact that the truncation level in the statement satisfies $2^{L_{n,m}} \le n$, so that $|\mathcal{I}_{n,m}| = 2^{L_{n,m}+1} - 1 \le 2n$.
Consequently, on the event $\mathcal{A} := \bigcap_{(\ell,k) \in \mathcal{I}_{n,m}}
\mathcal{A}_{\ell,k}$ of probability at least $1 - 1/n$, we have simultaneously for all
$(\ell,k) \in \mathcal{I}_{n,m}$,
\begin{equation}\label{eq:all coefficients}
    \bigl| \hat{f}_{\ell,k} - \langle f_0, \psi_{\ell,k} \rangle_2 \bigr|
    \;\leq\; C\left( 2^{\ell/2}\, q \sqrt{\log(1/q)} + \sqrt{\tfrac{\log n}{n}} \right).
\end{equation}

By the localisation property of the CDV basis, for any array of coefficients $(c_k)_k$,
\begin{equation}
\label{localisation CDV}
    \Bigl\| \sum_{k=0}^{2^\ell-1} c_k \psi_{\ell,k} \Bigr\|_\infty
    \;\leq\; \max_{0 \le k \le 2^\ell - 1} |c_k| \cdot
    \sup_{x \in [0,1]} \sum_{k=0}^{2^\ell-1} |\psi_{\ell,k}(x)|
    \;\leq\; C\, 2^{\ell/2} \max_{0 \le k \le 2^\ell - 1} |c_k| .
\end{equation}
Applying this at each level together with \eqref{eq:all coefficients}, we obtain, on $\mathcal{A}$,
\begin{equation}
\begin{aligned}
    \Bigl\| \sum_{\ell=0}^{L_{n,m}} \sum_{k=0}^{2^\ell-1}
    \bigl( \hat{f}_{\ell,k} - \langle f_0, \psi_{\ell,k} \rangle_2 \bigr) \psi_{\ell,k}
    \Bigr\|_\infty
    &\leq   C \sum_{\ell = 0}^{L_{n,m}} 2^{\ell/2}
    \left( 2^{\ell/2} q \sqrt{\log(1/q)} + \sqrt{\tfrac{\log n}{n}} \right)\\ 
    &\leq   C \left[ 2^{L_{n,m}} q \sqrt{\log(1/q)} + 2^{L_{n,m}/2} \sqrt{\tfrac{\log n}{n}} \right].
\end{aligned}
\end{equation}
For the tail, the assumption $f_0 \in B_{\infty,\infty}^\beta(M)$  together with \eqref{localisation CDV} imply that 
\begin{align*}
    \Bigl\| \sum_{\ell > L_{n,m}} \sum_{k=0}^{2^\ell-1}
    \langle f_0, \psi_{\ell,k} \rangle_2 \, \psi_{\ell,k} \Bigr\|_\infty
    &
    \leq  \sum_{\ell > L_{n,m}} M2^{-\ell(\beta+1/2)}\|\sum_{k=0}^{2^{\ell}-1} |\psi_{\ell,k}(x)| \|_{\infty}\\
    &\lesssim   \sum_{\ell > L_{n,m}} 2^{-\ell\beta}\lesssim 2^{-L_{n,m} \beta} .
\end{align*}
Combining the last two displays with the triangle inequality, on $\mathcal{A}$,
\begin{equation}\label{eq:three terms}
    \bigl\| \hat{f}_n - f_0 \bigr\|_\infty
    \;\leq\; C \left[ 2^{L_{n,m}} q \sqrt{\log(1/q)}
    \;+\; 2^{L_{n,m}/2} \sqrt{\tfrac{\log n}{n}}
    \;+\; \, 2^{-L_{n,m} \beta} \right].
\end{equation}
By definition of $L_{n,m}$ we have
\[
    2^{-L_{n,m}\beta} \leq 2^{\beta}\left[ \Bigl( q \sqrt{ \log(1/q)} \Bigr)^{\frac{\beta}{\beta+1}}
    + \Bigl(\frac{\log n}{n}\Bigr)^{\frac{\beta}{2\beta+1}} \right],
\]
and
\[
    2^{L_{n,m}} q \sqrt{ \log(1/q)} \leq \left ( q \sqrt{\log(1 /q)} \right )^{\frac{\beta}{\beta +1}} ,
    \qquad
    2^{L_{n,m}}\frac{\log n}{n} \leq \Bigl(\frac{\log n}{n}\Bigr)^{\frac{2\beta}{2\beta+1}}.
\]
Substituting these bounds into \eqref{eq:three terms} we get
\[
    \| \hat{f}_n - f_0 \bigr\|_\infty \lesssim \left ( q \sqrt{\log(1/q)} \right )^{\frac{\beta}{\beta +1}} + \left ( \frac{\log n}{n} \right )^{\frac{\beta}{2 \beta +1 }}.
\]
We split into two cases: if $\log n \leq m$, then $q \asymp m/n$, so that $q \sqrt{ \log(1/q)} \asymp m/n \sqrt{\log(n/m)}$. If $m < \log n$, then $q \lesssim \log n/n$ and $\log(1/q) \leq \log n$. It follows that
\[
    \left ( q \sqrt{\log(1/q)} \right )^{\frac{\beta}{\beta +1}} \lesssim \left ( \frac{(\log n)^{\frac{3}{2}}}{n} \right )^{\frac{\beta}{\beta +1}} = o \left ( \left (\frac{\log n}{n} \right )^{\frac{\beta}{2 \beta +1}} \right ).
\]
In both cases we obtain the stated bound. 
\end{proof}

\section{Technical results}

\subsection{Properties of Renyi divergence}\label{app : renyi}
Given $\lambda>0$, a statistical model $(P_f)_{f \in \cF}$ such that $P_f \ll P_{f_0}$ for all $f \in \cF$, we define the $\lambda$-R\'enyi divergence as
\begin{equation}\label{def : renyi}
    D_\lambda(P_f,P_{f_0}) := \frac{1}{\lambda-1} \log \int \left(\frac{dP_f}{dP_{f_0}}\right)^\lambda \, dP_{f_0}.
\end{equation}
One core property of the $D_\lambda$ divergence is the tensorization of product measures, for any $n \ge 1$
\[ D_\lambda(P_f^n,P_{f_0}^n) = nD_\lambda(P_f,P_{f_0}).\]
In the regression setting \eqref{def : model}, the R\'enyi divergence can easily be related to usual norms. In particular, any $D_\lambda$ contraction result can easily be converted into $L_2(P_X)$ contraction under some additional boundedness assumption using the following lemma (see \cite{JMLR:v26:24-0894}, Lemma 21).

\begin{lem}[Linking R\'enyi to $L_2$ and $L_\infty$-norms] \label{lem : L2 from renyi}
    Let $(P_f)$ be the statistical model from \eqref{def : model}, then, for any $\lambda >0$ and $f,f_0 \in \cF$
    \[D_{1+\lambda}(P_{f_0},P_{f}) \le \frac{1 + \lambda}{2 \si^2} \| f - f_0 \|^2_\infty. \]
    Assume further $\|f\|_\infty \vee \|f_0\|_\infty \le F$, for some $F>0$. Then, for any $\lambda \in (0,1)$,
    \[ D_\lambda(P_{f},P_{f_0}) \ge \frac{\lambda}2 e^{-2 F^2 \lambda(1-\lambda)}\| f-f_0\|_{2,P_X}^2.\]
\end{lem}

\subsection{The trimmed-mean estimator for sub-Gaussian variables}
\label{sec : trimmed mean lugosi}

In Corollary~\ref{corollary mendelson} and Theorem \ref{sup norm trimmed estimator} we derive error bounds for the proposed trimmed estimator of the Fourier and wavelet coefficients of the regression function, respectively. In the proofs, we invoke the upper bound, derived in \cite{robust_mean_estimation}, for the corresponding trimmed mean estimator. However, this result could not be applied directly in our setting and we had to adapt it. We show that data splitting is not required and the results also hold without it, improving the practical performance of the approach. We also slightly extend the result by not restricting ourselves to continuous random variables, an extension needed to deal with the wavelet coefficient estimators. For completeness we provide a self-contained proof following the proof of Theorem 1 of \cite{robust_mean_estimation}. For easier readability and clarity of presentation we did not aim to optimize the constant factors in the proof as they do not play an important role in our asymptotic results.

\subsubsection{Results for sub-Gaussian random variables}

Recall that a random variable $Z$ with mean $\mu = \mathbb{E} Z$ is
\emph{sub-Gaussian with variance proxy} $v^2$ if
\begin{equation}
\label{subgaussian}
    \mathbb{E}\, e^{\lambda \bar{Z}} \leq e^{\lambda^2 v^2 /2},
    \qquad \text{for all } \lambda \in \mathbb{R},
\end{equation}
where $\bar{Z} := Z - \mu$. Condition \eqref{subgaussian} implies in particular
that $\operatorname{Var}(Z) \leq v^2$ and that, for all $t \geq 0$,
\[
    \mathbb{P}(\bar{Z} \geq t) \leq e^{-t^2/(2 v^2)},
    \qquad
    \mathbb{P}(-\bar{Z} \geq t) \leq e^{-t^2/(2 v^2)}.
\]
Following \cite{robust_mean_estimation}, for $p \in (0,1)$ we write
$Q_p(Z) := \sup\{ M \in \mathbb{R} : \mathbb{P}(Z \geq M) \geq 1 - p \}$ for the $p$-quantile, so that $\mathbb{P}(Z <  Q_p(Z)) \leq p \leq \mathbb{P}(Z \leq  Q_p(Z)) $. We will use the following
elementary properties of sub-Gaussian random variables; see Lemma~$5.5$ in \cite{Vershynin_2012} for part~$1$,
Lemma~$1.3$ in \cite{rigollet2015highdimensional} for part~$2$, and the
remark on p.~$399$ in \cite{robust_mean_estimation} for part~$3$.

\begin{lem}
\label{lemma subgaussian property}
    Let $Z$ satisfy \eqref{subgaussian}. Then
    \begin{enumerate}
        \item  for $ p \geq 1$,
        $\left(\mathbb{E}|\bar{Z}|^p\right)^{1/p} \leq C_0\, v \sqrt{p}$;
        \item for $\rho \in (0,1/2]$,
        $|Q_\rho(Z) - \mu| \vee |Q_{1-\rho}(Z) - \mu|
        \leq v \sqrt{2 \log(1/\rho)}$;
        \item for $\rho \in (0,1/2]$,
        \[
            \mathbb{E}\!\left[|\bar{Z}|\,
            \mathbf{1}_{\{Z > Q_{1-\rho}(Z)\}}\right]
            \vee
            \mathbb{E}\!\left[|\bar{Z}|\,
            \mathbf{1}_{\{Z < Q_{\rho}(Z)\}}\right]
            \leq C_1\, v\, \rho \sqrt{\log(1/\rho)}.
        \]
    \end{enumerate}
\end{lem}

\subsubsection{Trimmed-mean estimator}

We briefly recall the construction of the trimmed-mean estimator of
\cite{robust_mean_estimation}. Let $Z_1,\dots,Z_n$ be a
sample of independent copies of $Z$. Suppose the
statistician receives corrupted data
$\widetilde{Z}_1,\dots,\widetilde{Z}_n$, obtained after an
adversary has changed at most $m$ observations in total across $Z_{1:n}$. Recall the clipping function $\phi_{a,b}$ from
\eqref{clipping function} and the trimming level $q$ defined as follows
\[
    q =q(m,n,\delta)= 8 \frac{m}{n} + 12\frac{\log(8/\delta)}{n}
\]
We denote by 
$\widetilde{a} = \widetilde{Z}_{(qn)}$ and
$\widetilde{b} = \widetilde{Z}_{(n - nq)}$ the $q$-th and $(n- nq)$-th order statistics of the corrupted sample, respectively. We define the estimator
\begin{equation}
\label{trimmed mean estimator}
    \hat{\mu} := \frac{1}{n} \sum_{i=1}^n
    \phi_{\widetilde{a},\widetilde{b}}(\widetilde{Z}_i).
\end{equation}

\begin{thm}{(Theorem~$1$ of \cite{robust_mean_estimation}, sub-Gaussian case)}
\label{lugosi and mendelson theorem}
Let $Z$ be a random variable with variance $s^2 = \text{Var}(Z)$ that satisfies \eqref{subgaussian} with mean $\mu$ and variance proxy $v^2$. Let
$\delta \in (0,1/2]$ and suppose that the contamination proportion $m/n$ and $\delta$ are such that $q \leq 1/4$. Then the
estimator \eqref{trimmed mean estimator} satisfies, with probability at least
$1-\delta$,
\begin{equation}
    |\hat{\mu} - \mu|
    \leq C  \left( v q \sqrt{\log(1/q)}
    + s \sqrt{\frac{\log(1/\delta)}{n}}\right).
\end{equation}
Consequently, there is an absolute constant $C'$ such that, with probability at
least $1-\delta$,
\begin{equation}
    |\hat{\mu} - \mu|
    \leq C'\, v \left( \frac{m}{n} \sqrt{\log \left( \frac{n}{m} \right) }
    + \sqrt{\frac{\log(1/\delta)}{n}}\right).
\end{equation}
\end{thm}

\begin{proof}
For convenience we introduce the notation $\eta=\eta_{n,m}:=m/n$.

\emph{Step 1: control of the data-driven thresholds $\widetilde{a},
\widetilde{b}$.}
Let us first consider the clean sample $Z_1,\dots,Z_n$ and recall that $\bar{Z} = Z - \mu$.  The indicators
$\mathbf{1}_{\{Z_i \geq \mu + Q_{1-2q}(\bar{Z})\}}$ are i.i.d.\ Bernoulli with
mean $\mathbb{P}(\bar{Z} \geq Q_{1-2q}(\bar{Z})) \geq 2q$. Their sum therefore dominates a $\text{Binom}(n,2q)$ random variable.
By Bernstein's inequality,
\[
    \bigl|\{ i : Z_i \geq \mu + Q_{1-2q}(\bar{Z})\}\bigr| \geq \tfrac{5}{4} q n
\]
with probability at least $1 - e^{- q n / 8}$. 
Applying a similar argument we can show that there exists an event $E$ that depends only on $Z_1,...,Z_n$ and satisfies, by a union bound, $\mathbb{P}(E) \geq 1 - \delta/2$, on which all four bounds 
\begin{equation}
\begin{aligned}
    &\bigl|\{ i : Z_i \geq \mu + Q_{1-2q}(\bar{Z})\}\bigr|
        \geq \tfrac{5}{4} q n,
    &&\bigl|\{ i : Z_i \leq \mu + Q_{1-q/4}(\bar{Z})\}\bigr|
        \geq \left (1 - \frac{3}{4} q \right ) n, \\
    &\bigl|\{ i : Z_i \leq \mu + Q_{2q}(\bar{Z})\}\bigr|
        \geq \tfrac{5}{4} q n,
    &&\bigl|\{ i : Z_i \geq \mu + Q_{q/4}(\bar{Z})\}\bigr|
        \geq \left (1 - \frac{3}{4} q \right ) n,
\end{aligned}
\end{equation}
hold simultaneously. Since the adversary changes at most
$\eta n$ of the $Z_i$, each of the four counts computed from
the corrupted $\widetilde{Z}_i$ changes by at most $ \eta n$. As $q > 8 \eta$, on $E$ we obtain 
\[
    \bigl|\{ i : \widetilde{Z}_i \geq \mu + Q_{1-2q}(\bar{Z})\}\bigr| > q n,
    \qquad
    \bigl|\{ i : \widetilde{Z}_i \leq \mu + Q_{1-q/4}(\bar{Z})\}\bigr|
        > (1 - q) n.
\]
The first bound forces the order statistic
$\widetilde{b} = \widetilde{Z}_{(n - nq)}$ to be at least
$\mu + Q_{1-2q}(\bar{Z})$, the second forces it to be at most
$\mu + Q_{1-q/4}(\bar{Z})$; the analogous bounds for
$\widetilde{a} = \widetilde{Z}_{(qn)}$ follow in the same way. In summary, on $E$,
\begin{equation}
\label{sandwich quantiles}
    Q_{1-2q}(\bar{Z}) \leq \widetilde{b} - \mu \leq Q_{1-q/4}(\bar{Z}),
    \qquad
    Q_{q/4}(\bar{Z}) \leq \widetilde{a} - \mu \leq Q_{2q}(\bar{Z}).
\end{equation}

\emph{Step 2: upper bound on the clean average.}
The function $\phi_{a,b}$, defined in \eqref{clipping function}, is nondecreasing in $(a,b)$ and, by
\eqref{sandwich quantiles}, $\widetilde{a} \leq \mu + Q_{2q}(\bar{Z}) =: r_1$ and
$\widetilde{b} \leq \mu + Q_{1-q/4}(\bar{Z}) =: r_2$ (on the event $E$). It follows that on $E$,
\[
    \phi_{\widetilde{a},\widetilde{b}}(x) \leq \phi_{r_1,r_2}(x) =: \psi(x)
    \qquad \text{for all } x .
\]
The variables $\psi(Z_i)$ are i.i.d.\ copies of $\psi(Z)$. Because
$\psi(x) - x = (r_1-x)_+ - (x-r_2)_+$,
\[
    \mathbb{E}\psi(Z) = \mu + \mathbb{E}(r_1 - Z)_+ - \mathbb{E}(Z - r_2)_+,
    \qquad
    |\mathbb{E}\psi(Z) - \mu| \leq \mathbb{E}(r_1 - Z)_+ \vee \mathbb{E}(Z - r_2)_+ .
\]
On $\{Z < r_1\} = \{\bar{Z} < Q_{2q}(\bar{Z})\}$ (event of probability at most $2q$), one
has $0 \leq r_1 - Z = Q_{2q}(\bar{Z}) - \bar{Z}
\leq |Q_{2q}(\bar{Z})| + |\bar{Z}|$, so
\[
    \mathbb{E}(r_1 - Z)_+
    \leq 2q\,|Q_{2q}(\bar{Z})|
    + \mathbb{E}\!\left[|\bar{Z}|\,
        \mathbf{1}_{\{\bar{Z} < Q_{2q}(\bar{Z})\}}\right]
    \leq C\, v\, q \sqrt{\log(1/q)},
\]
where the two terms are bounded by points~$2$ and $3$ of
Lemma~\ref{lemma subgaussian property}, both with $\rho = 2q$. The same argument
with $\rho = q/4$ gives $\mathbb{E}(Z - r_2)_+ \leq C\, v\, q \sqrt{\log(1/q)}$,
hence
\begin{equation}
\label{close to the mean}
    \mu - C\, v\, q \sqrt{\log(1/q)}
    \leq \mathbb{E}\psi(Z)
    \leq \mu + C\, v\, q \sqrt{\log(1/q)} .
\end{equation}
Since $\psi$ is $1$-Lipschitz, with an independent copy $Z'$ of $Z$,
\[
    \operatorname{Var}(\psi(Z))
    = \tfrac{1}{2}\mathbb{E}\bigl(\psi(Z) - \psi(Z')\bigr)^2
    \leq \tfrac{1}{2}\mathbb{E}(Z - Z')^2
    = \operatorname{Var}(Z) = s^2,
\]
and $\psi$ takes values in $[r_1,r_2]$, an interval of length
$r_2 - r_1 = Q_{1-q/4}(\bar{Z}) - Q_{2q}(\bar{Z}) \leq 2\kappa$, where
$\kappa := v\sqrt{2\log(4/q)}$ (Lemma~\ref{lemma subgaussian property},
point~$2$). Consequently $W_i := \psi(Z_i) - \mathbb{E}\psi(Z)$ satisfies
\[
    \mathbb{E} W_i = 0, \qquad |W_i| \leq 2\kappa, \qquad
    \sum_{i=1}^n \mathbb{E} W_i^2 = n\,\operatorname{Var}(\psi(Z)) \leq n s^2 .
\]
We use the following form of Bernstein's inequality: if $W_1,\dots,W_n$ are
independent, mean-zero, with $|W_i| \leq b_0$ and
$\sum_{i=1}^n \mathbb{E} W_i^2 \leq V$, then for every $x > 0$,
\[
    \mathbb{P}\!\left(\sum_{i=1}^n W_i
        \geq \sqrt{2 V x} + \tfrac{1}{3} b_0 x\right) \leq e^{-x} .
\]
Applying it with $V = n s^2$, $b_0 = 2\kappa$ and $x = \log(8/\delta)$, we get, with probability at least $1 - \delta/8$,
\begin{equation}
\label{berns bound}
    \frac{1}{n}\sum_{i=1}^n \psi(Z_i)
    \leq \mathbb{E}\psi(Z) + s \sqrt{\frac{2\log(8/\delta)}{n}}
    + \frac{2\kappa}{3}\frac{\log(8/\delta)}{n} .
\end{equation}
The last term in the sum of the above display is at most $C\, v\, q\sqrt{\log(1/q)}$, 
since
$\log(8/\delta)/n \leq q/12$ and we can bound $\kappa$ by $C' v \sqrt{\log(1 /q)}$ since $\log(4 /q) \leq 2 \log(1/q)$, under the assumption that $q \leq 1/4$. 

Combining \eqref{berns bound} with \eqref{close to the mean} and the pointwise bound $\phi_{\widetilde{a},\widetilde{b}} \leq \psi$ (on the event $E$),
\[
    \frac{1}{n}\sum_{i=1}^n \phi_{\widetilde{a},\widetilde{b}}(Z_i)
    \leq \frac{1}{n}\sum_{i=1}^n \psi(Z_i)
    \leq \mu + C\, v\, q\sqrt{\log(1/q)} + s \sqrt{\frac{2\log(8/\delta)}{n}} ,
\]
holds on the event of \eqref{berns bound} intersected with $E$.

\emph{Step 3: lower bound and corruption.}
On $E$ we also have, by \eqref{sandwich quantiles},
\[
    \phi_{\widetilde{a},\widetilde{b}}(x)
\geq \phi_{\mu + Q_{q/4}(\bar{Z}),\,\mu + Q_{1-2q}(\bar{Z})}(x) =: \psi_-(x),
\]
for all $x$, and repeating the argument for $\psi_-$ gives a second Bernstein event of probability at least $1 - \delta/8$ (we apply Bernstein as in \eqref{berns bound} but for $\psi_-$). On the intersection of this event with $E$,
\[
    \frac{1}{n}\sum_{i=1}^n \phi_{\widetilde{a},\widetilde{b}}(Z_i)
    \geq \mu - C\, v\, q\sqrt{\log(1/q)} - s \sqrt{\frac{2\log(8/\delta)}{n}} .
\]
The two Bernstein events hold simultaneously with probability at least $1 - \delta/4$. On their intersection with $E$,
\[
    \left| \frac{1}{n}\sum_{i=1}^n \phi_{\widetilde{a},\widetilde{b}}(Z_i)
    - \mu \right|
    \leq C\, v\, q\sqrt{\log(1/q)} + s \sqrt{\frac{2\log(8/\delta)}{n}} .
\]
Finally, the clean and corrupted trimmed averages differ at at most $m$ indices where $\widetilde{Z}_i \neq Z_i$, and at each such index both
$\phi_{\widetilde{a},\widetilde{b}}(\widetilde{Z}_i)$ and
$\phi_{\widetilde{a},\widetilde{b}}(Z_i)$ lie in $[\widetilde{a},\widetilde{b}]$,
an interval of length $\widetilde{b} - \widetilde{a} \leq 2\kappa$ on $E$. Hence, on
$E$,
\[
    \left| \frac{1}{n}\sum_{i=1}^n
        \phi_{\widetilde{a},\widetilde{b}}(\widetilde{Z}_i)
    - \frac{1}{n}\sum_{i=1}^n \phi_{\widetilde{a},\widetilde{b}}(Z_i) \right|
    \leq \frac{m}{n}\, 2\kappa
    \leq C\, v\, q\sqrt{\log(1/q)},
\]
using $m/n \leq q/8$. Combining the last two displays with the triangle
inequality, on the event $E$ intersected with the events of the above Bernstein inequalities, of probability at least $1-\delta/2-\delta/4=1- 3 \delta/4 \geq 1 - \delta$ ,
\[
    |\hat{\mu} - \mu|
    \leq C\, v\, q\sqrt{\log(1/q)} + s \sqrt{\frac{2\log(8/\delta)}{n}}
    \leq C \left( v q\sqrt{\log(1/q)} + s \sqrt{\frac{\log(1/\delta)}{n}}\right)
\]
which is the first claim.
\end{proof}

\subsection{On the map-formulation of adversarial contamination}
\label{sec : on the map formulation of adv cont}

Given two $n$-tuples $x_{1:n}$ and $y_{1:n}$, we denote by $d_{h}(x_{1:n}, y_{1:n})$ the Hamming distance between $x_{1:n}$ and $y_{1:n}$, namely
\[
    d_{h}(x_{1:n}, y_{1:n}) := | \{ i \in  [n] : x_i \neq y_i \} |,
\]
and by $B_m^{h}(x_{1:n}) := \{ \widetilde{x}_{1:n} : d_{h}(\widetilde{x}_{1:n},x_{1:n}) \leq m \}$ the ball of radius $m$ centered at $x_{1:n}$ with respect to the Hamming distance. Moreover, given a (standard) Borel space $\mathcal{X}$, we denote its Borel $\sigma$-field by $\cB(\mathcal{X})$.

Let $X_1,\dots,X_n$ be i.i.d. random variables with law $P_{\theta_0}$ on $(\cX, \cB(\cX))$, where $\theta_0\in\Theta$. In the main text of the paper, we defined adversarial contamination as follows. Given a probability space $(\cU, \cB(\cU), \nu)$, let $T: \cX^n \times \cU \to \cX^n$ be a measurable map such that $d_{h}(T(x_{1:n}, u), x_{1:n}) \leq m $ for all $(x_{1:n}, u) \in \cX^n \times \cU$. The contaminated data can be written as $(\widetilde{X}_1,...,\widetilde{X}_n) = T(X_{1:n}, U)$, where $U$ is a $\cU$-valued random variable with distribution $\nu$, independent of $X_{1:n}$. We denote by $\cT_m$ the collection of all such maps $T$ allowing the auxiliary probability space $(\cU,\cB(\cU),\nu)$ to vary across the maps. 

The law of the contaminated sample is then the push-forward through the measurable map $T$ of the probability measure $Q_{\theta_0} := P_{\theta_0}^{n} \otimes \nu$. We denote by $E_{Q_{\theta_0}}$ the expectation associated with $Q_{\theta_0}$. To study theoretical properties of statistical procedures under contamination, we consider the following quantity
\begin{equation}
\label{sup over maps}
 \sup_{T \in \cT_m} E_{Q_{\theta_0}}[ \Phi_{\theta_0}(T(X_{1:n},U))] = \sup_{T \in \cT_m} \widetilde{E}_{\theta_0}[ \Phi_{\theta_0}(\widetilde{X}_{1:n})],
\end{equation}
where $\Phi_{\theta_0}(\widetilde{X}_{1:n})$ denotes either $\mathbf{1}{ \{ L(\hat{\theta}_n(\widetilde{X}_{1:n}), \theta_0) > r(n,m) \} }$, for an estimator $\hat{\theta}_n$, a pseudo-metric $L(\cdot,\cdot)$ and a rate $r(n,m)$, or $G( S_{n,m}(\theta_0) \mid \widetilde{X}_{1:n})$ for a measurable set $S_{n,m}(\theta_0)\subset\Theta$ and a generalized posterior $G(\ \cdot \ \mid \widetilde{X}_{1:n})$ on $\Theta$.

Alternatively, we could define the contamination in terms of Markov kernels. We recall that a  (Markov) probability kernel $A : \cB(\cX^n) \times \cX^n \to [0,1]$ satisfies the following two properties: for each fixed $x_{1:n}$, the map $B\mapsto A(B , x_{1:n})$ is a probability measure on $(\cX^n,\cB(\cX^n))$ and, for each measurable set $B$, the map $x_{1:n} \mapsto A( B ,x_{1:n}) $ is $\cB(\cX^n)$-measurable. One can assume that $A$ represents a regular conditional distribution of the contaminated data $\widetilde{X}_{1:n}$ given the uncontaminated data $X_{1:n}$. To highlight this conditional structure we write $A(B|x_{1:n})=A(B,x_{1:n})$. In this case too, we have to constrain the number of contamination to be at most $m$, by requiring that, for all $x_{1:n} \in \cX^{n}$
\begin{equation}
\label{condition on the support of the kernel}
    A( B_m^{h}(x_{1:n}) \mid x_{1:n}) = 1.
\end{equation}
We denote by $\mathfrak{M}_m$ the collection of all such kernels and, as in equation \eqref{sup over maps}, we study
\begin{equation}
\label{sup over kernels}
    \sup_{A \in \mathfrak{M}_m} \widetilde{E}_{\theta_0, A} [ \Phi_{\theta_0}(\widetilde{X}_{1:n})],
\end{equation}
where $\widetilde{P}_{\theta_0,A}$ denotes the marginal law of $\widetilde{X}_{1:n}$ obtained by integrating the Markov kernel $A$ with respect to the law of the sample $X_{1:n}$ and $\widetilde{E}_{\theta_0,A}$ is the corresponding expectation operator.

The following lemma establishes that the two formulations lead to the same worst-case expectation, namely
\begin{equation}
    \sup_{T \in \cT_m}  E_{Q_{\theta_0}}[ \Phi_{\theta_0}(T(X_{1:n},U))] = \sup_{A \in \mathfrak{M}_m} \widetilde{E}_{\theta_0, A} [ \Phi_{\theta_0}(\widetilde{X}_{1:n})].
\end{equation}
\begin{lem}
\label{equivalence between sup over kernels and sup over maps}
 For any probability kernel $A \in \mathfrak{M}_m$, there exist a probability space $(\cU, \cB(\cU),\nu)$, a $\cU$-valued random variable $U \sim \nu$ independent of $X_{1:n}$ and a measurable map $T \in \cT_m$, such that $T(x_{1:n}, U) \sim A(\cdot \mid x_{1:n}) $ for all $x_{1:n}$. Conversely, given any probability space  $(\cU, \cB(\cU),\nu)$, a $\cU$-valued random variable $U \sim \nu$ and a measurable map $T \in \cT_m$, there exists a kernel $A \in \mathfrak{M}_m$, such that $T(x_{1:n}, U) \sim A(\ \cdot \ \mid x_{1:n})$ for all $x_{1:n}$. 
\end{lem}
\begin{proof}
    For the first statement, the randomization lemma (Lemma $3.2$ (vii) of \cite{Kallenberg2021},) provides a measurable map $T_0 : \cX^n \times [0,1] \to \cX^n$ such that $T_0(x_{1:n},U) \sim A(\cdot  \mid x_{1:n})$ for every $x_{1:n}$, where $U$ is uniformly distributed on $[0,1]$. Take $\cU = [0,1]$ and $\nu$ be the uniform probability measure. Since $A \in \mathfrak{M}_m$, we have that $\nu(\{ u: T_0(x_{1:n},u) \notin B_m^h(x_{1:n}) \}) = 0$ for all $x_{1:n}$. Thus we can take $T(x_{1:n},u)$ to be equal to $T_0(x_{1:n},u)$ outside this $\nu$-measure $0$ set and equal to $x_{1:n}$ otherwise, so that $T \in \cT_m$.
    
    To prove the other direction we construct the required kernel $A$ directly by setting, for each measurable set $B$ and all $x_{1:n}$
    \begin{equation}
    \label{first definition of A}
        A(B \mid x_{1:n}) := \nu( \{ u : T(x_{1:n}, u) \in B \}).
    \end{equation}
    Note that the map $u \mapsto  T(x_{1:n}, u)$ for a fixed $x_{1:n}$ is the section of a measurable map, so that $A(\ \cdot \ \mid x_{1:n})$ is the pushforward of $\nu$ through a measurable map, and hence a well-defined probability measure. Moreover, the definition in \eqref{first definition of A} is equivalent to
    \[
        A(B \mid x_{1:n}) = \int_{\cU} \mathbf{1}_{B}(T(x_{1:n},u)) \nu(du),
    \]
    so that the integral is measurable in $x_{1:n}$ by Fubini-Tonelli's theorem, the integrand being jointly measurable in $(x_{1:n},u)$ and nonnegative. Since $T \in \cT_m$, we have that $A$ satisfies the condition in Equation \eqref{condition on the support of the kernel}. 
\end{proof}

Note that in \eqref{sup over maps} the contamination map $T$ is deterministic and the supremum of the risk is taken, over all such maps $T$, outside of the expectation. The question naturally arises whether one could consider data driven contamination functions and taking the supremum inside of the expectation. This intuitively corresponds to a stronger poisoning setting, where the adversary, after seeing the exact realization of the data, may choose the contamination function for corrupting the data. This worst-case scenario can be formalized by considering the following risk 
\begin{equation}
\label{worst case risk}
    E_{\theta_0}\left [ \sup_{ \widetilde{x}_{1:n} \in  B_m^{h}(X_{1:n})} \Phi_{\theta_0}( \widetilde{x}_{1:n}) \right ],
\end{equation}

\begin{Remark}
We show below that under mild regularity conditions on the function $\Phi_{\theta_0}(\cdot)$ the above supremum is universally measurable and hence the above expectation has to be meant with respect to the completion of $P_{\theta_0}^n$. More precisely, one has to take the completion of $\mathcal{B}(\mathcal{X}^n)$ with respect to $P_{\theta_0}^n$. The completion of $P_{\theta_0}^n$, denoted as $\overline{P}_{\theta_0}^n$ (see the proof below), is then the unique measure that extends $P_{\theta_0}^n$ on this completed sigma-algebra.
\end{Remark}

In the next lemma, we follow and extend the arguments of Lemma B.2 in \cite{worstcaseramdas} by considering the supremum over data dependent maps required in our setting. We show that for Borel measurable maps the two formulations given in \eqref{worst case risk} and \eqref{sup over maps} are equivalent. Hence moving the supremum inside, which typically corresponds to worse case scenario, does not make the introduced contamination framework stronger.

\begin{lem}
\label{lem : worst case and maps}
    Let $\Phi_{\theta_0}: \cX^n \to [0,1]$ be a Borel measurable function. Then the map \[ x_{1:n} \mapsto s(x_{1:n}) := \sup_{\widetilde{x}_{1:n} \in B^{h}_m(x_{1:n}) } \Phi_{\theta_0}(\widetilde{x}_{1:n}) \] is universally measurable and
    \begin{equation}
    \label{statement pathwise}
        E_{\theta_0} [ s(X_{1:n})] = \sup_{ T \in \cT_{m}} E_{Q_{\theta_0}} \left [ \Phi_{\theta_0}(T(X_{1:n},U)) \right],
    \end{equation}
    where the expectation on the left hand side is taken with respect to the completion of $P_{\theta_0}^n$.
\end{lem}

\begin{proof}
    For brevity, throughout the proof we write $x$ and $\widetilde{x}$ for $x_{1:n}$ and $\widetilde{x}_{1:n}$, respectively. For $c \in \mathbb{R}$, consider the super-level sets $S_c := \{ x : s(x) > c\}$. Note that $x \in S_c$ if and only if there exists $\widetilde{x} \in \cX^n$ such that $\Phi_{\theta_0}(\widetilde x) \mathbf{1}_{\{ \widetilde x \in B^{h}_m(x) \} } > c$ or, equivalently 
    \begin{equation*}
        S_c = \pi_{\cX^n}( \underbrace{\{ (x,\widetilde x) :  \Phi_{\theta_0}(\widetilde x) \mathbf{1}_{\{ \widetilde x \in B^{h}_m(x) \}}  > c   \} }_{=: A_c}           ),
    \end{equation*}
    where $\pi_{\cX^n}(\cdot)$ denotes the projection onto the first $n$-coordinates. Since the map $(x, \widetilde x) \mapsto \Phi_{\theta_0}(\widetilde x) \mathbf{1}_{\{ \widetilde x \in B^{h}_m(x) \}}$ is Borel measurable on $\cX^{n} \times \cX^n$, being the product of Borel measurable maps, the set $A_c$ is Borel. The projection of a Borel set on a product space is analytic (see Exercise $14.3$ of \cite{Kechris1995}), and, thus universally measurable (Theorem 21.10 of \cite{Kechris1995}). It follows that the map $s$ is universally measurable. In particular, $s$ is measurable with respect to the completion of $\cB(\cX^n)$ under $P_{\theta_0}^n$  and, since $0 \leq s \leq 1$, the expectation
    \[
            E_{\theta_0}\big[ s(X_{1:n}) \big] := \int_{\cX^n} s \, d\overline{P}^n_{\theta_0}
    \]
    is well defined, where $\overline{P}_{\theta_0}^{\,n}$ denotes the completion of $P_{\theta_0}^n$. This proves the first statement. Now, since $T \in \cT_m$, we have that $T(x,u) \in B^{h}_m(x)$ and so $\Phi_{\theta_0}(T(x,u)) \leq s(x)$ for all $(x,u)$. Taking expectations and then taking the supremum over all $T \in \cT_m$, we see that the right hand side of \eqref{statement pathwise} is at most the left hand side. 

    Take a Borel map $s'$ that equals $s$ for $P_{\theta_0}^n$-almost every $x$ and such that $s' \leq s$ for all $x$. For completeness we provide such construction. Recall that $s$ is universally measurable, thus for any $j \in \mathbb{Q}$, there are two Borel sets $A_j,B_j$ such that $\{ s(x) > j \} \Delta A_j \subseteq B_j$ and $P_{\theta_0}^n(B_j) = 0$. Let us define the Borel set $B := \cup_{j} B_j$ such that $P_{\theta_0}^n(B) = 0$ and $s' = s \mathbf{1}_{B^c}$. We have that $\{s' \neq s\} \subseteq B$, $s' \leq s$ and for all $j \in \mathbb{Q}$ such that $j > 0$ we have $\{ s'(x) > j \} = \{ s > j \} \cap B^c = A_j \cap B^c $ which is Borel. If $j < 0$, since $s' \geq 0$, we have $\{ s' > j\} = \mathcal{X}^n$ that is also a Borel set.
    
    Then, for $\epsilon \in (0,1)$, let us define the following set
    \[
        C_\epsilon  := \left \{ (x,\widetilde x) : \Phi_{\theta_0}(\widetilde x) \mathbf{1}_{\{ \widetilde x \in B^{h}_m(x) \} } > s'(x)  - \epsilon \right \},
    \]
    which is Borel since $\Phi_{\theta_0}(\widetilde x) \mathbf{1}_{\{ \widetilde x \in B^{h}_m(x) \} }$ and $s'(x)$ are Borel. For any $x$, $s'(x) - \epsilon < s'(x) \leq s(x)$, so that, by definition of supremum, there exists some $\widetilde x \in \cX^n$ such that $(x,\widetilde x) \in C_\epsilon$. It follows that $\pi_{\cX^n}(C_\epsilon) = \cX^n$. By the Jankov--von Neumann theorem (Theorem $18.1$ of \cite{Kechris1995}), there is a universally measurable map $\psi_\epsilon : \cX^n \to \cX^n$ whose graph is contained in $C_\epsilon$,  namely
    \[
        \Phi_{\theta_0}(\psi_\epsilon(x)) \mathbf{1}_{\{ \psi_\epsilon(x) \in B^{h}_m(x) \} } > s'(x) - \epsilon.
    \]
    Since the map $\psi_\epsilon$ is only universally measurable, take a Borel measurable map $\psi_\epsilon': \cX^n \to \cX^n$ which equals $\psi_\epsilon$, $P_{\theta_0}^n$-almost surely and define, for all $x \in \cX^n$
    \[
      T_\epsilon(x) := \begin{cases} \psi_\epsilon'(x) & \text{if } d_{h}(\psi_\epsilon'(x), x) \leq m, \\
        x & \text{otherwise.}  \end{cases}
    \]
    The map $T_\epsilon$ is Borel measurable and, by construction, $T_\epsilon \in \cT_m$. Moreover, it satisfies, for $P_{\theta_0}^n$-almost every $x$
    \[
          \Phi_{\theta_0}(T_\epsilon(x)) \geq \Phi_{\theta_0}(\psi_\epsilon'(x)) \mathbf{1}_{\{ \psi_\epsilon'(x) \in B^{h}_m(x) \} } > s'(x) - \epsilon.
    \]
    Integrating both sides with respect to $P_{\theta_0}^n$ and upper bounding the left hand side by  the supremum over $\cT_m$, one gets
    \[
        \sup_{T \in \cT_m} E_{Q_{\theta_0}}[ \Phi_{\theta_0}(T(X,U)) ] \geq E_{\theta_0} [ \Phi_{\theta_0} (T_\epsilon(X))] \geq E_{\theta_0}[ s(X)] - \epsilon.
    \]
    Letting $\epsilon \to 0$, we obtain the result.
\end{proof}

\begin{Remark}

In principle one could alternatively define the contraction of the (generalized) posterior under contamination as follows. We say that a generalized posterior contracts around $\theta_0$ with respect to some loss $L$ at rate $\zeta_{n,m}$ if there exist a  large enough constant $M > 0$ and a measurable function $\Psi_{n,m} : \cX^n \to \mathbb{R}$ such that, for all $x_{1:n} \in \cX^{n}$ and all $\widetilde{x}_{1:n} \in B^{h}_{m}(x_{1:n})$
\begin{align}
    G( \{ \theta : L(\theta, \theta_0) > M \zeta_{n,m} \} \mid \widetilde{x}_{1:n}) \leq \Psi_{n,m}( x_{1:n}), \quad E_{\theta_0} \Psi_{n,m}(X_{1:n}) \to 0, \quad \text{ as } n \to \infty.\label{def:cont:alter}
\end{align}

Note that the inequality above holds also for the supremum over $\widetilde{x}_{1:n} \in B^{h}_{m}(x_{1:n})$. Then taking the expectation of this supremum under the completion of the clean data generating distribution $\overline{P}_{\theta_0}^n$ is upper bounded by $E_{\theta_0} \Psi_{n,m}(X_{1:n})$, which tends to zero. Hence this definition also implies the contraction under contamination as given in \eqref{worst case risk}.

However, the formulation with a map $T$ has several advantages. First, it provides a more natural connection with the lower bound results and the techniques applied to derive it. Furthermore, we believe, that this formulation is clearer and more intuitive than the one given in \eqref{def:cont:alter}. Therefore, despite the latter being also an equivalent alternative, we apply the formulation with the map $T$ throughout our paper.
\end{Remark}



\end{document}